\def\elsevier{0}
\def\arxiv{1}
\def\shortprep{0}
\def\longprep{0}
\def\longversion{1}
\def\ijwmip{0}

\documentclass[a4paper]{article}
\usepackage{amsmath,amssymb,amsthm}
\usepackage[utf8]{inputenc}

\if\elsevier0
\newcommand{\etal}{et al.}
\fi

\newcommand{\qspace}{\;\;\;\;}
\newcommand{\spaceafter}{\;\;\;\;}
\newcommand{\defterm}[1]{{\textit{#1}}}
\newcommand{\mraprepspace}{of }
\newcommand{\mraprepdomain}{on }
\newcommand{\projprep}{of }

\newcommand{\dualelem}[1]{{\tilde{{#1}}}}
\newcommand{\indexeddualelem}[2]{{\tilde{{#1}}_{#2}}}
\newcommand{\vectorstyle}[1]{{\mathbf{{#1}}}}
\newcommand{\indexedvectorstyle}[2]{{\mathbf{{#1}}_{#2}}}
\newcommand{\seqstyle}[1]{{\mathbf{{#1}}}}
\newcommand{\indexedseqstyle}[2]{{\mathbf{{#1}}_{#2}}}
\newcommand{\znstyle}[1]{{\mathbf{{#1}}}}

\newcommand{\znstyleprime}[1]{{{\mathbf{{#1}}}'}}

\newcommand{\firstznvarsymbol}{k}
\newcommand{\secondznvarsymbol}{\ell}
\newcommand{\thirdznvarsymbol}{m}
\newcommand{\firstznvar}{{\znstyle{\firstznvarsymbol}}}
\newcommand{\secondznvar}{{\znstyle{\secondznvarsymbol}}}
\newcommand{\thirdznvar}{{\znstyle{\thirdznvarsymbol}}}
\newcommand{\indexedfirstznvar}[1]{{\firstznvar_{#1}}}
\newcommand{\firstznvarp}{{{\znstyle{\firstznvarsymbol}}'}}
\newcommand{\zos}{{\znstyle{s}}}
\newcommand{\zot}{{\znstyle{t}}}
\newcommand{\zna}{{\znstyle{a}}}
\newcommand{\znb}{{\znstyle{b}}}
\newcommand{\zns}{{\znstyle{s}}}
\newcommand{\znt}{{\znstyle{t}}}

\newcommand{\zosp}{{\znstyleprime{s}}}

\newcommand{\rna}{{\vectorstyle{a}}}
\newcommand{\rnb}{{\vectorstyle{b}}}
\newcommand{\rnh}{{\vectorstyle{h}}}
\newcommand{\rns}{{\vectorstyle{s}}}
\newcommand{\rnt}{{\vectorstyle{t}}}
\newcommand{\rnx}{{\vectorstyle{x}}}
\newcommand{\rny}{{\vectorstyle{y}}}
\newcommand{\sqa}{\seqstyle{a}}
\newcommand{\sqb}{\seqstyle{b}}
\newcommand{\sqc}{\seqstyle{c}}
\newcommand{\sqd}{\seqstyle{d}}
\newcommand{\seqstyler}{\seqstyle{r}}
\newcommand{\sqs}{\seqstyle{s}}
\newcommand{\dualw}{{\tilde{w}}}
\newcommand{\dualf}{{\tilde{f}}}
\newcommand{\rnxo}{{\indexedvectorstyle{x}{0}}}
\newcommand{\rnxone}{{\indexedvectorstyle{x}{1}}}
\newcommand{\rnhone}{{\indexedvectorstyle{h}{1}}}

\DeclareMathOperator{\closop}{clos}
\DeclareMathOperator{\spanop}{span}

\DeclareMathOperator{\suppop}{supp}
\DeclareMathOperator{\borel}{bor}

\newcommand{\myepsilon}{\varepsilon}
\newcommand{\myphi}{\varphi}
\newcommand{\bphi}{\bar{\myphi}}
\newcommand{\br}{\bar{r}}

\newcommand{\bbigp}{\bar{P}}
\newcommand{\bbigq}{\bar{Q}}
\newcommand{\bbigv}{\bar{V}}
\newcommand{\balpha}{\bar{\alpha}}
\newcommand{\bbfh}{\mathbf{\bar{h}}}

\newcommand{\setsep}{:}
\newcommand{\realnumbers}{\mathbb{R}}
\newcommand{\complexnumbers}{\mathbb{C}}
\newcommand{\integernumbers}{\mathbb{Z}}
\newcommand{\naturalnumbers}{\mathbb{N}}
\newcommand{\positiverealnumbers}{\realnumbers_{+}}
\newcommand{\positiveintegers}{\integernumbers_{+}}
\newcommand{\nonnegrealnumbers}{\realnumbers_{0}}
\newcommand{\extreals}{{\realnumbers_*}}
\newcommand{\extnonnegreals}{{\realnumbers_{0*}}}
\newcommand{\intersection}{\cap}
\newcommand{\rn}{{\realnumbers^n}}
\newcommand{\zn}{{\integernumbers^n}}
\newcommand{\nn}{{\naturalnumbers^n}}
\newcommand{\dsum}{\dotplus}
\newcommand{\injtn}{\myepsilon}
\newcommand{\projtn}{\pi}
\newcommand{\ctp}{{\hat \otimes}}
\newcommand{\citp}{{\hat \otimes}_\injtn}
\newcommand{\cistp}{{\hat \otimes}_{\injtn^s}}
\newcommand{\cptp}{{\hat \otimes}_\pi}
\newcommand{\idfunc}{\mathrm{id}}
\newcommand{\idfunconspace}[1]{{\idfunc_{#1}}}
\newcommand{\convweakstar}{{\mathop{\to}\limits^{\textrm{w-*}}}}
\newcommand{\cbfuncsymbol}{C_b}
\newcommand{\vanishingfuncsymbol}{C_0}
\newcommand{\cscfuncsymbol}{C_\mathrm{com}}
\newcommand{\ucfuncsymbol}{C_\mathrm{u}}
\newcommand{\contfuncsymbol}{C}
\newcommand{\cospacesymbol}{c_0}
\newcommand{\onto}{\twoheadrightarrow}

\newcommand{\opcptp}{{\otimes_\projtn}}
\newcommand{\subsetset}{\subset}

\newcommand{\norm}[1]{{\left\Vert{}#1{}\right\Vert}}
\newcommand{\nsznorm}[1]{{\Vert{}#1{}\Vert}}
\newcommand{\abs}[1]{{\left\vert{#1}\right\vert}}
\newcommand{\topdual}[1]{#1^*}

\newcommand{\norminfty}[1]{\norm{#1}_\infty}
\newcommand{\nsznorminfty}[1]{\nsznorm{#1}_\infty}
\newcommand{\normtwo}[1]{\norm{#1}_2}
\newcommand{\normone}[1]{\norm{#1}_1}
\newcommand{\floor}[1]{\lfloor{}{#1}\rfloor}
\newcommand{\ceil}[1]{\lceil{}{#1}\rceil}
\newcommand{\szceil}[1]{\left\lceil{}{#1}\right\rceil}

\newcommand{\closedball}[3]{{\overline{B}}_{#1}({#2}{};{}{#3})}
\newcommand{\card}[1]{{\#{#1}}}
\newcommand{\contlinop}[2]{\mathcal{L}({#1},{#2})}
\newcommand{\contlinopautom}[1]{\contlinop{#1}{#1}}

\newcommand{\linopautom}[1]{{\contlinop{#1}{#1}}}
\newcommand{\dualappl}[2]{\langle #1 , #2 \rangle}
\newcommand{\szdualappl}[2]{\left\langle #1 , #2 \right\rangle}
\newcommand{\cbfunc}[2]{\cbfuncsymbol({#1}, {#2})}
\newcommand{\cbfunccv}[1]{\cbfuncsymbol({#1})}
\newcommand{\vanishingfunc}[2]{{\vanishingfuncsymbol({#1},{#2})}}
\newcommand{\vanishingfunccv}[1]{{\vanishingfuncsymbol({#1})}}
\newcommand{\cscfunc}[2]{{\cscfuncsymbol({#1},{#2})}}
\newcommand{\cscfunccv}[1]{{\cscfuncsymbol({#1})}}
\newcommand{\ucfunc}[2]{{\ucfuncsymbol({#1},{#2})}}
\newcommand{\ucfunccv}[1]{{\ucfuncsymbol({#1})}}

\newcommand{\contfunccv}[1]{{\contfuncsymbol({#1})}}
\newcommand{\vf}[1]{\vanishingfunc{\realnumbers}{#1}}
\newcommand{\dvf}[1]{\topdual{\vf{#1}}}
\newcommand{\vfn}[2]{\vanishingfunc{\realnumbers^{#1}}{#2}}
\newcommand{\dvfn}[2]{\topdual{\vfn{#1}{#2}}}
\newcommand{\biglp}[3]{L^{#1}({#2},{#3})}
\newcommand{\biglpcv}[2]{L^{#1}({#2})}
\newcommand{\littlelp}[3]{l^{#1}({#2},{#3})}
\newcommand{\littlelpcv}[2]{l^{#1}({#2})}
\newcommand{\szlittlelpcv}[2]{{l^{#1}\left({#2}\right)}}

\newcommand{\gencospace}[2]{{\cospacesymbol({#1},{#2})}}
\newcommand{\genmodcont}[4]{{\omega^{#1}_{#2}({#3};{#4})}}
\newcommand{\modcont}[2]{{\omega({#1};{#2})}}
\newcommand{\szmodcont}[2]{{\omega\left({#1};{#2}\right)}}
\newcommand{\funcdifffunc}[2]{{\Delta^{#1}_{#2}}}
\newcommand{\funcdiff}[3]{{\funcdifffunc{#1}{#2}}{#3}}
\newcommand{\funcdiffonefunc}[1]{{\Delta_{#1}}}
\newcommand{\funcdiffone}[2]{{{\funcdiffonefunc{#1}}{#2}}}
\newcommand{\restrictfunc}[2]{{{#1}\vert_{#2}}}

\newcommand{\kinterp}[4]{{({#3},{#4})_{{#1},{#2}}}}
\newcommand{\szkinterp}[4]{{\left({#3},{#4}\right)_{{#1},{#2}}}}
\newcommand{\borelfunc}[2]{{\borel({#1},{#2})}}
\newcommand{\schwartzdistr}[1]{{S'({#1})}}

\newcommand{\besovspace}[4]{B^{#1}_{{#2},{#3}}({#4})}

\newcommand{\holderspace}[2]{{C^{#1}({#2})}}
\newcommand{\zygmundspace}[2]{{\mathcal{Z}^{#1}({#2})}}
\newcommand{\contdiffspace}[2]{{C^{#1}({#2})}}
\newcommand{\contbesovnorm}[6]%
{{\norm{#6}^{(\textrm{c})}_{\besovspace{#1}{#2}{#3}{#4};{#5}}}}
\newcommand{\contbesovnormrl}[6]%
{{\norm{#6}^{(\textrm{c'})}_{\besovspace{#1}{#2}{#3}{#4};{#5}}}}
\newcommand{\lpcvnorm}[3]{{\norm{#3}_{\biglpcv{#1}{#2}}}}
\newcommand{\lpcvnormext}[3]{{\norm{#3}^{(\extreals)}_{\biglpcv{#1}{#2}}}}
\newcommand{\contbesovtail}[5]%
{{\lpcvnorm{#3}{]0,1[}{t \in ]0, 1[ \mapsto%
t^{-{#1}-\frac{1}{#3}} \genmodcont{#4}{#2}{#5}{t}}}}
\newcommand{\contbesovtailinfty}[5]%
{{\lpcvnormext{#3}{]0,1[}{t \in ]0, 1[ \mapsto%
t^{-{#1}} \genmodcont{#4}{#2}{#5}{t}}}}
\newcommand{\mdimderiv}[2]{{D^{#1}{#2}}}
\newcommand{\unitball}[1]{{B_{#1}}}

\newcommand{\indexedtensorproduct}[2]{\bigotimes_{#1}^{#2}}

\newcommand{\indexedctp}[3]{\sideset{}{_{#1}}{\mathop{\hat{\bigotimes}}}\limits_{#2}^{#3}}
\newcommand{\indexedopctp}[3]{\sideset{}{_{#1}}\bigotimes_{#2}^{#3}}

\newcommand{\indexedcitp}[2]{\indexedctp{\injtn}{#1}{#2}}
\newcommand{\indexedcptp}[2]{\indexedctp{\projtn}{#1}{#2}}

\newcommand{\indexedinhopctp}[4]{\indexedopctp{_{({#1},{#2})}}{#3}{#4}}

\newcommand{\indexedinhctp}[3]{\indexedctp{_{({#1})}}{#2}{#3}}

\newcommand{\indexedopcptp}[2]{\indexedopctp{\projtn}{#1}{#2}}

\newcommand{\indexeddirectsum}[2]{\mathop{\dot \sum}\limits_{#1}^{#2}}
\newcommand{\directsumoneindex}[1]{{\mathop{\dot \sum}\limits_{#1}}}

\newcommand{\setimage}[2]{{{#1}[{#2}]}}
\newcommand{\equalns}{=_{\mathrm{n.s.}}}
\newcommand{\defequalns}{:=_{\mathrm{n.s.}}}
\newcommand{\equaltvs}{=_{\mathrm{tvs}}}

\newcommand{\defequalset}{:=_{\mathrm{set}}}
\newcommand{\iisom}{\cong_1}
\newcommand{\closedsubspace}{\subset_{\textrm{c.s.}}}
\newcommand{\normedsubspace}{\subset_{\textrm{n.s.}}}
\newcommand{\tvsubspace}{\subset_{\textrm{tvs}}}
\newcommand{\isomemb}{\subset_1}
\newcommand{\setoneton}[1]{{Z({#1})}}
\newcommand{\setzeroton}[1]{{Z_0({#1})}}
\newcommand{\setplusminusn}[1]{{Z_{\pm}({#1})}}
\newcommand{\setplusminusnpower}[2]{{\left(\setplusminusn{#1}\right)^{#2}}}
\newcommand{\ncover}[1]{N_{\mathrm{cover}}({#1})}
\newcommand{\szncover}[1]{N_{\mathrm{cover}}\left({#1}\right)}
\newcommand{\suppopst}{\suppop\nolimits_{\mathrm{set}}}
\newcommand{\sequenceof}{\subset}
\newcommand{\seqelem}[2]{{{#1}[{#2}]}}
\newcommand{\szseqelem}[2]{{{#1}\!\left[{#2}\right]}}
\newcommand{\seqcombop}{{s_{\mathrm{comb}}}}
\newcommand{\seqprojop}{{s_{\mathrm{proj}}}}
\newcommand{\seqset}[2]{{{#2}^{#1}}}
\newcommand{\seqcomb}[2]{{\seqcombop({#1},\;{#2})}}
\newcommand{\seqproj}[2]{{\seqprojop({#1},{#2})}}
\newcommand{\cartprodelem}[2]{{\seqelem{#1}{#2}}}
\newcommand{\szcartprodelem}[2]{{\szseqelem{#1}{#2}}}
\newcommand{\gennatbasvec}[2]{{\mathbf{e}^{#1}_{#2}}}
\newcommand{\nnatbasvec}[1]{{\mathbf{e}_{#1}}}

\newcommand{\znatbasvec}[1]{{{\check{\mathbf{e}}}_{#1}}}
\newcommand{\zmnatbasvec}[1]{{{\check{\mathbf{e}}}_{#1}}}
\newcommand{\onetonnbv}[2]{{\mathbf{e}^{[{#1}]}_{#2}}}
\newcommand{\ntensornbv}[1]{{\mathbf{e}^{\otimes}_{#1}}}
\newcommand{\ztensornbv}[1]{{{\check{\mathbf{e}}}^{\otimes}_{#1}}}
\newcommand{\zeroset}{\{0\}}
\newcommand{\zerooneset}{\{0, 1\}}
\newcommand{\zeroonesetn}{{\zerooneset^n}}
\newcommand{\zeroonesetnprime}{{\zerooneset^{n'}}}
\newcommand{\cpzerooneset}[1]{{\zerooneset^{#1}}}
\newcommand{\finitezeroseq}[1]{{\mathbf{0}_{#1}}}
\newcommand{\finiteoneseq}[1]{{\mathbf{1}_{#1}}}
\newcommand{\sqordfunction}{{\sigma_{\mathrm{sq}}}}
\newcommand{\sqordfirstfunction}{{\sigma_{\mathrm{sq1}}}}
\newcommand{\sqordsecondfunction}{{\sigma_{\mathrm{sq2}}}}
\newcommand{\sqord}[1]{{\sqordfunction({#1})}}
\newcommand{\sqordfirst}[1]{{\sqordfirstfunction({#1})}}
\newcommand{\sqordsecond}[1]{{\sqordsecondfunction({#1})}}
\newcommand{\gensqordfunction}[1]{{\sigma_{\mathrm{sq}}^{[{#1}]}}}
\newcommand{\gensqordcompfunction}[2]{{\sigma_{\mathrm{sq}}^{({#1},{#2})}}}
\newcommand{\gensqord}[2]{{\gensqordfunction{#1}({#2})}}
\newcommand{\gensqordcomp}[3]{{\gensqordcompfunction{#1}{#2}({#3})}}

\newcommand{\cubeorderingfunction}[1]{{\sigma_{\mathrm{c}}^{[{#1}]}}}
\newcommand{\cubeordering}[2]{{\cubeorderingfunction{#1}({#2})}}
\newcommand{\cubediffordfunction}[2]{{\beta_{{#1},{#2}}}}
\newcommand{\cubedifford}[3]{{\cubediffordfunction{#1}{#2}({#3})}}
\newcommand{\constrcubeordfunction}[2]{{\alpha_{{#1},{#2}}}}

\newcommand{\rconstrcubeordfunction}[2]{{\alpha'_{{#1},{#2}}}}

\newcommand{\holdercoeff}[3]{{H({#1};{#2},{#3})}}
\newcommand{\indexsetequal}[2]{{I^{[{#1}]}_{\Sigma=}({#2})}}
\newcommand{\norminspace}[2]{\left\Vert{}{#1}\vert{}{#2}{}\right\Vert}
\newcommand{\sznorminspace}[2]{\left\Vert{}{#1}\biggl\vert{}{#2}{}\right\Vert}
\newcommand{\nsznorminspace}[2]{\Vert{}{#1}\vert{}{#2}{}\Vert}
\newcommand{\supportradius}[1]{{r_{\mathrm{supp}}({#1})}}
\newcommand{\existsdef}{{\exists_{\mathrm{def}}}}
\newcommand{\cinterp}[1]{{c_{\mathrm{interp}}({#1})}}
\newcommand{\zerocentredcube}[2]{{A_{{#1},{#2}}}}
\newcommand{\zerocentredcubediff}[2]{{B_{{#1},{#2}}}}
\newcommand{\rectangle}[2]{{I_{\mathrm{rect}}({#1},{#2})}}
\newcommand{\omegaseq}{{\seqstyle{a}}}
\newcommand{\szomegaseqelem}[1]{{\szseqelem{\omegaseq}{#1}}}
\newcommand{\itrans}[2]{{I_{\mathrm{trans}}({#1},{#2})}}
\newcommand{\znisomfunc}[1]{{\xi^{[#1]}}}
\newcommand{\znisom}[2]{{\znisomfunc{#1}({#2})}}

\newcommand{\phidual}{\tilde{\varphi}}
\newcommand{\psidual}{\tilde{\psi}}

\newcommand{\gdual}{\tilde{g}}
\newcommand{\onedimmothersf}{\myphi}

\newcommand{\onedimmotherdualsf}{\tilde{\myphi}}
\newcommand{\onedimmotherdualwavelet}{\tilde{\psi}}

\newcommand{\onedimvanvspace}[1]{{V^{(0)}_{1,{#1}}}}
\newcommand{\onedimvanwspace}[1]{{W^{(0)}_{1,{#1}}}}
\newcommand{\onedimvandualvspace}[1]{{\tilde{V}^{(0)}_{1,{#1}}}}

\newcommand{\onedimvanprojop}[1]{{P^{(0)}_{1,{#1}}}}
\newcommand{\onedimvandeltaop}[1]{{Q^{(0)}_{1,{#1}}}}
\newcommand{\onedimdualsf}[2]{{\onedimmotherdualsf_{{#1},{#2}}}}
\newcommand{\onedimdualwavelet}[2]{{\onedimmotherdualwavelet_{{#1},{#2}}}}

\newcommand{\mdimvanvspace}[2]{{V^{(0)}_{{#1},{#2}}}}
\newcommand{\mdimvanwspace}[2]{{W^{(0)}_{{#1},{#2}}}}
\newcommand{\mdimtensorvanvspace}[2]{{V^{\otimes}_{{#1},{#2}}}}
\newcommand{\mdimvanpartialwspace}[3]{{W^{(0)}_{{#1},{#2},{#3}}}}

\newcommand{\mdimvanpartialdualwspace}[3]{{\tilde{W}^{(0)}_{{#1},{#2},{#3}}}}

\newcommand{\mdimdualvspace}[2]{{\tilde{V}_{{#1},{#2}}}}
\newcommand{\mdimdualwspace}[2]{{\tilde{W}_{{#1},{#2}}}}
\newcommand{\mdimpartialdualwspace}[3]{{\tilde{W}_{{#1},{#2},{#3}}}}
\newcommand{\onedimdualvspace}[1]{{\tilde{V}_{1,{#1}}}}
\newcommand{\onedimdualvjiisomfunc}[1]{{{\tilde{\iota}}_{1,{#1}}}}

\newcommand{\nspace}[3]{{\tilde{N}_{{#1},{#2},{#3}}}}

\newcommand{\mdimvandeltaop}[2]{{Q^{(0)}_{{#1},{#2}}}}
\newcommand{\mdimvanpartialprojop}[3]{{Q^{(0)}_{{#1},{#2},{#3}}}}
\newcommand{\mdimvanprojop}[2]{{P^{(0)}_{{#1},{#2}}}}

\newcommand{\mdimvanpartialdualprojop}[3]{{\tilde{Q}^{(0)}_{{#1},{#2},{#3}}}}

\newcommand{\mdimpartialdualprojop}[3]{{\tilde{Q}_{{#1},{#2},{#3}}}}

\newcommand{\mdimuvspace}[2]{{V^{(\mathrm{u})}_{{#1},{#2}}}}
\newcommand{\mdimupartialwspace}[3]{{W^{(\mathrm{u})}_{{#1},{#2},{#3}}}}
\newcommand{\mdimuwspace}[2]{{W^{(\mathrm{u})}_{{#1},{#2}}}}

\newcommand{\mdimuprojop}[2]{P^{(\mathrm{u})}_{{#1},{#2}}}
\newcommand{\mdimupartialdeltaop}[3]{Q^{(\mathrm{u})}_{{#1},{#2},{#3}}}
\newcommand{\mdimudeltaop}[2]{Q^{(\mathrm{u})}_{{#1},{#2}}}

\newcommand{\mdimdualprojop}[2]{{\tilde{P}_{{#1},{#2}}}}
\newcommand{\mdimdualpartialdeltaop}[3]{\tilde{Q}_{{#1},{#2},{#3}}}
\newcommand{\mdimdualdeltaop}[2]{\tilde{Q}_{{#1},{#2}}}

\newcommand{\mdimmotherwavelet}[2]{\psi^{[{#1}]}_{#2}}
\newcommand{\mdimgenwavelet}[4]{\psi^{[{#1}]}_{{#2},{#3},{#4}}}
\newcommand{\mdimmotherdualwavelet}[2]{\tilde{\psi}^{[{#1}]}_{#2}}
\newcommand{\mdimgendualwavelet}[4]{\tilde{\psi}^{[{#1}]}_{{#2},{#3},{#4}}}
\newcommand{\mdimaltpartialdualprojop}[3]{{\tilde{R}_{{#1},{#3},{#2}}}}

\newcommand{\mdimdualvjiisomfunc}[2]{{{\tilde{\iota}}_{{#1},{#2}}}}
\newcommand{\mdimdualvjiisom}[3]{{\mdimdualvjiisomfunc{#1}{#2}({#3})}}
\newcommand{\mdimgenwaveletfilterelem}[3]{{g^{[{#1}]}_{{#2},{#3}}}}
\newcommand{\mdimgendualwaveletfilterelem}[3]{{\tilde{g}^{[{#1}]}_{{#2},{#3}}}}

\newcommand{\largevspace}[2]{{V_{{#1},{#2}}}}
\newcommand{\largevprojop}[2]{{P_{{#1},{#2}}}}
\newcommand{\genmotherwaveletonedim}[1]{\zeta_{#1}}
\newcommand{\genwaveletonedim}[3]{\zeta_{{#1},{#2},{#3}}}
\newcommand{\genmotherdualwaveletonedim}[1]{\tilde{\zeta}_{#1}}
\newcommand{\gendualwaveletonedim}[3]{\tilde{\zeta}_{{#1},{#2},{#3}}}

\newcommand{\tensorauxdualvspace}[2]{{\tilde{M}_{{#1},{#2}}}}
\newcommand{\mdimmothersf}[1]{{\myphi^{[{#1}]}}}
\newcommand{\mdimmotherdualsf}[1]{\tilde{\myphi}^{[{#1}]}}
\newcommand{\mdimsf}[3]{\myphi^{[{#1}]}_{{#2},{#3}}}
\newcommand{\mdimdualsf}[3]{\phidual^{[{#1}]}_{{#2},{#3}}}
\newcommand{\genprojoponedim}[2]{{\mdimvanpartialprojop{1}{#1}{#2}}}

\newcommand{\uspace}[2]{{\mdimvanpartialwspace{1}{#1}{#2}}}

\newcommand{\zeroonesetnnozero}[1]{{J_{+}({#1})}}

\newcommand{\vanvjtopisomfunc}[2]{{\iota^{(0)}_{{#1},{#2}}}}
\newcommand{\vanvjtopisom}[3]{{\vanvjtopisomfunc{#1}{#2}({#3})}}
\newcommand{\vanwjtopisomfunc}[2]{{\eta^{(0)}_{{#1},{#2}}}}
\newcommand{\vanwjtopisom}[3]{{\vanwjtopisomfunc{#1}{#2}({#3})}}
\newcommand{\vanwsjtopisomfunc}[3]{{\eta^{(0)}_{{#1},{#2},{#3}}}}
\newcommand{\vanwsjtopisom}[4]{{\vanwsjtopisomfunc{#1}{#2}{#3}({#4})}}

\newcommand{\dualwsjtopisomfunc}[3]{{{\tilde{\eta}}_{{#1},{#2},{#3}}}}
\newcommand{\dualwsjtopisom}[4]{{\dualwsjtopisomfunc{#1}{#2}{#3}({#4})}}
\newcommand{\uvjtopisomfunc}[2]{{\iota^{(\mathrm{u})}_{{#1},{#2}}}}
\newcommand{\uvjtopisom}[3]{{\uvjtopisomfunc{#1}{#2}({#3})}}
\newcommand{\uwjtopisomfunc}[2]{{\eta^{(\mathrm{u})}_{{#1},{#2}}}}
\newcommand{\uwjtopisom}[3]{{\uwjtopisomfunc{#1}{#2}({#3})}}
\newcommand{\uwsjtopisomfunc}[3]{{\eta^{(\mathrm{u})}_{{#1},{#2},{#3}}}}
\newcommand{\uwsjtopisom}[4]{{\uwsjtopisomfunc{#1}{#2}{#3}({#4})}}
\newcommand{\normintensorvp}[4]{\norm{#4}_{\mdimuvspace{#1}{#2}({#3})}}
\newcommand{\orthvp}[2]{\bbigv_{#2}({#1})}
\newcommand{\orthprojp}[2]{\bbigp^{({#1})}_{#2}}
\newcommand{\orthdeltaprojp}[2]{\bbigq^{({#1})}_{#2}}
\newcommand{\norminorthvp}[3]{\norm{#3}_{\orthvp{#1}{#2}}}
\newcommand{\bmiset}[2]{{\Omega({#1},{#2})}}
\newcommand{\tensorvp}[2]{V_{#2}({#1})}
\newcommand{\hfilterelem}[1]{{h_{#1}}}
\newcommand{\gfilterelem}[1]{{g_{#1}}}
\newcommand{\htfilterelem}[1]{{\tilde{h}_{#1}}}
\newcommand{\gtfilterelem}[1]{{\tilde{g}_{#1}}}
\newcommand{\onedimgenwaveletfilterelem}[2]{{g_{{#1},{#2}}}}
\newcommand{\onedimgendualwaveletfilterelem}[2]{{\tilde{g}_{{#1},{#2}}}}
\newcommand{\scfsum}[4]{{v_{{#1},{#2}}({#3},{#4})}}

\newcommand{\genvspace}[3]{{V^{#1}_{{#2},{#3}}}}
\newcommand{\genwspace}[3]{{W^{#1}_{{#2},{#3}}}}
\newcommand{\genpartialwspace}[4]{{W^{#1}_{{#2},{#3},{#4}}}}
\newcommand{\besovinterpcoeffnorm}[6]%
{{\norm{#6}^{(\textrm{w})}_{\besovspace{#3}{#4}{#5}{\realnumbers^{#1}};{#2}}}}
\newcommand{\besovorthwaveletnorm}[6]%
{{\norm{#6}^{(\textrm{o})}_{\besovspace{#3}{#4}{#5}{\realnumbers^{#1}};{#2}}}}
\newcommand{\besovinterpwaveletnorm}[6]%
{{\norm{#6}^{(\textrm{i})}_{\besovspace{#3}{#4}{#5}{\realnumbers^{#1}};{#2}}}}
\newcommand{\mdimuvpspace}[3]{{V^{(u)}_{{#1},{#3}}({#2})}}
\newcommand{\orthcoefffunctionalvalue}[3]%
{{\balpha_{{#1},{#2}}({#3})}}

\newcommand{\toperator}[3]{{T_{{#1},{#2},{#3}}}}
\newcommand{\dualtoperator}[3]{{\tilde{T}_{{#1},{#2},{#3}}}}
\newcommand{\aspace}[3]{{\tilde{A}_{{#1},{#2},{#3}}}}

\newcommand{\definek}{Let \(K = \realnumbers\) or \(K = \complexnumbers\). }
\newcommand{\definekn}
{Let \(K = \realnumbers\) or \(K = \complexnumbers\) and \(n \in \positiveintegers\). }
\newcommand{\givelemmawithoutproof}{}

\newcommand{\givelemmaswithoutproofsn}[1]{}
\newcommand{\mydef}{Definition }
\newcommand{\mydefs}{Definitions }
\newcommand{\mytheorem}{Theorem }
\newcommand{\mytheorems}{Theorems }
\newcommand{\mylemma}{Lemma }
\newcommand{\mylemmas}{Lemmas }
\newcommand{\mycorollary}{Corollary }
\newcommand{\myequation}{Equation }
\newcommand{\myequations}{Equations }

\if\elsevier1
\fi

\if\elsevier0
{
\newtheorem{theorem}{Theorem}[section]
\newtheorem{corollary}[theorem]{Corollary}
\newtheorem{lemma}[theorem]{Lemma}

\theoremstyle{remark}
\newtheorem{remark}[theorem]{Remark}

\theoremstyle{definition}
\newtheorem{definition}[theorem]{Definition}
}
\else
{
\newtheorem{theorem}{Theorem}
\newtheorem{lemma}[theorem]{Lemma}
\newdefinition{definition}{Definition}
\newproof{proof}{Proof}
}
\fi

\newenvironment{myenumerate}{\par\noindent}{\par}
\newcommand{\myitem}[1]{\par\noindent{}#1\;\;\;\;}

\begin{document}

\if\elsevier0
{


\if\longversion0
{
\title{Multiresolution Analysis for Compactly Supported 
  Interpolating Tensor Product Wavelets}
}
\else
{
\title{Multiresolution Analysis for Compactly Supported 
  Interpolating Tensor Product Wavelets (long version)}
}
\fi

\author{Tommi H\"oyn\"al\"anmaa \\
Physics, Tampere University, \\
P. O. Box 692, FI-33101 Tampere, Finland}

\if\arxiv1
\date{}
\else
\date{\today}
\fi

\maketitle

\begin{abstract}
We construct multidimensional interpolating tensor product MRA's
\mraprepspace the function spaces \(\vanishingfunc{\rn}{K}\)
, \(K = \realnumbers\) or
\(K = \complexnumbers\),
consisting of real or complex valued
functions on \(\rn\) vanishing at infinity
and the function spaces \(\ucfunc{\rn}{K}\) consisting
of bounded and uniformly continuous functions on \(\rn\).
We also construct an interpolating dual MRA for both of these spaces.
The theory of
the tensor products of Banach spaces is used. We generalize the
Besov space norm equivalence result from
Donoho (1992, Interpolating Wavelet Transforms) to our
\(n\)-dimensional construction.
\end{abstract}




}
\else
{
\begin{frontmatter}

	\title{Multiresolution Analysis for Compactly Supported %
  	Interpolating Tensor Product Wavelets}
	\author{Tommi H\"oyn\"al\"anmaa}
	\address{ %
	Institute of Physics, Tampere University of Technology, %
	P. O. Box 692, FI-33101 Tampere, Finland}
	\ead{tommi.hoynalanmaa@tut.fi}
	\begin{abstract}
	We construct multidimensional interpolating tensor product MRA's
	\mraprepspace the function spaces \(\vanishingfunc{\rn}{K}\),
	\(K = \realnumbers\) or \(K = \complexnumbers\),
	consisting of real or complex valued
	functions on \(\rn\) vanishing at infinity
	and the function spaces \(\ucfunc{\rn}{K}\) consisting
	of bounded and uniformly continuous functions on \(\rn\).
        We also construct an interpolating dual MRA for both of these spaces.
	The theory of
	the tensor products of Banach spaces is used. We generalize the
	Besov space norm equivalence result from
	Donoho (1992, Interpolating Wavelet Transforms)
	to our \(n\)-dimensional construction.
	\end{abstract}

\end{frontmatter}
}
\fi

\if\shortprep0
{
\tableofcontents
}
\fi

\section{Introduction}

Chui and Li \cite{cl1996} have constructed a one-dimensional MRA \mraprepspace 
function space
\(\ucfunccv{\realnumbers}\) (bounded and uniformly continuous complex valued
functions on \(\realnumbers\)) for
interpolating wavelets. Donoho \cite{donoho1992}
has derived convergence results for interpolating wavelets on space
\(\vanishingfunccv{\realnumbers}\).
Goedecker's book
\cite{goedecker1998} contains also
some material about multidimensional interpolating wavelets.
Goedecker \cite{goedecker1998} uses the term interpolating wavelets
to mean Deslauriers-Dubuc wavelets whereas some
other authors such as Chui and Li \cite{cl1996} use the term to mean roughly 
wavelet
families whose mother scaling function has the cardinal interpolation
property \(\myphi(k) = \delta_{k,0}\) for all \(k \in
\integernumbers\) (one-dimensional case). We follow the latter
convention in this article.  Chui and Li \cite{cl1994} have also constructed a
general framework for multivariate wavelets. However,
the theory in that article uses function space \(\biglpcv{2}{\rn}\)
and so it is unsuitable for the approach in this article.
Dubuc and Deslauriers \cite{dubuc1986,dd1989}
have investigated interpolation
processes related to the Deslauriers-Dubuc fundamental functions.
Han and Jia \cite{hj1998} have discussed fundamental functions
(see the article for a definition)
satisfying the cardinal interpolation property.
Numerical values for the wavelet filters of the Deslauriers-Dubuc
wavelets have been given in \cite{goedecker1998}
and \cite{ks2000}.

Theory for orthonormal wavelets has been developed e.g. by
Daubechies \cite{daubechies1988,daubechies1992}, Meyer \cite{meyer1992},
and Wojtaszczyk \cite{wojtaszczyk1997}.
Kova{\v c}evi{\'c} and Sweldens \cite{ks2000} 
have investigated the use of digital filters for multidimensional biorthogonal
wavelets.
Lewis \cite{lewis1994} has constructed a MRA \mraprepspace function space 
\(\biglpcv{2}{\realnumbers}\) for interpolating wavelets in one dimension.
Lawton \etal\ \cite{lls1997} have presented conditions for the refinement
mask of an orthonormal MRA.
Ji and Shen \cite{js1999} have given a condition for the refinement
mask to be interpolatory and a condition for a dual refinement mask.
Krivoshein and Skopina \cite{ks2011} have presented convergence
results for frame-like wavelets that are not frames.
Karakaz'yan \etal\ \cite{kst2008} have described symmetric
interpolatory masks generating dual compactly supported wavelet
systems and they have also given formulas for dual refinement masks.
Shui \etal\ \cite{sbz2001} have shown how to construct \(M\)-band
wavelets with all
the following properties: compact support, orthogonality,
linear-phase, regularity, and interpolation. This kind of scaling
functions exist only when \(M \geq 4\). Dahlke \etal\ \cite{dkm1999}
have developed a new method to construct higher dimensional scaling
functions. Their method is based on solutions to specific Lagrange
interpolation problems by polynomials.
Han \cite{han2002} has investigated symmetries of refinable functions.
\if11
{
DeVore \etal\ \cite{djl1992} have developed some
theory for
multidimensional orthonormal MRA \mraprepspace space
\(\biglpcv{p}{\realnumbers^d}\),
\(0 < p < \infty\), \(d \in \positiveintegers\).
}
\else
{
DeVore \etal\ \cite{djl1992} have developed some
theory for
multidimensional orthonormal MRA \mraprepspace space
\(\biglpcv{p}{\realnumbers^d}\),
\(0 < p < \infty\), \(d \in \positiveintegers\),
and \(\ucfunccv{\realnumbers^d}\),
\(d \in \positiveintegers\). }
\fi
He and Lai
\cite{hl1998} have
constructed nonseparable box spline wavelets on Sobolev spaces
\(H^s(\realnumbers^2)\).
Wavelets have also been discussed in \cite{hkpt}.

The theory about tensor products of Banach spaces has been presented e.g. in
\cite{ryan2002}. We use the notation from \cite{ryan2002} for Banach space
tensor products in this article.
Schaefer's book \cite{schaefer1966} contains some
material about tensor products of topological vector spaces.
Book \cite{lc1985} contains also an introduction to tensor products of
Banach spaces.
Doma\'nski \etal\ \cite{dl1997,dls1997} have done some research
on them.
Michor \cite{michor1978} has represented tensor products
of Banach spaces using category theory.
Grothendieck \cite{grothendieck1955}
gives a classical
presentation for tensor products of locally convex spaces.
Theory for tensor products of Banach spaces can also be found in book
\cite{df1993}.
Kustermans and Vaes \cite{kv2000} have developed some theory for space
\(C_0(G)\) where \(G\) is a compact or a locally compact group and
for the tensor products of \(C_0(G)\).
Daws \cite{daws2006} has presented some material on tensor products
of Banach algebras.
Gl\"ockner \cite{glockner2004} has shown that
the tensor products of topological vector spaces are not associative.
Reinov \cite{reinov2002} has investigated Banach spaces without the
approximation
property. Brodzki and Niblo \cite{bn2004} have done some research on
the rapid decay property of discrete groups
and the metric approximation property.

An introduction to Besov spaces can be found e.g. in
\cite{brenner1975}.
An introduction to Besov spaces using the
Fourier transform based definition of these spaces can be found in
\cite{triebel1983}.
DeVore and Popov \cite{dp1988} have investigated the connection of Besov spaces 
with the dyadic spline approximation
and interpolation of Besov spaces.
Kyriazis and Petrushev \cite{kp2002} have given 
a method for the construction of unconditional bases for
Triebel-Lizorkin and Besov spaces.
The relationship between orthonormal wavelets and Besov
spaces has been discussed by Meyer \cite{meyer1992}.
Almeida \cite{almeida2005} has investigated wavelet
bases in (generalized) Besov spaces.

Section \ref{sec:preliminaries} introduces some definitions
used in the rest of this article.
\if\shortprep0
{
Some results on function
spaces and sequences that are needed in the construction of the MRA's
are given in section \ref{sec:function-spaces-and-sequences}. Section
\ref{sec:tensor-products} contains results on Banach space tensor
products that are needed in the construction of the MRA's.
}
\fi
We give some general definitions needed by
the MRA's in section \ref{sec:gen-def-mra}.
A multivariate MRA \mraprepspace \(\ucfunc{\rn}{K}\),
\(K = \realnumbers\) or \(K = \complexnumbers\),
is constructed in section \ref{sec:uc-mra}.
A multivariate MRA \mraprepspace \(\vanishingfunc{\rn}{K}\),
\(K = \realnumbers\) or \(K = \complexnumbers\),
is constructed in section \ref{sec:van-mra}.
The interpolating dual MRA is presented in section \ref{sec:dual-mra}.
The Besov space norm equivalence from
Donoho \cite{donoho1992} is generalized
for the \(n\)-dimensional interpolating MRA's
in section \ref{sec:besov-norm-equivalence}.
\if\shortprep1
{
  A longer version of this article with more detailed proofs can be
  found at \cite{hoynalanmaa2010}.
}
\fi

\section{Preliminaries}
\label{sec:preliminaries}

\subsection{General}
\label{sec:preliminaries-general}

\if\shortprep0
{
When
\begin{math}
  n \in \naturalnumbers
\end{math},
\begin{math}
  n \geq 2
\end{math},
and
\begin{math}
  P_1, \ldots, P_n
\end{math}
are propositions we define
\begin{displaymath}
  P_1 \implies P_2 \implies \ldots \implies P_n  
\end{displaymath}
to mean
\begin{displaymath}
  (P_1 \implies P_2) \land (P_2 \implies P_3) \land \ldots \land
  (P_{n-1} \implies P_n) .
\end{displaymath}
}
\fi
\if\shortprep1
{
The set of all positive real numbers is denoted by
\begin{math}
  \positiverealnumbers
\end{math}
and the extended real line by
\begin{math}
  \extreals
\end{math}.
The set of all positive integer numbers is denoted by
\begin{math}
  \positiveintegers
\end{math}.
We define
\begin{math}
  \nonnegrealnumbers := \{ x \in \realnumbers \setsep x \geq 0 \}
\end{math}
and
\begin{math}
  \extnonnegreals := \{ x \in \extreals \setsep x \geq 0 \}
\end{math}.
}
\else
{
The set of all positive real numbers is denoted by
\begin{math}
  \positiverealnumbers
\end{math}
and the extended real line by
\begin{math}
  \extreals
\end{math}.
We define
\begin{math}
  \nonnegrealnumbers := \{ x \in \realnumbers \setsep x \geq 0 \}
\end{math},
\begin{math}
  \extnonnegreals := \{ x \in \extreals \setsep x \geq 0 \}
\end{math}, and
\begin{math}
  \positiveintegers
\end{math}
to be the set of positive integer numbers.
}
\fi
\if\shortprep0
{
If \(A\) is a set and \(P(x)\), \(x \in A\), is a proposition
then we define \(\existsdef x \in A : P(x)\) to mean that
\(\exists y \in A : P(y)\) and \(x\) is defined to
be some element of \(A\) for which \(P(x)\) is true.
}
\fi
When \(A\) and \(B\) are arbitrary sets, \(f\) is a function from
\(A\) into \(B\), and \(X \subset A\) the image of \(X\) under
\(f\) is denoted by \(\setimage{f}{X}\).
The set-theoretic support of a function
\(f : X \to \complexnumbers\) where \(X\) is
a set is denoted by \(\suppopst f\).
The topological support of a function
\(f : T \to \complexnumbers\) where \(T\) is a
topological space is denoted by \(\suppop f\).
When \(E\) is a metric space, \(x \in E\), and
\(r \in \positiverealnumbers\) the closed ball of radius \(r\) centred
at \(x\) is denoted by \(\closedball{E}{x}{r}\).
\if10
When \(A\) and \(B\) are some algebraic structures we may write
\(A \subsetset B\) to mean that \(B\) contains \(A\) as a subset.
\fi
\if\shortprep0
{
When \(A\) is a set we denote the cardinality of \(A\) by
\(\card{A}\).
}
\fi

\if\shortprep0
{
\begin{definition}
  \label{def:idfunc}
  When \(A\) is an arbitrary set define
  \(\idfunconspace{A} : A \to A\) to be
  the identity function on \(A\), i.e.
  \(\idfunconspace{A}(x) = x\) for all
  \(x \in A\).
\end{definition}
}
\fi

\begin{definition}
  \label{def:integerranges}
  When
  \begin{math}
    n \in \positiveintegers
  \end{math}
  define
  \begin{displaymath}
    \setoneton{n} := \{ k \in \positiveintegers \setsep k \leq n \} 
    .
  \end{displaymath}
  \if\shortprep0
  {
  When
  \begin{math}
    n \in \naturalnumbers
  \end{math}
  define
  \begin{displaymath}
    \setzeroton{n} := \{ k \in \naturalnumbers \setsep k \leq n \}
  \end{displaymath}
  and
  \begin{displaymath}
    \setplusminusn{n} := \{ k \in \integernumbers \setsep \abs{k} \leq n \} .
  \end{displaymath}
  }
  \fi
\end{definition}

\if\shortprep0
{
\begin{definition}
  When \(A\) is a set in which a commutative binary operation \(+\) 
  is defined,
  \(B \subset A\), and \(a \in A\)
  define
  \begin{displaymath}
  a + B := B + a := \{ a + x \setsep x \in B\} . 
  \end{displaymath}
\end{definition}
}
\fi

We define the differences \(\funcdifffunc{m}{\rnh}\) and
\(\funcdiffonefunc{\rnh}\)
as in \cite{triebel1983} and \cite[def. V.4.1]{bs1988}.

\begin{definition}
  Define
  \begin{displaymath}
    \left(\funcdiff{1}{\rnh}{f}\right)\left(\rnx\right)
    := (\funcdiffone{\rnh}{f})(\rnx)
    := f(\rnx + \rnh) - f(\rnx)
  \end{displaymath}
  for all
  \begin{math}
    \rnx \in \rn
  \end{math},
  \begin{math}
    \rnh \in \rn
  \end{math},
  and
  \begin{math}
    f \in \complexnumbers^{\rn}
  \end{math}
  and
  \begin{displaymath}
    \funcdiff{m}{\rnh}{f} := 
    \funcdiffone{\rnh}{(\funcdiff{m-1}{\rnh}{f})}
  \end{displaymath}
  for all
  \begin{math}
    m \in \naturalnumbers + 2
  \end{math},
  \begin{math}
    \rnh \in \rn
  \end{math},
  and
  \begin{math}
    f \in \complexnumbers^{\rn}
  \end{math}.
\end{definition}

\if\shortprep1
{
\begin{definition}
  \label{def:ncover}
  \definekn
  Let \(f\) be a function from \(\realnumbers^n\)
  into \(K\).
  Define
  \begin{displaymath}
    \ncover{f} := \max_{\rnx \in \realnumbers^n} \card{\left\{
    \firstznvar \in \integernumbers^n \setsep
    f(\rnx - \firstznvar) \neq 0 \right\}} .
  \end{displaymath}
\end{definition}
}
\fi

\if\shortprep1
{
\begin{definition}
  Let \(E\) and \(F\) be normed vector spaces.
  When \(f\) is a compactly supported function from \(E\)
  into \(F\) define
  \[
  \supportradius{f} :=
  \inf \{ r \in \nonnegrealnumbers \setsep \suppopst f \subset
  \closedball{E}{0}{r} \} .
  \]
\end{definition}
}
\fi

\if\shortprep1
{
\begin{definition}
  Let \(n \in \positiveintegers\).
  Define
  \begin{displaymath}
    \itrans{f}{\rnx} :=
    \{ \firstznvar \in \zn \setsep
    f( \rnx - \firstznvar ) \neq 0
    \}
  \end{displaymath}
  for all \(f \in \complexnumbers^{\realnumbers^n}\)
  and \(\rnx \in \rn\).
\end{definition}
}
\fi

\subsection{Sequences and Cartesian Products}
\label{sec:preliminaries-seq}

\if\shortprep1
{
\begin{definition}
  When \(n \in \positiveintegers\) define
  \begin{displaymath}
    \finitezeroseq{n} := (0)_{k \in \setoneton{n}} .
  \end{displaymath}
\end{definition}
}
\else
{
\begin{definition}
  When \(n \in \positiveintegers\) define
  \begin{displaymath}
    \finitezeroseq{n} := (0)_{k \in \setoneton{n}}
  \end{displaymath}
  and
  \begin{displaymath}
    \finiteoneseq{n} := (1)_{k \in \setoneton{n}} .
  \end{displaymath}
\end{definition}
}
\fi

\begin{definition}
  When \(n \in \positiveintegers\) define
  \(\zeroonesetnnozero{n} := \zeroonesetn \setminus \{ 
  \finitezeroseq{n} \}\).
\end{definition}

\if\shortprep0
{
\begin{definition}
  When \(I\) is a countable nonempty set define
  \begin{displaymath}
    \gennatbasvec{I}{k} := ( \delta_{j,k} )_{j \in I}
  \end{displaymath}
  for all \(k \in I\).
\end{definition}
}
\fi

\if\shortprep0
{
\begin{definition}

  Define
  \begin{align*}
    \nnatbasvec{k} & := \gennatbasvec{\naturalnumbers}{k}
    \spaceafter \textrm{for all} \; k \in \naturalnumbers \\
    \znatbasvec{k} & := \gennatbasvec{\integernumbers}{k}    
    \spaceafter \textrm{for all} \; k \in \integernumbers \\
    \zmnatbasvec{\firstznvar} & := \gennatbasvec{\zn}{\firstznvar}
    \spaceafter \textrm{for all} \; \firstznvar \in \zn .
  \end{align*}
  When \(n \in \positiveintegers\) define
  \begin{displaymath}
    \onetonnbv{n}{k} := \gennatbasvec{\setoneton{n}}{k}
    \spaceafter \textrm{for all} \; k \in \setoneton{n} .
  \end{displaymath}
\end{definition}
}
\fi

\if\shortprep1
{
  \begin{definition}
    Define
    \begin{math}
      \znatbasvec{k} := (\delta_{j,k})_{j \in \integernumbers}
    \end{math}
    for all \(k \in \integernumbers\).
    When \(n \in \positiveintegers\) define
    \begin{math}
      \onetonnbv{n}{k} := \gennatbasvec{\setoneton{n}}{k}
    \end{math}
    for all \(k \in \setoneton{n}\).
  \end{definition}
}
\fi

\if\shortprep0
{
\begin{definition}
  Let \(n, m \in \positiveintegers\), \(n \geq m\).
  Define
  \begin{displaymath}
    \seqproj{m}{\vectorstyle{s}}
    :=
    ( \seqelem{\vectorstyle{s}}{1}, \ldots, 
    \seqelem{\vectorstyle{s}}{m} )
    \in \complexnumbers^m
  \end{displaymath}
  where \(\vectorstyle{s} \in \complexnumbers^n\).
\end{definition}
}
\fi

\if\shortprep0
{
\begin{definition}
  Let \(n, m \in \positiveintegers\).
  Define
  \begin{displaymath}
    \seqcomb{\zns}{\znt} := ( \cartprodelem{\zns}{1}, \ldots,
      \cartprodelem{\zns}{n}, \cartprodelem{\znt}{1}, \ldots,
      \cartprodelem{\znt}{m} )
    \in \complexnumbers^{n+m}
  \end{displaymath}
  where \(\zns \in \complexnumbers^n\)
  and \(\znt \in \complexnumbers^m\).
\end{definition}
}
\fi

\subsection{Vector Spaces}
\label{sec:preliminaries-vectorspaces}

\if\shortprep0
{
Suppose that \(A\) and \(B\) are topological vector spaces.
We define
\begin{math}
  A \equaltvs B
\end{math}
to be true iff \(A\) and \(B\) are the same topological vector space, i.e.
\(A\) and \(B\) are the same vector space and have the same topology.
\(A\) and \(B\) are called \defterm{algebraically isomorphic} iff 
\(A\) and \(B\) are isomorphic as vector spaces. A function
\(\iota : A \onto B\) is called an
\defterm{algebraic isomorphism} iff \(\iota\) is an isomorphism from
vector space \(A\) onto vector space \(B\) (i.e. \(\iota\) is a linear
bijection from \(A\) onto \(B\)). A function \(\eta : A \onto B\)
is called a \defterm{topological isomorphism} iff \(\eta\) is an 
algebraic isomorphism from \(A\) onto \(B\) and a homeomorphism 
from topological
space \(A\) onto topological space \(B\) (i.e. \(\iota\) is an
algebraic isomorphism that preserves the topology). \(A\) and \(B\)
are called \defterm{topologically isomorphic} iff there exists a topological
isomorphism from \(A\) onto \(B\).
We define
\begin{math}
  \contlinop{A}{B}
\end{math}
to be the set of all continuous linear functions from \(A\) into
\(B\).

When \(B\) is a topological vector space we define
\begin{math}
  A \tvsubspace B
\end{math}
to be true iff \(A\) is a topological vector subspace of \(B\).
}
\fi

\if\shortprep1
{
\if10
{
Suppose that \(E\) and \(F\) are normed vector spaces. We define
\(E \equalns F\) to be true iff \(E\) and \(F\) are the same normed
vector space and \(E \iisom F\) to be true iff \(E\) and \(F\)
are isometrically isomorphic.
}
\fi
Suppose that \(A\) and \(B\) are topological vector spaces.
We define
\begin{math}
  \contlinop{A}{B}
\end{math}
to be the set of all continuous linear functions from \(A\) into
\(B\).
}
\else
{
\if\longprep1
{
Suppose that \(E\) and \(F\) are normed vector spaces. We define
\(E \equalns F\) to be true iff \(E\) and \(F\) are the same normed
vector space, i.e. \(E\) and \(F\) are the same vector space and they have
the same norm. We define \(E \iisom F\) to be true iff \(E\) and \(F\)
are isometrically isomorphic.
}
\else
{
Suppose that \(E\) and \(F\) are normed vector spaces. We define
\(E \equalns F\) to be true iff \(E\) and \(F\) are the same normed
vector space, i.e. \(E\) and \(F\) are the same vector space (they
consist of the same elements, have the same scalar field, and the
addition and scalar multiplication operations are the same) and have
the same norm. We define \(E \iisom F\) to be true iff \(E\) and \(F\)
are isometrically isomorphic.
}
\fi
}
\fi

The term \defterm{operator} shall mean a bounded linear
function from a normed
vector space into another one. The term \defterm{projection}
shall mean a linear
projection \projprep a vector space onto a vector space.
The topological dual of a Banach space \(A\) is denoted
by \(\topdual{A}\). 
Unless otherwise stated \(\topdual{A}\) is equipped with the norm
topology.
\if\shortprep0
{
The unit ball of a normed
vector space \(A\) is denoted by \(\unitball{A}\).
When \(A\) is a Banach space and \((\indexeddualelem{x}{k})_{k=0}^\infty\)
is a sequence in \(\topdual{A}\) we use the notation
\(\indexeddualelem{x}{k} \convweakstar \dualelem{x}\) to
mean that the sequence \((\indexeddualelem{x}{k})_{k=0}^\infty\) converges to
\(\dualelem{x}\) in the weak-* topology of
\(\topdual{A}\). Convergence in the norm topology
is denoted by \(\indexeddualelem{x}{k} \rightarrow \dualelem{x}\).
}
\fi

\if\shortprep0
{
When \(E\) is a Banach space and \(A\) and \(B\) open or closed
subspaces of \(E\) the topological direct sum of \(A\) and
\(B\) is denoted by \(A \dsum B\), i.e. \(A \dsum B\) is the algebraic
direct sum of vector spaces \(A\) and \(B\) and both \(A\) and \(B\)
are topologically complemented in the normed vector space
\(A \dsum B\).
}
\fi

When \(A\) and \(B\) are Banach spaces
we use the notation
\begin{math}
  A \isomemb B
\end{math}
to mean that \(A\) is isometrically embedded in \(B\), i.e. \(A \subset B\)
and the inclusion map is distance preserving.
\if\shortprep1
{
When \(A\) and \(B\) are topological vector spaces we define
\begin{math}
  A \equaltvs B
\end{math}
to be true iff \(A\) and \(B\) are equal topological vector spaces
}
\fi
When \(B\) is a Banach space
we use the notation
\begin{math}
  A \closedsubspace B
\end{math}
to mean that \(A\) is a closed subspace of \(B\).
When \(E\) is a Banach space, \(A \subsetset E\), and
\(B \subsetset \topdual{E}\) we define \(B \bot A\) to be true
iff
\begin{math}
  \dualappl{\dualelem{b}}{a} = 0
\end{math}
for all
\begin{math}
  a \in A
\end{math}
and
\begin{math}
\dualelem{b} \in B
\end{math}.
\if\shortprep1
\else
\begin{displaymath}
\forall a \in A, \dualelem{b} \in B :
\dualappl{\dualelem{b}}{a} = 0 .
\end{displaymath}
\fi

When \(E\) is a normed vector space and \(x \in E\) we may use the
notation \(\norminspace{x}{E}\) to mean the norm of \(x\) in \(E\).
When \(B\) is a normed vector space and we write
\begin{math}
  A \defequalns \{ x \in B \setsep P(x) \}
\end{math}
we assume that
\begin{math}
  \norminspace{x}{A} := \norminspace{x}{B}
\end{math}
for each
\begin{math}
  x \in A
\end{math}.
\if\shortprep0
{
When \(A\) is a vector space we may define a norm \(\norm{\cdot}_A\) on \(A\)
so that \(\norm{\cdot}_A\) is defined on a larger set \(E\) containing
\(A\) as a subset.
}
\fi

\if\shortprep1
{
When \(E\) is a Banach space and \(A\) and \(B\) open or closed
subspaces of \(E\) the topological direct sum of \(A\) and
\(B\) is denoted by \(A \dsum B\).
}
\fi

\begin{definition}
  Let \(E\) be a set and \(A\) be a vector space so that
  \(A \subsetset E\).
  Suppose that
  \begin{math}
    \norm{\cdot}_A : E \to \extnonnegreals 
  \end{math}
  is a norm on \(A\).
  We say that the norm
  \begin{math}
    \norm{\cdot}_A 
  \end{math}
  \defterm{characterizes \(A\) on \(E\)} iff
  \(\norm{x}_A < \infty\) for all \(x \in A\) and
  \(\norm{y}_A = \infty\) for all \(y \in E \setminus A\).
\end{definition}

\subsection{Function Spaces and Sequence Spaces}
\label{sec:preliminaries-functionspaces}


When \(p \in [1, \infty]\), \(B\) is a Banach
space, and \(I\) a denumerable set
we denote the \(l^p\) space consisting of a subset of \(B^I\) by
\begin{math}
  \littlelp{p}{I}{B}
\end{math}.
As usual,
\begin{math}
  \littlelpcv{p}{I} \defequalns \littlelp{p}{I}{\complexnumbers}
\end{math}.
When \(K = \realnumbers\) or \(K = \complexnumbers\),
\(n \in \positiveintegers\), and \(p \in [1, \infty]\)
we denote the \(L^p\) space consisting of a subset of the
Borel functions from \(\rn\) into \(K\)
by \(\biglp{p}{\rn}{K}\).
As usual,
\begin{math}
  \biglpcv{p}{\rn} \defequalns \biglp{p}{\rn}{\complexnumbers}
\end{math}.
When \(X\) and \(Y\) are Borel spaces
define
\begin{math}
  \borelfunc{X}{Y}
\end{math}
to be the set of all Borel functions from
\(X\) into \(Y\).

\if\shortprep1
{
Let 
\(K = \realnumbers\) or \(K = \complexnumbers\).
When \(T\) is a topological space we denote the space of \(K\)-valued
bounded and continuous functions on \(T\) by \(\cbfunc{T}{K}\).
When \(E\) is a metric space the space of \(K\)-valued bounded and
uniformly continuous functions on \(E\) is denoted by
\(\ucfunc{E}{K}\).
When \(T\) is a locally compact Hausdorff space the space of
\(K\)-valued bounded and continuous functions vanishing at infinity
\cite{lacey1974,donoho1992} is denoted by \(\vanishingfunc{T}{K}\).
The space of \(K\)-valued compactly supported and continuous functions
on \(T\) is denoted by \(\cscfunc{T}{K}\).
}
\fi

\if\shortprep0
{
\begin{definition}
Let \(T\) be a topological space and
\(K = \realnumbers\) or \(K = \complexnumbers\). Define
\begin{displaymath}
  \cbfunc{T}{K} := \left\{ f \setsep f \; \textrm{is a continuous and
      bounded function from}\; T \;\textrm{to}\; K \right\}
\end{displaymath}
where the norm of an element \(f\) is given by
\begin{displaymath}
  \norm{f} = \sup_{x \in T} \abs{f(x)}
\end{displaymath}
and
\begin{math}
  \cbfunccv{T} \defequalns \cbfunc{T}{\complexnumbers}
\end{math}.
\end{definition}

\if\longversion1
{
\(\cbfunc{T}{K}\) is a Banach space.
}
\fi

\if\longversion1
When \(f\) is a function from \(\realnumbers^n\) into
\(\complexnumbers\) and \(a \in \complexnumbers\)
\begin{displaymath}
  \lim_{\norm{x} \to \infty} f(x) = a
\end{displaymath}
is equivalent to
\begin{displaymath}
  \forall \myepsilon \in \positiverealnumbers :
  \exists r \in \positiverealnumbers :
  \forall \rnx \in \rn \setminus \closedball{\rn}{0}{r} : 
  \abs{f(\rnx) - a} < \myepsilon .
\end{displaymath}
\fi

\begin{definition}
  Let \(E\) be a metric space and
  \(K = \realnumbers\) or \(K = \complexnumbers\).
  Define
  \begin{displaymath}
    \ucfunc{E}{K}
    \defequalns
    \left\{
      f \in \cbfunc{E}{K}
      \setsep
      f \;\textrm{is uniformly continuous on}\; E
    \right\}
  \end{displaymath}
  and
  \begin{math}
    \ucfunccv{E} \defequalns \ucfunc{E}{\complexnumbers}
  \end{math}.
\end{definition}

Functions vanishing at infinity are defined
as in \cite{lacey1974}.

\begin{definition}
Let \(T\) be a locally compact Hausdorff space and
\(K = \realnumbers\) or \(K = \complexnumbers\). Function \(f \in
\cbfunc{T}{K}\) is said to \defterm{vanish at infinity} if the set
\(\{t \in T \setsep \abs{f(t)} \geq \myepsilon\}\) is compact for all
\(\myepsilon \in \positiverealnumbers\).
The Banach space of all functions vanishing at infinity is
defined by
\begin{displaymath}
  \vanishingfunc{T}{K} \defequalns \left\{ f \in \cbfunc{T}{K} \setsep f
    \;\textrm{vanishes at infinity} \right\}
\end{displaymath}                   
and
\begin{math}
  \vanishingfunccv{T} \defequalns \vanishingfunc{T}{\complexnumbers}
\end{math}.
\end{definition}

\begin{definition}
  Let \(T\) be a locally compact Hausdorff space and
  \(K = \realnumbers\) or \(K = \complexnumbers\).
  Define
  \begin{displaymath}
    \cscfunc{T}{K} \defequalns
    \left\{ f \in \cbfunc{T}{K} \setsep \suppop f \;\textrm{is compact}
    \right\}
  \end{displaymath}
  and
  \begin{math}
    \cscfunccv{T} \defequalns \cscfunc{T}{\complexnumbers}
  \end{math}.
\end{definition}

We have
\begin{math}
  \cscfunccv{\rn}
  \normedsubspace
  \vanishingfunccv{\rn}
  \closedsubspace
  \ucfunccv{\rn}
  \closedsubspace
  \cbfunccv{\rn}
\end{math}
for each \(n \in \positiveintegers\)
\cite[section 2]{donoho1992},
\cite{lacey1974}.
}
\fi

\if\shortprep1
{
\begin{definition}
  \label{def:gencospace}
  Let \(K = \realnumbers\) or \(K = \complexnumbers\) and
  \(I\) be a denumerable set.
  Define
  \begin{eqnarray*}
    & & \gencospace{I}{K} := \left\{ \left( a_\lambda \right)_{\lambda 
    \in
        I} \setsep \left( a_{\sigma(k)} \right)_{k=0}^\infty 
        \in
      \gencospace{\naturalnumbers}{K} \; \textrm{for some bijection}
      \; \sigma : \naturalnumbers \onto I
    \right\} \\
    & & \nsznorminspace{\sqa}{\gencospace{I}{K}}
    := \norminfty{\sqa} , \spaceafter \sqa \in \gencospace{I}{K} .
  \end{eqnarray*}
\end{definition}
}
\fi

\if\shortprep1
{
Banach space
\begin{math}
  \gencospace{I}{K}
\end{math}
is isometrically isomorphic to
\begin{math}
  \gencospace{\naturalnumbers}{K}
\end{math}.
}
\fi

\begin{definition}
  \label{def:distr-translation-and-dilatation}
  \definekn
  Let \(E\) be a closed subspace of
  \begin{math}
    \cbfunc{\rn}{K}
  \end{math}
  for which the following condition is true:
  \begin{displaymath}
    \forall f \in \cbfunc{\rn}{K}, c \in \positiverealnumbers ,
    \sqd \in \rn :
    f \in E \iff f(c \cdot - \sqd) \in E .
  \end{displaymath}
  Let \(a \in \positiverealnumbers\), \(\rnb \in \rn\),
  and \(\dualelem{f} \in \topdual{E}\).
  The \defterm{\(a\)-dilatation and \(\rnb\)-translation of 
  \(\dualelem{f}\)}, 
  denoted by
  \(\dualelem{f}(a \cdot - \rnb)\), is defined by \cite{cl1996}
  \begin{displaymath}
    \dualappl{\dualelem{f}(a \cdot - \rnb)}{f} := 
    \frac{1}{a^n}
    \szdualappl{\dualelem{f}}{f \left( \frac{\cdot + \rnb}{a} 
    \right)}
  \end{displaymath}
  for all
  \begin{math}
    f \in E
  \end{math}.
\end{definition}

\begin{definition}
  Let \(f \in \complexnumbers^\rn\),
  \(r_1 \in \positiverealnumbers\),
  and \(\dualf \in 
  \topdual{\contfunccv{\closedball{\rn}{0}{r_1}}}\).
  Let \(m \in \naturalnumbers\).
  We say that pair \((\dualf, f)\) \defterm{spans all
  the polynomials of degree at most \(m\)}
  iff
  \begin{displaymath}
    \sum_{\firstznvar \in \zn}
    \szdualappl{\dualf}
    {\restrictfunc{p(\cdot + \firstznvar)}{\closedball{\rn}{0}{r_1}}}
    f(\rnx - \firstznvar)
    = p(\rnx)
  \end{displaymath}
  for all \(\rnx \in \rn\) and for all polynomials
  \(p\) of \(n\) variables that are of degree at most \(m\).
\end{definition}

We define the modulus of continuity in the standard way
\cite[def. V.4.2]{bs1988}:

\begin{definition}
  \label{def:modcont}
  When \(n \in \positiveintegers\),
  \(m \in \positiveintegers\), and
  \(p \in [1, \infty]\)
  define
  \begin{displaymath}
    \genmodcont{m}{p}{f}{t}
    :=
    \sup
      \{
      \norminspace{\funcdiff{m}{\rnh}{f}}{\biglpcv{p}{\rn}}
      \setsep
      \rnh \in \closedball{\rn}{0}{t}
      \}
  \end{displaymath}
  and
  \begin{displaymath}
    \modcont{f}{t} := \genmodcont{1}{\infty}{f}{t}
  \end{displaymath}  
  for all
  \begin{math}
    f \in \borelfunc{\rn}{\complexnumbers}
  \end{math}
  and
  \begin{math}
    t \in \nonnegrealnumbers
  \end{math}.
\end{definition}

\if\shortprep0
{
We have
\begin{displaymath}
  \lim_{t \to 0} \modcont{f}{t} = 0
\end{displaymath}
for all \(f \in \ucfunccv{\rn}\) and
\(n \in \positiveintegers\).
}
\fi

\if\shortprep0
{
The Besov spaces on \(\realnumbers^n\) are denoted by
\(\besovspace{s}{p}{q}{\realnumbers^n}\).
There are two ways to define Besov spaces
and for certain combinations of parameters these yield
different spaces. One way to define these spaces is based on the Fourier
transform and the other on the modulus of continuity.
The former definition is used e.g. by Peetre \cite{peetre1975} and
the latter by DeVore and Popov \cite{dp1988}. Triebel
\cite{triebel1983} uses the  term ``Besov space'' for the
spaces defined by the latter definition and denotes them by
\begin{math}
  \Lambda^s_{p,q}
\end{math}
and he denotes the spaces defined by the former definition by
\begin{math}
  B^s_{p,q}  
\end{math}.
Both definitions yield the same spaces 
for the range of parameters \(s > \frac{n}{p}\) used in
this article \cite{donoho1992}.
The two definitions also yield the same spaces
\(\besovspace{s}{p}{q}{\realnumbers^n}\)
for \(s > 0\), \(1 \leq p < \infty\), and \(1 \leq q \leq \infty\)
\cite{triebel1983}.
The spaces \(\besovspace{s}{p}{q}{\realnumbers^n}\)
are quasi-Banach spaces for \(- \infty < s < \infty\), \(0 < p \leq
\infty\), \(0 < q \leq \infty\) and Banach spaces if 
\(1 \leq p \leq \infty\) and
\(1 \leq q \leq \infty\) \cite{triebel1983}.
}
\fi
\if10
According to DeVore and Popov \cite{dp1988} the definition based on
the modulus of continuity is more convenient in approximation 
theory.
\fi

\if\shortprep1
{
The Besov spaces on \(\realnumbers^n\) are denoted by
\(\besovspace{s}{p}{q}{\realnumbers^n}\).
}
\fi
We have
\begin{math}
  \zygmundspace{s}{\rn} \equaltvs
  \besovspace{s}{\infty}{\infty}{\realnumbers^n}
\end{math}
for all \(s \in \positiverealnumbers\) where
\(\zygmundspace{s}{\rn}\) are the H\"older-Zygmund spaces
\cite{triebel1983}.
We also have
\begin{math}
  \holderspace{s}{\rn} \equaltvs \zygmundspace{s}{\rn}
  \equaltvs \besovspace{s}{\infty}{\infty}{\rn}
\end{math}
when
\begin{math}
  s \in \positiverealnumbers \setminus \positiveintegers
\end{math}
where \(\holderspace{s}{\rn}\),
\(s \in \positiverealnumbers \setminus \positiveintegers\),
are the H\"older spaces on \(\realnumbers^n\) \cite{triebel1983}.
When \(m \in \naturalnumbers\) we denote the Banach space of
functions with bounded and uniformly continuous
partial derivatives up to \(m\)th degree
equipped with the usual norm by
\(\contdiffspace{m}{\rn}\),
see \cite{triebel1983}.
\if\shortprep0
{
When \(m \in \positiveintegers\) we have
\begin{math}
  \contdiffspace{m}{\rn} \subsetset \zygmundspace{m}{\rn}
\end{math}.
}
\fi


\begin{definition}
  \label{def:contbesovnorm}
  Let \(n \in \positiverealnumbers\),
  \(s \in \positiverealnumbers\),
  \(p \in [1, \infty]\), and
  \(q \in [1, \infty]\).
  Let \(m \in \positiveintegers\) and \(m > s\).
  Define
  \begin{displaymath}
    \contbesovnorm{s}{p}{q}{\rn}{m}{f}
    :=
    \lpcvnorm{p}{\rn}{f} + \contbesovtail{s}{p}{q}{m}{f}
  \end{displaymath}
  for all
  \begin{math}
    f \in \borelfunc{\rn}{\complexnumbers}
  \end{math}.
\end{definition}

\if01
{
Norm
\begin{math}
  \contbesovnorm{s}{p}{q}{\rn}{m}{\cdot}
\end{math}
is an equivalent norm on
\begin{math}
  \besovspace{s}{p}{q}{\rn}
\end{math}
for the range of parameters given in \mydef
\ref{def:contbesovnorm}.
Replacing the ranges \(]0,1[\) by
\(]0, \infty[\) in the definition above
results equivalent norms
\cite[def. V.4.3]{bs1988}.
Denote these norms by
\(\contbesovnormrl{s}{p}{q}{\rn}{m}{f}\).
Norms
\begin{math}
  \contbesovnorm{s}{p}{q}{\rn}{m}{\cdot}
\end{math}
and
\begin{math}
  \contbesovnormrl{s}{p}{q}{\rn}{m}{\cdot}
\end{math}
characterize the Besov space
\begin{math}
  \besovspace{s}{p}{q}{\rn}
\end{math}
on
\begin{math}
  \borelfunc{\rn}{\complexnumbers}
\end{math}
for the given range of parameters.
Note that
\begin{math}
  \biglpcv{p}{\rn} \subsetset \schwartzdistr{\rn}
\end{math}
for all
\(n \in \positiveintegers\)
and
\(p \in [1, \infty]\).
See \cite[def. V.4.3]{bs1988}.
}
\fi

\if\shortprep0
{
  Replacing the ranges \(]0,1[\) by
  \(]0, \infty[\) in the definition above
  results equivalent norms.
  Denote these norms by
  \begin{math}
    \contbesovnormrl{s}{p}{q}{\rn}{m}{\cdot}
  \end{math}.
  Norms
  \begin{math}
  \contbesovnorm{s}{p}{q}{\rn}{m}{\cdot}
  \end{math}
  and
  \begin{math}
    \contbesovnormrl{s}{p}{q}{\rn}{m}{\cdot}
  \end{math}
  are equivalent norms on
  \begin{math}
    \besovspace{s}{p}{q}{\rn}
  \end{math}
  for the range of parameters given in \mydef
  \ref{def:contbesovnorm}.
  Norms
  \begin{math}
    \contbesovnorm{s}{p}{q}{\rn}{m}{\cdot}
  \end{math}
  and
  \begin{math}
    \contbesovnormrl{s}{p}{q}{\rn}{m}{\cdot}
  \end{math}
  characterize the Besov space
  \begin{math}
    \besovspace{s}{p}{q}{\rn}
  \end{math}
  on
  \begin{math}
    \borelfunc{\rn}{\complexnumbers}
  \end{math}
  for the given range of parameters.
}
\else
{
  Norm
  \begin{math}
  \contbesovnorm{s}{p}{q}{\rn}{m}{\cdot}
  \end{math}
  is an equivalent norm on
  \begin{math}
    \besovspace{s}{p}{q}{\rn}
  \end{math}
  for the range of parameters given in \mydef
  \ref{def:contbesovnorm}.
  Norm
  \begin{math}
    \contbesovnorm{s}{p}{q}{\rn}{m}{\cdot}
  \end{math}
  characterizes the Besov space
  \begin{math}
    \besovspace{s}{p}{q}{\rn}
  \end{math}
  on
  \begin{math}
    \borelfunc{\rn}{\complexnumbers}
  \end{math}
  for the given range of parameters.
}
\fi
Note that
\begin{math}
  \biglpcv{p}{\rn} \subsetset \schwartzdistr{\rn}
\end{math}
for all
\(n \in \positiveintegers\)
and
\(p \in [1, \infty]\).
See \cite[def. V.4.3]{bs1988}.

\subsection{Tensor Products}
\label{sec:preliminaries-tensorproducts}

\begin{definition}
  When \(n \in \positiveintegers\) define
\if\shortprep0
{
  \begin{displaymath}
    \ntensornbv{\firstznvar}
    :=
    \indexedtensorproduct{l=1}{n} \nnatbasvec{\seqelem{\firstznvar}{l}}
  \end{displaymath}
  and
}
\fi
  \begin{displaymath}
    \ztensornbv{\firstznvar}
    :=
    \indexedtensorproduct{l=1}{n} \znatbasvec{\seqelem{\firstznvar}{l}} .
  \end{displaymath}
  for all \(\firstznvar \in \zn\).
\end{definition}

\if\shortprep0
{
The square ordering is defined as in \cite{ryan2002} and
\cite{singer1970} in the following.

\begin{definition}
  Define function \(\sqordfunction : \naturalnumbers \to \naturalnumbers^2\) by
  \begin{equation}
    \label{eq:sq-ordering-a}
    \sqord{0} = (0, 0)
  \end{equation}
  and
  \begin{equation}
    \label{eq:sq-ordering-b}
    \sqord{k} =
    \left\{
    \begin{array}{ll}
      (i, n) ; \;\; & k = n^2 + i \;\land\; i \in \setzeroton{n}
      \;\land\; n \in \naturalnumbers \\
      (n, n - i) ; \;\; & k = n^2 + n + i \;\land\; i \in \setoneton{n}
      \;\land\; n \in \naturalnumbers
    \end{array}
    \right.
  \end{equation}
  Function \(\sqordfunction\) is called the square ordering.
  Define also functions
  \(\sqordfirstfunction : \naturalnumbers \to \naturalnumbers\)
  and \(\sqordsecondfunction : \naturalnumbers \to \naturalnumbers\) by setting
  \begin{math}
    \sqordfirst{k} := \seqelem{\sqord{k}}{1}
  \end{math}
  and
  \begin{math}
    \sqordsecond{k} := \seqelem{\sqord{k}}{2}
  \end{math}
  for all \(k \in \naturalnumbers\).
\end{definition}

Function \(\sqordfunction\) is a bijection from \(\naturalnumbers\) onto
\(\naturalnumbers^2\). The tensor product basis is defined as in
\cite{ryan2002}:


\begin{definition}
  Let \(E\) and \(F\) be Banach spaces that both have Schauder bases.
  Let \((a_k)_{k=0}^\infty\) and \((b_k)_{k=0}^\infty\) be Schauder
  bases of \(E\) and \(F\), respectively. The sequence
  \begin{math}
    (a_{\sqordfirst{k}} \otimes b_{\sqordsecond{k}})_{k=0}^\infty
    \sequenceof E \otimes F
  \end{math}
  is called the \defterm{tensor product basis
  generated by Schauder bases
  \((a_k)_{k=0}^\infty\) and \((b_k)_{k=0}^\infty\)}.
\end{definition}

The tensor product basis is a Schauder basis for both \(E \ctp_\injtn
F\) and \(E \ctp_\projtn F\) \cite{ryan2002,singer1970}.
}
\fi

\if\shortprep0
{
When \(E\) and \(F\) are Banach spaces
consisting of complex valued functions
(in particular, real valued functions),
\(u \in E\), and \(v \in F\)
the tensor product \(u \otimes v\) is identified with a function on
the Cartesian product of the domains of \(u\) and \(v\) by setting
\begin{math}
  (u \otimes v)(x \otimes y) = u(x) v(y)
\end{math}
where \(x\) belongs to the domain of \(u\) and \(y\) belongs to the domain of
\(v\).
}
\fi

When \(n \in \positiveintegers\), \(\alpha\) is a uniform
crossnorm, and \(A_1, \ldots, A_n\) are Banach spaces then the
completed tensor product of several Banach spaces is defined by the
recursive formula
\begin{displaymath}
  A_1 \ctp_\alpha \cdots \ctp_\alpha A_n \defequalns
  ( A_1 \ctp_\alpha \cdots \ctp_\alpha A_{n-1} ) \ctp_\alpha A_n .
\end{displaymath}
The indexed completed tensor product is defined by
\begin{displaymath}
  \indexedctp{\alpha}{j=1}{n} A_j \defequalns
  \left( \indexedctp{\alpha}{j=1}{n-1} A_j \right) \ctp_\alpha A_n
\end{displaymath}
where \(n > 1\). When \(n = 1\) define
\begin{displaymath}
  \indexedctp{\alpha}{j=1}{1} A_j \defequalns A_1 .
\end{displaymath}

\if\shortprep0
{
The concept of norm in tensor product spaces is used in two slightly
different ways. When we speak about norm \(\alpha\) defined on space
\(E \otimes F\) we mean that \(\alpha\) is a function that maps each
element \(u \in E \otimes F\) to its norm \(\norm{u}\) in
\(\nonnegrealnumbers\). On the other hand uniform crossnorms and
tensor norms are assignments of a reasonable crossnorm to each pair of
Banach spaces.
}
\fi

\if\shortprep0
{
A tensor norm is defined to be a finitely generated uniform crossnorm \cite{ryan2002}. The Schatten dual of a tensor
norm \(\alpha\) is
denoted by \(\alpha^s\).
When \(\alpha\) is a tensor norm \(\alpha^s\) is a uniform crossnorm. 
}
\fi
Uniform crossnorms do not generally respect subspaces.
\if\shortprep0
{
I.e. if
\(X\) is a closed subspace of a Banach space \(E\) and \(Y\) is a 
closed subspace of a Banach space \(F\) the norm of an element
\(u \in X \ctp_\alpha Y\)
is not necessary equal in spaces \(X \ctp_\alpha Y\) and \(E
\ctp_\alpha F\) \cite{ryan2002}.
However, it is possible to use the norm 
inherited from \(E \ctp_\alpha F\) in vector space \(X \otimes Y\). 
Then the closure of
\(X \otimes Y\) with the inherited norm in \(E \ctp_\alpha F\) is a
closed subspace of \(E \ctp_\alpha F\) \cite[chapter 1]{lc1985}. 
}
\else
{
However, when \(X\) is a closed subspace of a Banach space \(E\),
\(Y\) is a closed subspace \(F\), and \(\alpha\) is a reasonable
crossnorm on \(E \otimes F\) it is possible to use the norm
inherited from \(E \otimes_\alpha F\) on vector space
\(X \otimes Y\) \cite[chapter 1]{lc1985}.
}
\fi
We give the following definition for this kind of construction.

\begin{definition}
  Let \(E\) and \(F\) be Banach spaces and \(\alpha\) a norm on \(E
  \otimes F\). Let \(X\) be a closed subspace of \(E\) and \(Y\) a
  closed subspace of \(F\). Define
  \begin{math}
    X \otimes_{\alpha ; E \ctp_\alpha F} Y
  \end{math}
  to be the normed vector space formed by using the norm inherited
  from \(E \ctp_\alpha F\) in the vector space \(X \otimes Y\), i.e.
  \begin{math}
    \norm{u}_{X \otimes_{\alpha ; E \ctp_\alpha F} Y} := \alpha_{E,F}(u)
  \end{math}
  for all \(u \in X \otimes Y\).
  Define
  \begin{math}
    X \ctp_{\alpha ; E \ctp_\alpha F} Y
  \end{math}
  to be the closure of \( X \otimes_{\alpha ; E \ctp_\alpha F} Y \) in
  space \(E \ctp_\alpha F\).
  The notations
  \begin{math}
    X \otimes_{\alpha ; E \otimes F} Y
  \end{math}
  and
  \begin{math}
    X \ctp_{\alpha ; E \otimes F} Y
  \end{math}
  may also be used for the aforementioned definitions.
\end{definition}

Now \(X \otimes_{\alpha ; E \ctp_\alpha F} Y\) is a normed subspace of
\(E \ctp_\alpha F\) and \(X \ctp_{\alpha ; E \ctp_\alpha F} Y\) is a
closed subspace of \(E \ctp_\alpha F\).

\if\shortprep1
{
\begin{definition}
  \label{def:indexed-inh-tp-alt}
  Let \(n \in \positiveintegers\). Let \(A_1, \ldots, A_n\) and
  \(B_1, \ldots, B_n\) be Banach spaces so that \(A_1\) is a closed subspace
  of \(B_1\) and \(B_k \otimes A_{k+1}\) is a linear subspace of
  \(B_{k+1}\) for \(k \in \setoneton{n-1}\).
  When \(k \geq 2\) define
  \begin{displaymath}
    \indexedinhctp{B_j}{j=1}{k} A_j \defequalns \closop_{B_k} T_k
  \end{displaymath}
  where
  \begin{displaymath}
    T_k \defequalns \left( \indexedinhctp{B_j}{j=1}{k - 1} A_j \right)
    \otimes_{(B_k)} A_k .
  \end{displaymath}
  When \(k = 1\) define
  \begin{displaymath}
    \indexedinhctp{B_j}{j=1}{k} A_j \defequalns A_1 .
  \end{displaymath}
\end{definition}
}
\fi

\if\shortprep1
{
When the assumptions of \mydef
\ref{def:indexed-inh-tp-alt} hold
\begin{displaymath}
  \indexedinhctp{B_j}{j=1}{k} A_j
\end{displaymath}
is a closed subspace of Banach space \(B_k\) for all \(k \in \setoneton{n}\).
}
\fi

\if\shortprep1
{
\begin{definition}
  \label{def:indexed-inh-op-tp}
  Let \(n \in \positiveintegers\).
  Let \(A_1, \ldots, A_n\), \(B_1, \ldots, B_n\),
  \(E_1, \ldots, E_n\), and \(F_1, \ldots, F_n\) be Banach spaces so
  that
  \begin{itemize}
  \item \(A_1\) is a closed subspace of \(E_1\) and
    \(E_k \otimes A_{k+1}\) is a linear subspace of \(E_{k+1}\) for
    all \(k = 1, \ldots, n - 1\).
  \item \(B_1\) is a closed subspace of \(F_1\) and
    \(F_k \otimes B_{k+1}\) is a linear subspace of \(F_{k+1}\) for
    all \(k = 1, \ldots, n - 1\).
  \end{itemize}
  Suppose that \(P_k : A_k \to B_k\), \(k = 1, \ldots, n\) are
  operators.
  Let \(S_1 := P_1\), \(T_1 := S_1 = P_1\), and
  \begin{math}
    S_k := T_{k-1} \otimes P_k
  \end{math}
  for \(k = 2, \ldots, n\).
  When \(k \in \{ 2, \ldots, n \}\) and \(S_k\) is continuous let
  \begin{displaymath}
    T_k : \indexedinhctp{E_l}{l=1}{k} A_l \to
    \indexedinhctp{F_l}{l=1}{k} B_l ,
  \end{displaymath}
  be the unique continuous linear extension of \(S_k\) to
  \begin{displaymath}
    \indexedinhctp{E_l}{l=1}{k} A_l .
  \end{displaymath}
  If \(k \in \{ 2, \ldots, n\}\) and \(S_k\) is not continuous let
  \(T_k = 0\).
  When all of the functions \(S_k\), \(T_k\), \(k = 1, \ldots, n\) are
  operators define
  \begin{displaymath}
    \indexedinhopctp{E_k}{F_k}{k=1}{n} P_k = T_n .
  \end{displaymath}
  If any of the functions \(S_k\), \(T_k\), \(k = 1, \ldots, n\) is
  not an operator then
  \begin{displaymath}
    \indexedinhopctp{E_k}{F_k}{k=1}{n} P_k
  \end{displaymath}
  is undefined.
\end{definition}
}
\fi

\if\shortprep0
{
If \(\alpha\) is a reasonable crossnorm, \(a \in \topdual{X}\), and \(b
\in \topdual{Y}\) then \(a \otimes b\) is a continuous linear
functional on \(X \otimes_\alpha Y\) and \(\norm{a \otimes b} =
\norm{a} \norm{b}\) \cite{ryan2002}. We have
\begin{equation}
  \label{eq:tensor-functional-a}
  (a \otimes b)(x \otimes y) = a(x)b(y)
\end{equation}
for all
\begin{math}
  x \in X
\end{math}
and
\begin{math}
y \in Y
\end{math}
and the elements of  \(\topdual{X} \otimes \topdual{Y}\) may be
interpreted as linear functionals on \(X \otimes_\alpha Y\).
There is an isometric embedding
\begin{equation}
  \label{eq:schatten-embedding}
  \topdual{X} \otimes_{\alpha^s} \topdual{Y} \isomemb \topdual{\left( X
      \otimes_\alpha Y \right)} .
\end{equation}
for tensor norms \(\alpha\) \cite{ryan2002}.
See also \cite{df1993}.
}
\else
{
There is an isometric embedding
\begin{equation}
  \label{eq:schatten-embedding}
  \topdual{X} \otimes_{\alpha^s} \topdual{Y} \isomemb \topdual{\left( X
      \otimes_\alpha Y \right)}
\end{equation}
for tensor norms \(\alpha\) and Banach spaces \(X\) and \(Y\) \cite{ryan2002},
see also \cite{df1993}.
}
\fi
\if\shortprep1
{
When \(K = \realnumbers\) or \(K = \complexnumbers\),
\(n \in \naturalnumbers + 2\), and \(X_1, \ldots, X_n\) are closed
subspaces of \(\cbfunc{\realnumbers}{K}\) the tensor
product \(X_1 \citp \cdots \citp X_n\) is a closed
subspace of \(\cbfunc{\realnumbers^n}{K}\).
}
\fi

\if\shortprep0
{

\section{Some Theorems on Function Spaces and Sequences}
\label{sec:function-spaces-and-sequences}

\if\shortprep0
{
\if\longversion0
{
\givelemmawithoutproof
}
\fi
}
\fi

\if\shortprep0
{
\begin{lemma}
  Let \(K = \realnumbers\) or \(K = \complexnumbers\). Let
  \(I\) be a denumerable set and
  \begin{math}
    \sqa \in \seqset{I}{K}
  \end{math}.
  Suppose that \(\sigma_1\) and \(\sigma_2\) are bijections from
  \(\naturalnumbers\) onto \(I\).
  Then
  \begin{displaymath}
    \left( \seqelem{\sqa}{\sigma_1(k)} \right)_{k=0}^\infty \in
    \gencospace{\naturalnumbers}{K}
    \iff
    \left( \seqelem{\sqa}{\sigma_2(k)} \right)_{k=0}^\infty \in
    \gencospace{\naturalnumbers}{K} .
  \end{displaymath}
\end{lemma}

\if\longversion1
{
\begin{proof}
  Suppose that
  \begin{math}
    ( \seqelem{\sqa}{\sigma_1(k)} )_{k=0}^\infty \in
    \gencospace{\naturalnumbers}{K}
  \end{math}.
  Let \(h > 0\). Then there exists \(n_1 \in \naturalnumbers\) so that
  \begin{math}
    \abs{\seqelem{\sqa}{\sigma_1(k)}} < h
  \end{math}
  for each \(k \in \naturalnumbers\), \(k > n_1\).
  Let
  \begin{displaymath}
    n_2 := \max \left\{ \sigma_2^{-1} \left( \sigma_1 \left( k
          \right) \right) \setsep k \in \setzeroton{n_1} \right\} .
  \end{displaymath}
  Suppose that \(k \in \naturalnumbers\), \(k > n_2\).
  Then
  \begin{align*}
    & k \neq \sigma_2^{-1} \left( \sigma_1 \left( l \right) \right)
    \qspace \forall l \in \setzeroton{n_1} \\
    & \implies
    \sigma_2 \left( k \right) \neq \sigma_1 \left( k \right)
    \qspace \forall l \in \setzeroton{n_1} \\
    & \implies
    \sigma_2(k) = \sigma_1 \left( l_0 \right)
    \; \textrm{for some} \; l_0 \in \naturalnumbers, l_0 > n_1 \\
    & \implies
    \abs{\seqelem{\sqa}{\sigma_2(k)}}
    = \abs{\seqelem{\sqa}{\sigma_1(l_0)}} < h .
  \end{align*}
  Hence
  \begin{displaymath}
    \lim_{k \to \infty} \seqelem{\sqa}{\sigma_2(k)} = 0 .
  \end{displaymath}
  and so
  \begin{math}
    ( \seqelem{\sqa}{\sigma_2(k)} )_{k=0}^\infty \in
    \gencospace{\naturalnumbers}{K}    
  \end{math}.
  Similarly, if
  \begin{math}
    ( \seqelem{\sqa}{\sigma_2(k)} )_{k=0}^\infty \in
    \gencospace{\naturalnumbers}{K}
  \end{math}
  it follows that
  \begin{math}
    ( \seqelem{\sqa}{\sigma_1(k)} )_{k=0}^\infty \in
    \gencospace{\naturalnumbers}{K}    
  \end{math}.
\end{proof}
}
\fi

Consequently the following definition makes sense.
}
\fi

\begin{definition}
  \label{def:gencospace}
  Let \(K = \realnumbers\) or \(K = \complexnumbers\) and
  \(I\) be a denumerable set.
  Define
  \begin{eqnarray*}
    & & \gencospace{I}{K} := \left\{ \left( a_\lambda \right)_{\lambda 
    \in
        I} \setsep \left( a_{\sigma(k)} \right)_{k=0}^\infty 
        \in
      \gencospace{\naturalnumbers}{K} \; \textrm{for some bijection}
      \; \sigma : \naturalnumbers \onto I
    \right\} \\
    & & \nsznorminspace{\sqa}{\gencospace{I}{K}}
    := \norminfty{\sqa} , \spaceafter \sqa \in \gencospace{I}{K} .
  \end{eqnarray*}
\end{definition}

Banach space
\begin{math}
  \gencospace{I}{K}
\end{math}
is isometrically isomorphic to
\begin{math}
  \gencospace{\naturalnumbers}{K}
\end{math}.

\begin{definition}
  \label{def:ncover}
  \definekn
  Let \(f\) be a function from \(\realnumbers^n\)
  into \(K\).
  Define
  \begin{displaymath}
    \ncover{f} := \max_{\rnx \in \realnumbers^n} \card{\left\{
    \firstznvar \in \integernumbers^n \setsep
    f(\rnx - \firstznvar) \neq 0 \right\}} .
  \end{displaymath}
\end{definition}

\if\longversion1
{
\(\ncover{f} < \infty\) for all compactly supported
functions \(f : \rn \to K\).
When \(a \in \positiverealnumbers\) function \(\rny \in \rn \mapsto a \rny\) 
is a bijection from \(\rn\) onto \(\rn\) and it follows that
\begin{equation}
  \label{eq:ncover-scaling}
  \max_{\rnx \in \realnumbers^n} \card{ \left\{
    \firstznvar \in \integernumbers^n
    \setsep f(a \rnx - \firstznvar) \neq 0 \right\}}
  = \ncover{f} .
\end{equation}
}
\fi

\if\shortprep0
{
\begin{definition}
  Let \(E\) and \(F\) be normed vector spaces.
  When \(f\) is a compactly supported function from \(E\)
  into \(F\) define
  \[
  \supportradius{f} :=
  \inf \{ r \in \nonnegrealnumbers \setsep \suppopst f \subset
  \closedball{E}{0}{r} \} .
  \]
\end{definition}
}
\fi

\if\shortprep0
{
\begin{definition}
	Let \(n \in \positiveintegers\).
	Define
	\begin{displaymath}
	  \itrans{f}{\rnx} :=
	  \{ \firstznvar \in \zn \setsep
	  f( \rnx - \firstznvar ) \neq 0
	  \}
	\end{displaymath}
	for all \(f \in \complexnumbers^{\realnumbers^n}\)
	and \(\rnx \in \rn\).
\end{definition}
}
\fi

\if\shortprep0
{
\begin{definition}
	Let \(\rna, \rnb \in \rn\) and \(n \in \positiveintegers\).
	Define
	\begin{displaymath}
	  \rectangle{\rna}{\rnb}
	  =
	  \mathop{\times}\limits_{k = 1}^{n}
	  I_1 \left(\cartprodelem{\rna}{k}, \cartprodelem{\rnb}{k}\right)
	\end{displaymath}
	where
	\begin{displaymath}
	  I_1(u, v)
		=
		[ \min \{ u, v \}, \max \{ u, v \} ] \intersection \integernumbers
	\end{displaymath}
	for all \(u, v\) in \(\realnumbers\).
\end{definition}
}
\fi

\if\longprep1
{
\givelemmaswithoutproofsn{five}
}
\fi

\if\shortprep0
{
\begin{lemma}
  \label{lem:grid-convergence-multidim}
  \definekn
  Let \(c \in \positiverealnumbers\), and
  \(f \in \cscfunc{\realnumbers^n}{K}\).
  When \(\seqstyle{a} \in \gencospace{\zn}{K}\) series
  \begin{displaymath}
    \sum_{\firstznvar \in \zn} \seqelem{\sqa}{\firstznvar}
    f(c \cdot - \firstznvar)
  \end{displaymath}
  converges unconditionally in \(\vanishingfunc{\realnumbers^n}{K}\)
  and
  \if\longversion0
  {
  \begin{displaymath}
    \norminfty{\sum_{\firstznvar \in \zn}
    \seqelem{\sqa}{\firstznvar} f(c \cdot - \firstznvar)} \leq
    \ncover{f} \norminfty{f} \norminfty{\sqa} .
  \end{displaymath}
  }
  \else
  {
  \begin{equation}
    \label{eq:grid-convergence-multidim-norm}
    \norminfty{\sum_{\firstznvar \in \zn}
    \seqelem{\sqa}{\firstznvar} f(c \cdot - \firstznvar)} \leq
    \ncover{f} \norminfty{f} \norminfty{\sqa} .
  \end{equation}
  }
  \fi
\end{lemma}
}
\fi

\if\longversion1
{
\begin{proof}
  Let \(\sqa \in \gencospace{\zn}{K}\).
  Let \(\sigma : \naturalnumbers \onto \zn\) be a bijection.
  Let
  \begin{displaymath}
    s_l = \sum_{i = 0}^l \seqelem{\sqa}{\sigma(i)}
    f \left( c \cdot - \sigma(i) \right)
  \end{displaymath}
  Let \(l, l' \in \naturalnumbers\), \(l' > l\). Now
  \begin{displaymath}
    s_{l'} - s_l = \sum_{i = l + 1}^{l'}
    \seqelem{\sqa}{\sigma(i)}
    f \left( c \cdot - \sigma(i) \right) .
  \end{displaymath}
  It follows from \myequation \eqref{eq:ncover-scaling} that
  \begin{displaymath}
    \abs{s_{l'}(\rnx) - s_l(\rnx)} =
    \abs{\sum_{i=l+1}^{l'}
    \seqelem{\sqa}{\sigma(i)}
      f \left( c \rnx - \sigma(i) \right)}
    \leq \ncover{f} b_{l,l'} m
  \end{displaymath}
  where
  \begin{math}
    b_{l,l'} = \max \left\{ \abs{\seqelem{\sqa}{\sigma(i)}}
    \setsep i = l+1, \ldots, l' \right\}
  \end{math}
  and
  \begin{math}
    m = \ncover{f} \norminfty{f}
  \end{math}.
  Let \(h > 0\). There exists \(i_0 \in \naturalnumbers\) so that
  \begin{displaymath}
    \abs{\seqelem{\sqa}{\sigma(i)}} < \frac{h}{m}
  \end{displaymath}
  for all
  \begin{math}
    i \in \naturalnumbers
  \end{math}
  and
  \begin{math}
    i > i_0
  \end{math}.
  Now
  \begin{math}
    b_{l,l'} < \frac{h}{m}
  \end{math}
  for all
  \begin{math}
    \forall l, l' > i_0
  \end{math}
  and
  \begin{math}
    l' > l
  \end{math}.
  We also have
  \begin{displaymath}
    \abs{s_{l'}(\rnx) - s_l(\rnx)} < \frac{h}{m} \cdot m = h
  \end{displaymath}
  for all \(\rnx \in \realnumbers\) and for all \(l, l' \in
  \naturalnumbers\), \(l, l' > i_0\), \(l' > l\).
  Hence \((s_l)_{l = 0}^\infty\) is a Cauchy sequence in
  \(\vanishingfunc{\rn}{K}\) and it converges
  in \(\vanishingfunc{\rn}{K}\).
  Therefore
  \begin{displaymath}
    \sum_{\firstznvar \in \zn}
    \seqelem{\sqa}{\firstznvar} f(c \cdot - \firstznvar)
  \end{displaymath}
  converges unconditionally in \(\vanishingfunc{\rn}{K}\).
  It follows from \myequation \eqref{eq:ncover-scaling} that
  \begin{displaymath}
    \abs{\sum_{\firstznvar \in \zn}
    \seqelem{\sqa}{\firstznvar}
    f(c \rnx - \firstznvar)}
    \leq \ncover{f} \norminfty{f} \norminfty{\sqa}
  \end{displaymath}
  for all \(\rnx \in \rn\). Hence inequality
  \eqref{eq:grid-convergence-multidim-norm} is true.
\end{proof}
}
\fi

\if\shortprep0
{
\begin{lemma}
  \label{lem:single-res-norm-equiv-multidim}
  \definekn
  Let \(c \in \positiverealnumbers\).
  Let \(f \in \cscfunc{\realnumbers^n}{K}\)
  so that
  \begin{equation}
    \label{eq:card-a}
    f(\firstznvar + \vectorstyle{b}) =
    \delta_{\firstznvar,0}
  \end{equation}
  for all
  \begin{math}
    \firstznvar \in \integernumbers^n
  \end{math}
  for some constant \(\vectorstyle{b} \in \realnumbers^n\).
  Then
  \if\longversion0
  {
  \begin{displaymath}
    \norminfty{\sqa}
    \leq
    \norminfty{\sum_{\firstznvar \in \zn}
      \seqelem{\sqa}{\firstznvar} f( c \cdot - \firstznvar )}
    \leq
    \ncover{f} \norminfty{f} \norminfty{\sqa}
  \end{displaymath}
  for all \(\sqa \in \gencospace{\zn}{K}\).
  }
  \else
  {
  \begin{equation}
    \label{eq:single-res-norm-inequality}
    \norminfty{\sqa}
    \leq
    \norminfty{\sum_{\firstznvar \in \zn}
      \seqelem{\sqa}{\firstznvar} f( c \cdot - \firstznvar )}
    \leq
    \ncover{f} \norminfty{f} \norminfty{\sqa}
  \end{equation}
  for all \(\sqa \in \gencospace{\zn}{K}\).
  }
  \fi
\end{lemma}
}
\fi

\if\longversion1
{
\begin{proof}
  Let \(\sqa \in \gencospace{\zn}{K}\).
  By \mylemma \ref{lem:grid-convergence-multidim} series
  \begin{displaymath}
    \sum_{\firstznvar \in \zn} \seqelem{\sqa}{\firstznvar}
    f( c \cdot - \firstznvar )
  \end{displaymath}
  converges unconditionally in \(\vanishingfunc{\realnumbers^n}{K}\)
  and the right-hand side of inequality
  \eqref{eq:single-res-norm-inequality} is true.
  Let \(\sigma : \naturalnumbers \onto \zn\) be a bijection and
  \(m := \norminfty{\sqa}\). Since
  \begin{displaymath}
    \lim_{l \to \infty} \seqelem{\sqa}{\sigma(l)} = 0
  \end{displaymath}
  there exists \(l_0 \in \naturalnumbers\) so that
  \(\abs{\seqelem{\sqa}{\sigma(l_0)}} = m\). Now
  \begin{align*}
    \norminfty{\sum_{\firstznvar \in \zn}
    \seqelem{\sqa}{\firstznvar}
    f( c \cdot - \firstznvar )}
    & \geq
    \abs{\sum_{\firstznvar \in \zn}
    \seqelem{\sqa}{\firstznvar}
    f \left( c \frac{\sigma(l_0) + \vectorstyle{b}}{c}
        - \firstznvar \right)} \\
    & =
    \abs{\sum_{\firstznvar \in \zn} \seqelem{\sqa}{\firstznvar}
    f \left( \sigma(l_0) + \vectorstyle{b} - \firstznvar \right)} \\
    & =
    \abs{\seqelem{\sqa}{\sigma(l_0)}} = \norminfty{\sqa} .
  \end{align*}
\end{proof}
}
\fi

\if\shortprep0
{
It follows from \myequation \eqref{eq:card-a} that
\(\norminfty{f} \geq 1\) and \(\ncover{f} \geq 1\) under the conditions of 
\mylemma \ref{lem:single-res-norm-equiv-multidim}.
}
\fi

\if\shortprep0
{
\begin{lemma}
  \label{lem:cospace-top-isomorphism}
  \definek
  Let \(c \in \positiverealnumbers\).
  Let \(f \in \cscfunc{\rn}{K}\)
  so that
  \begin{math}
    f(\firstznvar + \rnb) = \delta_{\firstznvar,0}
  \end{math}
  for all \(\firstznvar \in \integernumbers^n\)
  for some \(\rnb \in \realnumbers^n\).
  Let
\if\longversion1
{
  \begin{eqnarray}
    \nonumber
    & & V := \left\{ \sum_{\firstznvar \in \zn}
      \seqelem{\sqa}{\firstznvar} f(c \cdot - \firstznvar)
      \setsep \seqstyle{a} \in \gencospace{\zn}{K} \right\} \\
    \label{eq:cospace-isom-a}
    & & \norminspace{g}{V} := \norminfty{g}, \spaceafter g \in V .
  \end{eqnarray}
}
\else
{
  \begin{eqnarray*}
    & & V := \left\{ \sum_{\firstznvar \in \zn}
      \seqelem{\sqa}{\firstznvar} f(c \cdot - \firstznvar)
      \setsep \seqstyle{a} \in \gencospace{\zn}{K} \right\} \\
    & & \norminspace{g}{V} := \norminfty{g}, \spaceafter g \in V .
  \end{eqnarray*}
}
\fi
  Define function \(\iota : \gencospace{\zn}{K} \to V\) by
  \begin{displaymath}
    \iota(\seqstyle{a}) :=
    \sum_{\firstznvar \in \zn}
    \seqelem{\sqa}{\firstznvar} f(c \cdot - \firstznvar)
  \end{displaymath}
  for all
  \begin{math}
    \seqstyle{a} \in \gencospace{\zn}{K}
  \end{math}.
  Then
  \begin{itemize}
  \item[(i)]
    \(V\) is a closed subspace of \(\vfn{n}{K}\).
  \item[(ii)]
    Function \(\iota\) is a topological isomorphism from 
    \(\gencospace{\zn}{K}\) onto
    \(V\) and
    \if\longversion1
    {
    \begin{equation}
      \label{eq:cospace-norm-equiv}
      \norminfty{\seqstyle{a}}
      \leq \norminfty{\iota(\seqstyle{a})}
      \leq \ncover{f} \norminfty{f} \norminfty{\seqstyle{a}}
    \end{equation}
    }
    \else
    {
    \begin{displaymath}
      \norminfty{\seqstyle{a}}
      \leq \norminfty{\iota(\seqstyle{a})}
      \leq \ncover{f} \norminfty{f} \norminfty{\seqstyle{a}}
    \end{displaymath}
    }
    \fi
    for all
    \begin{math}
      \seqstyle{a} \in \gencospace{\zn}{K}
    \end{math}.
  \end{itemize}
\end{lemma}
}
\fi

\if\longversion1
{
\begin{proof}
  The series in \myequation \eqref{eq:cospace-isom-a} converges
  unconditionally by \mylemma \ref{lem:grid-convergence-multidim}.
  Function \(\iota\) is linear and a surjection onto \(V\).
  By \mylemma \ref{lem:single-res-norm-equiv-multidim} inequality
  \eqref{eq:cospace-norm-equiv} is true.

  By the definition of \(\iota\) we have
  \begin{math}
    \norminfty{\seqstyle{a}}
    \leq
    \norminfty{\iota(\seqstyle{a})}
  \end{math}
  for all
  \begin{math}
    \seqstyle{a} \in \gencospace{\zn}{K}
  \end{math}.
  Suppose that \(\seqstyle{w}, \seqstyle{z} \in 
  \gencospace{\zn}{K}\) and
  \(\seqstyle{w} \neq \seqstyle{z}\).
  Now
  \begin{math}
    \norminfty{\iota(\seqstyle{w} - \seqstyle{z})}
    \geq \norminfty{\seqstyle{w} - \seqstyle{z}}
    > 0
  \end{math}.
  Hence \(\iota\) is an injection. It follows that \(\iota\) is a
  linear bijection and vector spaces \(V\) and 
  \(\gencospace{\zn}{K}\) are
  algebraically isomorphic.
  We have
  \begin{math}
    \norminfty{\seqstyle{a}}
    \leq \norminfty{\iota(\seqstyle{a})}
    \leq \ncover{f} \norminfty{f} \norminfty{\seqstyle{a}}    
  \end{math}
  for all \(\seqstyle{a} \in \gencospace{\zn}{K}\). Hence \(\iota\) 
  and \(\iota^{-1}\) are continuous.
  Thus (ii) is true.

  Since \(V\) and \(\gencospace{\zn}{K}\) are topologically 
  isomorphic and
  \(\gencospace{\zn}{K}\) is a Banach space it follows that \(V\) is
  also a Banach space and consequently a closed subspace of
  \(\vfn{n}{K}\). Thus (i) is true.
\end{proof}
}
\fi

\if\longversion1
\givelemmawithoutproof
\fi

\if\shortprep0
{
\begin{lemma}
  \label{lem:linftyspace-top-isomorphism}
  \definek
  Let \(c \in \positiverealnumbers\).
  Let \(f \in \cscfunc{\rn}{K}\)
  and
  \begin{math}
    f(\firstznvar + \rnb) = \delta_{\firstznvar,0}
  \end{math}
  for all \(\firstznvar \in \integernumbers^n\)
  for some \(\rnb \in \realnumbers^n\).
  Let
  \begin{eqnarray*}
    & & V := \left\{ \rnx \in \rn \mapsto
      \sum_{\firstznvar \in \zn} \seqelem{\sqa}{\firstznvar}
      f(c \rnx - \firstznvar)
      \setsep \seqstyle{a} \in \littlelp{\infty}{\zn}{K} \right\} \\
    & & \norminspace{g}{V} := \norminfty{g} .
  \end{eqnarray*}
  Define function \(\iota : \littlelp{\infty}{\zn}{K} \to V\) by
  \begin{displaymath}
    \iota(\seqstyle{a}) :=
    \rnx \in \rn \mapsto \sum_{\firstznvar \in \zn} 
    \seqelem{\sqa}{\firstznvar} f(c \rnx - \firstznvar)
  \end{displaymath}
  for all
  \begin{math}
    \seqstyle{a} \in \littlelp{\infty}{\zn}{K}
  \end{math}.
  Then
  \begin{itemize}
  \item[(i)]
    \(V\) is a closed subspace of \(\ucfunc{\rn}{K}\).
  \item[(ii)]
    Function \(\iota\) is a topological isomorphism from
    \(\littlelp{\infty}{\zn}{K}\) onto
    \(V\) and
    \begin{displaymath}
      \norminfty{\seqstyle{a}}
      \leq \norminfty{\iota(\seqstyle{a})}
      \leq \ncover{f} \norminfty{f} \norminfty{\seqstyle{a}}
    \end{displaymath}
    for all
    \begin{math}
      \seqstyle{a} \in \littlelp{\infty}{\zn}{K}
    \end{math}.
  \end{itemize}
\end{lemma}
}
\fi

\if\shortprep0
{
\begin{lemma}
  \label{lem:dual-gen-conv}
  \definekn
  Let
  \(E \equalns \ucfunc{\realnumbers^n}{K}\)
  or
  \(E \equalns \vanishingfunc{\realnumbers^n}{K}\).
  Let
  \(c \in \positiverealnumbers\),
  \(\dualelem{f} \in \topdual{E}\), and
  \begin{math}
    \sqd \in \littlelp{1}{\zn}{K}
  \end{math}.
  Then the series
  \begin{displaymath}
    \sum_{\firstznvar \in \integernumbers^n}
    \seqelem{\sqd}{\firstznvar} \dualelem{f} (c \cdot - \firstznvar)
  \end{displaymath}
  converges absolutely in
  \(\topdual{E}\)
  and
  \begin{displaymath}
    \sznorminspace{\sum_{\firstznvar \in \integernumbers^n}
    \seqelem{\sqd}{\firstznvar} \dualelem{f} (c \cdot
      - \firstznvar)}{\topdual{E}}
    \leq \frac{\norm{\dualelem{f}}}{c^n} 
    \normone{\sqd}
  \end{displaymath}
\end{lemma}
}
\fi

\if\longversion1
{
\begin{proof}
  Now
  \begin{displaymath}
    \norm{\dualelem{f} (c \cdot - \firstznvar)} = 
    \frac{1}{c^n}
    \norm{\dualelem{f}}
  \end{displaymath}
  for all
  \begin{math}
    \firstznvar \in \integernumbers^n
  \end{math}
  and
  \begin{displaymath}
    \sum_{\firstznvar \in \integernumbers^n}
    \norm{\seqelem{\sqd}{\firstznvar}
    \dualelem{f} (c \cdot - \firstznvar)}
    =
    \frac{\norm{\dualelem{f}}}{c^n}
    \sum_{\firstznvar \in \integernumbers^n}
    \abs{\seqelem{\sqd}{\firstznvar}}
  \end{displaymath}
  that converges because
  \(\sqd \in \littlelp{1}{\zn}{K}\).
  Consequently
  \begin{displaymath}
    \norm{\sum_{\firstznvar \in \integernumbers^n}
    \seqelem{\sqd}{\firstznvar} \dualelem{f} (c \cdot
      - \firstznvar)}
    \leq
    \sum_{\firstznvar \in \integernumbers^n}
    \norm{\seqelem{\sqd}{\firstznvar}
    \dualelem{f} (c \cdot - \firstznvar)}
    = \frac{\norm{\dualelem{f}}}{c^n} \normone{\sqd} .
  \end{displaymath}
\end{proof}
}
\fi

\if\shortprep0
{
Absolute convergence of a sequence in a Banach space implies unconditional
convergence.
}
\fi

\if\shortprep0
{
\begin{lemma}
  \label{lem:dual-norm}
  \definekn
  Let
  \(E \equalns \ucfunc{\realnumbers^n}{K}\)
  or
  \(E \equalns \vanishingfunc{\realnumbers^n}{K}\).
  Let \(c \in \positiverealnumbers\).
  Then
  \begin{displaymath}
    \sznorminspace{\sum_{\firstznvar \in \integernumbers^n}
      \seqelem{\sqd}{\firstznvar}
      \delta(c \cdot - \firstznvar)}{\topdual{E}}
    = \frac{1}{c^n} \normone{\sqd}
  \end{displaymath}
  for all \(\sqd \in \littlelp{1}{\zn}{K}\).
\end{lemma}
}
\fi

\if\shortprep0
{
\begin{proof}
  See also the proof of \cite[theorem 2.6]{cl1996}.
  By \mylemma \ref{lem:dual-gen-conv}
  \begin{displaymath}
    \sznorminspace{\sum_{\firstznvar \in \integernumbers^n}
      \seqelem{\sqd}{\firstznvar}
      \delta(c \cdot - \firstznvar)}
      {\topdual{E}}
    \leq \frac{1}{c^n} \normone{\sqd}
  \end{displaymath}
  and the series on the left-hand side converges absolutely.
  Let
  \begin{displaymath}
    h(x) :=
    \left\{
    \begin{array}{ll}
      4x + 1 ; & \textrm{if}\; x \in [-\frac{1}{4}, 0[ \\
      -4x + 1 ; & \textrm{if}\; x \in [0, \frac{1}{4}[ \\
      0 ; & \textrm{otherwise}
    \end{array}
    \right.
  \end{displaymath}
  where \(x \in \realnumbers\) and
  \(h^{[n]} : \realnumbers^n \to \realnumbers\),
  \begin{displaymath}
    h^{[n]} = \bigotimes_{k=1}^n h .
  \end{displaymath}
  Let \(\sigma : \naturalnumbers \to \integernumbers^n\) be a bijection.
  Suppose that
  \(\sqa \in  \gencospace{\zn}{K}\).
  By \mylemma \ref{lem:grid-convergence-multidim} series
  \begin{math}
    \sum_{\firstznvar \in \zn} \seqelem{\sqa}{\firstznvar}
    h(c \cdot - \firstznvar)
  \end{math}
  converges unconditionally in 
  \(\vanishingfunc{\realnumbers^n}{K}\).
  By \mylemma \ref{lem:dual-gen-conv} the series
  \begin{math}
    \sum_{\firstznvar \in \integernumbers^n}
    \seqelem{\sqd}{\firstznvar}
    \delta(c \cdot - \firstznvar)
  \end{math}
  converges absolutely.
\if\longversion1
{
  Furthermore,
  \begin{displaymath}
    \dualappl{\delta(c \cdot - \secondznvar)}{h^{[n]}(c \cdot -
      \firstznvar)}
    = \frac{1}{c^n} \szdualappl{\delta \left( \cdot -
        \frac{\secondznvar}{c} \right) }{h^{[n]} \left( c \cdot -
        \firstznvar \right)}
    = \delta_{\secondznvar,\firstznvar}
  \end{displaymath}
  for all
  \begin{math}
    \firstznvar, \secondznvar \in \zn
  \end{math}
  and
  \begin{displaymath}
    \szdualappl{\sum_{\secondznvar \in \zn}
    \seqelem{\sqd}{\secondznvar}
    \delta(c \cdot - \secondznvar)}
    {\sum_{\firstznvar \in \zn} \seqelem{\sqa}{\firstznvar}
    h^{[n]}(c \cdot - \firstznvar)}
    = \frac{1}{c^n} \sum_{\firstznvar \in \zn}
    \seqelem{\sqd}{\firstznvar}
    \seqelem{\sqa}{\firstznvar}
  \end{displaymath}
  and
  \begin{displaymath}
    \norminfty{\sum_{\firstznvar \in \zn}
    \seqelem{\sqa}{\firstznvar}
    h^{[n]}(c \cdot - \firstznvar)}
    = \norminfty{\sqa} .
  \end{displaymath}
  So
  \begin{eqnarray*}
    \frac{1}{c^n} \abs{\sum_{\firstznvar \in \zn}
    \seqelem{\sqd}{\firstznvar}
    \seqelem{\sqa}{\firstznvar}}
    & = &
    \abs{\szdualappl{\sum_{\secondznvar \in \zn}
    \seqelem{\sqd}{\secondznvar}
    \delta(c \cdot - \secondznvar)}
    {\sum_{\firstznvar \in \zn}
    \seqelem{\sqa}{\firstznvar}
    h^{[n]}(c \cdot - \firstznvar)}} \\
    & \leq & \sznorminspace{\sum_{\secondznvar \in \zn}
    \seqelem{\sqd}{\firstznvar}
    \delta(c \cdot - \secondznvar)}
    {\topdual{E}}
    \norminfty{\sqa} .
  \end{eqnarray*}
}
\else
{
  Furthermore,
  \begin{displaymath}
    \szdualappl{\sum_{\secondznvar \in \zn}
    \seqelem{\sqd}{\secondznvar}
    \delta(c \cdot - \secondznvar)}
    {\sum_{\firstznvar \in \zn}
    \seqelem{\sqa}{\firstznvar}
    h^{[n]}(c \cdot - \firstznvar)}
    = \frac{1}{c^n} \sum_{\firstznvar \in \zn}
    \seqelem{\sqd}{\firstznvar}
    \seqelem{\sqa}{\firstznvar}
  \end{displaymath}
  and
  \begin{displaymath}
    \norminfty{\sum_{\firstznvar \in \zn}
    \seqelem{\sqa}{\firstznvar}
    h^{[n]}(c \cdot - \firstznvar)}
    = \norminfty{\sqa} .
  \end{displaymath}
  So
  \begin{displaymath}
    \frac{1}{c^n} \abs{\sum_{\firstznvar \in \zn}
    \seqelem{\sqd}{\firstznvar}
    \seqelem{\sqa}{\firstznvar}}    
    \leq \norm{\sum_{\secondznvar \in \zn}
    \seqelem{\sqd}{\firstznvar}
    \delta(c \cdot - \secondznvar)}
    \norminfty{\sqa} .
  \end{displaymath}
}
\fi
  Hence
\if\elsevier1
{
  \begin{align*}
    \sznorminspace{\sum_{\firstznvar \in \integernumbers^n}
    \seqelem{\sqd}{\firstznvar} \delta(c \cdot - \firstznvar)}
    {\topdual{E}}
    & \geq
    \frac{1}{c^n} \sup \left\{
       \abs{\sum_{\firstznvar \in \zn}
       \seqelem{\sqd}{\firstznvar}
       \seqelem{\sqa}{\firstznvar}}
      \setsep \seqstyle{a} \in \gencospace{\zn}{K}
    \land \norminfty{\sqa} \leq 1 \right\} \\
    & = \frac{1}{c^n} \normone{\sqd}
  \end{align*}
}
\else
{
  \begin{eqnarray*}
    & & \sznorminspace{\sum_{\firstznvar \in \integernumbers^n}
    \seqelem{\sqd}{\firstznvar} \delta(c \cdot - \firstznvar)}
    {\topdual{E}} \\
    & & \geq
    \frac{1}{c^n} \sup \left\{
       \abs{\sum_{\firstznvar \in \zn}
       \seqelem{\sqd}{\firstznvar}
       \seqelem{\sqa}{\firstznvar}}
      \setsep \seqstyle{a} \in \gencospace{\zn}{K}
    \land \norminfty{\sqa} \leq 1 \right\} \\
    & & = \frac{1}{c^n} \normone{\sqd}
  \end{eqnarray*}
}
\fi  
  where the last equality follows from
  \begin{math}
    \topdual{\gencospace{\zn}{K}} \equalns \littlelp{1}{\zn}{K}
  \end{math}.
\end{proof}
}
\fi

\if\longprep1
{
\givelemmaswithoutproofsn{two}
}
\fi

\if\shortprep0
{
\begin{lemma}
  \label{lem:dual-appl-in-l-one}
  \definekn
  Let
  \(E \equalns \ucfunc{\realnumbers^n}{K}\)
  or
  \(E \equalns \vanishingfunc{\realnumbers^n}{K}\).
  Let
  \(c \in \positiverealnumbers\),
  \(\dualelem{f} \in \topdual{E}\),
  \(f \in \cscfunc{\rn}{K}\), and
  \begin{math}
     \seqelem{\sqd}{\firstznvar} =
     \dualappl{\dualelem{f}}
     {f(c \cdot - \firstznvar)}
  \end{math}
  for all \(\firstznvar \in \integernumbers^n\).
  Then
  \begin{math}
    \sqd \in \littlelp{1}{\zn}{K}
  \end{math}
  and
  \begin{math}
    \normone{\sqd} \leq \ncover{f} \norminfty{f} 
    \norm{\dualelem{f}}
  \end{math}.
\end{lemma}
}
\fi

\if\longversion1
{
\begin{proof}
  Suppose that
  \(\sqa \in \gencospace{\zn}{K}\).
  Let
  \begin{displaymath}
    g := \sum_{\firstznvar \in \zn}
    \seqelem{\sqa}{\firstznvar}
    f(c \cdot - \firstznvar) .
  \end{displaymath}
  By \mylemma \ref{lem:grid-convergence-multidim} we have
  \(g \in \vanishingfunc{\realnumbers^n}{K}\) and
  \(\norminfty{g} \leq \ncover{f} \norminfty{f} 
  \norminfty{\sqa}\).
  Now
  \begin{displaymath}
    \dualappl{\dualelem{f}}{g}
    = \sum_{\firstznvar \in \zn}
      \seqelem{\sqa}{\firstznvar}
      \dualappl{\dualelem{f}}{f(c \cdot - \firstznvar)}
    = \sum_{\firstznvar \in \zn}
      \seqelem{\sqa}{\firstznvar}
      \seqelem{\sqd}{\firstznvar}
    \in K
  \end{displaymath}
  where the series converge absolutely.
  By \mylemma \ref{lem:grid-convergence-multidim} we have
  \begin{equation}
    \label{eq:dual-appl-in-l-one-a}
    \abs{\sum_{\firstznvar \in \zn}
    \seqelem{\sqa}{\firstznvar}
    \seqelem{\sqd}{\firstznvar}}
    =\abs{\dualappl{\dualelem{f}}{g}}
    \leq \nsznorm{\dualelem{f}} \norm{g}
    \leq \ncover{f} \norminfty{f} \nsznorm{\dualelem{f}} 
    \norminfty{\sqa} .
  \end{equation}
  As \(\sqa \in \gencospace{\zn}{K}\) was arbitrary it follows
  that
  \begin{math}
    \sqd \in \topdual{\gencospace{\zn}{K}} \equalns
    \littlelp{1}{\zn}{K}
  \end{math}.
  It follows from \myequation \eqref{eq:dual-appl-in-l-one-a} that
  \begin{math}
    \normone{\sqd} \leq
    \ncover{f} \norminfty{f} \nsznorm{\dualelem{f}}
  \end{math}.
\end{proof}
}
\fi

\if\shortprep0
{
\begin{lemma}
  \label{lem:mdim-vanishing-samples}
  \definekn
  Let
  \(c \in \positiverealnumbers\) and
  \(f \in \vanishingfunc{\realnumbers^n}{K}\).
  Then
  \begin{math}
    \left( f \left( c \firstznvar \right) \right)_{\firstznvar \in \zn}
    \in \gencospace{\zn}{K}
  \end{math}
  and
  \begin{math}
    \norminfty{\left( f \left( c \firstznvar \right)
      \right)_{\firstznvar \in \zn}} \leq \norminfty{f}
  \end{math}.
\end{lemma}
}
\fi

\if\longversion1
{
\begin{proof}
  Let \(\sigma : \naturalnumbers \onto \zn\)
  be a bijection.
  Let \(h > 0\).
  There exists \(r > 0\) so that
  \begin{math}
    \abs{f(\rnx)} < h
  \end{math}
  for all
  \begin{math}
    \rnx \in \rn
  \end{math},
  \begin{math}
    \norm{x} > r
  \end{math}.
  The set
  \begin{math}
    \left\{ \firstznvar \in \integernumbers^n \setsep
      \normtwo{c \firstznvar} \leq r \right\}
  \end{math}
  is finite and hence there exists \(l_0 \in \naturalnumbers\) so 
  that
  \(\normtwo{c \sigma(l)} > r\) for all \(l > l_0\). Now
  \begin{math}
     \abs{f (c \sigma(l))} < h
  \end{math}
  for all
  \begin{math}
    l > l_0
  \end{math}.
  Consequently
  \begin{math}
    \left( f \left( c \firstznvar \right) \right)_{\firstznvar \in 
    \zn}
    \in \gencospace{\zn}{K}
  \end{math}.
  It follows from the definition of the supremum norm that
  \begin{math}
    \norminfty{\left( f \left( c \firstznvar \right)
      \right)_{\firstznvar \in \zn}} \leq \norminfty{f}
  \end{math}.  
\end{proof}
}
\fi

\if\shortprep0
{
\begin{lemma}
  \label{lem:dual-single-res-norm-equiv-multidim}
  \definekn
  Let
  \(E \equalns \ucfunc{\realnumbers^n}{K}\)
  or
  \(E \equalns \vanishingfunc{\realnumbers^n}{K}\).
  Let \(c \in \positiverealnumbers\),
  \(\dualelem{f} \in \topdual{E}\),
  and \(f \in \cscfunc{\rn}{K}\).
  Suppose that
  \if\shortprep1
  {
  \begin{displaymath}
    \forall \firstznvar \in \integernumbers^n :
    \dualappl{\dualelem{f}}{f(\cdot - \firstznvar)}
    = \delta_{\firstznvar,0} .
  \end{displaymath}
  }
  \else
  {
  \begin{equation}
    \label{eq:card-b}
    \forall \firstznvar \in \integernumbers^n :
    \dualappl{\dualelem{f}}{f(\cdot - \firstznvar)}
    = \delta_{\firstznvar,0} .
  \end{equation}
  }
  \fi
  Then
  \begin{displaymath}
    \frac{1}{c^n} \frac{1}{\ncover{f} \norminfty{f}} 
    \normone{\sqd}
    \leq \sznorminspace{\sum_{\firstznvar \in \integernumbers^n}
    \seqelem{\sqd}{\firstznvar} \dualelem{f}(c \cdot - \firstznvar)}
    {\topdual{E}}
    \leq \frac{1}{c^n} \norm{\dualelem{f}} \normone{\sqd}
  \end{displaymath}
  for all
  \begin{math}
    \seqstyle{d} \in \littlelp{1}{\zn}{K}
  \end{math}.
\end{lemma}
}
\fi

\if\shortprep0
{
\begin{proof}
  Let
  \begin{math}
    \sqd \in \littlelp{1}{\zn}{K}
  \end{math}.
  Let
  \begin{displaymath}
    \dualelem{g} :=
    \sum_{\firstznvar \in \integernumbers^n}
    \seqelem{\sqd}{\firstznvar}
    \dualelem{f}(c \cdot - \firstznvar) .
  \end{displaymath}
  By \mylemma \ref{lem:dual-gen-conv} the series on the right-hand side
  converges absolutely and
  \begin{displaymath}
    \norminspace{\dualelem{g}}{\topdual{E}} \leq \frac{1}{c^n}
    \norminspace{\dualelem{f}}{\topdual{E}} 
    \normone{\sqd} .
  \end{displaymath}
  Suppose that
  \(\sqa \in \gencospace{\zn}{K}\).
  Let
  \begin{displaymath}
    g := \sum_{\firstznvar \in \zn} \seqelem{\sqa}{\firstznvar}
    f(c \cdot - \firstznvar)
  \end{displaymath}
  By \mylemma \ref{lem:grid-convergence-multidim} the series on the
  right-hand side converges unconditionally and
  \begin{math}
    \norminfty{g} \leq \ncover{f} \norminfty{f} 
    \norminfty{\sqa}    
  \end{math}.
\if\longversion1
{
  By \mydef \ref{def:distr-translation-and-dilatation} we have
  \begin{displaymath}
    \dualappl{\dualelem{f}(c \cdot - \firstznvar)}
    {f(c \cdot - \secondznvar)}
    = \frac{1}{c^n} \dualappl{\dualelem{f}}
    {f(\cdot + \firstznvar - \secondznvar)}
  \end{displaymath}
  and hence by \myequation \eqref{eq:card-b}
  \begin{displaymath}
    \dualappl{\dualelem{g}}{g}
      = \frac{1}{c^n} \sum_{\firstznvar \in \integernumbers^n}
      \sum_{\secondznvar \in \zn}
      \seqelem{\sqd}{\firstznvar}
      \seqelem{\sqa}{\secondznvar}
      \delta_{\firstznvar, \secondznvar}
      = \frac{1}{c^n} \sum_{\firstznvar \in \zn}
      \seqelem{\sqd}{\firstznvar}
      \seqelem{\sqa}{\firstznvar}
  \end{displaymath}
  and it follows that
  \begin{displaymath}
    \frac{1}{c^n} \abs{\sum_{\firstznvar \in \zn}
      \seqelem{\sqd}{\firstznvar}
      \seqelem{\sqa}{\firstznvar}}
    \leq \norm{\dualelem{g}} \norminfty{g}
    \leq \norm{\dualelem{g}} \ncover{f} \norminfty{f} 
    \norminfty{\sqa} .
  \end{displaymath}
}
\else
{
  By \mydef \ref{def:distr-translation-and-dilatation}
  and \myequation \eqref{eq:card-b} we have
  \begin{displaymath}
    \dualappl{\dualelem{g}}{g}
      = \frac{1}{c^n} \sum_{\firstznvar \in \zn}
      \seqelem{\sqd}{\firstznvar}
      \seqelem{\sqa}{\firstznvar}
  \end{displaymath}
  and it follows that
  \begin{displaymath}
    \frac{1}{c^n} \abs{\sum_{\firstznvar \in \zn}
      \seqelem{\sqd}{\firstznvar}
      \seqelem{\sqa}{\firstznvar}}
    \leq \norm{\dualelem{g}} \norminfty{g}
    \leq \norm{\dualelem{g}} \ncover{f} \norminfty{f} 
    \norminfty{\seqstyle{a}} .
  \end{displaymath}
}
\fi
  Thus
  \begin{equation}
    \label{eq:dual-single-res-inequality-a}
    \norm{\dualelem{g}}
    \geq
    \frac{1}{c^n} \frac{1}{\ncover{f} \norminfty{f}
      \norminfty{\seqstyle{a}}}
    \abs{\sum_{\firstznvar \in \zn}
      \seqelem{\sqd}{\firstznvar}
      \seqelem{\sqa}{\firstznvar}} .
  \end{equation}
  Let \(\sigma : \naturalnumbers \onto \zn\) be a bijection.
  When \(m \in \naturalnumbers\) define sequence 
  \(\indexedseqstyle{a}{m} \in \gencospace{\zn}{K}\) by
  \begin{displaymath}
    \seqelem{\indexedseqstyle{a}{m}}{\firstznvar} := \left\{
      \begin{array}{ll}
        e^{-i \arg \seqelem{\sqd}{\firstznvar}}
        ;
        & \textrm{if}\; \sigma^{-1}(\firstznvar) \leq m \\
        0 ; & \textrm{otherwise.}
      \end{array}
    \right.
  \end{displaymath}
  where \(\firstznvar \in \zn\).
  Sequence \(\seqstyle{a} \in \gencospace{\zn}{K}\) was arbitrary 
  and it follows from
  \eqref{eq:dual-single-res-inequality-a} that
  \begin{eqnarray*}
    \norm{\dualelem{g}}
    & \geq &
    \frac{1}{c^n} \frac{1}{\ncover{f} \norminfty{f}
      \norminfty{\indexedseqstyle{a}{m}}}
    \abs{\sum_{\firstznvar \in \zn}
    \seqelem{\sqd}{\firstznvar}
    \seqelem{\indexedseqstyle{a}{m}}{\firstznvar}} \\
    & = &
    \frac{1}{c^n} \frac{1}{\ncover{f} \norminfty{f}
      \norminfty{\indexedseqstyle{a}{m}}}
    \sum_{u=0}^m \abs{\seqelem{\sqd}{\sigma(u)}}
  \end{eqnarray*}
  for all \(m \in \naturalnumbers\).
\if\longversion0
{
  Now \(\norminfty{\indexedseqstyle{a}{m}} \leq 1\) for all
  \(m \in \naturalnumbers\) and consequently
  \begin{displaymath}
    \norminspace{\dualelem{g}}{\topdual{E}}
    \geq
    \frac{1}{c^n} \frac{1}{\ncover{f} \norminfty{f}}
    \sum_{\firstznvar \in \zn} \abs{\seqelem{\sqd}{\firstznvar}}
    =
    \frac{1}{c^n} \frac{1}{\ncover{f} \norminfty{f}}
    \normone{\sqd} .
  \end{displaymath}
}
\else
{
  Now \(\norminfty{\indexedseqstyle{a}{m}} \leq 1\) for all
  \(m \in \naturalnumbers\) and hence
  \begin{displaymath}
    \norminspace{\dualelem{g}}{\topdual{E}}
    \geq
    \frac{1}{c^n} \frac{1}{\ncover{f} \norminfty{f}}
    \sum_{u=0}^m \abs{\seqelem{\sqd}{\sigma(u)}}
  \end{displaymath}
  for all \(m \in \naturalnumbers\).
  Consequently
  \begin{displaymath}
    \norminspace{\dualelem{g}}{\topdual{E}}
    \geq
    \frac{1}{c^n} \frac{1}{\ncover{f} \norminfty{f}}
    \sum_{\firstznvar \in \zn} \abs{\seqelem{\sqd}{\firstznvar}}
    =
    \frac{1}{c^n} \frac{1}{\ncover{f} \norminfty{f}}
    \normone{\sqd} .
  \end{displaymath}
}
\fi
\end{proof}
}
\fi

\if\longprep1
{
\givelemmaswithoutproofsn{six}
}
\fi

\if\shortprep0
{
\begin{lemma}
  \label{lem:multidim-dual-isomorphism}
  \definekn
  Let
  \(E \equalns \ucfunc{\realnumbers^n}{K}\)
  or
  \(E \equalns \vanishingfunc{\realnumbers^n}{K}\).
  Let \(c \in \positiverealnumbers\),
  \(\dualelem{f} \in \topdual{E}\),
  and \(f \in \cscfunc{\rn}{K}\).
  Suppose that
  \begin{displaymath}
    \dualappl{\dualelem{f}}{f(\cdot - \firstznvar)}
    = \delta_{\firstznvar,0}
  \end{displaymath}
  for all
  \begin{math}
    \firstznvar \in \integernumbers^n
  \end{math}.
  Define
  \if\longversion1
  {
  \begin{eqnarray}
    \nonumber
    & & M := \left\{ c^n \sum_{\firstznvar \in \integernumbers^n}
      \seqelem{\sqd}{\firstznvar} \dualelem{f}
      \left( c \cdot - \firstznvar \right) \setsep
      \sqd \in
      \littlelp{1}{\zn}{K}
    \right\} \\
    \label{eq:m-space-a}
    & &
    \norminspace{\dualelem{g}}{M} := \norminspace{\dualelem{g}}{\topdual{E}}
    .
  \end{eqnarray}
  }
  \else
  {
  \begin{eqnarray*}
    & & M := \left\{ c^n \sum_{\firstznvar \in \integernumbers^n}
      \seqelem{\sqd}{\firstznvar} \dualelem{f}
      \left( c \cdot - \firstznvar \right) \setsep
      \sqd \in
      \littlelp{1}{\zn}{K}
    \right\} \\
    & &
    \norminspace{\dualelem{g}}{M} := \norminspace{\dualelem{g}}{\topdual{E}}
    .
  \end{eqnarray*}
  }
  \fi
  Define function \(\iota : \littlelp{1}{\zn}{K} \to M\) by
  \begin{displaymath}
    \iota(\sqd) :=
     c^n \sum_{\firstznvar \in \integernumbers^n}
    \seqelem{\sqd}{\firstznvar}
    \dualelem{f} \left( c \cdot - \firstznvar \right)
  \end{displaymath}
  for all
  \begin{math}
    \sqd \in \littlelp{1}{\zn}{K}
  \end{math}.
  Then
  \begin{itemize}
  \item[(i)]
    \(M\) is a closed subspace of \(\topdual{E}\).
  \item[(ii)]
    Function \(\iota\) is a topological isomorphism from
    \(\littlelp{1}{\zn}{K}\) onto \(M\)
    and
    \begin{equation}
      \label{eq:dual-norm-equiv}
      \frac{1}{\ncover{f} \norminfty{f}} \normone{\sqd}
      \leq \norm{\iota(\sqd)}
      \leq \norm{\dualelem{f}} \normone{\sqd}
    \end{equation}
    for all
    \begin{math}
      \sqd \in \littlelp{1}{\zn}{K}
    \end{math}.
  \end{itemize}
\end{lemma}
}
\fi

\if\longversion1
{
\begin{proof}
  The series in \myequation \eqref{eq:m-space-a} converges absolutely by
  \mylemma \ref{lem:dual-gen-conv}.
  Function \(\iota\) is linear and a surjection onto \(M\).
  By \mylemma \ref{lem:dual-single-res-norm-equiv-multidim}
  \myequation \eqref{eq:dual-norm-equiv} is true.
  
  Suppose that
  \begin{math}
    \seqstyle{w}, \seqstyle{z} \in \littlelp{1}{\zn}{K}
  \end{math}
  and \(\seqstyle{w} \neq \seqstyle{z}\).
  Now
  \begin{math}
    \normone{\seqstyle{w} - \seqstyle{z}} \neq 0
  \end{math}
  so
  \begin{displaymath}
    \norm{\iota(\seqstyle{w} - \seqstyle{z})}
    \geq
    \frac{1}{\ncover{f} \norminfty{f}}
    \normone{\seqstyle{w} - \seqstyle{z}}
    > 0 .
  \end{displaymath}
  It follows that
  \begin{math}
    \iota(\seqstyle{w} - \seqstyle{z}) \neq 0
  \end{math}.
  Hence \(\iota\) is an injection.
  It follows that \(\iota\) is a linear bijection
  and vector spaces \(M\) and \(\littlelp{1}{\zn}{K}\) are
  algebraically isomorphic.
  By \myequation \eqref{eq:dual-norm-equiv}
  \(\iota\) and \(\iota^{-1}\) are continuous. Thus (ii) is true.
  
  Since \(M\) and \(\littlelp{1}{\zn}{K}\) are topologically
  isomorphic and \(\littlelp{1}{\zn}{K}\) is a Banach space it
  follows that \(M\) is also a Banach space and consequently a 
  closed subspace of \(\topdual{E}\). Thus (i) is true.
\end{proof}
}
\fi

\if\shortprep0
{
\begin{lemma}
  \label{lem:multidim-delta-isomorphism}
  \definekn
  Let
  \(E \equalns \ucfunc{\realnumbers^n}{K}\)
  or
  \(E \equalns \vanishingfunc{\realnumbers^n}{K}\).
  Let \(c \in \positiverealnumbers\).
  Define
  \if\longversion1
  {
  \begin{eqnarray}
    \nonumber
    & & M := \left\{ \sum_{\firstznvar \in \integernumbers^n}
      \seqelem{\sqd}{\firstznvar} \delta
      \left( \cdot - \frac{\firstznvar}{c} \right) \setsep
      \sqd \in
      \littlelp{1}{\zn}{K}
    \right\} \\
    \label{eq:m-space-b}
    & &
    \norminspace{\dualelem{g}}{M} := \norminspace{\dualelem{g}}{\topdual{E}}
  \end{eqnarray}
  }
  \else
  {
  \begin{eqnarray*}
    & & M := \left\{ \sum_{\firstznvar \in \integernumbers^n}
      \seqelem{\sqd}{\firstznvar} \delta
      \left( \cdot - \frac{\firstznvar}{c} \right) \setsep
      \sqd \in
      \littlelp{1}{\zn}{K}
    \right\} \\
    & &
    \norminspace{\dualelem{g}}{M} := \norminspace{\dualelem{g}}{\topdual{E}}
  \end{eqnarray*}
  }
  \fi
  Define function \(\iota : \littlelp{1}{\zn}{K} \to M\) by
  \begin{displaymath}
    \iota(\sqd) := \sum_{\firstznvar \in \integernumbers^n}
    \seqelem{\sqd}{\firstznvar}
    \delta \left( \cdot - \frac{\firstznvar}{c} \right)
  \end{displaymath}
  for all
  \begin{math}
    \sqd \in \littlelp{1}{\zn}{K}
  \end{math}.
  Then \(M\) is a closed subspace of
  \(\topdual{E}\) and function
  \(\iota\) is an isometric isomorphism from
  \(\littlelp{1}{\zn}{K}\) onto \(M\).
\end{lemma}
}
\fi

\if\longversion1
{
\begin{proof}
  See also \cite[theorem 2.6]{goedecker1998}.
  The series in \myequation \eqref{eq:m-space-b} converges absolutely by
  \mylemma \ref{lem:dual-gen-conv}. Function \(\iota\) is linear and a
  surjection onto \(M\). By \mylemma \ref{lem:dual-norm}
  \(\norm{\iota(\sqd)} = \normone{\sqd}\)
  for all \(\sqd \in \littlelp{1}{\zn}{K}\). It follows
  that \(\iota\) is injective.
  Hence \(\iota\) is an isometric isomorphism from
  \(\littlelp{1}{\zn}{K}\) onto \(M\) and \(M\) is a closed
  subspace of \(\topdual{E}\).
\end{proof}
}
\fi

\if\longversion1
{
\givelemmaswithoutproofsn{five}
}
\else
{
  See also \cite[theorem 2.6]{cl1996}.
}
\fi

\if\longversion1
{
\begin{lemma}
  \label{lem:double-series}
  Let \(m \in \naturalnumbers\).
  Let \(\sqa \in \seqset{\setzeroton{m} \times \naturalnumbers}
    {\nonnegrealnumbers}\)
  and
  \begin{math}
    \indexedseqstyle{b}{k} :=
    (\seqelem{\sqa}{k,l})_{l=0}^\infty
  \end{math}
  for each \(k \in \setzeroton{m}\).
  Then
  \begin{itemize}
    \item[(i)]
      Series
      \begin{displaymath}
        \sum_{\lambda \in \setzeroton{m} \times \naturalnumbers}
        \seqelem{\sqa}{\lambda}
      \end{displaymath}
      converges if and only if the series
      \begin{displaymath}
        \sum_{l \in \naturalnumbers}
        \seqelem{\indexedseqstyle{b}{k}}{l}
      \end{displaymath}
      converges for all \(k \in \setzeroton{m}\).
    \item[(ii)]
      \begin{displaymath}
        \sum_{\lambda \in \setzeroton{m} \times \naturalnumbers}
        \seqelem{\sqa}{\lambda}
			  =
			  \sum_{k \in \setzeroton{m}}
			  \sum_{l \in \naturalnumbers}
			  \seqelem{\indexedseqstyle{b}{k}}{l}
      \end{displaymath}
    \item[(iii)]
      When \(p \in [1, \infty]\)
      \begin{displaymath}
        \norm{\sqa}_p = \norm{\left(
          \norm{\indexedseqstyle{b}{k}}_p \right)_{k \in 
          \setzeroton{m}}}_p .
      \end{displaymath}
  \end{itemize}
\end{lemma}
}
\fi

\if\shortprep0
{
\begin{lemma}
  \label{lem:change-int-and-sum-order}
  Let \(n \in \positiveintegers\),
  \(\sqa \in l^1\),
  \(v_k \in \cscfunccv{\rn}\) for each \(k \in \naturalnumbers\),
  and
  \begin{equation}
    \label{eq:series-a}
    f(\rnx) := \sum_{k \in \naturalnumbers}
    \seqelem{\sqa}{k} v_k(\rnx)
  \end{equation}
  for all \(\rnx \in \rn\).
  Suppose also that
  \begin{itemize}
    \item[(i)]
      The series in \eqref{eq:series-a}
      converges absolutely for each \(\rnx \in \rn\).
    \item[(ii)]
      There exists \(c_1 \in \positiverealnumbers\)
      so that
      \begin{math}
        \normone{v_k} \leq c_1
      \end{math}
      for all \(k \in \naturalnumbers\).
  \end{itemize}
  Then
  \begin{displaymath}
    \int_{\rnx \in \rn} \left(
      \sum_{k \in \naturalnumbers}
      \abs{\seqelem{\sqa}{k}}
      \abs{v_k(\rnx)}
    \right)
    d \tau
    =
    \sum_{k \in \naturalnumbers}
    \abs{\seqelem{\sqa}{k}}
    \int_{\rnx \in \rn}    
    \abs{v_k(\rnx)}
    d \tau .
  \end{displaymath}
\end{lemma}
}
\fi

\if\shortprep0
{
\begin{lemma}
  \label{lem:uc-convergence}
  Let \(n \in \positiveintegers\),
  \(f \in \cscfunccv{\rn}\),
  \(c \in \positiverealnumbers\),
  and
  \(\sqa \in \littlelpcv{\infty}{\zn}\).
  Then the series
  \begin{displaymath}
    g(\rnx) := \sum_{\firstznvar \in \zn}
    \seqelem{\sqa}{\firstznvar}
    f(c \rnx - \firstznvar)
  \end{displaymath}
  converges absolutely for each \(\rnx \in \rn\).
  We also have
  \begin{math}
    g \in \ucfunccv{\rn}
  \end{math}.
\end{lemma}
}
\fi

\if\shortprep0
{
\begin{lemma}
  \label{lem:delta-and-pointwise-conv}
  Let \(n \in \positiveintegers\) and
  \(I\) be a countably infinite set.
  Let
  \begin{math}
    v_\alpha \in \cscfunccv{\rn}
  \end{math}
  and
  \begin{math}
    \indexeddualelem{v}{\alpha} \in \cscfunccv{\rn}
  \end{math}
  for all \(\alpha \in I\).
  When \(\alpha \in I\) define
  \begin{math}
    J(\alpha) := \{ \beta \in I
    \setsep
    \suppop \indexeddualelem{v}{\alpha}
    \intersection
    \suppop v_\beta \neq \emptyset
    \}
  \end{math}.
  Assume that \(J(\alpha)\) is finite
  for each \(\alpha \in I\).
  Let \(\sqa \in \littlelpcv{\infty}{I}\) and
  define
  \begin{equation}
    \label{eq:def-f}
    f(\rnx) := \sum_{\alpha \in I}
    \seqelem{\sqa}{\alpha} v_\alpha (\rnx)
  \end{equation}
  for all \(\rnx \in \rn\).
  Assume that the series in \myequation \eqref{eq:def-f}
  converges absolutely for each \(\rnx \in \rn\).
  Then
  \begin{displaymath}
    \dualappl{\indexeddualelem{v}{\alpha}}{f}
    =
    \sum_{\beta \in I}
    \seqelem{\sqa}{\beta}
    \dualappl{\indexeddualelem{v}{\alpha}}{v_\beta} .
  \end{displaymath}
\end{lemma}
}
\fi

\if\shortprep0
{
\begin{lemma}
  \label{lem:conv-in-different-norms}
  Let \(n \in \positiveintegers\), \(p \in [1, \infty[\),
  \(x_k \in \biglpcv{p}{\rn} \intersection \biglpcv{\infty}{\rn}\),
  \(y \in \biglpcv{p}{\rn}\), and
  \(z \in \biglpcv{\infty}{\rn}\).
  Suppose that
  \begin{math}
    \norm{x_k - y}_{\biglpcv{p}{\rn}} \to 0
  \end{math}
  as \(k \to \infty\)
  and
  \begin{math}
    \norm{x_k - z}_{\biglpcv{\infty}{\rn}} \to 0
  \end{math}
  as \(k \to \infty\).
  Then \(y = z\) almost everywhere.
\end{lemma}
}
\fi

\if\shortprep0
{
\begin{lemma}
  \label{lem:lp-oper-interp}
  Let \((R, \mu)\) and \((S, \nu)\) be totally
  \(\sigma\)-finite measure spaces.
  Suppose that
  \(A_1 \closedsubspace L^1(R, \mu)\),
  \(A_\infty \closedsubspace L^\infty(R, \mu)\),
  \(T_1\) is an operator from \(A_1\) into
  \(L^1(S,\nu)\), and
  \(T_\infty\) is an operator from \(A_\infty\) into
  \(L^\infty(S,\nu)\).
  When \(p \in ] 1, \infty [\) define
  \begin{math}
    T_p := \kinterp{1-\frac{1}{p}}{p}{T_1}{T_\infty}
  \end{math},
  \begin{math}
    A_p \defequalset \kinterp{1-\frac{1}{p}}{p}{A_1}{A_\infty}
  \end{math}, and
  \begin{math}
    \norminspace{f}{A_p} := \norminspace{f}{L^p(R, \mu)}
  \end{math}
  for all \(f \in A_p\).
  Then
  \begin{displaymath}
    \norminspace{T_p}{\contlinop{A_p}{L^p(S, \nu)}}
    \leq
    \frac{p}{p-1}
    \norm{T_1}^{\frac{1}{p}}
    \norm{T_\infty}^{1-\frac{1}{p}}
  \end{displaymath}
  for all \(p \in ] 1, \infty [\).
\end{lemma}
}
\fi

\if\shortprep0
{
\begin{proof}
  Use definitions IV.4.4 and V.1.7,
  lemma IV.4.5, and
  theorems V.1.6 and V.1.12
  in \cite{bs1988}.
\end{proof}
}
\fi

\if10
{
\givelemmawithoutproof

\begin{lemma}
  \label{lem:delta-lin-comb}
  Let \(n, m \in \positiveintegers\).
  Suppose that
  \(a_k \in \complexnumbers\),
  \(\indexedseqstyle{b}{k} \in \rn\),
  and
  \(\indexedseqstyle{b}{k} \neq
  \indexedseqstyle{b}{l}\)
  for all \(k, l \in \setoneton{m}\)
  and \(k \neq l\).
  Define
  \begin{displaymath}
    \dualelem{f} :=
    \sum_{k = 1}^m a_k \delta(\cdot - \indexedseqstyle{b}{k}) . 
  \end{displaymath}
  Then
  \begin{displaymath}
    \norminspace{\dualelem{f}}{\topdual{\vanishingfunccv{\rn}}}
    =
    \norminspace{\dualelem{f}}{\topdual{\ucfunccv{\rn}}} .
  \end{displaymath}
\end{lemma}
}
\fi

\if\shortprep0
{
\begin{lemma}
  Let \(f \in \complexnumbers\),
  \(m \in \naturalnumbers\),
  \(r_1 \in \positiverealnumbers\),
  and \(\dualf \in 
  \topdual{\contfunccv{\closedball{\realnumbers}{0}{r_1}}}\).
  Let
  \begin{displaymath}
    f^{[n]} := \indexedtensorproduct{i=0}{n} f
  \end{displaymath}
  and
  \begin{displaymath}
    \dualf^{[n]} := \indexedtensorproduct{i=0}{n} \dualf .
  \end{displaymath}
  If \((\dualf, f)\) spans all polynomials of degree at
  most \(m\) 
  then \((\dualf^{[n]}, f^{[n]})\) spans all polynomials of degree
  at most \(m\).
\end{lemma}
}
\fi

\if\shortprep0
{
\begin{proof}
  Let
  \begin{displaymath}
    p(\rnx) :=
    \prod_{i=1}^n p_i(\seqelem{\rnx}{i})
    =
    \prod_{i=1}^n \left( \seqelem{\rnx}{i} \right)^{\seqelem{\sqa}{i}}
  \end{displaymath}
  where
  \begin{math}
    \seqelem{\sqa}{i} \in \setzeroton{m}
  \end{math},
  \begin{math}
    i = 1, \ldots, n
  \end{math}.
  Define \(r_2 \in \positiverealnumbers\) so that
  \begin{math}
    \left( \closedball{\realnumbers}{0}{r_1} \right)^ n \subset
    \closedball{\rn}{0}{r_2}
  \end{math}.
  Now
  \begin{eqnarray*}
    & & \sum_{\firstznvar \in \zn}
    \szdualappl{\dualf^{[n]}(\cdot -
      \firstznvar)}{\restrictfunc{p}{\closedball{\rn}{\firstznvar}{r_2}}}
    f^{[n]}(\rnx-\firstznvar) \\
    & = &
    \sum_{k_1 \in \integernumbers} \cdots \sum_{k_n \in
      \integernumbers}
    \prod_{i=1}^n 
    \szdualappl{\dualf(\cdot -
      \seqelem{\firstznvar}{i})}{\restrictfunc{p_i}
      {\closedball{\realnumbers}{\seqelem{\firstznvar}{i}}{r_1}}}
    f\left( \seqelem{\rnx}{i} - \seqelem{\firstznvar}{i} \right)
    \\
    & = &
    \prod_{i=1}^n \sum_{k \in \integernumbers}
    \szdualappl{\dualf(\cdot-k)}
               {\restrictfunc{p_i}
                 {\closedball{\realnumbers}{k}{r_1}}}
    f\left( \seqelem{\rnx}{i} - k \right)
    \\
    & = &
    \prod_{i=1}^n p_i(\seqelem{\rnx}{i}) = p(\rnx)
  \end{eqnarray*}
  for all
  \begin{math}
    \rnx \in \rn
  \end{math}.
\end{proof}
}
\fi

\if01
{
\begin{definition}
  Let \(n \in \positiveintegers\).
  Let \(\zna, \znb \in \zn\).
  We say that \(\zna\) and \(\znb\) are
  \defterm{neighbours} iff
  \begin{math}
    \norminfty{\zna-\znb} = 1
  \end{math}.
\end{definition}
}
\fi

\if\shortprep0
{
\begin{definition}
  Let \(n \in \positiveintegers\), \(I \subset \naturalnumbers\),
  \(I \neq \emptyset\).
  Let \(\alpha \in (\zn)^I\).
  We say that \(\alpha\)
  \defterm{preserves neighbours}
  iff \(k, k+1 \in I\) implies
  \begin{math}
    \nsznorminfty{\alpha(k+1)-\alpha(k)} = 1
  \end{math}.
\end{definition}
}
\fi

\if\shortprep0
{
\begin{definition}
  When \(n \in \positiveintegers\)
  define
  \begin{math}
    \zerocentredcube{n}{k}
    :=
    \left(
      \setplusminusn{k}
    \right)^n
  \end{math}
  for all \(k \in \naturalnumbers\)
  and
  \begin{math}
    \zerocentredcubediff{n}{k}
    :=
    \zerocentredcube{n}{k}
    \setminus
    \zerocentredcube{n}{k-1}
  \end{math}  
  for all \(k \in \positiveintegers\).
  Define also
  \begin{math}
    \zerocentredcubediff{n}{0} := \{\finitezeroseq{n}\}
  \end{math}.
\end{definition}
}
\fi

\if\shortprep0
{
\begin{definition}
  \label{def:rectangle-ordering}
  Let \(k \in \positiveintegers\). Define function
  \begin{math}
    \cubediffordfunction{2}{k} :
    \setzeroton{\card{\zerocentredcubediff{n}{k}}-1} \to
    \zerocentredcubediff{n}{k}
  \end{math}
  by
  \begin{displaymath}
    \cubedifford{2}{k}{m} :=
    \left\{
    \begin{array}{ll}
      (-k+1+c,-k) ; & m = c \in \setzeroton{2k-1} \\
      (k,-k+c+1) ; & m = 2k+c, c \in \setzeroton{2k-1} \\
      (-k+1+c,-k) ; & m = 4k+c, c \in \setzeroton{2k-1} \\
      (-k+1+c,-k) ; & m = 6k+c, c \in \setzeroton{2k-1}
    \end{array}
    \right.
  \end{displaymath}
  Define also
  \begin{math}
    \cubediffordfunction{2}{0} :
    \setzeroton{\card{\zerocentredcubediff{n}{0}}-1} \to
    \zerocentredcubediff{n}{0}
  \end{math}
  by
  \begin{math}
    \cubedifford{2}{0}{0} := (0,0)
  \end{math}.
\end{definition}
}
\fi

\if\shortprep0
{
  Define
  $
  \alpha_{2,k}(m) := \beta_{2,j}(m-\card{A_{2,j-1}})
  $
  where $\card{A_{2,j-1}} \leq m < \card{A_{2,j}}$ and
  $m \in \setzeroton{\card{A_{2,k}}-1}$, $j \in \setzeroton{k}$ and
  $k \in \positiveintegers$. When $k=0$ let $\alpha_{2,0}(0) = (0,0)$.
  Let also
  $\alpha'_{2,k}(m) := \alpha_{2,k}(\card{A_{2,k}}-1-m)$ for 
  $m \in \setzeroton{\card{A_{2,k}}-1}$.
  
  Now \(\cubediffordfunction{2}{k}\) is a bijection from
  \(\setzeroton{\card{\zerocentredcubediff{n}{k}}-1}\) onto
  \(\zerocentredcubediff{n}{k}\) and it preserves neighbours for each
  \(k \in \naturalnumbers\). Functions \(\constrcubeordfunction{2}{k}\)
  and \(\rconstrcubeordfunction{2}{k}\)
  are bijections from \(\setzeroton{\card{\zerocentredcube{n}{k}}-1}\)
  onto \(\zerocentredcube{n}{k}\) and they preserve neighbours for
  each \(k \in \naturalnumbers\).
  
  If $k < j$ we define $\sum_{i=j}^k x_i = 0$.
}
\fi

\if\shortprep0
{
  \begin{lemma}
    \label{lem:ordering-bijections}
    Let $n \in \naturalnumbers$, $n \geq 3$. There exist bijections
    $\beta_{n,k} : \setzeroton{\card{B_{n,k}}-1} \onto B_{n,k}$,
    $k \in \naturalnumbers$, so that $\beta_{n,k}$ preserves neighbours,
    $\beta_{n,k}(0) = (k, \finitezeroseq{n-1})$, and
    $\beta_{n,k}(\card{B_{n,k}}-1) = (-k, \finitezeroseq{n-1})$ for all
    $k \in \naturalnumbers$.
  \end{lemma}
}
\fi

\if\shortprep0
{
  \begin{proof}
    We prove case \(n = 3\) first. Suppose that $k \in \positiveintegers$.
    Let $$C_0 := \{(k,\mathbf{x}) \setsep \mathbf{x} \in A_{2,k}\}$$ and
    $$p_0(m) := (k,\alpha_{2,k}(m))$$ for $m \in
    \setzeroton{\card{A_{2,k}}-1}$.
    Let $$C_s := \{(k-s,\mathbf{x}) \setsep \mathbf{x} \in B_{2,k}\}$$ and
    $$p_s(m) := (k-s,\beta_{2,k}(m))$$ for
    $m \in \setzeroton{\card{B_{2,k}}-1}$ and $s = 1, \ldots, 2k-1$.
    Let $$C_{2k} := \{(-k,\mathbf{x}) \setsep \mathbf{x} \in A_{2,k}\}$$ and
    $$p_{2k}(m) := (-k,\alpha'_{2,k}(m))$$ for $m \in
    \setzeroton{\card{A_{2,k}}-1}$. Define
    $$\beta_{3,k}(m) := p_s\left(m - \sum_{s'=0}^{s-1}\card{C_{s'}}\right)$$
    when
    $$
    \sum_{s'=0}^{s-1} \card{C_{s'}} \leq m < \sum_{s'=0}^s
    \card{C_{s'}}
    $$
    and $m \in \setzeroton{\card{B_{3,k}}-1}$
    Now functions $\beta_{3,k}$ satisfy the conditions of the lemma.
    
    Suppose that the
    lemma is true for some $n \in \naturalnumbers$, $n \geq 3$. Define
    $\beta^1_{n,k} := \beta_{n,k}$ for all $k \in \naturalnumbers$ and
    $\beta^{-1}_{n,k}(m) := \beta_{n,k}(\card{B_{n,k}}-1-m)$ for all
    $m \in \setzeroton{\card{B_{n,k}}-1}$, $k \in \naturalnumbers$.
    Define also
    $
    \alpha_{n,k}(m) := \beta^{(-1)^{k-j}}_{n,j}(m-\card{A_{n,j-1}})
    $
    where $\card{A_{n,j-1}} \leq m < \card{A_{n,j}}$ and
    $m \in \setzeroton{\card{A_{n,k}}-1}$, $j \in \setzeroton{k}$ and
    $k \in \positiveintegers$. When $k=0$ let
    $\alpha_{n,0}(0) = \finitezeroseq{n}$.
    Let also
    $$
    \alpha'_{n,k}(m) := \beta^{(-1)^j}_{n,k-j}\left( m - \sum_{j'=0}^{j-1}
    \card{B_{n,k-j'}} \right)
    $$
    for $m \in \setzeroton{\card{A_{n,k}}-1}$ and $j \in \setzeroton{k}$
    where
    $$
    \sum_{j'=0}^{j-1} \card{B_{n,k-j'}} \leq m < \sum_{j'=0}^j
    \card{B_{n,k-j'}} .
    $$
    Let $k \in \positiveintegers$.
    Let $$C_0 := \{(k,\mathbf{x}) \setsep \mathbf{x} \in A_{n,k}\}$$ and
    $$p_0(m) := (k,\alpha_{n,k}(m))$$ for
    $m \in \setzeroton{\card{A_{n,k}}-1}$.
    Let $$C_s := \{(k-s,\mathbf{x}) \setsep \mathbf{x} \in B_{n,k}\}$$ and
    $$
    p_s(m) := (k-s, \beta^{(-1)^s}_{n,k}(m))
    $$ for
    $m \in \setzeroton{\card{B_{n,k}}-1}$ and $a = 1, \ldots, 2k-1$.
    Let $$C_{2k} := \{(-k,\mathbf{x}) \setsep \mathbf{x} \in A_{n,k}\}$$ and
    $$p_{2k}(m) := (-k,\alpha'_{n,k}(m))$$ for
    $m \in \setzeroton{\card{A_{n,k}}-1}$.
    Define
    $$\beta_{n+1,k}(m) := p_s\left(m - \sum_{s'=0}^{s-1}\card{C_{s'}}\right)$$
    when
    $$
    \sum_{s'=0}^{s-1} \card{C_{s'}} \leq m < \sum_{s'=0}^s
    \card{C_{s'}} .
    $$
    and $m \in \setzeroton{\card{B_{n+1,k}}-1}$.
    Now functions $\beta_{n+1,k}$ satisfy the conditions of the lemma.
  \end{proof}
}
\fi

\if\shortprep0
{
  \begin{definition}
    Define
    \begin{displaymath}
      \cubeordering{2}{m} := \beta_{2,k}\left(m-\card{A_{2,k-1}}\right)
    \end{displaymath}
    where $\card{A_{2,k-1}} \leq m < \card{A_{2,k}}$ and
    $m, k \in \naturalnumbers$.
    When $n \geq 3$ define
    \begin{displaymath}
      \cubeordering{n}{m} := \beta^{(-1)^k}_{n,k}\left( m - \card{A_{n,k-1}}\right)
    \end{displaymath}
    where $\card{A_{n,k-1}} \leq m < \card{A_{n,k}}$ and
    $m, k \in \naturalnumbers$.
  \end{definition}
  
  Now \(\cubeorderingfunction{n}\) is a bijection from \(\naturalnumbers\)
  onto \(\zn\) preserving neighbours for each \(n \in \naturalnumbers\),
  \(n \geq 2\). Case \(n \geq 3\) follows from Lemma
  \ref{lem:ordering-bijections}.
}
\fi

\section{On Tensor Product Spaces}
\label{sec:tensor-products}

\if\shortprep0
{
\begin{definition}
  When \(n \in \positiveintegers\) define function
  \begin{math}
    \gensqordfunction{n} : \naturalnumbers \to \nn
  \end{math}
  by
  \begin{displaymath}
    \gensqord{n}{k}
    :=
    \left\{
      \begin{array}{ll}
        k ; & \;\; n = 1 \\
        \seqcomb{\gensqord{n-1}{\sqordfirst{k}}}{\sqordsecond{k}} ;
        & \;\; n \geq 2
      \end{array}
    \right.
  \end{displaymath}
  for all \(k \in \naturalnumbers\).
  When \(n \in \positiveintegers\) and \(j \in \setoneton{n}\)
  define function
  \begin{math}
    \gensqordcompfunction{n}{j} : \naturalnumbers \to \naturalnumbers
  \end{math}
  by
  \begin{displaymath}
    \gensqordcomp{n}{j}{k}
    :=
    \cartprodelem{\gensqord{n}{k}}{j}
  \end{displaymath}
  for all \(k \in \naturalnumbers\).
\end{definition}
}
\fi

\if\shortprep0
{
It can be shown by induction that
function \(\gensqordfunction{n}\) is a bijection from \(\naturalnumbers\) onto
\(\nn\) for each \(n \in \positiveintegers\).
}
\fi

\if\shortprep0
{
\begin{lemma}
  \label{lem:tp-schauder-basis}
  Let \(n \in \positiveintegers\) and
  \(\alpha = \projtn\) or \(\alpha = \injtn\).
  Suppose that \(A_k\) is a Banach space for each
  \(k \in \setoneton{n}\).
  Suppose also that
  \begin{math}
    (a_{k,j})_{j=0}^\infty
  \end{math}
  is a Schauder basis of Banach space \(A_k\)
  for each
  \(k \in \setoneton{n}\).
  Then the sequence
  \begin{displaymath}
    \left(
      \indexedtensorproduct{k=1}{n} a_{k,\gensqordcomp{n}{k}{j}}
    \right)_{j=0}^\infty
  \end{displaymath}
  is a Schauder basis of Banach space
  \begin{displaymath}
    \indexedctp{\alpha}{k=1}{n} A_k .
  \end{displaymath}
\end{lemma}

\begin{proof}
  Use the two dimensional case and induction by \(n\).
\end{proof}
}
\fi

\if\shortprep0
{
We give the following definition to be used with inherited tensor
products.
}
\fi

\if\shortprep0
{
\begin{definition}
  \label{def:preserve-operator-continuity}
  Let \(A\) and \(B\) be Banach spaces and \(\alpha\) a reasonable
  crossnorm on \(A \otimes B\). We say that \(\alpha\)
  \defterm{preserves operator continuity} iff
  \(S \otimes T\) is an operator from \(A \otimes_\alpha B\) into
  \(A \otimes_\alpha B\) whenever \(S : A \to A\) and \(T : B \to B\)
  are operators.
\end{definition}
}
\fi

\if\shortprep0
{
If \(\alpha\) is a uniform crossnorm \(\alpha_{A,B}\) preserves operator
continuity for all Banach spaces \(A\) and \(B\).
The distributive law for completed Banach space tensor
products and direct sums is
\if\shortprep1
given
\else
proved
\fi
next.
Relationship between tensor products and direct sums
has also been investigated by
Grecu and Ryan \cite{gr2005} and
Ansemil and Floret \cite{af2005}.
Some lemmas are
\if\shortprep1
given
\else
proved
\fi
first.
\if\longprep1
\givelemmawithoutproof
\fi
}
\fi

\if\shortprep0
{
\begin{lemma}
  \label{lem:dsum-closed}
  Let \(E\) be a Banach space.
  Let \(B\) and \(C\) be closed subspaces of \(E\) so that
  \(B \cap C = \{ 0 \}\).
  Define operators
  \(P_B : B \dsum C \to B\) and \(P_C : B \dsum C \to C\) by
  \begin{displaymath}
    P_B (b + c) := b , \spaceafter x = b + c \in B \dsum C
  \end{displaymath}
  and
  \begin{displaymath}
    P_C (b + c) := c , \spaceafter x = b + c \in B \dsum C .
  \end{displaymath}
  Then \(B \dsum C\) is a closed subspace of \(E\).
\end{lemma}
}
\fi

\if\longversion1
{
\begin{proof}
  The direct sum \(B \dsum C\) is a normed subspace of \(E\).
  Let \(M = \max \{ \norm{P_B}, \norm{P_C} \}\).
  Let \((x_k)_{k=0}^\infty \subset B \dsum C\) be a Cauchy sequence
  and \(x_k = b_k + c_k\) where \(b_k \in B\) and \(c_k \in C\) for
  all \(k \in \naturalnumbers\).
  Let \(h > 0\). There exists \(N \in \naturalnumbers\) so that
  \(\norm{x_{k'} - x_k} < \frac{h}{M}\) for all \(k, k' \in
  \naturalnumbers\) so that \(k, k' \geq N\).
  Now \(\norm{b_{k'} - b_k} = \norm{P_B(x_{k'}-x_k)} \leq \norm{P_B}
  \norm{x_{k'} - x_k} < h\) when \(k, k' \geq N\).
  Similarly \(\norm{c_{k'} - c_k} < h\) when \(k, k' \geq N\).
  Hence \((b_k)_{k=0}^\infty\) is a Cauchy sequence in \(B\) and
  \((c_l)_{l=0}^\infty\) is a Cauchy sequence in \(C\). So
  \(b_k \to b \in B\) as \(k \to \infty\) and
  \(c_l \to c \in C\) as \(l \to \infty\).
  Furthermore, \(x_k = b_k + c_k \to b + c \in B \dsum C\) as \(k \to
  \infty\).
  Hence every Cauchy sequence in \(B \dsum C\) converges in that space
  and \(B \dsum C\) is a closed subspace of \(E\).
\end{proof}
}
\fi

\if\shortprep0
{
\begin{lemma}
  \label{lem:disjoint-tp}
  Let \(A\) and \(E\) be
  Banach spaces. Let \(B\) and \(C\) be closed subspaces of \(E\) so
  that \(B \cap C = \{ 0 \}\).
  Let \(\alpha\) be a reasonable crossnorm on \(A \otimes (B \dsum C)\)
  and \(\beta\) a reasonable crossnorm on \((B \dsum C) \otimes A\).
  Suppose that \(\alpha\) and \(\beta\) preserve operator continuity.
  Define \(P_B\) and \(P_C\) as in \mylemma
  \ref{lem:dsum-closed} and suppose that \(P_B\) and
  \(P_C\) are operators.
  Then
  \begin{math}
    A \ctp_{\alpha; A \otimes (B \dsum C)} B \cap A \ctp_{\alpha; A
      \otimes (B \dsum C)} C = \{ 0 \}
  \end{math}
  and
  \begin{math}
    B \ctp_{\beta; A \otimes (B \dsum C)} A \cap C \ctp_{\beta; A
      \otimes (B \dsum C)} A = \{ 0 \}
  \end{math}.
\end{lemma}
}
\fi

\if\shortprep0
{
\begin{proof}
  Suppose that \(P_B\) is continuous.

  Function \(\idfunconspace{A}\) is an operator. 
  Define function \(P_1 : A \otimes_\alpha (B \dsum C) \to A
  \otimes_{\alpha; A \otimes (B \dsum C)} B\) by \( P_1 = 
  \idfunconspace{A} \otimes P_B \). Function
  \(P_1\) is an operator and (see \cite[chapter 6.1]{ryan2002}) it
  has a unique continuous linear extension
  \(P'_1 : A \ctp_{\alpha; A \otimes (B \dsum C)} (B \dsum C)
  \to A \ctp_{\alpha; A \otimes (B \dsum C)} B\).
  \(P_1\) and \(P'_1\) have the same norm.

  Let \(f \in A \ctp_{\alpha; A \otimes (B \dsum C)} B\). There exists
  a sequence
  \((f_k)_{k=0}^\infty\) in the noncompleted tensor product \(A
  \otimes_{\alpha; A \otimes (B \dsum C)} B\) so that
  \(\norm{f - f_k}_\alpha \to 0\) as \(k
  \to \infty\). Now \(P'_1 f_k = P_1 f_k = f_k\) for all \(k \in
  \naturalnumbers\) and hence \(P'_1 f_k \to f\) as \(k \to
  \infty\). On the other hand, \(P'_1\) is continuous and \(P'_1 f_k
  \to P'_1 f\) as \(k \to \infty\). Hence \(P'_1 f = f\). Therefore
  \(P'_1 f = f\) for all \(f \in A \ctp_{\alpha; A \otimes (B \dsum C)} B\).

  Let \(g \in A \ctp_{\alpha; A \otimes (B \dsum C)} C\). There exists
  a sequence
  \((g_k)_{k=0}^\infty\) in the noncompleted tensor product \(A
  \otimes_{\alpha; A \otimes (B \dsum C)} C\) so that \(\norm{g -
    g_k}_\alpha \to 0\) as \(k
  \to \infty\). Now
  \begin{displaymath}
    g_k = \sum_{i=1}^{m(k)} a_{k,i} \otimes c_{k,i}
  \end{displaymath}
  where \(a_{k,i} \in A\), \(c_{k,i} \in C\),
  \(m(k) \in \naturalnumbers\), and \(k \in \naturalnumbers\).
  Since \(P_B c = 0\) for all \(c \in C\)
  \begin{displaymath}
    P_B g_k = \sum_{i=1}^{m(k)} a_{k,i} \otimes \left( P_B
    c_{k,i} \right) = 0 \qspace \forall k \in \naturalnumbers .
  \end{displaymath}
  Hence \(P'_1 g = 0\).
  If \(g \in A \ctp_{\alpha; A \otimes (B \dsum C)} C\) and
  \(g \neq 0\) then \(P'_1 g = 0\)
  and \(P'_1 g \neq g\) and hence
  \(g \not\in A \ctp_{\alpha; A \otimes (B \dsum C)} B\). So
  \begin{math}
    A \ctp_{\alpha; A \otimes (B \dsum C)} B \cap A \ctp_{\alpha; A
      \otimes (B \dsum C)} C = \{ 0 \}
  \end{math}.
  Proof of equation
  \begin{math}
    B \ctp_{\beta; A \otimes (B \dsum C)} A \cap C \ctp_{\beta; A
      \otimes (B \dsum C)} A = \{ 0 \}
  \end{math}
  is similar.
\end{proof}
}
\fi

\if\shortprep0
{
\begin{theorem}
  \label{th:distr-law}
  Let \(A\) and \(E\) be
  Banach spaces. Let \(B\) and \(C\) be closed subspaces of \(E\) so
  that \(B \cap C = \{ 0 \}\).
  Define \(P_B\) and \(P_C\) as in \mylemma
  \ref{lem:dsum-closed} and suppose that \(P_B\) and
  \(P_C\) are operators.
  Let \(\alpha\) be a reasonable crossnorm on \(A \otimes (B \dsum C)\)
  and \(\beta\) a reasonable crossnorm on \((B \dsum C) \otimes A\).
  Suppose that \(\alpha\) and \(\beta\) preserve operator continuity.
  Now
  \begin{math}
    A \ctp_\alpha \left( B \dsum C \right) \equalns
    A \ctp_{\alpha ; A \otimes (B \dsum C)} B \dsum
    A \ctp_{\alpha ; A \otimes (B \dsum C)} C
  \end{math}
  and
  \begin{math}
    \left( B \dsum C \right) \ctp_\beta A \equalns
    B \ctp_{\beta ; (B \dsum C) \otimes A} A
    \dsum C \ctp_{\beta ; (B \dsum C) \otimes A} A
  \end{math}.
\end{theorem}
}
\fi

\if10
{
Relationship between tensor products and direct sums
has also been investigated by
Grecu and Ryan \cite{gr2005} and
Ansemil and Floret \cite{af2005}.
}
\fi

\if\shortprep0
{
\begin{proof}
  By \mylemma \ref{lem:disjoint-tp}
  \begin{displaymath}
    A \ctp_{\alpha; A \otimes (B \dsum C)} B \cap A \ctp_{\alpha; A
      \otimes (B \dsum C)} C = \{ 0 \} .
  \end{displaymath}
  Let
  \begin{displaymath}
    M \defequalns A \ctp_{\alpha; A \otimes (B \dsum C)} B \dsum A \ctp_{\alpha; A
      \otimes (B \dsum C)} C .
  \end{displaymath}
  Now \(M\) is a normed subspace of \(A \ctp_\alpha (B \dsum C)\).
  Define operators
  \(P_1 : A \otimes_\alpha (B \dsum C) \to A \otimes_\alpha B\)
  and
  \(P_2 : A \otimes_\alpha (B \dsum C) \to A \otimes_\alpha C\)
  by \(P_1 = \idfunconspace{A} \otimes P_B\) and
  \(P_2 = \idfunconspace{A} \otimes P_C\).
  Operator \(P_1\) has a unique continuous linear
  extension \(P_1' : A \ctp_\alpha (B \dsum C) \to A \ctp_\alpha B\).
  We have \(\norm{P_1} = \norm{P_1'}\).
  Similarly, operator \(P_2\) has a unique continuous linear extension
  \(P_2' : A \ctp_\alpha (B \dsum C) \to A \ctp_\alpha C\)
  and we have \(\norm{P_2} = \norm{P_2'}\).
  Let \(x \in M\). Now \(x = u + v\) where
  \(u \in A \ctp_{\alpha; A \otimes (B \dsum C)} B\)
  and
  \(v \in A \ctp_{\alpha; A \otimes (B \dsum C)} C\).
  There exists
  \((u_k)_{k=0}^\infty \subset
  A \otimes_{\alpha; A \otimes (B \dsum C)} B\)
  so that \(u_k \to u\) as \(k \to \infty\) and
  \begin{displaymath}
    u_k = \sum_{j=1}^{m_1(k)} w_{j,k} \otimes b_{j,k}
  \end{displaymath}
  where \(m_1(k) \in \positiveintegers\), \(w_{j,k} \in A\), and
  \(b_{j,k} \in B\).
  There exists
  \((v_k)_{k=0}^\infty \subset
  A \otimes_{\alpha; A \otimes (B \dsum C)} C\)
  so that \(v_k \to v\) as \(k \to \infty\) and
  \begin{displaymath}
    v_k = \sum_{j=1}^{m_2(k)} z_{j,k} \otimes c_{j,k}
  \end{displaymath}
  where \(m_2(k) \in \positiveintegers\), \(z_{j,k} \in A\), and
  \(c_{j,k} \in C\).
  Now \(P_1 u_k = u_k\) and \(P_1 v_k = 0\) for all
  \(k \in \naturalnumbers\) and hence \(P_1' u = u\) and
  \(P_1' v = 0\).
  Similarly, \(P_2 u_k = 0\) and \(P_2 v_k = v_k\) for all
  \(k \in \naturalnumbers\) and hence \(P_2' u = 0\) and
  \(P_2' v = v\).
  So \(P_1' x = u\) and \(P_2' x = v\).
  By \mylemma \ref{lem:dsum-closed} the normed vector space \(M\) is a
  closed subspace of \(A \ctp_\alpha (B \dsum C)\).

  Let \(t \in A \ctp_\alpha (B \dsum C)\). Let \(h > 0\). There exists
  \(y \in A \otimes_\alpha (B \dsum C)\) so that \(\norm{y - t} < h\).
  Now
  \begin{displaymath}
    y = \sum_{j=1}^n a_j \otimes (b'_j + c'_j) =
    \sum_{j=1}^n a_j \otimes b'_j + \sum_{j=1}^n a_j \otimes c'_j
  \end{displaymath}
  where \(n \in \naturalnumbers\) and \(a_j \in A\), \(b'_j \in B\),
  \(c'_j \in C\) for \(j \in \setoneton{n}\).
  Therefore
  \(y \in A \otimes_{\alpha; A \otimes (B \dsum C)} B \dsum
  A \otimes_{\alpha; A \otimes (B \dsum C)} C \subset M\).
  Number \(h > 0\) was arbitrary so \(t \in {\overline M} = M\).
  Hence
  \begin{displaymath}
    A \ctp_\alpha \left( B \dsum C \right) \equalns
    A \ctp_{\alpha ; A \otimes (B \dsum C)} B \dsum
    A \ctp_{\alpha ; A \otimes (B \dsum C)} C .
  \end{displaymath}
  Proof of
  \begin{displaymath}
    \left( B \dsum C \right) \ctp_\beta A \equalns
    B \ctp_{\beta ; (B \dsum C) \otimes A} A
    \dsum C \ctp_{\beta ; (B \dsum C) \otimes A} A
  \end{displaymath}
  is similar.
\end{proof}
}
\fi

\if\shortprep0
{
\begin{corollary}
  \label{cor:distr-law-ucn}
  Let \(\alpha\) be an injective uniform crossnorm. Let \(A\) and \(E\) be
  Banach spaces. Let \(B\) and \(C\) be closed subspaces of \(E\) so
  that \(B \cap C = \{ 0 \}\).
  Define \(P_B\) and \(P_C\) as in \mylemma
  \ref{lem:dsum-closed} and suppose that \(P_B\) and
  \(P_C\) are operators.
  Then
  \begin{align*}
    A \ctp_\alpha B \cap A \ctp_\alpha C & = \{ 0 \} , \\
    B \ctp_\alpha A \cap C \ctp_\alpha A & = \{ 0 \} ,
  \end{align*}
  and
  \begin{align*}
    A \ctp_\alpha \left( B \dsum C \right) & \equalns
    A \ctp_\alpha B \dsum A \ctp_\alpha C , \\
    \left( B \dsum C \right) \ctp_\alpha A & \equalns
    B \ctp_\alpha A \dsum C \ctp_\alpha A .
  \end{align*}
\end{corollary}
}
\fi

\if\shortprep0
{
\begin{lemma}
\label{lem:inj-and-sup-norm-noncompleted}
\definek
Let \(n, m \in \positiveintegers\),
\(X\) be a closed subspace of \(\cbfunc{\realnumbers^n}{K}\) and
\(Y\) be a closed subspace of \(\cbfunc{\realnumbers^m}{K}\).
Then
\begin{displaymath}
  \injtn(u) = \sup_{\vectorstyle{t} \in \realnumbers^{n+m}}
    \abs{\sum_{i=1}^k
      x_i(\cartprodelem{\vectorstyle{t}}{1},
        \ldots, \cartprodelem{\vectorstyle{t}}{n})
      y_i(\cartprodelem{\vectorstyle{t}}{n+1},
        \ldots, \cartprodelem{\vectorstyle{t}}{n+m})}
\end{displaymath}
for all \(u \in X \otimes Y\)
where
\begin{displaymath}
  u = \sum_{i=1}^k x_i \otimes y_i ,
\end{displaymath}
\(k \in \positiveintegers\), and \(x_i \in X\),
\(y_i \in Y\) for \(i \in \setoneton{k}\).
I.e. the injective tensor norm of \(u\) equals to the norm of \(u\) as
an element of \(\cbfunc{\realnumbers^{n+m}}{K}\).
\end{lemma}
}
\fi

\if\shortprep0
{
\begin{proof}
Let \(u \in X \otimes Y\),
\begin{displaymath}
  u = \sum_{i=1}^k x_i \otimes y_i ,
\end{displaymath}
where \(k \in \positiveintegers\) and \(x_i \in X\), \(y_i \in Y\) for
\(i \in \setoneton{k}\).
Now
\begin{displaymath}
  \injtn(u) = \sup \left\{ \abs{\sum_{i=1}^k \myphi(x_i) \psi(y_i)}
      \setsep \myphi \in B_{\topdual{X}}, \psi \in B_{\topdual{Y}}
    \right\}
\end{displaymath}
Set
\begin{math}
  A_1 = \left\{ \delta_\rnt \in \topdual{X}
    \setsep \rnt \in \realnumbers^n \right\}
\end{math}
is a norming set of \(X\), set
\begin{math}
  A_2 = \left\{ \delta_\rnt \in \topdual{Y}
    \setsep \rnt \in \realnumbers^m \right\}
\end{math}
is a norming set of \(Y\), and consequently
\cite[chapter 3.1]{ryan2002}
\if\elsevier0
{
\begin{align*}
  \injtn(u) & = \sup \left\{ \abs{\sum_{i=1}^k \myphi(x_i) \psi(y_i)}
      \setsep \myphi \in A_1, \psi \in A_2
    \right\} \\
  & = \sup \left\{ \abs{\sum_{i=1}^k 
    \delta_{\indexedvectorstyle{t}{1}}(x_i)
    \delta_{\indexedvectorstyle{t}{2}}(y_i)}
      \setsep \indexedvectorstyle{t}{1} \in \realnumbers^n,
      \indexedvectorstyle{t}{2} \in \realnumbers^m
    \right\} \\
  & = \sup_{t \in \realnumbers^{n+m}}
  \abs{\sum_{i=1}^k 
    x_i(\cartprodelem{\vectorstyle{t}}{1},
      \ldots, \cartprodelem{\vectorstyle{t}}{n})
    y_i(\cartprodelem{\vectorstyle{t}}{n+1},
    \ldots, \cartprodelem{\vectorstyle{t}}{n+m})} .
\end{align*}
}
\else
{
\begin{align*}
  \injtn(u) & = \sup \left\{ \abs{\sum_{i=1}^k \myphi(x_i) \psi(y_i)}
      \setsep \myphi \in A_1, \psi \in A_2
    \right\} \\
  & = \sup \left\{ \abs{\sum_{i=1}^k 
    \delta_{\indexedvectorstyle{t}{1}}(x_i)
    \delta_{\indexedvectorstyle{t}{2}}(y_i)}
      \setsep \indexedvectorstyle{t}{1} \in \realnumbers^n,
      \indexedvectorstyle{t}{2} \in \realnumbers^m
    \right\} \\
  & = \sup_{t \in \realnumbers^{n+m}}
  \abs{\sum_{i=1}^k 
    x_i(\cartprodelem{\vectorstyle{t}}{1},
      \ldots, \cartprodelem{\vectorstyle{t}}{n})
    y_i(\cartprodelem{\vectorstyle{t}}{n+1},
    \ldots, \cartprodelem{\vectorstyle{t}}{n+m})} .
\end{align*}
}
\fi
\end{proof}
}
\fi

\if\longversion1
{
\begin{theorem}
\label{th:inj-and-sup-norm}
\definek
Let \(n, m \in \positiveintegers\),
\(X\) be a closed subspace of \(\cbfunc{\realnumbers^n}{K}\), and
\(Y\) be a closed subspace of \(\cbfunc{\realnumbers^m}{K}\). The tensor
product \(X \citp Y\) is a closed subspace of
\(\cbfunc{\realnumbers^{n + m}}{K}\)
and consequently
\begin{displaymath}
  \injtn(u) = \norminfty{u} \qspace \forall u \in X \citp Y
\end{displaymath}
where \(\norminfty{u}\) is the supremum norm of \(u\) as an element of
\(\cbfunc{\realnumbers^{n + m}}{K}\).
\end{theorem}
}
\else
{
\if\longprep1
{
\begin{corollary}
\label{cor:inj-and-sup-norm}
\definek
Let \(n, m \in \positiveintegers\),
\(X\) be a closed subspace of \(\cbfunc{\realnumbers^n}{K}\), and
\(Y\) be a closed subspace of \(\cbfunc{\realnumbers^m}{K}\). The tensor
product \(X \citp Y\) is a closed subspace of
\(\cbfunc{\realnumbers^{n + m}}{K}\)
and consequently
\begin{displaymath}
  \injtn(u) = \norminfty{u} \qspace \forall u \in X \citp Y
\end{displaymath}
where \(\norminfty{u}\) is the supremum norm of \(u\) as an element of
\(\cbfunc{\realnumbers^{n + m}}{K}\).
\end{corollary}
}
\fi
}
\fi

\if\longversion1
{
\begin{proof}
  The algebraic tensor product \(X \otimes Y\) is isomorphic to a
  subspace of \(\cbfunc{\realnumbers^{n+m}}{K}\) as a vector space. By
  \mylemma \ref{lem:inj-and-sup-norm-noncompleted} the norm in \(X
  \otimes_\injtn Y\) is equal to the supremum norm. Therefore the
  completion \(X \citp Y\) is the closure of
  \(X \otimes_\injtn Y\) in \(\cbfunc{\realnumbers^{n+m}}{K}\). Hence
  the norm of an element \(u \in X \citp Y\) equals the norm of \(u\)
  in \(\cbfunc{\realnumbers^{n+m}}{K}\).
\end{proof}
}
\fi

\if\shortprep1
{
\begin{theorem}
\label{th:inj-and-sup-norm-general}
Let \(K = \realnumbers\) or \(K = \complexnumbers\),
\(n \in \naturalnumbers + 2\), and \(X_1, \ldots, X_n\) be closed
subspaces of \(\cbfunc{\realnumbers}{K}\). Then the tensor
product \(X_1 \citp \cdots \citp X_n\) is a closed
subspace of \(\cbfunc{\realnumbers^n}{K}\).
\end{theorem}
}
\else
{
\begin{corollary}
\label{cor:inj-and-sup-norm-general}
Let \(K = \realnumbers\) or \(K = \complexnumbers\),
\(n \in \naturalnumbers + 2\), and \(X_1, \ldots, X_n\) be closed
subspaces of \(\cbfunc{\realnumbers}{K}\). Then the tensor
product \(X_1 \citp \cdots \citp X_n\) is a closed
subspace of \(\cbfunc{\realnumbers^n}{K}\).
\end{corollary}
}
\fi

\if\longprep1
{
\givelemmaswithoutproofsn{two}
}
\fi

\if\shortprep0
{
\begin{lemma}
  \label{lem:function-tp-vanishing}
  Let \(K = \realnumbers\) or \(K = \complexnumbers\). Let
  \(n, m \in \positiveintegers\),
  \(f \in \vanishingfunc{\realnumbers^n}{K}\) and
  \(g \in \vanishingfunc{\realnumbers^m}{K}\).
  Then
  \begin{math}
    f \otimes g \in \vanishingfunc{\realnumbers^{n+m}}{K}
  \end{math}.
\end{lemma}
}
\fi

\if\longversion1
{
\begin{proof}
  Define operator
  \(A_1 : \realnumbers^{n+m} \to \realnumbers^n\) by
  \begin{displaymath}
    A_1 \vectorstyle{x} := (\cartprodelem{\vectorstyle{x}}{1},
      \ldots, \cartprodelem{\vectorstyle{x}}{n}),
    \spaceafter \vectorstyle{x} \in \realnumbers^{n+m}
  \end{displaymath}
  and operator
  \(A_2 : \realnumbers^{n+m} \to \realnumbers^m\) by
  \begin{displaymath}
    A_2 \vectorstyle{x}
    := (\cartprodelem{\vectorstyle{x}}{n+1},
      \ldots, \cartprodelem{\vectorstyle{x}}{n+m}),
    \spaceafter \vectorstyle{x} \in \realnumbers^{n+m} . 
  \end{displaymath}
  Now
  \begin{math}
    (f \otimes g)(\vectorstyle{x})
    = f(A_1 \vectorstyle{x}) g(A_2 \vectorstyle{x})
  \end{math}
  for all
  \begin{math}
    \vectorstyle{x} \in \realnumbers^{n+m}
  \end{math}
  and hence \(f \otimes g\) is continuous.
  Since
  \begin{math}
    \abs{(f \otimes g)(\vectorstyle{x})}
    = \abs{f(A_1 \vectorstyle{x}) g(A_2 \vectorstyle{x})}
    \leq \norminfty{f} \norminfty{g}
  \end{math}
  for all
  \begin{math}
    \vectorstyle{x} \in \realnumbers^{n+m}
  \end{math}
  function \(f \otimes g\) is bounded.

  Let \(c \in ] 0, 1 [\).
  There exists \(r_1 \in \positiverealnumbers\) so that
  \begin{displaymath}
    \abs{f(\vectorstyle{x})} < \frac{c}{\norminfty{g} + 1}
    \qspace \forall \vectorstyle{x}
    \in \realnumbers^{n+m} \setminus \closedball{\rn}{0}{r_1} .
  \end{displaymath}
  Similarly, there exists \(r_2 \in \positiverealnumbers\) so that
  \begin{displaymath}
    \abs{g(\vectorstyle{x})} < \frac{c}{\norminfty{f} + 1}
    \qspace \forall \vectorstyle{x}
    \in \realnumbers^{n+m} \setminus \closedball{\rn}{0}{r_1} .
  \end{displaymath}
  Let \(r := 2 \max \{ r_1, r_2 \}\).
  Suppose that
  \begin{math}
    \vectorstyle{y} \in \realnumbers^{n+m} \setminus \closedball{\rn}{0}{r}
  \end{math}.
  Now
  \begin{displaymath}
    \norm{\vectorstyle{y}}
    = \sqrt{\sum_{k=1}^{n+m} \left( \cartprodelem{\vectorstyle{y}}{k} \right)^2}
    = \sqrt{\norm{A_1 \vectorstyle{y}}^2 + \norm{A_2 \vectorstyle{y}}^2} > r .
  \end{displaymath}
  Consequently
  \begin{eqnarray*}
    \norm{A_1 \vectorstyle{y}}^2 + \norm{A_2 \vectorstyle{y}}^2 >  r^2
    & \implies &
    \max \{ \norm{A_1 \vectorstyle{y}}^2, \norm{A_2 \vectorstyle{y}}^2 \}
    > \frac{r^2}{2}
    \\
    & \implies & \norm{A_1 \vectorstyle{y}} > \max \{ r_1, r_2 \} \geq r_1 \\
    & & \lor \;
    \norm{A_2 \vectorstyle{y}} > \max \{ r_1, r_2 \} \geq r_2 .
  \end{eqnarray*}
  In case \(\norm{A_1 \vectorstyle{y}} > r_1\)
  \begin{displaymath}
    \abs{(f \otimes g)(\vectorstyle{y})}
    = \abs{f(A_1 \vectorstyle{y})} \abs{g(A_2 \vectorstyle{y})}
    \leq \frac{c}{\norminfty{g} + 1} \norminfty{g} < c
  \end{displaymath}
  and in case \(\norm{A_2 \vectorstyle{y}} > r_2\)
  \begin{displaymath}
    \abs{(f \otimes g)(\vectorstyle{y})}
    = \abs{f(A_1 \vectorstyle{y})} \abs{g(A_2 \vectorstyle{y})}
    \leq \norminfty{f} \frac{c}{\norminfty{f} + 1} < c .
  \end{displaymath}
  Hence \(\abs{(f \otimes g)(\vectorstyle{y})} < c\).
  Consequently
  \begin{displaymath}
    \lim_{\norm{\vectorstyle{x}} \to \infty} (f \otimes g)(\vectorstyle{y}) = 0
  \end{displaymath}
  and \(f \otimes g \in \vanishingfunc{\realnumbers^{n+m}}{K}\).
\end{proof}
}
\fi

\if\shortprep0
{
\begin{lemma}
  \label{lem:function-tp-uniform}
  Let \(K = \realnumbers\) or \(K = \complexnumbers\). Let
  \(n, m \in \positiveintegers\),
  \(f \in \ucfunc{\realnumbers^n}{K}\) and
  \(g \in \ucfunc{\realnumbers^m}{K}\).
  Then
  \begin{math}
    f \otimes g \in \ucfunc{\realnumbers^{n+m}}{K}
  \end{math}.
\end{lemma}
}
\fi

\if\shortprep0
{
\begin{lemma}
  \label{lem:uniform-tp-inclusion}
  \definekn
  Then
  \begin{displaymath}
    \indexedcitp{j=1}{n} \ucfunc{\realnumbers}{K} \closedsubspace
    \ucfunc{\realnumbers^n}{K}
  \end{displaymath}
\end{lemma}
}
\fi

\if\shortprep0
{
\begin{proof}
  Define
  \begin{displaymath}
    F^{[n']} \defequalns \indexedcitp{j=1}{n'} \ucfunc{\realnumbers}{K}
  \end{displaymath}
	for all \(n' \in \integernumbers\).
	By \mycorollary \ref{cor:inj-and-sup-norm-general}
	\begin{math}
	  F^{[n']} \closedsubspace \cbfunc{\rn}{K}
	\end{math}
	for each \(n' \in \positiveintegers\).
	When \(n = 1\) the lemma is true.
	Suppose that the lemma is true for some \(n' \in \positiveintegers\).
	By \mylemma \ref{lem:function-tp-uniform}
	\begin{math}
	  F^{[n']} \otimes_\injtn \ucfunc{\realnumbers}{K}
	  \normedsubspace \ucfunc{\realnumbers^{n'+1}}{K}
	\end{math}.
	As \(\ucfunc{\realnumbers^{n'+1}}{K}\) is complete
	it follows that
	\begin{math}
	  F^{[n'+1]} \equalns F^{[n']} \citp \ucfunc{\realnumbers}{K}
	  \closedsubspace \ucfunc{\realnumbers^{n'+1}}{K}
	\end{math}.
\end{proof}
}
\fi

\if\shortprep0
{
\begin{lemma}
  \label{lem:vanishing-tp-inclusion}
  \definekn
  Then
  \begin{displaymath}
    \indexedcitp{j=1}{n} \vanishingfunc{\realnumbers}{K} \closedsubspace
    \vanishingfunc{\realnumbers^n}{K} .
  \end{displaymath}
\end{lemma}
}
\fi

\if\longprep1
{
\begin{proof}
  The proof is similar to the proof of \mylemma
  \ref{lem:uniform-tp-inclusion}.
\end{proof}
}
\else
{
\if\longversion1
{
\begin{proof}
  Define
  \begin{displaymath}
    F^{[n']} \defequalns \indexedcitp{j=1}{n'} \vanishingfunc{\realnumbers}{K}
  \end{displaymath}
	for all \(n' \in \integernumbers\).
	By \mycorollary \ref{cor:inj-and-sup-norm-general}
	\begin{math}
	  F^{[n']} \closedsubspace \cbfunc{\rn}{K}
	\end{math}
	for each \(n' \in \positiveintegers\).
	When \(n = 1\) the lemma is true.
	Suppose that the lemma is true for some \(n' \in \positiveintegers\).
	By \mylemma \ref{lem:function-tp-vanishing}
	\begin{math}
	  F^{[n']} \otimes_\injtn \vanishingfunc{\realnumbers}{K}
	  \normedsubspace \vanishingfunc{\realnumbers^{n'+1}}{K}
	\end{math}.
	As \(\vanishingfunc{\realnumbers^{n'+1}}{K}\) is complete
	it follows that
	\begin{math}
	  F^{[n'+1]} \equalns F^{[n']} \citp \vanishingfunc{\realnumbers}{K}
	  \closedsubspace \vanishingfunc{\realnumbers^{n'+1}}{K}
	\end{math}.
\end{proof}
}
\fi
}
\fi

\if\longprep1
{
  \givelemmawithoutproof
}
\fi

\if\longversion1
{
\begin{lemma}
  \label{lem:operator-norm-tp-operator-continuity}
  \definek
  Let \(\alpha\) be a uniform crossnorm.
  Let \(E_1\), \(E_2\), \(F_1\), and \(F_2\) be Banach spaces.
  Let \(A_1\) be a closed subspace of \(E_1\), \(A_2\) a closed
  subspace of \(E_2\), \(B_1\) a closed subspace of \(F_1\), and
  \(B_2\) a closed subspace of \(F_2\).
  Suppose that \(A_1\) is topologically complemented in \(E_1\) and
  \(A_2\) is topologically complemented in \(E_2\). 
  Let \(S : A_1 \to B_1\) and \(T : A_2 \to B_2\) be
  operators. Then \(S \otimes T\) is an operator from
  \(A_1 \otimes_{\alpha ; E_1 \otimes E_2} A_2\) into
  \(B_1 \otimes_{\alpha ; F_1 \otimes F_2} B_2\).
\end{lemma}
}
\fi

\if\longversion1
{
\begin{proof}
  Let \(P_1\) be a continuous projection \projprep
  \(E_1\) onto \(A_1\) and
  \(P_2\) be a continuous projection \projprep
  \(E_2\) onto \(A_2\).
  Let \(X = E_1 \otimes_\alpha E_2\) and
  \(Y = F_1 \otimes_\alpha F_2\).
  Let \(R = (S \circ P_1) \otimes (T \circ P_2)\). Now \(R\) is an
  operator from \(X\) into \(Y\) and
  \begin{math}
    \setimage{R}{X} = (S \otimes T) \left[
      A_1 \otimes_{\alpha ; E_1 \otimes E_2} A_2 \right]
    \subset B_1 \otimes_{\alpha ; F_1 \otimes F_2} B_2
  \end{math}.
  Hence \(S \otimes T = R \vert \left( A_1 \otimes_{\alpha ; E_1 \otimes E_2}
  A_2 \right)\) is an operator and
  \begin{displaymath}
    \norm{S \otimes T} \leq \norm{R} = \norm{S \circ P_1} \norm{T
      \circ P_2}
    \leq \norm{S} \norm{T} \norm{P_1} \norm{P_2} .
  \end{displaymath}
\end{proof}
}
\fi

\if\shortprep0
{
\begin{definition}
  \label{def:inh-tp-alt}
  Let \(E_1\), \(E_2\), and \(F\) Banach spaces so that \(E_1 \otimes
  E_2\) is a linear subspace of \(F\). Define \(E_1 \otimes_{(F)}
  E_2\) to be the vector space \(E_1 \otimes E_2\) equipped with the
  norm inherited from \(F\) and
  \begin{displaymath}
    E_1 \ctp_{(F)} E_2 \defequalns \closop_F \left( E_1 \otimes_{(F)} E_2
    \right) .
  \end{displaymath}
\end{definition}
}
\fi

\begin{definition}
  \label{def:indexed-inh-tp-alt}
  Let \(n \in \positiveintegers\). Let \(A_1, \ldots, A_n\) and
  \(B_1, \ldots, B_n\) be Banach spaces so that \(A_1\) is a closed subspace
  of \(B_1\) and \(B_k \otimes A_{k+1}\) is a linear subspace of
  \(B_{k+1}\) for \(k \in \setoneton{n-1}\).
  When \(k \geq 2\) define
  \begin{displaymath}
    \indexedinhctp{B_j}{j=1}{k} A_j \defequalns \closop_{B_k} T_k
  \end{displaymath}
  where
  \begin{displaymath}
    T_k \defequalns \left( \indexedinhctp{B_j}{j=1}{k - 1} A_j \right)
    \otimes_{(B_k)} A_k .
  \end{displaymath}
  When \(k = 1\) define
  \begin{displaymath}
    \indexedinhctp{B_j}{j=1}{k} A_j \defequalns A_1 .
  \end{displaymath}
\end{definition}

When the assumptions of \mydef
\ref{def:indexed-inh-tp-alt} hold
\begin{displaymath}
  \indexedinhctp{B_j}{j=1}{k} A_j
\end{displaymath}
is a closed subspace of Banach space \(B_k\) for all \(k \in \setoneton{n}\).

\if\shortprep0
{
\begin{definition}
  \label{def:inh-op-tp}
  Let \(A_1, A_2, B_1, B_2, E\), and \(F\) be Banach spaces
  so that \(A_1 \otimes A_2\) is a linear subspace of \(E\) and
  \(B_1 \otimes B_2\) is a linear subspace of \(F\).
  Let \(P_1 : A_1 \to B_1\) and \(P_2 : A_2 \to B_2\) be operators
  so that \(P_1 \otimes P_2\) is an
  operator from \(A_1 \otimes_{(E)} A_2\) into
  \(B_1 \otimes_{(F)} B_2\).
  Define operator
  \begin{math}
    P_1 \otimes_{(E, F)} P_2 : A_1 \ctp_{(E)} A_2 \to B_1
    \ctp_{(F)} B_2
  \end{math}
  to be the unique continuous linear extension of \(P_1 \otimes P_2\) to
  \(A_1 \ctp_{(E)} A_2\).
\end{definition}
}
\fi

\if\shortprep0
{
The extension in the previous definion always exists.
See \cite[theorem 3.4.4]{hpc2005} and
\cite[chapter 6.1]{ryan2002}.
}
\fi

\begin{definition}
  \label{def:indexed-inh-op-tp}
  Let \(n \in \positiveintegers\).
  Let \(A_1, \ldots, A_n\), \(B_1, \ldots, B_n\),
  \(E_1, \ldots, E_n\), and \(F_1, \ldots, F_n\) be Banach spaces so
  that
  \begin{itemize}
  \item \(A_1\) is a closed subspace of \(E_1\) and
    \(E_k \otimes A_{k+1}\) is a linear subspace of \(E_{k+1}\) for
    all \(k = 1, \ldots, n - 1\).
  \item \(B_1\) is a closed subspace of \(F_1\) and
    \(F_k \otimes B_{k+1}\) is a linear subspace of \(F_{k+1}\) for
    all \(k = 1, \ldots, n - 1\).
  \end{itemize}
  Suppose that \(P_k : A_k \to B_k\), \(k = 1, \ldots, n\) are
  operators.
  Let \(S_1 := P_1\), \(T_1 := S_1 = P_1\), and
  \begin{math}
    S_k := T_{k-1} \otimes P_k
  \end{math}
  for \(k = 2, \ldots, n\).
  When \(k \in \{ 2, \ldots, n \}\) and \(S_k\) is continuous let
  \begin{displaymath}
    T_k : \indexedinhctp{E_l}{l=1}{k} A_l \to
    \indexedinhctp{F_l}{l=1}{k} B_l ,
  \end{displaymath}
  be the unique continuous linear extension of \(S_k\) to
  \begin{displaymath}
    \indexedinhctp{E_l}{l=1}{k} A_l .
  \end{displaymath}
  If \(k \in \{ 2, \ldots, n\}\) and \(S_k\) is not continuous let
  \(T_k = 0\).
  When all of the functions \(S_k\), \(T_k\), \(k = 1, \ldots, n\) are
  operators define
  \begin{displaymath}
    \indexedinhopctp{E_k}{F_k}{k=1}{n} P_k = T_n .
  \end{displaymath}
  If any of the functions \(S_k\), \(T_k\), \(k = 1, \ldots, n\) is
  not an operator then
  \begin{displaymath}
    \indexedinhopctp{E_k}{F_k}{k=1}{n} P_k
  \end{displaymath}
  is undefined.
\end{definition}

\if\longprep1
{
\givelemmaswithoutproofsn{three}
}
\fi

\if\shortprep0
{
\begin{lemma}
  \label{lem:tp-projection}
  Let \(n \in \positiveintegers\). Let \(A_1, \ldots, A_n\) and
  \(B_1, \ldots, B_n\) be Banach spaces so that \(B_k\) is a closed
  subspace of \(A_k\) for each \(k \in \setoneton{n}\).
  Let \(\alpha\) be a uniform crossnorm. Let
  \(P_k : A_k \to B_k\) be a continuous projection
  \projprep \(A_k\)
  onto \(B_k\)
  for each \(k \in \setoneton{n}\).
  Let
  \begin{displaymath}
    E_n \defequalns \indexedctp{\alpha}{k=1}{n} A_k ,
    \;
    F_n \defequalns \indexedctp{\alpha}{k=1}{n} B_k ,
    \;\;\textrm{and}\;\;
    Q_n \defequalns \indexedopctp{\alpha}{k=1}{n} P_k .
  \end{displaymath}
  Suppose that \(F_n \closedsubspace E_n\).
  Then \(Q_n\) is a continuous projection \projprep \(E_n\) onto
  \(F_n\).
\end{lemma}
}
\fi

\if\longversion1
{
\begin{proof}
  \(Q_n\) is a continuous linear function from \(E_n\) into \(F_n\).
  We shall prove by induction that \(Q_n\) is a projection onto
  \(F_n\).
  Suppose that the proposition is true for some
  \(n' \in \setoneton{n - 1} \) (induction assumption).
  Then \(Q_{n'+1} = Q_{n'} \otimes_\alpha P_{n'+1}\),
  \(E_{n'+1} = E_{n'} \ctp_\alpha A_{n'+1}\), and
  \(F_{n'+1} = F_{n'} \ctp_\alpha B_{n'+1}\).
  Let \(x \in F_{n'+1}\). There exists
  \((x_k)_{k=0}^\infty \subset F_{n'} \otimes_\alpha B_{n'+1}\) so
  that \(x_k \to x\) as \(k \to \infty\).
  Let \(k_1 \in \naturalnumbers\). Then
  \begin{displaymath}
    x_{k_1} = \sum_{j=1}^m f_j \otimes b_j
  \end{displaymath}
  where \(m \in \positiveintegers\) and \(f_j \in F_{n'}\),
  \(b_j \in B_{n'+1}\) for each \(j \in \setoneton{m}\).
  Now
  \begin{displaymath}
    Q_{n'+1} x_{k_1} = \sum_{j=1}^m (Q_{n'} f_j) \otimes (P_{n'+1} b_j)
    = \sum_{j=1}^m f_j \otimes b_j  = x_{k_1}
  \end{displaymath}
  where the second equality follows the induction assumption.
  Since \(k_1 \in \naturalnumbers\) was arbitrary and \(Q_{n'+1}\) is
  continuous \(Q_{n'+1} x = x\). Therefore the proposition is true for
  \(n'+1\) and hence for all \(n \in \positiveintegers\).
\end{proof}
}
\fi

\if\shortprep0
{
Uniform crossnorms behave well with isometric and topological
isomorphisms. The following two lemmas are related to this.
}
\fi

\if\shortprep0
{
\begin{lemma}
  \label{lem:isom-tp}
  Let \(\alpha\) be a uniform crossnorm.
  Let \(A, B, A'\), and \(B'\) be Banach spaces so that
  \(A \iisom A'\) and \(B \iisom B'\). Let \(\iota_A : A \to A'\) and
  \(\iota_B : B \to B'\) be isometric isomorphisms.
  Then \(\iota_A \otimes_\alpha \iota_B\) is an isometric isomorphism
  from \(A \ctp_\alpha B\) onto \(A' \ctp_\alpha B'\).
\end{lemma}
}
\fi

\if\longversion1
{
\begin{proof}
  Function \(\iota_A \otimes_\alpha \iota_B\) is a linear and
  continuous function from \(A \ctp_\alpha B\) into
  \(A' \ctp_\alpha B'\).

  We prove first that \(\iota_A \otimes_\alpha \iota_B\) is a
  surjection onto \(A' \ctp_\alpha B'\).
  Let \(v \in A' \ctp_\alpha B'\). There exists a sequence
  \((v_k)_{k=0}^\infty \subset A' \otimes_\alpha B'\) so that
  \(v_k \to v\) as \(k \to \infty\). Now
  \begin{displaymath}
    v_k = \sum_{j=1}^{m_k} x_{k,j} \otimes y_{k,j} \qspace \forall k
    \in \naturalnumbers
  \end{displaymath}
  where \(m_k \in \positiveintegers\) for all
  \(k \in \naturalnumbers\) and \(x_{k,j} \in A'\),
  \(y_{k,j} \in B'\) for
  all \(k \in \positiveintegers\) and \(j \in \setoneton{m_k}\).
  Let \(a_{k,j} = \iota_A^{-1}(x_{k,j})\) and
  \(b_{k,j} = \iota_B^{-1}(y_{k,j})\) for
  all \(k \in \positiveintegers\) and \(j \in \setoneton{m_k}\).
  Define
  \begin{displaymath}
    u_k := \sum_{j=1}^{m_k} a_{k,j} \otimes b_{k,j} \qspace \forall k
    \in \naturalnumbers
  \end{displaymath}
  Now \((\iota_A \otimes \iota_B)(u_k) = v_k\) and
  \((\iota_A^{-1} \otimes \iota_B^{-1})(v_k) = u_k\) for all
  \(k \in \naturalnumbers\).
  Since \(\iota_A^{-1}\) is an isometric isomorphism from \(A'\) onto
  \(A\) and \(\iota_B^{-1}\) is an isometric isomorphism from \(B'\) onto
  \(B\) it follows that \(\iota_A^{-1} \ctp_\alpha \iota_B^{-1}\) is a
  continuous linear function from \(A' \ctp_\alpha B'\) into
  \(A \ctp_\alpha B\). Therefore
  \(u_k = (\iota_A^{-1} \otimes_\alpha \iota_B^{-1})(v_k) \to (\iota_A^{-1}
  \otimes_\alpha \iota_B^{-1})(v) \in A \ctp_\alpha B\) as \(k \to \infty\).
  Furthermore, \((\iota_A \otimes \iota_B)(u_k) = v_k \to v\) as
  \(k \to \infty\) and it follows that
  \((\iota_A \otimes_\alpha \iota_B)((\iota_A^{-1} \otimes_\alpha \iota_B^{-1})(v))
  = v\). Hence \(\iota_A \otimes_\alpha \iota_B\) is a surjection onto
  \(A' \ctp_\alpha B'\).

  We then prove that \(\iota_A \otimes_\alpha \iota_B\) is distance
  preserving.  Let \(v \in A \ctp_\alpha B\). Since \(\alpha\) is a
  uniform crossnorm it follows that \(\norm{\iota_A \otimes_\alpha
    \iota_B} \leq \norm{\iota_A} \norm{\iota_B} = 1\) and
  \(\norm{\iota_A^{-1} \otimes_\alpha \iota_B^{-1}} \leq
  \norm{\iota_A^{-1}} \norm{\iota_B^{-1}} = 1\). Now
  \begin{displaymath}
    \alpha_{A',B'} \left( \left( \iota_A \otimes_\alpha \iota_B
      \right) \left( v \right) \right) \leq
    \norm{\iota_A \otimes_\alpha \iota_B} \alpha_{A,B}(v)
    \leq \alpha_{A,B}(v)
  \end{displaymath}
  and
  \begin{eqnarray*}
    \alpha_{A,B} (v) & = & \alpha_{A,B} \left( \left( \iota_A^{-1}
        \otimes_\alpha \iota_B^{-1} \right) \left( \left( \iota_A
          \otimes_\alpha \iota_B
        \right) \left( v \right) \right) \right) \\
    & \leq &
    \norm{\iota_A^{-1} \otimes_\alpha \iota_B^{-1}}
    \alpha_{A',B'} \left( \left( \iota_A
        \otimes_\alpha \iota_B
      \right) \left( v \right) \right) \\
    & \leq &
    \alpha_{A',B'} \left( \left( \iota_A
        \otimes_\alpha \iota_B
      \right) \left( v \right) \right) .
  \end{eqnarray*}
  Consequently
  \(\alpha_{A,B} (v) = \alpha_{A',B'} (( \iota_A
  \otimes_\alpha \iota_B) ( v ) )\). Thus function
  \(\iota_A \otimes_\alpha \iota_B\) is distance preserving and
  an isometric isomorphism from \(A \ctp_\alpha B\) onto
  \(A' \ctp_\alpha B'\).
\end{proof}
}
\fi

\if\shortprep0
{
\begin{lemma}
  \label{lem:top-isom-tp}
  Let \(\alpha\) be a uniform crossnorm.
  Let \(A, B, A'\), and \(B'\) be Banach spaces. Let
  \(\iota_A : A \to A'\) and
  \(\iota_B : B \to B'\) be topological isomorphisms.
  Then \(\iota_A \otimes_\alpha \iota_B\) is a topological isomorphism
  from \(A \ctp_\alpha B\) onto \(A' \ctp_\alpha B'\).
\end{lemma}
}
\fi

\if\longversion1
{
\begin{proof}
  Function \(\iota_A \otimes_\alpha \iota_B\) is a linear and
  continuous function from \(A \ctp_\alpha B\) into
  \(A' \ctp_\alpha B'\).
  The proof of the surjectivity of
  \begin{math}
    \iota_A \otimes_\alpha \iota_B
  \end{math}
  is similar to that in the proof of \mylemma \ref{lem:isom-tp}.

  We then prove that \(\iota_A \otimes_\alpha \iota_B\) defines a norm
  equivalence between spaces \(A \ctp_\alpha B\) and
  \(A' \ctp_\alpha B'\).
  We may assume that none of the spaces \(A\), \(B\), \(A'\), or
  \(B'\) is \(\zeroset\).
  We have \(\norm{\iota_A \otimes_\alpha \iota_B} > 0\) and
  \(\norm{\iota_A^{-1} \otimes_\alpha \iota_B^{-1}} > 0\)
  as neither of these two mappings is identically zero.
  Now
  \begin{math}
    \alpha_{A',B'} ( ( \iota_A \otimes_\alpha \iota_B ) ( v ) ) \leq
    \norm{\iota_A \otimes_\alpha \iota_B} \alpha_{A,B}(v)
  \end{math}
  and
  \begin{eqnarray*}
    \alpha_{A,B} (v) & = & \alpha_{A,B} \left( \left( \iota_A^{-1}
        \otimes_\alpha \iota_B^{-1} \right) \left( \left( \iota_A
          \otimes_\alpha \iota_B
        \right) \left( v \right) \right) \right) \\
    & \leq &
    \norm{\iota_A^{-1} \otimes_\alpha \iota_B^{-1}}
    \alpha_{A',B'} \left( \left( \iota_A
        \otimes_\alpha \iota_B
      \right) \left( v \right) \right)
  \end{eqnarray*}
  from which it follows that
  \begin{displaymath}
    \frac{1}{\norm{\iota_A^{-1} \otimes_\alpha \iota_B^{-1}}}
    \alpha_{A,B} (v)
    \leq
    \alpha_{A',B'} ( ( \iota_A
        \otimes_\alpha \iota_B ) ( v ) ) .
  \end{displaymath}
  Consequently
  \begin{eqnarray*}
    \frac{1}{\norm{\iota_A^{-1} \otimes_\alpha \iota_B^{-1}}}
    \alpha_{A,B} (v)
    & \leq &
    \alpha_{A',B'} \left( \left( \iota_A
        \otimes_\alpha \iota_B
      \right) \left( v \right) \right) \\
    & \leq &
    \norm{\iota_A^{-1} \otimes_\alpha \iota_B^{-1}}
    \alpha_{A',B'} \left( \left( \iota_A
        \otimes_\alpha \iota_B
      \right) \left( v \right) \right) .
  \end{eqnarray*}
  Hence function
  \(\iota_A \otimes_\alpha \iota_B\) is defines a norm equivalence
  between spaces \(A \ctp_\alpha B\) and
  \(A' \ctp_\alpha B'\).
  If
  \begin{math}
    x, y \in A \ctp_\alpha B
  \end{math},
  \begin{math}
    x \neq y
  \end{math}
  it follows that
  \begin{eqnarray*}
    \alpha_{A',B'} \left( \left( \iota_A
        \otimes_\alpha \iota_B
      \right) \left( x \right) - 
        \left( \iota_A
          \otimes_\alpha \iota_B
        \right) \left( y \right) \right)
    & = & \alpha_{A',B'} \left( \left( \iota_A
        \otimes_\alpha \iota_B
      \right) \left(x - y \right) \right) \\
    & \geq &
    \frac{1}{\norm{\iota_A^{-1} \otimes_\alpha \iota_B^{-1}}}
    \alpha_{A,B} \left( x - y \right)
    > 0
  \end{eqnarray*}
  and consequently
  \begin{math}
    ( \iota_A \otimes_\alpha \iota_B ) ( x ) \neq ( \iota_A \otimes_\alpha \iota_B ) ( y )
  \end{math}.
  Hence function
  \begin{math}
    \iota_A \otimes_\alpha \iota_B
  \end{math}
  is a topological isomorphism from
  \begin{math}
    A \ctp_\alpha B
  \end{math}
  onto
  \begin{math}
    A' \ctp_\alpha B'
  \end{math}.
\end{proof}
}
\fi

\if\shortprep0
{
Book \cite{ryan2002} has an exercise for proving isometric isomorphism
\(l^1 \cptp l^1 \iisom l^1\). However, we also need to construct the
isometric isomorphism for this article so we give the following
lemma without proof.
}
\fi

\if\shortprep0
{
\begin{lemma}
  \label{lem:lone-isomorphism}
  \definek
  Define function
  \begin{displaymath}
    \iota : \littlelp{1}{\naturalnumbers}{K}
    \to
    \littlelp{1}{\naturalnumbers}{K} \cptp \littlelp{1}{\naturalnumbers}{K}
  \end{displaymath}
  by
  \begin{displaymath}
    \iota(\nnatbasvec{k}) := \ntensornbv{\sqord{k}} ,
    \spaceafter k \in \naturalnumbers ,
  \end{displaymath}
  and extending by linearity and continuity onto whole
  \begin{math}
    \littlelp{1}{\naturalnumbers}{K}
  \end{math}.
  Then \(\iota\) is an isometric isomorphism from Banach space
  \begin{math}
    \littlelp{1}{\naturalnumbers}{K}
  \end{math}
  onto Banach space
  \begin{math}
    \littlelp{1}{\naturalnumbers}{K} \cptp \littlelp{1}{\naturalnumbers}{K}
  \end{math}.  
\end{lemma}
}
\fi

\if\shortprep0
{
This result can be generalized for the completed projective tensor product
of more than two 
\begin{math}
  \littlelp{1}{\naturalnumbers}{K}
\end{math}:
}
\fi

\if\shortprep0
{
\begin{lemma}
  \label{lem:lone-gen-isomorphism}
  \definekn
  Define function
  \begin{displaymath}
    \eta^{[n]} : \littlelp{1}{\naturalnumbers}{K}
    \to
    \indexedcptp{k=1}{n} \littlelp{1}{\naturalnumbers}{K}
  \end{displaymath}
  by
  \begin{displaymath}
    \eta^{[n]}(\nnatbasvec{k}) := \ntensornbv{\gensqord{n}{k}} ,
    \spaceafter k \in \naturalnumbers ,
  \end{displaymath}
  and extending by linearity and continuity onto whole
  \begin{math}
    \littlelp{1}{\naturalnumbers}{K}
  \end{math}.
  Then \(\eta\) is an isometric isomorphism from Banach space
  \begin{math}
    \littlelp{1}{\naturalnumbers}{K}
  \end{math}
  onto Banach space
  \begin{displaymath}
    \indexedcptp{k=1}{n} \littlelp{1}{\naturalnumbers}{K} .
  \end{displaymath}  
\end{lemma}
}
\fi

\if\shortprep0
{
\begin{proof}
  Case \(n = 1\) is trivial and case \(n = 2\) is true by \mylemma
  \ref{lem:lone-isomorphism}. The lemma can be proven using induction by
  \(n\) and equation
  \begin{displaymath}
    \forall n \in \naturalnumbers + 2 :
    \eta^{[n+1]}
    =
    (\eta^{[n]} \opcptp
    \idfunconspace{\littlelp{1}{\naturalnumbers}{K}})
    \circ
    \eta^{[2]} .
  \end{displaymath}
\end{proof}
}
\fi

\if\longprep1
{
  \givelemmawithoutproof
}
\fi

\if\shortprep0
{
\begin{lemma}
  \label{lem:lone-zn-isomorphism}
  \definekn
  Define function
  \begin{displaymath}
    \znisomfunc{n} : 
    \indexedcptp{k=1}{n} \littlelp{1}{\integernumbers}{K}
    \to
    \littlelp{1}{\zn}{K}
  \end{displaymath}
  by
  \begin{displaymath}
    \znisom{n}{\ztensornbv{\firstznvar}} := \znatbasvec{\firstznvar} ,
    \spaceafter \firstznvar \in \zn ,
  \end{displaymath}
  and extending by linearity and continuity onto whole
  \begin{math}
    \indexedcptp{k=1}{n} \littlelp{1}{\integernumbers}{K}
  \end{math}.
  Then \(\znisomfunc{n}\) is well-defined and it is an isometric isomorphism from
  Banach space
  \begin{math}
    \indexedcptp{k=1}{n} \littlelp{1}{\integernumbers}{K}
  \end{math}
  onto Banach space
  \begin{math}
    \littlelp{1}{\zn}{K} .
  \end{math}      
\end{lemma}
}
\fi

\if\longversion1
{
\begin{proof}
  Let \(\alpha\) be some bijection from
  \(\integernumbers\) onto \(\naturalnumbers\).
  Define function
  \begin{math}
    \beta : \littlelp{1}{\integernumbers}{K} \to
    \littlelp{1}{\naturalnumbers}{K}
  \end{math}
  by
  \begin{displaymath}
    \beta(\znatbasvec{k}) := \nnatbasvec{\alpha(k)}
  \end{displaymath}
  for all
  \begin{math}
    k \in \integernumbers
  \end{math}.
  and extending by linearity and continuity onto whole
  \begin{math}
    \littlelp{1}{\integernumbers}{K}
  \end{math}.
  Function \(\beta\) is an isometric isomorphism from
  Banach space
  \begin{math}
    \littlelp{1}{\integernumbers}{K}
  \end{math}
  onto Banach space
  \begin{math}
    \littlelp{1}{\naturalnumbers}{K}
  \end{math}.
  Define
  \begin{displaymath}
    \beta^{[n]} := \indexedopcptp{k=1}{n} \beta .
  \end{displaymath}
  Now \(\beta^{[n]}\) is an isometric isomorphism
  from Banach space
  \begin{displaymath}
    \indexedcptp{k=1}{n} \littlelp{1}{\integernumbers}{K}
  \end{displaymath}
  onto Banach space
  \begin{displaymath}
    \indexedcptp{k=1}{n} \littlelp{1}{\naturalnumbers}{K} .
  \end{displaymath}
  Let
  \begin{displaymath}
    \alpha^{[n]}(\secondznvar)
    :=
    \sum_{k=1}^n \alpha(\cartprodelem{\secondznvar}{k}) \onetonnbv{n}{k}
  \end{displaymath}
  for all
  \begin{math}
    \secondznvar \in \zn
  \end{math}.
  It can be shown that \(\alpha^{[n]}\) is a bijection from
  \(\zn\) onto \(\nn\).
  Let
  \begin{math}
    \rho := ( \alpha^{[n]} )^{-1} \circ \gensqordfunction{n}
  \end{math}.
  Now \(\rho\) is a bijection from \(\naturalnumbers\) onto \(\zn\).
  Define function
  \begin{math}
    \gamma : \littlelp{1}{\naturalnumbers}{K} \to
    \littlelp{1}{\zn}{K}
  \end{math}
  by
  \begin{displaymath}
    \gamma(\nnatbasvec{k}) := \zmnatbasvec{\rho(k)} ,
    \spaceafter k \in \naturalnumbers ,
  \end{displaymath}
  and extending by linearity and continuity onto whole
  \begin{math}
    \littlelp{1}{\naturalnumbers}{K}
  \end{math}.
  Function \(\gamma\) is an isometric isomorphism from Banach space
  \begin{math}
    \littlelp{1}{\naturalnumbers}{K}
  \end{math}
  onto Banach space
  \begin{math}
    \littlelp{1}{\zn}{K}
  \end{math}.
  Define function \(\eta^{[n]}\) as in \mylemma 
  \ref{lem:lone-gen-isomorphism}.
  Let
  \begin{math}
    \znisomfunc{n} := \gamma \circ (\eta^{[n]})^{-1} \circ \beta^{[n]}
  \end{math}
  Now function \(\znisomfunc{n}\) is an isometric isomorphism from Banach space
  \begin{displaymath}
    \indexedcptp{k=1}{n} \littlelp{1}{\integernumbers}{K}
  \end{displaymath}
  onto Banach space
  \begin{math}
    \littlelp{1}{\zn}{K}
  \end{math}.
\end{proof}
}
\fi

\if10
{
Book \cite{ryan2002} has an exercise for proving isometric isomorphism
\(l^1 \cptp l^1 \iisom l^1\).
}
\fi
\if\shortprep0
{
We give then a lemma about multiplication of series in spaces
topologically isomorphic to \(c_0\) or \(l^1\).
}
\fi

\if\shortprep0
{
\begin{lemma}
  \label{lem:tensor-product-of-sums}
  \definekn
  Let \(F\) be a Banach space with scalar field \(K\)
  and \(\alpha\) be a uniform crossnorm.
  Let \(X \equalns \gencospace{\naturalnumbers}{K}\) or
  \(X \equalns \littlelp{1}{\naturalnumbers}{K}\).
  Let \(\eta : F \to X\) be a
  topological isomorphism and suppose that
  \begin{math}
    \eta(f_j) = \nnatbasvec{j}
  \end{math}
  where
  \begin{math}
    f_j \in F
  \end{math}
  for each
  \begin{math}
    j \in \naturalnumbers
  \end{math}.
  Let \(\indexedseqstyle{a}{k} \in X\)
  for all
  \begin{math}
    k \in \setoneton{n}
  \end{math}.
  Then
  \begin{eqnarray*}
    \indexedtensorproduct{l=1}{n} \left( \sum_{j \in \naturalnumbers}
    \seqelem{\indexedseqstyle{a}{l}}{j} f_j \right)
    & = & \sum_{\firstznvar \in \nn} \left( \prod_{l=1}^n
      \seqelem{\indexedseqstyle{a}{l}}
      {\cartprodelem{\firstznvar}{l}}
      \right)
      f_{\cartprodelem{\firstznvar}{1}} \otimes \cdots
      \otimes f_{\cartprodelem{\firstznvar}{n}} \\
    & = & \sum_{\firstznvar \in \nn} \indexedtensorproduct{l=1}{n}
    \seqelem{\indexedseqstyle{a}{l}}
    {{\cartprodelem{\firstznvar}{l}}}
    f_{\cartprodelem{\firstznvar}{{l}}}
  \end{eqnarray*}
  where the second and third series are computed
  on space
  \begin{math}
    \indexedctp{\alpha}{k=1}{n} F
  \end{math}.
\end{lemma}
}
\fi

\if10
{
\begin{proof}
Use \mylemma \ref{lem:top-isom-tp}.
\end{proof}
}
\else
\if\shortprep0
{
\begin{proof}
  Let
  \begin{displaymath}
    \eta^{[n]} := \indexedopctp{\alpha}{k=1}{n} \eta , \;\;
    F^{[n]} \defequalns \indexedctp{\alpha}{k=1}{n} F , \;\; \textrm{and} \;
    X^{[n]} \defequalns \indexedctp{\alpha}{k=1}{n} X .
  \end{displaymath}
  By \mylemma \ref{lem:top-isom-tp}
  \begin{math}
    \eta^{[n]}
  \end{math}
  is a topological isomorphism from
  \begin{math}
    F^{[n]}
  \end{math}
  onto
  \begin{math}
    X^{[n]}
  \end{math}.
  
  Let
  \begin{displaymath}
    t := \indexedtensorproduct{k=1}{n} \left( \sum_{j \in \naturalnumbers}
    \seqelem{\indexedseqstyle{a}{k}}{j} f_j \right) .
  \end{displaymath}
  Now
\if\longversion0
{
  \if\elsevier1
  {
  	\begin{align*}
    	\eta^{[n]} \left( t \right)
    	& = \indexedtensorproduct{l=1}{n} \left(
     	   \eta \left( \sum_{j \in \naturalnumbers}
      	    \seqelem{\indexedseqstyle{a}{l}}{j} f_j \right)
      	\right) \\
    	& = \sum_{m_1 \in \naturalnumbers} \cdots
      	\sum_{m_n \in \naturalnumbers}
      	\left( \prod_{l=1}^n \seqelem{\indexedseqstyle{a}{l}}{m_l} \right)
      	\nnatbasvec{m_1} \otimes \cdots \otimes \nnatbasvec{m_n}
  	\end{align*}
  }
  \else
  {
  	\begin{displaymath}
    	\eta^{[n]} \left( t \right)
    	= \indexedtensorproduct{l=1}{n} \left(
     	   \eta \left( \sum_{j \in \naturalnumbers}
      	    \seqelem{\indexedseqstyle{a}{l}}{j} f_j \right)
      	\right)
    	= \sum_{m_1 \in \naturalnumbers} \cdots
      	\sum_{m_n \in \naturalnumbers}
      	\left( \prod_{l=1}^n \seqelem{\indexedseqstyle{a}{l}}{m_l} \right)
      	\nnatbasvec{m_1} \otimes \cdots \otimes \nnatbasvec{m_n}
  	\end{displaymath}
  }
  \fi
}
\else
{
  \begin{eqnarray*}
    \eta^{[n]} \left( t \right)
    & = & \indexedtensorproduct{l=1}{n} \left(
        \eta \left( \sum_{j \in \naturalnumbers}
          \seqelem{\indexedseqstyle{a}{l}}{j} f_j \right)
      \right)
    = \indexedtensorproduct{l=1}{n} \left(
        \sum_{j \in \naturalnumbers}
          \seqelem{\indexedseqstyle{a}{l}}{j}
          \nnatbasvec{j}
      \right) \\
    & = & \sum_{m_1 \in \naturalnumbers} \cdots
      \sum_{m_n \in \naturalnumbers}
      \left( \prod_{l=1}^n \seqelem{\indexedseqstyle{a}{l}}{m_l} \right)
      \nnatbasvec{m_1} \otimes \cdots \otimes \nnatbasvec{m_n}
  \end{eqnarray*}
}
\fi
  and
  \begin{displaymath}
    \szseqelem{\left( \eta^{[n]} \left( t \right) \right)}{\firstznvar}
    = \prod_{l=1}^n
      \szseqelem{\indexedseqstyle{a}{l}}
      {\szcartprodelem{\firstznvar}{l}}
  \end{displaymath}
  for each
  \begin{math}
    \firstznvar \in \nn
  \end{math}.
  Consequently
  \begin{displaymath}
    \eta^{[n]} \left( t \right)
    = \sum_{\firstznvar \in \nn}
      \left( \prod_{l=1}^n 
        \seqelem{\indexedseqstyle{a}{l}}
        {\cartprodelem{\firstznvar}{l}}
      \right)
      \nnatbasvec{\cartprodelem{\firstznvar}{1}}
      \otimes \cdots \otimes
      \nnatbasvec{\cartprodelem{\firstznvar}{n}} .
  \end{displaymath}
  It follows that
  \if\elsevier0
  {
  \begin{eqnarray*}
    t & = & \left( \eta^{[n]} \right)^{-1} \left( \eta^{[n]} \left( t
      \right) \right)
    = \sum_{\lambda \in \nn}
      \left( \prod_{l=1}^n
        \seqelem{\indexedseqstyle{a}{l}}
        {\cartprodelem{\firstznvar}{l}}
        \right)
      f_{\cartprodelem{\firstznvar}{1}}
      \otimes \cdots \otimes
      f_{\cartprodelem{\firstznvar}{n}} \\
    & = & \sum_{\firstznvar \in \nn} \indexedtensorproduct{l=1}{n}
        \seqelem{\indexedseqstyle{a}{l}}
        {\cartprodelem{\firstznvar}{l}}      
        f_{\cartprodelem{\firstznvar}{l}} .
  \end{eqnarray*}
  }
  \else
  {
  \begin{align*}
    t & = \left( \eta^{[n]} \right)^{-1} \left( \eta^{[n]} \left( t
      \right) \right)
    = \sum_{\lambda \in \nn}
      \left( \prod_{l=1}^n
        \seqelem{\indexedseqstyle{a}{l}}
        {\cartprodelem{\firstznvar}{l}}
        \right)
      f_{\cartprodelem{\firstznvar}{1}}
      \otimes \cdots \otimes
      f_{\cartprodelem{\firstznvar}{n}} \\
    & = \sum_{\firstznvar \in \nn} \indexedtensorproduct{l=1}{n}
        \seqelem{\indexedseqstyle{a}{l}}
        {\cartprodelem{\firstznvar}{l}}      
        f_{\cartprodelem{\firstznvar}{l}} .
  \end{align*}
  }
  \fi
\end{proof}
}
\fi
\fi

}
\fi

\section{General Definitions for a Compactly Supported Interpolating MRA}
\label{sec:gen-def-mra}

\subsection{Mother Scaling Function}

\begin{definition}
  \definekn
  A \defterm{compactly supported interpolating mother scaling function}  
  is a function \(\myphi \in \cscfunc{\rn}{K}\) satisfying the following
  conditions:
\begin{myenumerate}
  \myitem{(MSF.1)}
    \begin{displaymath}
 		\forall \firstznvar \in \zn : \myphi(\firstznvar) = 
 		\delta_{\firstznvar,0}
 	\end{displaymath}
  \myitem{(MSF.2)}
    \begin{displaymath}
      \forall \rnx \in \rn : \myphi ( \rnx )
      = \sum_{\firstznvar \in \zn} 
      \myphi \left( \frac{\firstznvar}{2} \right)
      \myphi ( 2 \rnx - \firstznvar ) .
    \end{displaymath}
\end{myenumerate}
\end{definition}

\if10
\begin{definition}
\label{def:msf}
\definekn
Let \(\myphi\) be a compactly supported interpolating mother scaling function.
Define
\begin{eqnarray*}
  & & \largevspace{n}{j} :=
  \left\{
  	\rnx \in \rn \mapsto
  	\sum_{\firstznvar \in \zn} \seqelem{\seqstyle{a}}{\firstznvar}
  	\myphi(2^j \rnx - \firstznvar)
  	\setsep \seqstyle{a} \in \seqset{\zn}{K}
  \right\} \\
  & & \norminspace{f}{\largevspace{n}{j}} := \norminfty{f},
  \spaceafter f \in \largevspace{n}{j}
  .
\end{eqnarray*}
\end{definition}
\fi

\if10
\begin{definition}
  \label{def:genvjspace}
	\definekn
	Let \(\myphi\) be a compactly supported interpolating mother scaling
        function.
	Define
	\begin{displaymath}
	  (\largevprojop{n}{j} f)(\rnx) :=
	  \sum_{\firstznvar \in \zn} f \left( \frac{\firstznvar}{2^j} \right)
	  \myphi(2^j \rnx - \firstznvar)
	\end{displaymath}
	for all \(j \in \integernumbers\),
	\(\rnx \in \rn\),
	and \(f \in K^{\rn}\).
\end{definition}
\fi

The Deslauriers-Dubuc fundamental functions satisfy these conditions
\cite{donoho1992}.

\subsection{General Definitions for the Univariate MRA's}

We shall denote the function space for which the MRA is defined by \(E\) in this
section. We have either
\(E = \ucfunc{\realnumbers}{K}\)
or
\(E = \vanishingfunc{\realnumbers}{K}\)
where \(K = \realnumbers\) or \(K = \complexnumbers\).
We shall assume that \(\myphi \in \cscfunc{\realnumbers}{K}\)
is a compactly supported interpolating mother scaling function
throughout this subsection.

\begin{definition}
  Define function \(\psi \in \cscfunc{\realnumbers}{K}\) by
  \begin{displaymath}
    \psi(x) := \myphi(2 x - 1)
  \end{displaymath}
  for all \(x \in \realnumbers\).
  Function \(\psi\) is called \defterm{the mother wavelet}.
\end{definition}

\begin{definition}
  Define
  \begin{displaymath}
    \myphi_{j,k} := \myphi(2^j \cdot - k)
  \end{displaymath}
  and
  \begin{displaymath}
    \psi_{j,k} := \psi(2^j \cdot - k)
  \end{displaymath}
  for all \(j \in \integernumbers\) and
  \(k \in \integernumbers\).
\end{definition}

\if\shortprep1
{
\begin{definition}
  Define
  \begin{displaymath}
    \genmotherwaveletonedim{s} :=
    \left\{
      \begin{array}{ll}
        \myphi ; \;\; & s = 0 \\
        \psi ; \;\; & s = 1 
      \end{array}
    \right.
  \end{displaymath}
  where \(s \in \{ 0, 1 \}\).
\end{definition}
}
\else
{
\begin{definition}
  Define
  \begin{displaymath}
    \genmotherwaveletonedim{s} :=
    \left\{
      \begin{array}{ll}
        \myphi ; \;\; & s = 0 \\
        \psi ; \;\; & s = 1 
      \end{array}
    \right.
  \end{displaymath}
  and
  \begin{displaymath}
    \genwaveletonedim{s}{j}{k} :=
    \left\{
      \begin{array}{ll}
        \myphi_{j,k} ; \;\; & s = 0 \\
        \psi_{j,k} ; \;\; & s = 1 
      \end{array}
    \right.
  \end{displaymath}
  where \(s \in \{ 0, 1 \}\),
  \(j \in \integernumbers\),
  and \(k \in \integernumbers\).
\end{definition}
}
\fi

\begin{definition}
  \label{def:onedimwaveletfilters-a}
  When \(k \in \integernumbers\) define
  \begin{eqnarray*}
    \hfilterelem{k} & := & \myphi \left( \frac{k}{2} \right) \\
    \gfilterelem{k} & := & \delta_{k,1} \\
    \htfilterelem{k} & := & \delta_{k,0} \\
    \gtfilterelem{k} & := & (-1)^{k-1} h_{1-k} .
  \end{eqnarray*}
\end{definition}

\if\shortprep0
{
\begin{definition}
  \label{def:onedimwaveletfilters}
	When \(t \in \{ 0, 1\}\) and \(k \in \integernumbers\)
	define
	\begin{displaymath}
	  \onedimgenwaveletfilterelem{t}{k} :=
	    \left\{
	    \begin{array}{ll}
	      \hfilterelem{k} ; \;\; & t = 0 \\
	      \gfilterelem{k} ; \;\; & t = 1
	    \end{array}
	    \right.
	\end{displaymath}
	and
	\begin{displaymath}
	  \onedimgendualwaveletfilterelem{t}{k} :=
	    \left\{
	    \begin{array}{ll}
	      \htfilterelem{k} ; \;\; & t = 0 \\
	      \gtfilterelem{k} ; \;\; & t = 1
	    \end{array}
	    \right. .
	\end{displaymath}
\end{definition}
}
\fi

\begin{definition}
  Define
  \begin{math}
    \onedimmotherdualsf := \delta \in \topdual{E}
  \end{math}
  and
  \begin{math}
    \onedimmotherdualwavelet \in \topdual{E}
  \end{math}
  by
  \begin{equation}
    \label{eq:onedimmotherdualwavelet}
    \onedimmotherdualwavelet :=
    2 \sum_{k \in \integernumbers} \gtfilterelem{k}
    \onedimmotherdualsf(2 \cdot - k) .
  \end{equation}
  Define
  \begin{math}
    \onedimdualsf{j}{k} := 2^j \onedimmotherdualsf(2^j \cdot - k)
  \end{math}
  and
  \begin{math}
    \onedimdualwavelet{j}{k}
    := 2^j \onedimmotherdualwavelet(2^j \cdot - k)
  \end{math}
  where \(j, k \in \integernumbers\).
\end{definition}

As \(\onedimmothersf\) is compactly supported only a finite number
of numbers \(\hfilterelem{k}\), \(k \in \integernumbers\), and the other
three filters defined by \mydef \ref{def:onedimwaveletfilters-a} are nonzero.
Consequently the series in \myequation \eqref{eq:onedimmotherdualwavelet}
has a finite number of nonzero terms.

\if\shortprep1
{
\begin{definition}
  \label{def:genwaveletonedim}
  Define
  \begin{displaymath}
    \genmotherdualwaveletonedim{s} :=
    \left\{
      \begin{array}{ll}
        \onedimmotherdualsf ; \;\; & s = 0 \\
        \onedimmotherdualwavelet ; \;\; & s = 1
      \end{array}
    \right.
  \end{displaymath}
  where \(s \in \zerooneset\).
\end{definition}
}
\else
{
\begin{definition}
  \label{def:genwaveletonedim}
  Define
  \begin{displaymath}
    \genmotherdualwaveletonedim{s} :=
    \left\{
      \begin{array}{ll}
        \onedimmotherdualsf ; \;\; & s = 0 \\
        \onedimmotherdualwavelet ; \;\; & s = 1
      \end{array}
    \right.
  \end{displaymath}
  and
  \begin{displaymath}
    \gendualwaveletonedim{s}{j}{k} :=
    \left\{
      \begin{array}{ll}
        \onedimdualsf{j}{k} ; \;\; & s = 0 \\
        \onedimdualwavelet{j}{k} ; \;\; & s = 1
      \end{array}
    \right.
  \end{displaymath}
  where \(j, k \in \integernumbers\) and \(s \in \zerooneset\).
\end{definition}
}
\fi

\if\shortprep0
{
\begin{lemma}
  \label{lem:orthogonality-one-dim}
  Let \(k, l \in \integernumbers\). Then
  \begin{itemize}
  \item[(i)]
    \begin{displaymath}
      \dualappl{\phidual_{j,k}}{\myphi_{j,l}} =
      \delta_{k,l}
    \end{displaymath}
  \item[(ii)]
    \begin{displaymath}
      \dualappl{\phidual_{j,k}}{\psi_{j,l}} =
      0
    \end{displaymath}
  \item[(iii)]
    \begin{displaymath}
      \dualappl{\psidual_{j,k}}{\myphi_{j,l}} =
      0
    \end{displaymath}
  \item[(iv)]
    \begin{displaymath}
      \dualappl{\psidual_{j,k}}{\psi_{j,l}} =
      \delta_{k,l}
    \end{displaymath}
  \end{itemize}
\end{lemma}
}
\fi

\if\shortprep0
{
\begin{proof}
  See also \cite[section 2]{cl1996}
  and \cite{goedecker1998}.
  \begin{itemize}
  \item[(i)]
    \begin{displaymath}
      \dualappl{\delta(2^j \cdot - k)}{\myphi(2^j \cdot - l)}
      =
      \dualappl{\frac{1}{2^j} \delta \left( \cdot -
          \frac{k}{2^j} \right) }{\myphi(2^j \cdot - l)}
      = \frac{1}{2^j} \myphi(k - l)
      = \frac{1}{2^j} \delta_{k,l}
    \end{displaymath}
    The last equality follows from (MSF.1).
  \item[(ii)]
    As \(2k - 2l - 1 \neq 0\) it follows from (MSF.1) that 
    \begin{displaymath}
      \dualappl{\delta(2^j \cdot - k)}{\psi(2^j \cdot - l)}
      = \frac{1}{2^j} \psi(k - l)
      = \frac{1}{2^j} \myphi(2k - 2l - 1)
      = 0 .
    \end{displaymath}
  \item[(iii)]
    \begin{eqnarray*}
      \dualappl{\psidual(2^j \cdot - k)}{\myphi(2^j \cdot - l)}
      & = &
      2 \sum_{\nu \in \integernumbers} \gdual_\nu
      \frac{1}{2^{j+1}} \dualappl{\delta(\cdot - \frac{2k +
          \nu}{2^{j+1}})}{\myphi(2^j \cdot - l)} \\
      & = & \frac{1}{2^j} \sum_{\nu \in \integernumbers} \gdual_\nu
      \myphi \left( k - l + \frac{\nu}{2} \right) \\
      & = & \frac{1}{2^j} \sum_{\nu \in \integernumbers} \gdual_\nu
      \sum_{\mu \in \integernumbers}
      \myphi \left( \frac{\mu}{2} \right)
      \myphi \left( 2k - 2l + \nu - \mu \right) \\
      & = &
      \frac{1}{2^j}
      \sum_{\nu \in \integernumbers}
      \sum_{\mu \in \integernumbers}
      (-1)^{\nu-1} \myphi \left( \frac{1-\nu}{2} \right)
      \myphi \left( \frac{\mu}{2} \right)
      \delta_{2l-2k,\nu-\mu} \\
      & = &
      \frac{1}{2^j}
      \sum_{\nu \in \integernumbers}
      (-1)^{\nu-1} \myphi \left( \frac{1-\nu}{2} \right)
      \myphi \left( \frac{\nu-2l+2k}{2} \right) \\
      & = &
      \frac{1}{2^j}
      \biggl(
        \myphi \left( \frac{1-2l+2k}{2} \right) \\
        & & - \sum_{\nu \in \integernumbers}
        \myphi \left( \frac{1-2 \nu}{2} \right)
        \myphi \left( 2 \nu' - 2l + 2k \right)
      \biggl) \\
      & = &
      \frac{1}{2^j}
      \left(
        \myphi \left( \frac{1-2l+2k}{2} \right)
        - \myphi \left( \frac{1-2l+2k}{2} \right)
      \right) = 0
    \end{eqnarray*}
    The 6th equality follows from (MSF.1).
  \item[(iv)]
    \if\elsevier0
    {
    \begin{align*}
      \dualappl{\psidual(2^j \cdot - k)}{\psi(2^j \cdot - l)}
      & = \szdualappl{2 \sum_{\nu \in \integernumbers} \gdual_\nu
        \delta(2^{j+1} \cdot - 2k - \nu)}{\myphi(2^{j+1} \cdot - 2l -
        1)} \\
      & = \frac{1}{2^j}
      \sum_{\nu \in \integernumbers} \gdual_\nu \myphi(2k + \nu - 2l -
      1) \\
      & = \frac{1}{2^j}
      \sum_{\nu \in \integernumbers} \gdual_\nu \delta_{\nu,2l-2k+1}
      =
      \frac{1}{2^j} \delta_{k,l}
    \end{align*}
    }
    \else
    {
    \begin{align*}
      \dualappl{\psidual(2^j \cdot - k)}{\psi(2^j \cdot - l)}
      & = \szdualappl{2 \sum_{\nu \in \integernumbers} \gdual_\nu
        \delta(2^{j+1} \cdot - 2k - \nu)}{\myphi(2^{j+1} \cdot - 2l -
        1)} \\
      & = \frac{1}{2^j}
      \sum_{\nu \in \integernumbers} \gdual_\nu \myphi(2k + \nu - 2l -
      1)
      = \frac{1}{2^j}
      \sum_{\nu \in \integernumbers} \gdual_\nu \delta_{\nu,2l-2k+1} \\
      & =
      \frac{1}{2^j} \delta_{k,l}
    \end{align*}
    }
    \fi
    The last equality follows from the fact that \(\gdual_\nu =
    0\) for all \(\nu \in \integernumbers\), \(\nu\) odd, and \(\nu
    \neq 1\).
  \end{itemize}
\end{proof}
}
\fi

\subsection{General Definitions for Multivariate MRA's}

We will assume that \(n \in \positiveintegers\) and either
\(K = \realnumbers\) or \(K = \complexnumbers\) throughout
this subsection. We will also assume that
\(E \equalns \ucfunc{\realnumbers}{K}\) or
\(E \equalns \vanishingfunc{\realnumbers}{K}\).
Furthermore,
\(\onedimmothersf \in \cscfunc{\realnumbers}{K}\)
shall be a compactly supported interpolating mother scaling function
throughout this subsection.

We set
\begin{displaymath}
	F \defequalns \left\{
	  \begin{array}{ll}
	    \ucfunc{\rn}{K} ; & E = \ucfunc{\realnumbers}{K} \\
	    \vanishingfunc{\rn}{K} ; & E = \vanishingfunc{\realnumbers}{K} \\
	  \end{array}
	\right.
\end{displaymath}

\begin{definition}
  \label{def:mdimsf}
  Define function \(\mdimmothersf{n} \in \cscfunc{\rn}{K}\) by
  \begin{displaymath}
    \mdimmothersf{n} := \indexedtensorproduct{k=1}{n} \onedimmothersf .
  \end{displaymath}
  Function \(\mdimmothersf{n}\) is called
  \defterm{an \(n\)-dimensional tensor product mother scaling
  function generated by \(\onedimmothersf\)}.
  Define also
  \if\shortprep1
  {
  \begin{math}
    \mdimsf{n}{j}{\firstznvar} :=
    \mdimmothersf{n}(2^j \cdot - \firstznvar)
  \end{math}
  }
  \else
  {
  \begin{displaymath}
    \mdimsf{n}{j}{\firstznvar} :=
    \mdimmothersf{n}(2^j \cdot - \firstznvar)
  \end{displaymath}
  }
  \fi
  where \(j \in \integernumbers\) and
  \(\firstznvar \in \zn\).
\end{definition}

Function
\begin{math}
  \mdimmothersf{n}
\end{math}
is a compactly supported interpolating mother scaling function on
\begin{math}
  \rn
\end{math}.

\if\shortprep1
{
  \begin{definition}
    When
    \begin{math}
      \zos \in \zeroonesetn
    \end{math}
    define function
    \begin{math}
      \mdimmotherwavelet{n}{\zos} \in \cscfunc{\rn}{K}
    \end{math}
    by
    \begin{displaymath}
      \mdimmotherwavelet{n}{\zos} :=
      \indexedtensorproduct{k=1}{n}
      \genmotherwaveletonedim{\cartprodelem{\zos}{k}}
    \end{displaymath}
    and
    \begin{math}
      \mdimgenwavelet{n}{\zos}{j}{\firstznvar} :=
      \mdimmotherwavelet{n}{\zos}(2^j \cdot - \firstznvar)
    \end{math}.
  \end{definition}
}
\else
{
  \begin{definition}
    Define function
    \begin{math}
      \mdimmotherwavelet{n}{\zos} \in \cscfunc{\rn}{K}
    \end{math}
    by
    \begin{displaymath}
      \mdimmotherwavelet{n}{\zos} :=
      \indexedtensorproduct{k=1}{n}
      \genmotherwaveletonedim{\cartprodelem{\zos}{k}}
    \end{displaymath}
    for all
    \begin{math}
      \zos \in \zeroonesetn
    \end{math}.
    Define functions
    \begin{math}
      \mdimgenwavelet{n}{\zos}{j}{\firstznvar} \in \cscfunc{\rn}{K}
    \end{math}
    by
    \begin{displaymath}
      \mdimgenwavelet{n}{\zos}{j}{\firstznvar} :=
      \indexedtensorproduct{l=1}{n}
      \genwaveletonedim{\cartprodelem{\zos}{l}}{j}{\cartprodelem{\firstznvar}{l}}
    \end{displaymath}
    for all
    \begin{math}
      \zos \in \zeroonesetn
    \end{math},
    \begin{math}
      j \in \integernumbers
    \end{math},
    and
    \begin{math}
      \firstznvar \in \zn
    \end{math}.
  \end{definition}
}
\fi

\if\shortprep0
{
We have
\begin{equation}
  \label{eq:gen-wavelet-card-interp}
  \mdimmotherwavelet{n}{\zos} \left( \firstznvar + \frac{1}{2} \zos \right)
  = \delta_{\firstznvar,0}
\end{equation}
for all \(\zos \in \zeroonesetn\) and
\(\firstznvar \in \zn\).
}
\fi

\if\shortprep0
{
\begin{definition}
  \label{def:mdimwaveletfilters}
  Let \(n \in \positiveintegers\),
  \(\zos \in \zeroonesetn\), and
  \(\znt \in \zn\).
  Define
  \begin{displaymath}
    \mdimgenwaveletfilterelem{n}{\zos}{\znt}
    :=
    \prod_{l=1}^n
    \onedimgenwaveletfilterelem{\cartprodelem{\zos}{l}}
    {\cartprodelem{\znt}{l}}
  \end{displaymath}
  and
  \begin{displaymath}
    \mdimgendualwaveletfilterelem{n}{\zos}{\znt}
    :=
    \prod_{l=1}^n
    \onedimgendualwaveletfilterelem{\cartprodelem{\zos}{l}}
    {\cartprodelem{\znt}{l}} .
  \end{displaymath}
\end{definition}
}
\fi

\if\longprep1
{
\givelemmawithoutproof
}
\fi

\if\shortprep0
{
\begin{lemma}
  \label{lem:mdim-wavelet-formulas}
  Suppose that \(m \in \positiveintegers\) so that
  \begin{displaymath}
    \forall k \in \integernumbers :
    \abs{k} \geq m \implies \hfilterelem{k} = 0 .
  \end{displaymath}
  Then
  \if\shortprep1
  {
  \begin{eqnarray*}
    \mdimmotherwavelet{n}{\zos}
    & = &
    \sum_{\secondznvar \in \setplusminusnpower{m}{n}}
    \mdimgenwaveletfilterelem{n}{\zos}
    {\secondznvar}
    \mdimsf{n}{1}{\secondznvar} \\
    \mdimgenwavelet{n}{\zos}{j}{\firstznvar}
    & = &
    \sum_{\secondznvar \in \setplusminusnpower{m}{n}}
    \mdimgenwaveletfilterelem{n}{\zos}
    {\secondznvar}
    \mdimsf{n}{j+1}{2 \firstznvar + \secondznvar}
    =
    \sum_{\thirdznvar \in \zn}
    \mdimgenwaveletfilterelem{n}{\zos}
    {\thirdznvar - 2 \firstznvar}
    \mdimsf{n}{j+1}{\thirdznvar}
  \end{eqnarray*}
  }
  \else
  {
  \begin{eqnarray*}
    \mdimmotherwavelet{n}{\zos}
    & = &
    \sum_{\secondznvar \in \setplusminusnpower{m}{n}}
    \mdimgenwaveletfilterelem{n}{\zos}
    {\secondznvar}
    \mdimsf{n}{1}{\secondznvar} \\
    \mdimgenwavelet{n}{\zos}{j}{\firstznvar}
    & = &
    \sum_{\secondznvar \in \setplusminusnpower{m}{n}}
    \mdimgenwaveletfilterelem{n}{\zos}
    {\secondznvar}
    \mdimsf{n}{j+1}{2 \firstznvar + \secondznvar}
    =
    \sum_{\thirdznvar \in \zn}
    \mdimgenwaveletfilterelem{n}{\zos}
    {\thirdznvar - 2 \firstznvar}
    \mdimsf{n}{j+1}{\thirdznvar} \\
    \mdimmothersf{n} & = &
    \sum_{\secondznvar \in \setplusminusnpower{m}{n}}
    \mdimgenwaveletfilterelem{n}{\finitezeroseq{n}}{\secondznvar}
    \mdimsf{n}{1}{\secondznvar} \\
    \mdimsf{n}{j}{\firstznvar} & = &
    \sum_{\secondznvar \in \setplusminusnpower{m}{n}}
    \mdimgenwaveletfilterelem{n}{\finitezeroseq{n}}{\secondznvar}
    \mdimsf{n}{j+1}{2 \firstznvar + \secondznvar}
    =
    \sum_{\thirdznvar \in \zn}
    \mdimgenwaveletfilterelem{n}{\finitezeroseq{n}}
    {\thirdznvar - 2 \firstznvar}
    \mdimsf{n}{j+1}{\thirdznvar} ,
  \end{eqnarray*}
  }
  \fi
  for all
  \(\zos \in \zeroonesetn\),
  \(j \in \integernumbers\),
  and \(\firstznvar \in \zn\).
\end{lemma}
}
\fi

The domain of the Dirac \(\delta\) functional varies in this article.
I.e. we may keep \(\delta\) as an element of different dual spaces \(A^*\).
When \(z_1, \ldots, z_m \in \realnumbers\) we will identify
\begin{math}
	\delta(\cdot - z_1) \otimes \ldots \otimes 
	\delta(\cdot - z_m)
\end{math}
with
\begin{math}
  \delta(\cdot - (z_1, \ldots, z_m))
\end{math}.

\if\shortprep1
{
  \begin{definition}
    Define \(\mdimmotherdualsf{n} \in \topdual{F}\) by
    \begin{displaymath}
      \mdimmotherdualsf{n} := \indexedtensorproduct{l=1}{n} \onedimmotherdualsf
    \end{displaymath}
    and
    \(\mdimdualsf{n}{j}{\firstznvar} \in \topdual{F}\)
    by
    \begin{math}
      \mdimdualsf{n}{j}{\firstznvar} :=
      2^j \mdimmotherdualsf{n}(2^j \cdot - \firstznvar)
    \end{math}
    where \(j \in \integernumbers\) and
    \(\firstznvar \in \zn\).
    Define also
    \(\mdimmotherdualwavelet{n}{\zos} \in \topdual{F}\)
    by
    \begin{displaymath}
      \mdimmotherdualwavelet{n}{\zos} :=
      \indexedtensorproduct{l=1}{n} \genmotherdualwaveletonedim{\cartprodelem{\zos}{l}}
    \end{displaymath}
    where \(\zos \in \zeroonesetn\)
    and
    \(\mdimgendualwavelet{n}{\zos}{j}{\firstznvar} \in \topdual{F}\)
    by
    \begin{math}
      \mdimgendualwavelet{n}{\zos}{j}{\firstznvar} :=
      2^j \mdimmotherdualwavelet{n}{\zos}(2^j \cdot - \firstznvar)
    \end{math}
    where \(\zos \in \zeroonesetn\),
    \(j \in \integernumbers\),
    and
    \(\firstznvar \in \zn\).
  \end{definition}
}
\else
{
  \begin{definition}
    Define \(\mdimmotherdualsf{n} \in \topdual{F}\) by
    \begin{displaymath}
      \mdimmotherdualsf{n} := \indexedtensorproduct{l=1}{n} \onedimmotherdualsf
    \end{displaymath}
    and
    \(\mdimdualsf{n}{j}{\firstznvar} \in \topdual{F}\)
    by
    \begin{displaymath}
      \mdimdualsf{n}{j}{\firstznvar} :=
      \indexedtensorproduct{l=1}{n} \onedimdualsf{j}{\cartprodelem{\firstznvar}{l}}
    \end{displaymath}
    where \(j \in \integernumbers\) and
    \(\firstznvar \in \zn\).
    Define also
    \(\mdimmotherdualwavelet{n}{\zos} \in \topdual{F}\)
    by
    \begin{displaymath}
      \mdimmotherdualwavelet{n}{\zos} :=
      \indexedtensorproduct{l=1}{n} \genmotherdualwaveletonedim{\cartprodelem{\zos}{l}}
    \end{displaymath}
    where \(\zos \in \zeroonesetn\)
    and
    \(\mdimgendualwavelet{n}{\zos}{j}{\firstznvar} \in \topdual{F}\)
    by
    \begin{displaymath}
      \mdimgendualwavelet{n}{\zos}{j}{\firstznvar} :=
      \indexedtensorproduct{l=1}{n}
      \gendualwaveletonedim{\cartprodelem{\zos}{l}}{j}{\cartprodelem{\firstznvar}{l}}
    \end{displaymath}
    where \(\zos \in \zeroonesetn\),
    \(j \in \integernumbers\),
    and
    \(\firstznvar \in \zn\).
  \end{definition}
}
\fi

\if\shortprep0
{
Note that \(\mdimmotherdualsf{n} = \delta\) (the Dirac delta functional),
\begin{displaymath}
  \mdimdualsf{n}{j}{\firstznvar} = 2^{nj} \delta(2^j \cdot - \firstznvar)
  = \delta \left( \cdot - \frac{\firstznvar}{2^j} \right) ,
\end{displaymath}
and
\if\elsevier1
{
\begin{eqnarray}
  \label{eq:mdim-dual-wavelet-formula}
  \mdimmotherdualwavelet{n}{\zos}
  & = &
  \sum_{\secondznvar \in \zn}
  \mdimgendualwaveletfilterelem{n}{\zos}{\secondznvar}
  \mdimdualsf{n}{1}{\secondznvar} \\
  \nonumber
  \mdimgendualwavelet{n}{\zos}{j}{\firstznvar}
  & = & 2^{nj} \mdimmotherdualwavelet{n}{\zos} \left(
    2^j \cdot - \firstznvar \right)
  = \sum_{\secondznvar \in \zn} 
  \mdimgendualwaveletfilterelem{n}{\zos}{\secondznvar}
  \delta \left( \cdot - \frac{2 \firstznvar + 
  \secondznvar}{2^{j+1}}\right) \\
  \label{eq:mdim-dual-wavelet-formula-b}
  & = &
  \sum_{\thirdznvar \in \zn}
  \mdimgendualwaveletfilterelem{n}{\zos}
  {\thirdznvar - 2 \firstznvar}
  \mdimdualsf{n}{j+1}{\thirdznvar}
\end{eqnarray}
}
\else
{
\begin{eqnarray}
  \label{eq:mdim-dual-wavelet-formula}
  \mdimmotherdualwavelet{n}{\zos}
  & = &
  \sum_{\secondznvar \in \zn}
  \mdimgendualwaveletfilterelem{n}{\zos}{\secondznvar}
  \mdimdualsf{n}{1}{\secondznvar} \\
  \nonumber
  \mdimgendualwavelet{n}{\zos}{j}{\firstznvar}
  & = & 2^{nj} \mdimmotherdualwavelet{n}{\zos} \left(
    2^j \cdot - \firstznvar \right)
  = \sum_{\secondznvar \in \zn} 
  \mdimgendualwaveletfilterelem{n}{\zos}{\secondznvar}
  \delta \left( \cdot - \frac{2 \firstznvar + 
  \secondznvar}{2^{j+1}}\right) \\
  \label{eq:mdim-dual-wavelet-formula-b}
  & = &
  \sum_{\thirdznvar \in \zn}
  \mdimgendualwaveletfilterelem{n}{\zos}
  {\thirdznvar - 2 \firstznvar}
  \mdimdualsf{n}{j+1}{\thirdznvar}
\end{eqnarray}
}
\fi
for all
\(j \in \integernumbers\),
\(\zos \in \zeroonesetn\),
and \(\firstznvar \in \zn\).
We also have
\begin{eqnarray}
  \label{eq:mdim-wavelet-orthogonality}
  \szdualappl{\mdimgendualwavelet{n}{\zos}{j}{\secondznvar}}
    {\mdimgenwavelet{n}{\zot}{j}{\firstznvar}}
  & = & \delta_{\zos,\zot}
  \delta_{\secondznvar,\firstznvar} \\
  \label{eq:mdim-wavelet-grid-value}
  \szdualappl{\mdimdualsf{n}{j+1}{\secondznvar}}
    {\mdimgenwavelet{n}{\zos}{j}{\firstznvar}}
  & = & \mdimgenwavelet{n}{\zos}{j}{\firstznvar}
    \left( \frac{\secondznvar}{2^{j+1}} \right)
  = \mdimgenwaveletfilterelem{n}{\zos}{\secondznvar - 2 \firstznvar} \\
  \label{eq:mdim-dual-wavelet-grid-value}
  \szdualappl{\mdimgendualwavelet{n}{\zos}{j}
  {\secondznvar}}
  {\mdimsf{n}{j+1}{\firstznvar}}
  & = & \mdimgendualwaveletfilterelem{n}{\zos}
  {\firstznvar - 2 \secondznvar}
\end{eqnarray}
for all \(j \in \integernumbers\),
\(\firstznvar, \secondznvar \in \zn\), and
\(\zos, \zot \in \zeroonesetn\).
}
\fi

\if\shortprep0
{
\begin{lemma}
  \label{lem:mdim-filter-formula}
  Let \(\firstznvar, \secondznvar \in \zn\).
  Then
  \begin{displaymath}
    \sum_{\zos \in \zeroonesetn}
    \sum_{\thirdznvar \in \zn}
    \mdimgendualwaveletfilterelem{n}{\zos}{\firstznvar - 2 \thirdznvar}
    \mdimgenwaveletfilterelem{n}{\zos}{\secondznvar - 2 \thirdznvar}
    = \delta_{\firstznvar,\secondznvar} .
  \end{displaymath}
\end{lemma}
}
\fi

\if\shortprep0
{
\begin{proof}
  Let
  \begin{displaymath}
    a_{i,j,t} := \onedimgendualwaveletfilterelem{t}
    {\cartprodelem{\firstznvar}{i}-2j}
    \onedimgenwaveletfilterelem{t}
    {\cartprodelem{\secondznvar}{i}-2j} .
  \end{displaymath}
  Now
  \begin{displaymath}
    \sum_{\zos \in \zeroonesetn}
    \sum_{\thirdznvar \in \zn}
    \mdimgendualwaveletfilterelem{n}{\zos}{\firstznvar - 2 
    \thirdznvar}
    \mdimgenwaveletfilterelem{n}{\zos}{\secondznvar - 2 \thirdznvar}
    =
    \sum_{\zos \in \zeroonesetn}
    \sum_{\thirdznvar \in \zn}
	  \prod_{i=1}^n 
	  a_{i,\cartprodelem{\thirdznvar}{i},\cartprodelem{\zos}{i}}
	  = \prod_{i=1}^n \sum_{j \in \integernumbers}
	    \sum_{t=0}^1 a_{i,j,t}
  \end{displaymath}
  Let
  \begin{displaymath}
    b_i := \sum_{j \in \integernumbers}
    \sum_{t=0}^1 a_{i,j,t}
  \end{displaymath}
  for all \(i \in \integernumbers\).
  
  Suppose that \(i_0 \in \setoneton{n}\) and
  \(\cartprodelem{\secondznvar}{i_0}\) is even. Then
  \(\cartprodelem{\secondznvar}{i_0} = 2w\) for some
  \(w \in \integernumbers\).
  It follows that
  \begin{math}
    a_{i,j,0} = \delta_{\cartprodelem{\firstznvar}{i},2j}
      \delta_{w,j}
  \end{math}
  and
  \begin{math}
    a_{i,j,1} = \gtfilterelem{\cartprodelem{\firstznvar}{i_0} - 2j}
      \gfilterelem{2w-2j}
      = 0
  \end{math}
  Consequently
  \begin{equation}
    \label{eq:case-a}
    b_{i_0} = \sum_{j \in \integernumbers} 
    \delta_{\cartprodelem{\firstznvar}{i_0},2j}
      \delta_{w,j}
      = \delta_{\cartprodelem{\firstznvar}{i_0},
      \cartprodelem{\secondznvar}{i_0}} .
  \end{equation}
  
  Assume then that \(i_0 \in \setoneton{n}\) and
  \(\cartprodelem{\secondznvar}{i_0}\) is odd.
  Now
  \(\cartprodelem{\secondznvar}{i_0} = 2w + 1\)
  for some \(w \in \integernumbers\),
  \begin{math}
    a_{i_0,j,0} = \delta_{\cartprodelem{\firstznvar}{i_0},2j}
    \hfilterelem{2w+1-2j}
  \end{math}, and
  \begin{math}
    a_{i_0,j,1} =
    \gtfilterelem{\cartprodelem{\firstznvar}{i_0} - 2j}
    \gfilterelem{2w+1-2j}
    = (-1)^{\cartprodelem{\firstznvar}{i_0}-2w-1}
      \hfilterelem{2w+1-\cartprodelem{\firstznvar}{i_0}}
      \delta_{2w+1-2j,1}
  \end{math}.
  Assume that
  \(\cartprodelem{\firstznvar}{i_0}\) is even.
  Now
  \(\cartprodelem{\firstznvar}{i_0} = 2z\) for some
  \(z \in \integernumbers\),
  \begin{math}
    a_{i_0,j,0} = \delta_{z,j} \hfilterelem{2w+1-2z} ,
  \end{math}
  and
  \begin{math}
    a_{i_0,j,1} = - \hfilterelem{2w+1-2z} \delta_{w,j}
  \end{math}.
  It follows that
  \begin{equation}
    \label{eq:case-b}
    b_{i_0} =
    \sum_{j \in \integernumbers} \delta_{z,j} \hfilterelem{2w+1-2z}
    - \sum_{j \in \integernumbers} \hfilterelem{2w+1-2z} 
    \delta_{w,j}
    = 0 = \delta_{\cartprodelem{\firstznvar}{i_0},
    \cartprodelem{\secondznvar}{i_0}}.
  \end{equation}
  Assume finally that
  \(\cartprodelem{\firstznvar}{i_0}\) is odd.
  Now \(\cartprodelem{\firstznvar}{i_0} = 2z+1\) for
  some \(z \in \integernumbers\).
  We have
  \begin{math}
    a_{i_0,j,0} = \delta_{2z+1,2l} \hfilterelem{2w+1-2l}
    = 0
  \end{math},
  \begin{math}
    a_{i_0,j,1} = (-1)^{2z+1-2w-1} \hfilterelem{2w+1-2z-1}
      \delta_{2w+1-2j,1}
      = \delta_{w,z} \delta_{w,j}
  \end{math}.
  It follows that
  \begin{equation}
    \label{eq:case-c}
    b_{i_0} =
    \sum_{j \in \integernumbers} a_{i_0,j,0}
    + \sum_{j \in \integernumbers} a_{i_0,j,1}
	  = \delta_{w,z}
	  = \delta_{\cartprodelem{\firstznvar}{i_0},
	  \cartprodelem{\secondznvar}{i_0}} .
  \end{equation}
  By \myequations \eqref{eq:case-a}, \eqref{eq:case-b}, and \eqref{eq:case-c}
  the lemma is true.
\end{proof}
}
\fi

Goedecker \cite{goedecker1998} gives formulas for
wavelet filters, too.

\section{Compactly Supported Interpolating MRA \mraprepspace 
$\ucfunc{\rn}{K}$}

\label{sec:uc-mra}

We will assume that \(n \in \positiveintegers\) and
\(K = \realnumbers\) or \(K = \complexnumbers\) throughout this
section. We will also assume that
\(\mdimmothersf{n} \in \cscfunc{\rn}{K}\) is
an \(n\)-dimensional tensor product mother
scaling function in this section.
Chui and Li \cite{cl1996} have developed a MRA in the
univariate case \(\ucfunccv{\realnumbers}\).

\begin{definition}
  \label{def:uvjspace}
  Define
  \begin{eqnarray}
    \label{eq:uvjspace}
    & & \mdimuvspace{n}{j} :=
    \left\{
      \rnx \in \rn \mapsto
      \sum_{\firstznvar \in \zn}
      \seqelem{\seqstyle{a}}{\firstznvar}
      \mdimmothersf{n} ( 2^j \rnx - \firstznvar )
      \setsep
      \seqstyle{a} \in \littlelp{\infty}{\zn}{K}
    \right\} \\
    \nonumber
    & & \norminspace{f}{\mdimuvspace{n}{j}} := \norminfty{f} ,
    \spaceafter f \in \mdimuvspace{n}{j}
  \end{eqnarray}
  for all \(j \in \integernumbers\).
\end{definition}

\begin{definition}
Let spaces
\(\mdimuvspace{n}{j}\), \(j \in \integernumbers\),
be defined by \mydef \ref{def:uvjspace}.
We call
\(\{\mdimuvspace{n}{j} \setsep j \in \integernumbers\}\)
an \defterm{interpolating tensor product MRA} \mraprepspace 
\(\ucfunc{\rn}{K}\)
generated by \(\mdimmothersf{n}\) provided that the following
conditions are satisfied:
\begin{myenumerate}
  \myitem{(MRA1.1)}
    \begin{math}
  	  \forall j \in \integernumbers :
  	  \mdimuvspace{n}{j} \subset \mdimuvspace{n}{j+1}
  	\end{math}
  \myitem{(MRA1.2)}
    \begin{math}
      \overline{\bigcup_{j \in \integernumbers} \mdimuvspace{n}{j}}
      = \ucfunc{\rn}{K}
    \end{math}
  \myitem{(MRA1.3)}
    \begin{math}
      \bigcap_{j \in \integernumbers} \mdimuvspace{n}{j} = K
    \end{math}
  \myitem{(MRA1.4)}
    \begin{math}
      \forall j \in \integernumbers, f \in K^{\rn} :
      f \in \mdimuvspace{n}{j} \iff f(2 \cdot) \in \mdimuvspace{n}{j+1}
    \end{math}
  \myitem{(MRA1.5)}
    \begin{math}
      \forall j \in \integernumbers, \firstznvar \in \zn,
      f \in K^{\rn} :
      f \in \mdimuvspace{n}{j} \iff
      f(\cdot - 2^{-j} \firstznvar) \in \mdimuvspace{n}{j}
    \end{math}
  \myitem{(MRA1.6)}
    \begin{math}
      \forall \firstznvar \in \zn :
      \mdimmothersf{n}(\firstznvar) = \delta_{\firstznvar, 0} 
    \end{math}
\end{myenumerate}
\end{definition}

Our requirements for the definition of interpolating
multiresolution analysis are stricter and simpler than those in \cite{cl1996}.
Condition (MRA1.6) is replaced by a weaker condition
for \(\myphi\) in \cite{cl1996} but it is possible to construct function
\(\myphi_L\)
that satisfies condition (MRA1.6) and generates the
same subspaces \(\mdimuvspace{1}{j}\) as function \(\myphi\).

It has been proved in \cite{jm1991} that
\begin{displaymath}
  \sum_{\firstznvar \in \zn} \mdimmothersf{n}(\rny - \firstznvar)
  = \widehat{\mdimmothersf{n}}(0)
\end{displaymath}
for all \(\rny \in \rn\).
By the cardinal interpolation property (MSF.1) we have
\begin{displaymath}
  \sum_{\firstznvar \in \zn} \mdimmothersf{n}(0 - \firstznvar)
  = 1
\end{displaymath}
from which it follows that
\begin{equation}
  \label{eq:sf-sum}
  \sum_{\firstznvar \in \zn} \mdimmothersf{n}(\rny - \firstznvar)
  = 1
\end{equation}
for all \(\rny \in \rn\).

\begin{remark}
  The series in \myequation \eqref{eq:uvjspace} converges (pointwise) for all
  \begin{math}
    \seqstyle{a} \in \littlelp{\infty}{\zn}{K}
  \end{math}
  but the series
  \begin{math}
    \sum_{\firstznvar \in \zn}
    \seqelem{\seqstyle{a}}{\firstznvar}
    \mdimmothersf{n} ( 2^j \cdot - \firstznvar )
  \end{math}
  need not converge in the norm topology of
  \begin{math}
    \ucfunc{\rn}{K}
  \end{math}.
  For example,
  consider sequence \(\seqelem{\seqstyle{a}}{\firstznvar} =
  1\), \(\firstznvar \in \zn\). Now
  \begin{math}
    \sum_{\firstznvar \in \zn}
    \seqelem{\seqstyle{a}}{\firstznvar}
    \mdimmothersf{n} ( 2^j \rnx - \firstznvar )
    = 1
  \end{math}
  for all \(\rnx \in \rn\) but the series
  \begin{math}
    \sum_{\firstznvar \in \zn}
    \seqelem{\seqstyle{a}}{\firstznvar}
    \mdimmothersf{n} ( 2^j \cdot - \firstznvar )
  \end{math}
  does not converge in the norm topology of
  \begin{math}
    \ucfunc{\rn}{K}
  \end{math}.
\end{remark}

\if\longprep1
{
\givelemmawithoutproof
}
\fi

\if\shortprep1
{
  Under the conditions given at the start of this section
  (\(\mdimmothersf{n}\) is an \(n\)-dimensional tensor
  product mother scaling function)
  the conditions (MRA1.1),
  (MRA1.4), (MRA1.5), and (MRA1.6) are true.
  Spaces \(\mdimuvspace{n}{j}\), \(j \in \integernumbers\),
  are topologically isomorphic to
  \(\littlelp{\infty}{\zn}{K}\).
}
\else
{
\begin{lemma}
  \label{lem:uv-mra-properties}
  Under the conditions given at the start of this section
  (\(\mdimmothersf{n}\) is an \(n\)-dimensional tensor
  product mother scaling function)
  the conditions (MRA1.1),
  (MRA1.4), (MRA1.5), and (MRA1.6) are true.
\end{lemma}
}
\fi

We will show later that conditions (MRA1.2) and (MRA1.3) are true,
too.

\if\shortprep0
{
\begin{proof}
  Use \mylemmas \ref{lem:linftyspace-top-isomorphism}
  and \ref{lem:mdim-wavelet-formulas}.
\end{proof}
}
\fi

\if\shortprep0
{
Define function
\begin{math}
  \uvjtopisomfunc{n}{j} : \littlelp{\infty}{\zn}{K} \to \mdimuvspace{n}{j}
\end{math}
by
\begin{displaymath}
  \uvjtopisom{n}{j}{\seqstyle{a}} :=
  \rnx \in \rn \mapsto \sum_{\firstznvar \in \zn} \seqelem{\sqa}{\firstznvar}
    \mdimsf{n}{j}{\firstznvar}(\rnx)
\end{displaymath}
for all
\begin{math}
  \seqstyle{a} \in \littlelp{\infty}{\zn}{K}
\end{math}.
By \mylemma \ref{lem:linftyspace-top-isomorphism}
function
\(\uvjtopisomfunc{n}{j}\) is a topological isomorphism from
\(\littlelp{\infty}{\zn}{K}\) onto
\(\mdimuvspace{n}{j}\) and
\begin{equation}
  \label{eq:linftyspace-norm-equiv}
  \norminfty{\seqstyle{a}}
  \leq \norminfty{\uvjtopisom{n}{j}{\seqstyle{a}}}
  \leq \ncover{\mdimmothersf{n}} \norminfty{\mdimmothersf{n}}
  \norminfty{\seqstyle{a}}
\end{equation}
for all
\begin{math}
  \seqstyle{a} \in \littlelp{\infty}{\zn}{K}
\end{math}
and
\(j \in \integernumbers\).
}
\fi

\begin{theorem}
  \label{th:intersection}
  We have
  \begin{displaymath}
    \bigcap_{j \in \integernumbers} \mdimuvspace{n}{j} = K .
  \end{displaymath}
\end{theorem}

\begin{proof}
  The proof is similar to a part of the proof of
  \cite[theorem 3.2]{cl1996}.
\if\longversion1
{
  Suppose that \(g \in \mdimuvspace{n}{j}\) for all \(j \in \integernumbers\).
  Now
  \begin{displaymath}
    \forall \rnx \in \rn, j \in \integernumbers :
    g(\rnx) = \sum_{\firstznvar \in \zn}
    \seqelem{\indexedseqstyle{a}{j}}{\firstznvar}
    \mdimmothersf{n}(2^j \rnx - \firstznvar)
  \end{displaymath}
  where
  \begin{math}
    \indexedseqstyle{a}{j}
    = (a_{j,\firstznvar})_{\firstznvar \in \zn}
    \in \littlelp{\infty}{\zn}{K}
  \end{math}
  for all \(j \in \integernumbers\).
  It follows from condition (MRA1.6) that
  \begin{displaymath}
    \norminfty{g} \geq
    \abs{g \left( \frac{\firstznvar}{2^j} \right)} =
    \abs{\sum_{\secondznvar \in \integernumbers^n} a_{j,\secondznvar} 
    \mdimmothersf{n}(\firstznvar -
      \secondznvar)}
    = \abs{a_{j,\firstznvar}}
  \end{displaymath}
  for all \(\firstznvar \in \integernumbers^n\) and
  \(j \in \integernumbers\). Hence
  \(\norminfty{\indexedseqstyle{a}{j}} \leq \norminfty{g}\)
  for all \(j \in \integernumbers\).
  Function \(\mdimmothersf{n}\) is compactly supported so there exists
  \(r \in \positiverealnumbers\) so that
  \begin{math}
    \suppop \mdimmothersf{n} \subset \closedball{\realnumbers^n}{0}{r}
  \end{math}.
  
  Let \(\rnx, \rny \in \realnumbers^n\). Now
  \begin{subequations}
    \begin{align}
      \label{eq:dist-a}
      \abs{g(\rnx) - g(\rny)} & =
      \abs{\sum_{\firstznvar \in \integernumbers^n} a_{j,\firstznvar}
        \mdimmothersf{n}(2^j \rnx -
        \firstznvar)
        - \sum_{\secondznvar \in \integernumbers^n} a_{j,\secondznvar}
        \mdimmothersf{n}(2^j \rny -
        \secondznvar)} \\
      & =
      \label{eq:dist-b}
      \abs{\sum_{\firstznvar \in \integernumbers^n} a_{j,\firstznvar}
        \left( \mdimmothersf{n}(2^j \rnx - \firstznvar)
          - \mdimmothersf{n}(2^j \rny - \firstznvar) 
        \right)}
    \end{align}
  \end{subequations}
  for all \(j \in \integernumbers\).
  Since \(f\) is compactly supported all the series in \myequations
  \eqref{eq:dist-a} and \eqref{eq:dist-b} contain only finite number
  of nonzero terms.
  Therefore
  \begin{align*}
    \abs{g(\rnx) - g(\rny)} & \leq
    \norminfty{\indexedseqstyle{a}{j}}
    \sum_{\firstznvar \in \integernumbers^n}
    \abs{\mdimmothersf{n}(2^j \rnx - \firstznvar)
      - \mdimmothersf{n}(2^j \rny - \firstznvar)} \\
    &
    \leq \norminfty{g}
    \sum_{\firstznvar \in \integernumbers^n}
    \abs{\mdimmothersf{n}(2^j \rnx - \firstznvar)
      - \mdimmothersf{n}(2^j \rny - \firstznvar)}
  \end{align*}
  for all \(j \in \integernumbers\)
  where the series contain only finite number of nonzero terms.
  Function \(\mdimmothersf{n}\) is uniformly continuous and hence
  \begin{math}
    \abs{\mdimmothersf{n}(2^j \rnx - \firstznvar)
    - \mdimmothersf{n}(2^j \rny - \firstznvar)} \to 0 \;\textrm{as}\;
    j \to -\infty
  \end{math}
  for each \(\firstznvar \in \integernumbers^n\).
  
  Let
  \begin{math}
    m = \max \{ \norm{\rnx}, \norm{\rny} \}
  \end{math}.
  Suppose that \(j \in \integernumbers\), \(j \leq 0\),
  \(\thirdznvar \in \integernumbers^n\), and
  \(\normtwo{\thirdznvar} > r + m\).
  Now
  \begin{math}
    2^j \norm{\rnx} \leq 2^j m \leq m
  \end{math}
  and
  \begin{math}
    \norm{2^j \rnx - \thirdznvar} \geq \abs{2^j \norm{\rnx} -
      \normtwo{\thirdznvar}}
    = \normtwo{\thirdznvar} - 2^j \norm{\rnx}
    > r + m - m = r
  \end{math}.
  Hence
  \begin{math}
    \mdimmothersf{n}(2^j \rnx - \thirdznvar) = 0
  \end{math}.
  We also get
  \begin{math}
    \mdimmothersf{n}(2^j \rny - \thirdznvar) = 0
  \end{math}
  similarly.
  Let
  \begin{displaymath}
    \seqelem{\sqd}{\firstznvar} :=
    \left\{
      \begin{array}{ll}
        2 \norminfty{\mdimmothersf{n}} ; & \textrm{if} \; 
        \normtwo{\firstznvar} \leq r + m \\
        0; & \textrm{otherwise}
      \end{array}
    \right.
  \end{displaymath}
  for all \(\firstznvar \in \integernumbers^n\).
  Now
  \begin{math}
    \abs{\mdimmothersf{n}(2^j \rnx - \firstznvar)
    - \mdimmothersf{n}(2^j \rny - \firstznvar)}
    \leq 
    \seqelem{\sqd}{\firstznvar}
  \end{math}
  for all \(\firstznvar \in \integernumbers^n\),
  \(j \in \integernumbers\), \(j \leq 0\)
  and
  \begin{displaymath}
    \sum_{\firstznvar \in \integernumbers^n} \seqelem{\sqd}{\firstznvar}
    < \infty .
  \end{displaymath}
  Thus by the Dominated Convergence Theorem for Series 
  \cite{krieger2005}
  \begin{displaymath}
    \sum_{\firstznvar \in \integernumbers^n}
    \abs{\mdimmothersf{n}(2^j \rnx - \firstznvar)
    - \mdimmothersf{n}(2^j \rny - \firstznvar)}
    \to 0  \;\textrm{as}\; j \to -\infty .
  \end{displaymath}
  Hence \(\abs{g(\rnx) - g(\rny)} = 0\).
  Vectors \(\rnx, \rny \in \realnumbers^n\) were arbitrary so \(g\) 
  is a constant function.
  Consequently
  \begin{displaymath}
    \bigcap_{j \in \integernumbers} \mdimuvspace{n}{j} = K .
  \end{displaymath}
}
\fi
\end{proof}

\begin{definition}
  \label{def:mdimpartialuwspace}
  When \(j \in \integernumbers\) and
  \(\zos \in \zeroonesetn\) define
  \begin{eqnarray*}
    & & \mdimupartialwspace{n}{\zos}{j}
    :=
    \left\{
    	\rnx \in \rn \mapsto
    	\sum_{\firstznvar \in \zn}
    	\seqelem{\sqa}{\firstznvar}
    	\mdimgenwavelet{n}{\zos}{j}{\firstznvar}(\rnx)
    	\setsep
    	\sqa \in \littlelp{\infty}{\zn}{K}
    \right\} \\
    & & \norminspace{f}{\mdimupartialwspace{n}{\zos}{j}}
    :=
    \norminfty{f} ,
    \spaceafter
    f \in \mdimupartialwspace{n}{\zos}{j}
    .
  \end{eqnarray*}
\end{definition}

\if\shortprep1
{
  Spaces \(\mdimupartialwspace{n}{\zos}{j}\),
  \(\zos \in \zeroonesetn\),
  \(j \in \integernumbers\),
  are topologically isomorphic to
  \(\littlelp{\infty}{\zn}{K}\).
}
\else
{
When \(j \in \integernumbers\) we have
\begin{math}
  \mdimuvspace{n}{j} \equalns
  \mdimupartialwspace{n}{\finitezeroseq{n}}{j}
\end{math}.
By \mylemma \ref{lem:linftyspace-top-isomorphism}
function
\begin{math}
  \uwsjtopisomfunc{n}{\zos}{j} :
  \littlelp{\infty}{\zn}{K} \to 
  \mdimupartialwspace{n}{\zos}{j}
\end{math}
defined by
\begin{displaymath}
  \uwsjtopisom{n}{\zos}{j}{\seqstyle{a}}
  :=
  \sum_{\firstznvar \in \integernumbers} 
  \seqelem{\seqstyle{a}}{\firstznvar}
  \mdimgenwavelet{n}{\zos}{j}{\firstznvar}
\end{displaymath}
for all
\begin{math}
  \seqstyle{a} \in \littlelp{\infty}{\zn}{K}
\end{math}
is a topological isomorphism from
\begin{math}
  \littlelp{\infty}{\zn}{K}
\end{math}
onto
\begin{math}
  \mdimupartialwspace{n}{\zos}{j}
\end{math}
and
\begin{equation}
  \label{eq:mdimpartialuwspace-norm-equiv}
  \norminfty{\seqstyle{a}}
  \leq \norminfty{\uwsjtopisom{n}{\zos}{j}{\seqstyle{a}}}
  \leq \szncover{\mdimmotherwavelet{n}{\zos}}
  \norminfty{\mdimmotherwavelet{n}{\zos}}
  \norminfty{\seqstyle{a}}
\end{equation}
for all
\begin{math}
  \seqstyle{a} \in \littlelp{\infty}{\zn}{K}
\end{math},
\begin{math}
  \zos \in \zeroonesetn
\end{math}, and
\begin{math}
  j \in \integernumbers
\end{math}.
}
\fi

\begin{definition}
  When \(\zos \in \zeroonesetn\)
  and \(j \in \integernumbers\)
  define
  \begin{displaymath}
    \left(\mdimupartialdeltaop{n}{\zos}{j} f\right)
    \left( \rnx \right) :=
    \sum_{\firstznvar \in \zn}
    \szdualappl{\mdimgendualwavelet{n}{\zos}{j}{\firstznvar}}{f}
    \mdimgenwavelet{n}{\zos}{j}{\firstznvar}
    \left( \rnx \right)    
  \end{displaymath}
  for all \(\rnx \in \rn\)
  and \(f \in \ucfunc{\rn}{K}\).
  When \(j \in \integernumbers\) define
  \begin{math}
    \mdimuprojop{n}{j} := 
    \mdimupartialdeltaop{n}{\finitezeroseq{n}}{j}
  \end{math}.
\end{definition}

Operator \(\mdimupartialdeltaop{n}{\zos}{j}\)
is a continuous projection \projprep
\(\ucfunc{\rn}{K}\) onto
\(\mdimupartialwspace{n}{\zos}{j}\)
for each
\(\zos \in \zeroonesetn\)
and \(j \in \integernumbers\).
\if\shortprep0
{
When \(\zos \in \zeroonesetn\)
operators \(\mdimupartialdeltaop{n}{\zos}{j}\),
\(j \in \integernumbers\),
are uniformly bounded by
\begin{displaymath}
  \norm{\mdimupartialdeltaop{n}{\zos}{j}}
  \leq
  \norm{\mdimmotherdualwavelet{n}{\zos}}
  \norminfty{\mdimmotherwavelet{n}{\zos}}
  \ncover{\mdimmotherwavelet{n}{\zos}}
\end{displaymath}
for all \(j \in \integernumbers\).
Operators \(\mdimuprojop{n}{j}\),
\(j \in \integernumbers\),
are uniformly bounded by
\begin{equation}
  \label{eq:mdimuprojop-uniform-bound}
  \norm{\mdimuprojop{n}{j}}
  \leq
  \norminfty{\mdimmothersf{n}}
  \ncover{\mdimmothersf{n}}
\end{equation}
for all \(j \in \integernumbers\).
}
\fi

\if\shortprep0
{
\begin{lemma}
  \label{lem:functional-and-proj-op}
  Let \(j \in \integernumbers\),
  \(\zos \in \zeroonesetn\),
  \(\firstznvar \in \zn\), and
  \(f \in \ucfunccv{\rn}\).
  Then
  \begin{math}
    \dualappl{\mdimgendualwavelet{n}{\zos}{j}{\firstznvar}}{f}
    = 
    \dualappl{\mdimgendualwavelet{n}{\zos}{j}{\firstznvar}}
    	{\mdimuprojop{n}{j+1} f}
  \end{math}.
\end{lemma}

\begin{proof}
  Use \myequation \eqref{eq:mdim-dual-wavelet-formula}.
\end{proof}
}
\fi

\if\shortprep0
{
\begin{lemma}
  \label{lem:mdimpartialdeltaop-circ}
  Let \(j, j' \in \integernumbers\),
  \( j < j'\),
  and
  \(\zos \in \zeroonesetn\).
  Then
  \begin{math}
    \mdimupartialdeltaop{n}{\zos}{j} =
    \mdimupartialdeltaop{n}{\zos}{j} \circ \mdimuprojop{n}{j'}
  \end{math}.
\end{lemma}

\begin{proof}
  Use \mylemma \ref{lem:functional-and-proj-op}.
\end{proof}
}
\fi

\if\longprep1
{
\givelemmawithoutproof
}
\fi

\if\shortprep0
{
\begin{lemma}
  \label{lem:partialuwspace-intersection}
  Let \(j \in \integernumbers\).
  Let \(\zos, \zot \in \zeroonesetn\) and \(\zos \neq \zot\).
  Then
  \begin{math}
    \mdimupartialwspace{n}{\zos}{j} \intersection
    \mdimupartialwspace{n}{\zot}{j}
    = \{ 0 \} .
  \end{math}
\end{lemma}
}
\fi

\if\shortprep0
{
\begin{lemma}
  \label{lem:udeltaop-grid-value}
  Let \(j \in \integernumbers\),
  \(\zos \in \zeroonesetn\), and
  \(\secondznvar \in \zn\).
  Let \(\sqa \in \littlelpcv{\infty}{\zn}\) and
  \begin{displaymath}
    f(\rnx) =
    \sum_{\firstznvar \in \zn}
      \seqelem{\sqa}{\firstznvar}
      \mdimsf{n}{j+1}{\firstznvar}(\rnx)  
  \end{displaymath}
  for all \(\rnx \in \rn\).
  Then
  \begin{displaymath}
    \left(\mdimupartialdeltaop{n}{\zos}{j} f\right)
    \left( \frac{\secondznvar}{2^{j+1}} \right)
    = \sum_{\firstznvar \in \zn}
    \seqelem{\sqa}{\firstznvar}
    \sum_{\seqstyle{z} \in \zn}
    \mdimgendualwaveletfilterelem{n}{\zos}{\firstznvar - 2 \seqstyle{z}}
    \mdimgenwaveletfilterelem{n}{\zos}{\secondznvar - 2 \seqstyle{z}} .
  \end{displaymath}
\end{lemma}

\begin{proof}
  Use \myequations \eqref{eq:mdim-wavelet-grid-value}
  and \eqref{eq:mdim-dual-wavelet-grid-value}.
\end{proof}
}
\fi

\begin{definition}
  \label{def:mdimuwspace}
  When \(j \in \integernumbers\) define
  \begin{displaymath}
    \mdimuwspace{n}{j} \defequalns
      \directsumoneindex{\zos \in \zeroonesetnnozero{n}}
      \mdimupartialwspace{n}{\zos}{j}
  \end{displaymath}
  and
  \begin{displaymath}
    \mdimudeltaop{n}{j} :=
      \sum_{\zos \in \zeroonesetnnozero{n}}
      \mdimupartialdeltaop{n}{\zos}{j} .
  \end{displaymath}
\end{definition}

\if\shortprep1
{
  Spaces \(\mdimuwspace{n}{j}\), \(j \in \integernumbers\),
  are topologically isomorphic to
  \(\littlelp{\infty}{\zeroonesetnnozero{n} \times \zn}{K}\).
}
\fi

\if\shortprep0
{
\begin{lemma}
  \label{lem:mdimpartialuprojop-orthogonality}
  Let \(j \in \integernumbers\). Then
  \begin{itemize}
  \item[(i)]
    \begin{math}
      \forall \zos \in \zeroonesetn,
      \zot \in \zeroonesetn :
      \zos \neq \zot \implies
      \forall f \in 
      \mdimupartialwspace{n}{\zos}{j} :
      \mdimupartialdeltaop{n}{\zot}{j} f = 0
    \end{math}
  \item[(ii)]
    \begin{math}
      \forall j' \in \integernumbers,
      f \in \mdimuvspace{n}{j'} : j' \leq j
      \implies \mdimudeltaop{n}{j} f = 0
    \end{math}
  \item[(iii)]
    \begin{math}
      \forall j' \in \integernumbers,
      f \in \mdimuwspace{n}{j'} : j' \neq j
      \implies \mdimudeltaop{n}{j} f = 0
    \end{math}
  \item[(iv)]
    \begin{math}
      \forall j' \in \integernumbers,
      f \in \mdimuwspace{n}{j'} : j' \geq j
      \implies \mdimuprojop{n}{j} f = 0
    \end{math}.
  \end{itemize}
\end{lemma}
}
\fi

\if\shortprep0
{
\begin{proof}
  \mbox{ }
  \begin{itemize}
  \item[(i)]
    Let
    \begin{math}
      \zos, \zot \in \zeroonesetn
    \end{math}
    and
    \begin{math}
      \zos \neq \zot
    \end{math}.
    Let \(f \in 
    \mdimupartialwspace{n}{\zos}{j}\).
    Now
    \begin{displaymath}
      f(\rnx) = \sum_{\firstznvar \in \zn}
      \seqelem{\sqa}{\firstznvar}
      \mdimgenwavelet{n}{\zos}{j}{\firstznvar}
      (\rnx)
    \end{displaymath}
    for all \(\rnx \in \rn\),
    where \(\sqa \in \littlelpcv{\infty}{\zn}\).
    By \mylemma \ref{lem:delta-and-pointwise-conv}
    we have
    \begin{displaymath}
      \dualappl{\mdimgendualwavelet{n}
      {\zot}{j}{\secondznvar}}{f}
      = \sum_{\firstznvar \in \zn}
      \seqelem{\sqa}{\firstznvar}
      \dualappl{\mdimgendualwavelet{n}{\zot}{j}
      {\secondznvar}}
      {\mdimgenwavelet{n}{\zos}{j}
      {\firstznvar}}
      = 0
    \end{displaymath}
    for all
    \(\secondznvar \in \zn\).
    Hence
    \begin{math}
      \mdimupartialdeltaop{n}{\zot}{j} f = 0
    \end{math}.
  \item[(ii)]
    If \(j' < j\) then \(\mdimuvspace{n}{j'} \subset
    \mdimuvspace{n}{j} = 
    \mdimupartialwspace{n}{\finitezeroseq{n}}{j}\) and
    proposition (ii) follows from proposition (i).
  \item[(iii)]
    If \(j < j'\) and \(g \in \mdimuwspace{n}{j'}\) then
    \begin{math}
      \mdimudeltaop{n}{j} g
      = \mdimudeltaop{n}{j} \mdimuprojop{n}{j'} g
    \end{math}.
    Since \(\mdimuprojop{n}{j'} =
    \mdimupartialdeltaop{n}{\finitezeroseq{n}}{j'}\) it
    follows that \(\mdimudeltaop{n}{j} g = 0\) by proposition
    (i).
    If \(j > j'\) and \(g \in \mdimuwspace{n}{j'}\) then
    \(\mdimuwspace{n}{j'} \subset \mdimuvspace{n}{j}\)
    and \(\mdimudeltaop{n}{j} g = 0\) by proposition (ii).
  \item[(iv)]
    Let \(j' \in \integernumbers\), \(j' \geq j\).
    Suppose that \(w \in \mdimuwspace{n}{j'}\).
    Since \(\mdimuprojop{n}{j'} =
    \mdimupartialdeltaop{n}{\finitezeroseq{n}}{j'}\) it
    follows that
    \(\mdimuprojop{n}{j} w
    = \mdimuprojop{n}{j} \mdimuprojop{n}{j'} w = 0\)
    by proposition (i).
  \end{itemize}
\end{proof}
}
\fi

\if\shortprep1
{
  When \(j \in \integernumbers\)
  we have
  \begin{displaymath}
    \mdimuvspace{n}{j+1}
    \equalns
    \mdimuvspace{n}{j} \dsum \mdimuwspace{n}{j}
    \equalns
    \mdimuvspace{n}{j} \dsum
    \directsumoneindex{\zos \in \zeroonesetnnozero{n}}
    \mdimupartialwspace{n}{\zos}{j}
    \equalns
    \directsumoneindex{\zos \in \zeroonesetn}
    \mdimupartialwspace{n}{\zos}{j} .
  \end{displaymath}
}
\else
{
\begin{lemma}
  \label{lem:mdimuvspace-direct-sum}
  Let \(j \in \integernumbers\).
  Then
  \begin{displaymath}
    \mdimuvspace{n}{j+1}
    \equalns
    \mdimuvspace{n}{j} \dsum \mdimuwspace{n}{j}
    \equalns
    \mdimuvspace{n}{j} \dsum
    \directsumoneindex{\zos \in \zeroonesetnnozero{n}}
    \mdimupartialwspace{n}{\zos}{j}
    \equalns
    \directsumoneindex{\zos \in \zeroonesetn}
    \mdimupartialwspace{n}{\zos}{j} .
  \end{displaymath}
\end{lemma}
}
\fi

\if\shortprep1
{
\if10
{
  \begin{proof}
    By \mylemma \ref{lem:partialuwspace-intersection}
    the sums in the lemma are direct.
    Use \mylemmas \ref{lem:udeltaop-grid-value} and
    \ref{lem:mdim-filter-formula} for the proof.
  \end{proof}
}
\fi
}
\else
{
\begin{proof}
  By \mylemma \ref{lem:partialuwspace-intersection}
  the sums in the lemma are direct.
  
  Suppose that
  \begin{math}
    f \in \mdimuvspace{n}{j+1}
  \end{math}.
  Now
  \begin{displaymath}
  	f(\rnx) = \sum_{\firstznvar \in \zn} \seqelem{\sqa}{\firstznvar}
  		\mdimsf{n}{j+1}{\firstznvar}(\rnx)
  \end{displaymath}
  for all \(\rnx \in \rn\) where
  \(\sqa \in \littlelp{\infty}{\zn}{K}\).
  Let
  \begin{math}
    f_{\zos} := \mdimupartialdeltaop{n}{\zos}{j} f
    \in \mdimupartialwspace{n}{\zos}{j}
  \end{math}
  for each \(\zos \in \zeroonesetn\).
  Suppose that \(\secondznvar \in \zn\).
  By \mylemma \ref{lem:udeltaop-grid-value}
  \begin{eqnarray*}
    \sum_{\zos \in \zeroonesetn}
    f_\zos \left( \frac{\secondznvar}{2^{j+1}} \right)
    & = &
    \sum_{\zos \in \zeroonesetn}
    \sum_{\firstznvar \in \zn}
    \seqelem{\sqa}{\firstznvar}
    \sum_{\thirdznvar \in \zn}
    \mdimgendualwaveletfilterelem{n}{\zos}{\firstznvar - 2 
    \thirdznvar}
    \mdimgenwaveletfilterelem{n}{\zos}{\secondznvar - 2 
    \thirdznvar} \\
    & = &
    \sum_{\firstznvar \in \zn}
    \seqelem{\sqa}{\firstznvar}
    \sum_{\zos \in \zeroonesetn}
    \sum_{\thirdznvar \in \zn}
    \mdimgendualwaveletfilterelem{n}{\zos}{\firstznvar - 2 
    \thirdznvar}
    \mdimgenwaveletfilterelem{n}{\zos}{\secondznvar - 2 
    \thirdznvar} .
  \end{eqnarray*}
  By \mylemma \ref{lem:mdim-filter-formula}
  \begin{displaymath}
    \sum_{\zos \in \zeroonesetn}
    f_\zos \left( \frac{\secondznvar}{2^{j+1}} \right)
    =
    \sum_{\firstznvar \in \zn}
    \seqelem{\sqa}{\firstznvar}
    \delta_{\firstznvar,\secondznvar}
    =
    \seqelem{\sqa}{\secondznvar}
  \end{displaymath}
  Consequently
  \begin{eqnarray*}
    \left( \sum_{\zos \in \zeroonesetn}
    f_\zos \right) 
    \left( \rnx \right)
    & = &
    \sum_{\secondznvar \in \zn}
      \left( \sum_{\zos \in \zeroonesetn} f_\zos \right)
      \left( \frac{\secondznvar}{2^{j+1}} \right)
      \mdimsf{n}{j+1}{\secondznvar}(\rnx) \\
    & = &
    \sum_{\secondznvar \in \zn}
    \seqelem{\sqa}{\secondznvar}
    \mdimsf{n}{j+1}{\secondznvar}(\rnx) .
  \end{eqnarray*}
  Hence
  \begin{displaymath}
    f = \sum_{\zos \in \zeroonesetn} f_\zos .
  \end{displaymath}
\end{proof}
}
\fi

\if\shortprep1
{
  Operator \(\mdimudeltaop{n}{j}\) is a continuous projection
  \projprep \(\ucfunc{\rn}{K}\) onto
  \(\mdimuwspace{n}{j}\)
  for each \(j \in \integernumbers\).
  We also have
  \( \mdimudeltaop{n}{j} = \mdimuprojop{n}{j+1} - \mdimuprojop{n}{j}\)
  for all \(j \in \integernumbers\).
}
\fi

\if\shortprep0
{
\begin{lemma}
  \label{lem:mdimudeltaop-properties}
  \mbox{ }
  \begin{itemize}
    \item[(i)]
      \(\mdimudeltaop{n}{j}\) is a continuous projection
      \projprep \(\ucfunc{\rn}{K}\) onto
      \(\mdimuwspace{n}{j}\)
      for each \(j \in \integernumbers\).
    \item[(ii)]
      We have
      \( \mdimudeltaop{n}{j} = \mdimuprojop{n}{j+1} - \mdimuprojop{n}{j}\)
      for all \(j \in \integernumbers\).
    \item[(iii)]
      Operators \(\mdimudeltaop{n}{j}\), \(j \in \integernumbers\),
      are uniformly bounded by
      \begin{displaymath}
        \norm{\mdimudeltaop{n}{j}}
        \leq
        2 \ncover{\mdimmothersf{n}}
        \norminfty{\mdimmothersf{n}} .
      \end{displaymath}
  \end{itemize}
\end{lemma}
}
\fi

\if\shortprep0
{
\begin{proof}
  Use \mylemmas
  \ref{lem:mdimpartialdeltaop-circ},
  \ref{lem:mdimpartialuprojop-orthogonality}, and
  \ref{lem:mdimuvspace-direct-sum}
  and \myequation
  \eqref{eq:mdimuprojop-uniform-bound}.
\end{proof}
}
\fi

\if\shortprep0
{
\begin{definition}
  When
  \begin{math}
    j \in \integernumbers
  \end{math}
  define
  \begin{displaymath}
    \uwjtopisom{n}{j}{\sqa} := \rnx \in \rn \mapsto
    \sum_{\zos \in \zeroonesetnnozero{n}}
    \sum_{\firstznvar \in \zn}
    \seqelem{\sqa}{\zos,\firstznvar}
    \mdimgenwavelet{n}{\zos}{j}{\firstznvar}(\rnx)
  \end{displaymath}
  for all
  \begin{math}
    \sqa \in \littlelp{\infty}{\zeroonesetnnozero{n} \times \zn}{K}
  \end{math}.
\end{definition}
}
\fi

\if\shortprep0
{
\begin{lemma}
  \label{lem:uw-top-isomorphism}
  \mbox{ }
  \begin{itemize}
    \item[(i)]
      \(\mdimuwspace{n}{j} \closedsubspace \ucfunc{\rn}{K}\)
      for all \(j \in \integernumbers\).
    \item[(ii)]
      Function \(\uwjtopisomfunc{n}{j}\) is a topological
      isomorphism from
      \(\littlelp{\infty}{\zeroonesetnnozero{n} \times \zn}{K}\)
      onto
      \(\mdimuwspace{n}{j}\).
    \item[(iii)]
      There exists \(c_1 \in \positiverealnumbers\) so that
      \(\norm{\uwjtopisomfunc{n}{j}} \leq c_1\) for all
      \(j \in \integernumbers\).
    \item[(iv)]
      There exists \(c_2 \in \positiverealnumbers\) so that
      \(\norm{(\uwjtopisomfunc{n}{j})^{-1}} \leq c_2\) for all
      \(j \in \integernumbers\).
  \end{itemize}
\end{lemma}
}
\fi

\if\shortprep0
{
\begin{proof}
  Use \mylemma \ref{lem:uc-convergence}
  and \mydefs
  \ref{def:mdimpartialuwspace} and
  \ref{def:mdimuwspace}.
\end{proof}
}
\fi

\if\shortprep1
{
  There exists
  \(c_1 \in \positiverealnumbers\)
  and
  \(c_2 \in \positiverealnumbers\)
  so that
  \begin{math}
    \norminfty{f - \mdimuprojop{n}{j} f}
    \leq
    c_1 \modcont{f}{2^{-j} c_2}
  \end{math}
  for all
  \(f \in \ucfunccv{\rn}\)
  and
  \(j \in \integernumbers\).
  When \(f \in \ucfunccv{\rn}\) we have
  \begin{displaymath}
    \lim_{j \to \infty} \norminfty{f - \mdimuprojop{n}{j} f} = 0 .
  \end{displaymath}
  Consequently
  \begin{displaymath}
    \overline{\bigcup_{j \in \integernumbers} \mdimuvspace{n}{j}} =
    \ucfunc{\rn}{K} .
  \end{displaymath}
}
\fi

\begin{theorem}
  \label{th:upj-ineq}
  There exist
  \(c_1 \in \positiverealnumbers\)
  and
  \(c_2 \in \positiverealnumbers\)
  so that
  \begin{displaymath}
    \norminfty{f - \mdimuprojop{n}{j} f}
    \leq
    c_1 \modcont{f}{2^{-j} c_2}
  \end{displaymath}
  for all
  \(f \in \ucfunccv{\rn}\)
  and
  \(j \in \integernumbers\).
\end{theorem}

\begin{proof}
  Let
  \(f \in \ucfunccv{\rn}\)
  and
  \(j \in \integernumbers\).
  Define \(r_1 := \supportradius{\mdimmothersf{n}}\).
  Let \(\rnx \in \rn\).
  Define
  \begin{math}
    I := \itrans{\mdimmothersf{n}}{2^j \rnx}
  \end{math}.
  Now by \myequation \eqref{eq:sf-sum}
  we have
  \begin{eqnarray*}
    \left( f - \mdimuprojop{n}{j} f \right) \left(\rnx\right)
    & = &
    \sum_{\firstznvar \in \zn}
    \left( f\left(\rnx\right)
      - f\left(\frac{\firstznvar}{2^j}\right)
    \right)
    \mdimmothersf{n}\left(2^j \rnx - \firstznvar\right) \\
    & = &
    \sum_{\firstznvar \in I}
    \left( f\left(\rnx\right)
      - f\left(\frac{\firstznvar}{2^j}\right)
    \right)
    \mdimmothersf{n}\left(2^j \rnx - \firstznvar\right)
  \end{eqnarray*}
  Let \(\secondznvar \in I\).
  Now
  \begin{math}
    2^j \rnx - \secondznvar
    \in \closedball{\rn}{0}{r_1}
  \end{math}
  from which it follows that
  \begin{math}
    \rnx - 2^{-j} \secondznvar
    \in \closedball{\rn}{0}{2^{-j} r_1}
  \end{math}.
  Hence
  \begin{equation}
    \label{eq:upj-compl-ineq}
    \abs{f\left(\rnx\right)
    - f\left(\frac{\secondznvar}{2^j}\right)}
    \leq
    \modcont{f}{2^{-j} r_1} .
  \end{equation}
  Let
  \begin{math}
    c_1 := \ncover{\mdimmothersf{n}}
    \norminfty{\mdimmothersf{n}}
  \end{math}
  and
  \begin{math}
    c_2 := r_1 = \supportradius{\mdimmothersf{n}}
  \end{math}.  
  By \myequation \eqref{eq:upj-compl-ineq} we have
  \begin{displaymath}
    \abs{\left( f - \mdimuprojop{n}{j} f \right) \left(\rnx\right)}
    \leq
    c_1 \modcont{f}{2^{-j} c_2}.
  \end{displaymath}
\end{proof}

\if\shortprep0
{
\begin{theorem}
  \label{th:upj-convergence}
  Let \(f \in \ucfunccv{\rn}\). Now
  \begin{displaymath}
    \lim_{j \to \infty} \norminfty{f - \mdimuprojop{n}{j} f} = 0 .
  \end{displaymath}
\end{theorem}
}
\fi

\if\shortprep0
{
\begin{proof}
  \if10
  {
  See \cite[theorem 3.2]{cl1996}
  and \cite[theorem 2.4]{donoho1992}.
  }
  \else
  {
  See also \cite[theorem 3.2]{cl1996}
  and \cite[theorem 2.4]{donoho1992}.
  Let \(j \in \integernumbers\).
  By \mytheorem \ref{th:upj-ineq} we have
  \begin{math}
    \nsznorminfty{f - \mdimuprojop{n}{j} f}
    \leq c_1 \omega(f ; 2^{-j} c_2) 
  \end{math}
  where \(c_1\) and \(c_2\) do not depend
  on \(j\) or \(f\).
  Since \(f\) is uniformly continuous
  \begin{displaymath}
    \lim_{t \to 0} \omega( f; t ) = 0
  \end{displaymath}
  and hence
  \begin{displaymath}
    \lim_{j \to \infty} \norminfty{f - \mdimuprojop{n}{j} f} = 0 .
  \end{displaymath}
  }
  \fi
\end{proof}
}
\fi

\if\shortprep1
{
\if10
{
\begin{corollary}
  \label{cor:mdimuv-union}
  \begin{displaymath}
    \overline{\bigcup_{j \in \integernumbers} \mdimuvspace{n}{j}} =
    \ucfunc{\rn}{K} .
  \end{displaymath}
\end{corollary}
}
\fi
}
\else
{
\begin{theorem}
  \label{th:mdimuv-union}
  \begin{displaymath}
    \overline{\bigcup_{j \in \integernumbers} \mdimuvspace{n}{j}} =
    \ucfunc{\rn}{K} .
  \end{displaymath}
\end{theorem}
}
\fi

\if\shortprep0
{
\begin{proof}
  Since \(\mdimuvspace{n}{j} \in \ucfunc{\rn}{K}\)
  for all
  \begin{math}
    j \in \integernumbers
  \end{math}
  it follows that
  \begin{math}
    \overline{\bigcup_{j \in \integernumbers} \mdimuvspace{n}{j}} \subset
    \ucfunc{\rn}{K}
  \end{math}.
  Let \(f \in \ucfunc{\rn}{K}\). Now
  \((\mdimuprojop{n}{j} f)_{j=0}^\infty\) is a sequence in the set
  \begin{math}
    \bigcup_{j \in \integernumbers} \mdimuvspace{n}{j}
  \end{math}.
  It follows from \mytheorem \ref{th:upj-convergence} 
  that \(\mdimuprojop{n}{j} f \to f\) as \(j \to \infty\). 
  Consequently
  \begin{math}
    f \in \overline{\bigcup_{j \in \integernumbers} \mdimuvspace{n}{j}}
  \end{math}.
\end{proof}
}
\fi

\if\shortprep0
{
By \mytheorem \ref{th:upj-convergence}
we have
\begin{equation}
  \label{eq:ucfunction-expansion-mdim}
  f = \mdimuprojop{n}{j_0} f + \sum_{j=j_0}^\infty 
  \mdimudeltaop{n}{j} f
\end{equation}
for all \(f \in \ucfunc{\realnumbers^n}{K}\).
If
\begin{displaymath}
  f = v + \sum_{j=j_0}^\infty w_j
\end{displaymath}
where \(v \in \mdimuvspace{n}{j_0}\) and
\(w_j \in \mdimuwspace{n}{j}\) for all \(j \in \integernumbers\),
\(j \geq j_0\),
it follows from \mylemmas
\ref{lem:mdimudeltaop-properties} (ii) and
\ref{lem:mdimpartialuprojop-orthogonality} 
that \(v = \mdimuprojop{n}{j_0} f\) and
\(w_j = \mdimudeltaop{n}{j} f\) for all
\(j \in \integernumbers\),
\(j \geq j_0\).
}
\fi

\if\shortprep0
{
\begin{definition}
  Let \(n \in \naturalnumbers + 2\).
  When \(j \in \integernumbers\) and
  \(l \in \naturalnumbers\)
  define
  \begin{displaymath}
    \scfsum{j}{l}{\rnx}{\rnh}
    :=
    \sum_{p=0}^l
    \left(
      \mdimmothersf{n}(2^j \rnx
      - \cubeordering{n}{p} - \floor{2^j \rnx}
      + 2^j \rnh)
      -
      \mdimmothersf{n}(2^j \rnx
      - \cubeordering{n}{p} - \floor{2^j \rnx})
    \right)
  \end{displaymath}
  for all
  \(\rnx, \rnh\) in \(\rn\).
\end{definition}
}
\fi

\if\shortprep0
{
\begin{lemma}
  \label{lem:scf-sum}
  Let \(n \in \naturalnumbers + 2\) and
  \(t \in \positiverealnumbers\).
  Then there exists \(c_1 \in \positiverealnumbers\),
  which may depend on \(t\),
  so that
  \begin{displaymath}
    \sum_{l=0}^\infty
      \abs{\scfsum{j}{l}{\rnx}{\rnh}}
    \leq c_1
  \end{displaymath}
  for all
  \(\rnx \in \rn\),
  \(\rnh \in \closedball{\rn}{0}{2^{-j}t}\),
  and
  \(j \in \integernumbers\).
\end{lemma}
}
\fi

\if\shortprep0
{
\begin{proof}
  Let \(\rnxone \in \rn\).
  Let
  \begin{displaymath}
    I := \{
      \firstznvar \in \zn
      \setsep
      \exists \rny \in \closedball{\rn}{2^j \rnxone}{t} :
      \mdimmothersf{n}(\rny - \firstznvar) \neq 0
    \} .
  \end{displaymath}
  There exists \(r_1 \in \positiverealnumbers\) so that
  \begin{math}
    \suppop \mdimmothersf{n} \subset
    \closedball{\rn}{0}{r_1}
  \end{math}.
  Let \(\indexedfirstznvar{1} \in I\).
  Now
  \(\mdimmothersf{n}(\rny - \indexedfirstznvar{1}) \neq 0\)
  for some
  \(\rny \in \closedball{\rn}{2^j \rnxone}{t}\).
  Furthermore,
  \begin{eqnarray*}
    \norm{\indexedfirstznvar{1} - \floor{2^j \rnxone}}
    & \leq &
    \norm{\indexedfirstznvar{1} - \rny}
    + \norm{\rny - 2^j \rnxone}
    + \norm{2^j \rnxone - \floor{2^j \rnxone}} \\
    & \leq & r_2 := r_1 + t + \sqrt{n} .
  \end{eqnarray*}
  Now
  \begin{math}
   I \subset
   K := \rectangle{\floor{2^j \rnxone} - \ceil{r_2}}
   {\floor{2^j \rnxone} + \ceil{r_2}}
  \end{math}
  and \(m := \card{K} = (2 \ceil{r_2} + 1)^n\).
  As \(\cubeorderingfunction{n}\) increases along
  zero-centred cubes we have
  \begin{math}
    \setimage{\cubeorderingfunction{n}}
    {\setzeroton{\card{K}-1}}
    +
    \floor{2^j \rnxone}
    = K
  \end{math}.
  Using \myequation \eqref{eq:sf-sum} we get
  \begin{eqnarray*}
    \scfsum{j}{l_1}{\rnxone}{\rnh}
    & = & \sum_{p=0}^{l_1}
    \biggl(
      \mdimmothersf{n}\left( 2^j \rnxone + 2^j \rnh - \cubeordering{n}{p}
      - \floor{2^j \rnxone} \right) \\
      & & -
      \mdimmothersf{n}\left( 2^j \rnxone - \cubeordering{n}{p}
      - \floor{2^j \rnxone} \right)
    \biggl) \\
    & = & 1 - 1 = 0
  \end{eqnarray*}
  for all \(\rnh \in \closedball{\rn}{0}{2^{-j}t}\) and
  \(l_1 \in \naturalnumbers + m\).
  Hence
  \begin{eqnarray*}
    \sum_{l=0}^\infty \abs{\scfsum{j}{l}{\rnxone}{\rnh}}
    & = &
    \sum_{l=0}^{m-1} \abs{\scfsum{j}{l}{\rnxone}{\rnh}} \\
    & \leq &
    \sum_{l=0}^{m-1}
    \sum_{p=0}^l    
    \biggl\vert
      \mdimmothersf{n}\left( 2^j \rnxone + 2^j \rnh - 
      \cubeordering{n}{p} - \floor{2^j \rnxone} \right) \\
      & & -
      \mdimmothersf{n}\left( 2^j \rnxone - \cubeordering{n}{p}
      - \floor{2^j \rnxone} \right)
    \biggl\vert \\
    & \leq & c_1 := 2 m^2 \norminfty{\mdimmothersf{n}}
  \end{eqnarray*}
  for all \(\rnh \in \closedball{\rn}{0}{2^{-j}t}\).
\end{proof}
}
\fi

\begin{theorem}
  \label{th:omega-pj-inequality}
  Let \(n \in \positiveintegers\) and \(t \in \positiverealnumbers\).
  There exists a constant
  \begin{math}
    c_1 \in \positiverealnumbers
  \end{math},
  which may depend on \(t\),
  so that
  \begin{displaymath}
    \forall f \in \ucfunccv{\realnumbers^n},
    j \in \integernumbers : \;
    \modcont{\mdimuprojop{n}{j} f}{2^{-j} t}
    \leq c_1 \modcont{f}{2^{-j} \sqrt{n}} .
  \end{displaymath}
\end{theorem}

\if\shortprep1
{
\begin{proof}
  The proof of case \(n = 1\) is similar to the
  proof in \cite[section 7.1]{donoho1992}.
  We will assume that \(n \geq 2\) in the sequel.
  Define bijection
  \begin{math}
    \cubeorderingfunction{n} : \naturalnumbers \onto \zn
  \end{math}
  so that
  \begin{math}
    \norminfty{\cubeordering{n}{p+1} - \cubeordering{n}{p}}
    = 1
  \end{math}
  for all
  \begin{math}
    p \in \naturalnumbers
  \end{math}.
  It can be proved that this kind of bijection always exists.
  Define
  \begin{displaymath}
    \scfsum{j}{l}{\rnx}{\rnh}
    :=
    \sum_{p=0}^l
    \left(
      \mdimmothersf{n}(2^j \rnx
      - \cubeordering{n}{p} - \floor{2^j \rnx}
      + 2^j \rnh)
      -
      \mdimmothersf{n}(2^j \rnx
      - \cubeordering{n}{p} - \floor{2^j \rnx})
    \right)
  \end{displaymath}
  where \(j, l \in \naturalnumbers\),
  \(\rnx \in \rn\), and \(\rnh \in \closedball{\rn}{0}{2^{-j}t}\)
  There exists \(c_1 \in \positiverealnumbers\) so that
  \begin{displaymath}
    \sum_{l=0}^\infty \abs{\scfsum{j}{l}{\rnx}{\rnh}} \leq c_1
  \end{displaymath}
  for all \(j \in \naturalnumbers\),
  \(\rnx \in \rn\), and \(\rnh \in \closedball{\rn}{0}{2^{-j}t}\).
  Let \(f \in \ucfunccv{\rn}\),
  \(j_1 \in \integernumbers\),
  \(\rnxone \in \rn\), and
  \(\rnhone \in \closedball{\rn}{0}{2^{-j}t}\).
  and
  Let
  \begin{math}
    \omegaseq
    := (f(2^{-j}\firstznvar))_{\firstznvar \in \zn}
  \end{math}.
  Now
  \begin{displaymath}
    z := 
    \left( \mdimuprojop{n}{j} f \right)
    \left( \rnxone + \rnhone \right)
    -
    \left( \mdimuprojop{n}{j} f \right)
    \left( \rnxone \right)
    =
    \sum_{l=0}^\infty
    \left(
      \szomegaseqelem{\cubeordering{n}{l}}
      -
      \szomegaseqelem{\cubeordering{n}{l+1}}
    \right)
    \scfsum{j_1}{l}{\rnxone}{\rnhone}
  \end{displaymath}
  from which it follows that
  \begin{displaymath}
    \abs{z}
    \leq
    \left(
      \sup_{l \in \naturalnumbers}
      \abs{\szomegaseqelem{\cubeordering{n}{l}}
        - \szomegaseqelem{\cubeordering{n}{l+1}}}
    \right)
    \cdot
    \sum_{l=0}^\infty
    \abs{\scfsum{j_1}{l}{\rnxone}{\rnhone}}
    \leq 
    \modcont{f}{2^{-j}\sqrt{n}} \cdot c_1 .
  \end{displaymath}
\end{proof}
}
\else
{
\begin{proof}
   The proof of case \(n = 1\) is similar to the
   proof in \cite[section 7.1]{donoho1992}.
   We will assume that \(n \geq 2\) in the sequel.
   
   Let \(f \in \ucfunccv{\rn}\) and
   \(j \in \integernumbers\).
   Let
   \begin{math}
     \omegaseq
     := (f(2^{-j}\firstznvar))_{\firstznvar \in \zn}
   \end{math}.
   \if\shortprep0
   {
   Now
   \begin{displaymath}
     \left( \mdimuprojop{n}{j} f \right)
     \left( \rnx + \rnh \right)
     -
     \left( \mdimuprojop{n}{j} f \right)
     \left( \rnx \right)
     =
     \sum_{\firstznvar \in \zn}
     \szomegaseqelem{\firstznvar}
     \left(
       \mdimmothersf{n}\left(2^j \rnx + 2^j \rnh
       - \firstznvar \right)
       -
       \mdimmothersf{n}\left(2^j \rnx - \firstznvar \right)
     \right)
   \end{displaymath}
   for all \(\rnx, \rnh \in \rn\).
   }
   \fi
   Assume that \(\rnxone \in \rn\) and
   \(\rnhone \in \closedball{\rn}{0}{2^{-j}t}\).
   Let
   \begin{math}
     \eta(p) := \cubeordering{n}{p} + \floor{2^j \rnxone}
   \end{math}
   for all
   \begin{math}
     p \in \naturalnumbers
   \end{math}.
   Define
   \begin{math}
     v_l := \scfsum{j}{l}{\rnxone}{\rnhone}
   \end{math}
   and
   \begin{displaymath}
     w_l := \szomegaseqelem{\cubeordering{n}{l}}
     \left(
       \mdimmothersf{n}\left(2^j \rnxone + 2^j \rnhone - 
       \eta(l)\right)
       -
       \mdimmothersf{n}\left(2^j \rnxone - 
       \eta(l)\right)
     \right)
   \end{displaymath}
   for all \(l \in \naturalnumbers\).
   Now
   \begin{eqnarray*}
     w_l & = & \szomegaseqelem{\cubeordering{n}{l}}
     \left(
       \sum_{p=0}^l
         \left(
           \mdimmothersf{n}\left(2^j \rnxone + 2^j \rnhone
           - \eta(p)\right) -
           \mdimmothersf{n}\left(2^j \rnxone 
           - \eta(p)\right)
         \right)
     \right. \\
       & & -
     \left.
       \sum_{p=0}^{l-1}
         \left(
           \mdimmothersf{n}\left(2^j \rnxone + 2^j \rnhone
           - \eta(p)\right) -
           \mdimmothersf{n}\left(2^j \rnxone 
           - \eta(p)\right)
         \right)
     \right) \\
     & = &
     \szomegaseqelem{\cubeordering{n}{l}} v_l
     - \szomegaseqelem{\cubeordering{n}{l}} v_{l-1}
   \end{eqnarray*}
   for all
   \begin{math}
     l \in \positiveintegers
   \end{math}.
   \if\shortprep0
   {
   Furthermore,
   \begin{eqnarray*}
     \sum_{l=1}^\infty w_l
     & = &
     \sum_{l=1}^\infty \szomegaseqelem{\cubeordering{n}{l}} v_l
     -
     \sum_{l=1}^\infty \szomegaseqelem{\cubeordering{n}{l}} v_{l-1} \\
     & = &
     \sum_{l=1}^\infty \szomegaseqelem{\cubeordering{n}{l}} v_l
     -
     \sum_{l=0}^\infty \szomegaseqelem{\cubeordering{n}{l+1}} v_l \\
     & = &
     - \szomegaseqelem{\cubeordering{n}{1}} v_0
     + \sum_{l=1}^\infty
     \left(
       \szomegaseqelem{\cubeordering{n}{l}}
       -
       \szomegaseqelem{\cubeordering{n}{l+1}}
     \right)
     v_l .
   \end{eqnarray*}
   We also have
   \begin{math}
     w_0 = \szomegaseqelem{\cubeordering{n}{0}} v_0
   \end{math}
   from which it follows that
   \begin{eqnarray*}
     z & := &
     \sum_{l=0}^\infty w_l \\
     & = &
     \szomegaseqelem{\cubeordering{n}{0}} v_0
     - \szomegaseqelem{\cubeordering{n}{1}} v_0
     + \sum_{l=1}^\infty
     \left(
       \szomegaseqelem{\cubeordering{n}{l}}
       -
       \szomegaseqelem{\cubeordering{n}{l+1}}
     \right)
     v_l \\
     & = &
     \sum_{l=0}^\infty
     \left(
       \szomegaseqelem{\cubeordering{n}{l}}
       -
       \szomegaseqelem{\cubeordering{n}{l+1}}
     \right)
     v_l .
   \end{eqnarray*}
   Consequently
   \begin{displaymath}
     \abs{z}
     \leq
     \left(
       \sup_{l \in \naturalnumbers}
       \abs{\szomegaseqelem{\cubeordering{n}{l}}
         - \szomegaseqelem{\cubeordering{n}{l+1}}}
     \right)
     \cdot
     \sum_{l=0}^\infty \abs{v_l}
   \end{displaymath}
   Function \(\cubeorderingfunction{n}\)
   preserves neighbours from which it follows
   that
   \begin{eqnarray*}
     & & \norminfty{\cubeordering{n}{l} - \cubeordering{n}{l+1}}
     = 1 \\
     & & \implies
     \normtwo{\cubeordering{n}{l} - \cubeordering{n}{l+1}}
     \leq \sqrt{n} \\
     & & \implies
     \normtwo{2^{-j} \cubeordering{n}{l}
     - 2^{-j} \cubeordering{n}{l+1}}
     \leq 2^{-j} \sqrt{n} .
   \end{eqnarray*}
   for all
   \begin{math}
     l \in \naturalnumbers
   \end{math}.
   Consequently
   \begin{equation}
     \label{eq:coeff-bound-b}
     \sup_{l \in \naturalnumbers}
     \abs{\szomegaseqelem{\cubeordering{n}{l}}
       - \szomegaseqelem{\cubeordering{n}{l+1}}}
     \leq
     \modcont{f}{2^{-j}\sqrt{n}} .
   \end{equation}
   It follows from \mylemma \ref{lem:scf-sum}
   and \myequation \eqref{eq:coeff-bound-b}
   that
   \begin{displaymath}
     \abs{\left( \mdimuprojop{n}{j} f \right)
     \left( \rnxone + \rnhone \right)
     -
     \left( \mdimuprojop{n}{j} f \right)
     \left( \rnxone \right)}
     = \abs{z}
     \leq
     c_1 \modcont{f}{2^{-j}\sqrt{n}}
   \end{displaymath}
   }
   \else
   {
   \if10
   {
   Furthermore,
   \begin{eqnarray*}
     & &
     \left( \mdimuprojop{n}{j} f \right)
     \left( \rnxone + \rnhone \right)
     -
     \left( \mdimuprojop{n}{j} f \right)
     \left( \rnxone \right)
     =
     \sum_{l=0}^\infty w_l
     =
     w_0 + \sum_{l=1} \szomegaseqelem{\cubeordering{n}{l}} v_l
     -
     \sum_{l=1} \szomegaseqelem{\cubeordering{n}{l}} v_{l-1} \\
     & & =
     \szomegaseqelem{\cubeordering{n}{0}} v_0
     -
     \szomegaseqelem{\cubeordering{n}{1}} v_0
     +
     \sum_{l=1}
     \left(
       \szomegaseqelem{\cubeordering{n}{l}}
       -
       \szomegaseqelem{\cubeordering{n}{l+1}}       
     \right) v_l \\
     & & =
     \sum_{l=0}
     \left(
       \szomegaseqelem{\cubeordering{n}{l}}
       -
       \szomegaseqelem{\cubeordering{n}{l+1}}       
     \right) v_l
   \end{eqnarray*}
   from which the theorem follows.
   }
   \fi
   }
   \fi
\end{proof}
}
\fi

\begin{theorem}
  Let \(j_0 \in \integernumbers\).
  Suppose that the mother scaling function \(\myphi\) is Lipschitz 
  continuous and let \(\mdimmothersf{n}\) be
  the tensor product mother scaling function generated
  by \(\myphi\).
  Then
  \begin{displaymath}
    \mdimuvspace{n}{j_0}
    \dsum \indexeddirectsum{j=j_0}{\infty} \mdimuwspace{n}{j}
    \neq
    \ucfunc{\rn}{K} .
  \end{displaymath}
\end{theorem}

\begin{proof}
Let
\begin{displaymath}
  A := \mdimuvspace{n}{j_0}
  \dsum \indexeddirectsum{j=j_0}{\infty}
  \mdimuwspace{n}{j} .
\end{displaymath}
Define function \(f \in \ucfunc{\realnumbers}{K}\) by
\begin{displaymath}
  f(x) := \left\{
    \begin{array}{ll}
      \sqrt{x} ; & x \in [ 0, 1 ] \\
      -x + 2 ; & x \in [ 1, 2 ] \\
      0 ; & x \leq 0 \lor x \geq 2
    \end{array}
    \right.
\end{displaymath}
and function \(f^{[n]} \in \ucfunc{\rn}{K}\) by
\begin{displaymath}
  f^{[n]} := \bigotimes_{k=1}^n f .
\end{displaymath}
Functions \(f\) and \(f^{[n]}\) are not Lipschitz continuous.
As \(\myphi\) is Lipschitz
continuous
all the functions \(\mdimuprojop{n}{j} g\)
are Lipschitz continuous for each
\(g \in \ucfunc{\realnumbers^n}{K}\).
\(A\) is a locally convex space with an inductive limit
topology. It follows from
\cite[theorem 6.2]{schaefer1966} that the space \(A\) is complete.

Suppose that \(A\) would be equal to \(\ucfunc{\rn}{K}\)
as a set.
Then we would have \(f^{[n]} \in A\). It follows from the definition of
the locally convex direct sum that
\begin{displaymath}
  f^{[n]} \in \mdimuvspace{n}{j_0} \dsum 
  \indexeddirectsum{j=j_0}{j_1}
  \mdimuwspace{n}{j}
\end{displaymath}
for some \(j_1 \in \integernumbers\), \(j_1 \geq j_0\).
Now \(f^{[n]} = \mdimuprojop{n}{j_1+1} f^{[n]}\). It follows that
\(f^{[n]}\) is Lipschitz continuous, which is a contradiction. Hence
\(A\) is not equal to
\(\ucfunc{\rn}{K}\) as a set.
\end{proof}

The definition of locally convex direct sums from
\cite{schaefer1966} is used here.
An example of a Lipschitz continuous mother scaling function is
a Deslauriers-Dubuc fundamental function with Hölder regularity greater than 1.
By \mylemma \ref{lem:mdimuvspace-direct-sum} and \mytheorem
\ref{th:mdimuv-union} we have
\begin{displaymath}
  \ucfunc{\realnumbers^{n}}{K}
  \equalns
  \closop \left( \bigcup_{l=j_0}^\infty \left( 
  \mdimuvspace{n}{j_0} \dsum
    \indexeddirectsum{j=j_0}{l} \mdimuwspace{n}{j} \right) 
    \right)
\end{displaymath}
for all \(j_0 \in \integernumbers\).

\if\shortprep0
{

\section{Compactly Supported Interpolating MRA
  \mraprepspace $\vanishingfunc{\realnumbers^n}{K}$}
\label{sec:van-mra}

We shall assume that
\(K = \realnumbers\) or \(K = \complexnumbers\) throughout this section.
We shall also assume that
\begin{math}
  \myphi \in \cscfunc{\realnumbers}{K}
\end{math}
is a function for which conditions (MSF.1) and
(MSF.2) hold.
Unless otherwise stated, we
shall assume that 
the same values of \(K\) and \(\myphi\) are used throughout this
section.

\subsection{General}

Donoho \cite{donoho1992} constructs
projection operators and gives convergence results for interpolating
wavelets on \(\vanishingfunccv{\realnumbers}\).

\begin{definition}
  \label{def:vj-mdim}
  \definekn
  Let \(\myphi \in \cscfunc{\realnumbers}{K}\) be a 
  function for which conditions (MSF.1) and (MSF.2)
  hold and let \(\mdimmothersf{n}\) be defined by
  \mydef \ref{def:mdimsf}.
  Define
  \begin{eqnarray}
    \nonumber
    & & \mdimvanvspace{n}{j} := \left\{ \sum_{\firstznvar \in \zn}
     \seqelem{\seqstyle{a}}{\firstznvar}
      \mdimmothersf{n}(2^j \cdot - \firstznvar)
      \setsep \seqstyle{a} \in \gencospace{\zn}{K}
      \right\} \\
    \label{eq:vj-mdim}
    & & \norminspace{f}{\mdimvanvspace{n}{j}} := \norminfty{f} ,
    \spaceafter f \in \mdimvanvspace{n}{j}
  \end{eqnarray}
  for each
  \(j \in \integernumbers\).
\end{definition}

\if\shortprep0
{
The infinite sum in \myequation \eqref{eq:vj-mdim} converges
unconditionally by \mylemma \ref{lem:grid-convergence-multidim}.
}
\fi

\begin{definition}
Let spaces
\(\mdimvanvspace{n}{j}\), \(j \in \integernumbers\),
be defined by \mydef \ref{def:vj-mdim}.
We call
\(\{\mdimvanvspace{n}{j} \setsep j \in \integernumbers\}\)
an interpolating tensor product
MRA \mraprepspace \(\vanishingfunc{\rn}{K}\)
generated by
\(\mdimmothersf{n}\) provided that the following conditions are
satisfied:
\begin{myenumerate}
  \myitem{(MRA2.1)}
    \begin{math}
    \forall j \in \integernumbers :
      \mdimvanvspace{n}{j} \subset \mdimvanvspace{n}{j+1}    
    \end{math}
  \myitem{(MRA2.2)}
    \begin{math}
      \overline{\bigcup_{j \in \integernumbers} 
      \mdimvanvspace{n}{j}} =
      \vanishingfunc{\realnumbers^n}{K}
    \end{math}
  \myitem{(MRA2.3)}
    \begin{math}
      \bigcap_{j \in \integernumbers} \mdimvanvspace{n}{j} = \{ 0 \}
    \end{math}
  \myitem{(MRA2.4)}
    \begin{math}
      \forall j \in \integernumbers, f \in K^{\rn} :
      f \in \mdimvanvspace{n}{j} \iff f(2 \cdot) \in 
      \mdimvanvspace{n}{j+1}
    \end{math}
  \myitem{(MRA2.5)}
    \begin{math}
      \forall j \in \integernumbers, \firstznvar \in \zn,
      f \in K^{\rn} :
      f \in \mdimvanvspace{n}{j} \iff
      f(\cdot - 2^{-j} \firstznvar) \in \mdimvanvspace{n}{j}
    \end{math}
  \myitem{(MRA2.6)}
    \begin{math}
      \forall \firstznvar \in \zn :
      \mdimmothersf{n}(\firstznvar) = \delta_{\firstznvar, 0}
    \end{math}
\end{myenumerate}
\end{definition}

Note that in this definition the intersection of spaces \(\mdimvanvspace{n}{j}\) is
\(\{ 0 \}\) instead of \(\complexnumbers\) (all complex valued
constant functions on \(\rn\)) as in \cite{cl1996} since here the MRA is
constructed for functions vanishing at infinity.
Donoho \cite{donoho1992} does not require the mother scaling function
\(\myphi\) to be compactly supported but \(\myphi\) has to be of rapid decay
in his construction. He also includes requirements for the regularity and
polynomial span of \(\myphi\), which are needed for the norm equivalences
to Besov and Triebel-Lizorkin spaces, into the definition of the MRA.
\if\longprep1
{
\givelemmawithoutproof
}
\fi
\if\shortprep1
{
  Under the condition that \(\mdimmothersf{n}\)
  is an \(n\)-dimensional tensor product
  mother scaling function
  the conditions
  (MRA2.1), (MRA2.4), (MRA2.5), and (MRA2.6) are true.
}
\else
{

\begin{lemma}
  Under the condition that \(\mdimmothersf{n}\)
  is an \(n\)-dimensional tensor product
  mother scaling function
  the conditions
  (MRA2.1), (MRA2.2), (MRA2.3), (MRA2.4), (MRA2.5), and (MRA2.6) are
  true.
\end{lemma}

See also \mylemma \ref{lem:uv-mra-properties}.
}
\fi

\if11
\subsection{Basic Definitions for the Univariate MRA \mraprepspace 
$\vanishingfunc{\realnumbers}{K}$}
\else
\subsection{Univariate MRA \mraprepspace 
\(\vanishingfunc{\realnumbers}{K}\)}
\fi

\if\shortprep0
{
By \mylemma \ref{lem:cospace-top-isomorphism} spaces
\(\onedimvanvspace{j}\), \(j \in \integernumbers\), are closed subspaces of
\(\vanishingfunc{\realnumbers}{K}\).
}
\fi

\begin{definition}
  \label{def:proj-op-one-dim}
  When \(j \in \integernumbers\) define operator \( \onedimvanprojop{j} :
  \vanishingfunc{\realnumbers}{K} \to \onedimvanvspace{j} \) by
  \begin{equation}
    \label{eq:proj-op-one-dim}
    \onedimvanprojop{j} f := \sum_{k \in \integernumbers}
      \dualappl{\phidual_{j,k}}{f} \myphi_{j,k} .
  \end{equation}
\end{definition}

\if\shortprep1
{
Operator \(\onedimvanprojop{j}\) is a continuous
projection onto Banach space \(\onedimvanvspace{j}\) for each
\(j \in \integernumbers\).
}
\else
{
The infinite sum in \eqref{eq:proj-op-one-dim} converges unconditionally
by \mylemma \ref{lem:grid-convergence-multidim}.
By \mydef \ref{def:vj-mdim}
\begin{displaymath}
    \sum_{k \in \integernumbers} f \left( \frac{k}{2^j}
    \right) \myphi \left( 2^j \cdot - k \right) \in \onedimvanvspace{j} 
\end{displaymath}
for all \(f \in \vanishingfunc{\realnumbers}{K}\) and
\(j \in \integernumbers\).
Hence functions \(\onedimvanprojop{j}\), \(j \in \integernumbers\), are well defined.
\if\shortprep1
By \mylemma \ref{lem:grid-convergence-multidim} functions
\(\onedimvanprojop{j}\), \(j \in \integernumbers\), are operators.
\else
Functions \(\onedimvanprojop{j}\), \(j \in \integernumbers\),  are linear by their definition.
We have
\begin{displaymath}
  \norminfty{\left( f \left( \frac{k}{2^j} \right) \right)_{k \in
      \integernumbers}} \leq \norminfty{f}
\end{displaymath}
for all \(f \in \vanishingfunc{\realnumbers}{K}\).
Hence by \mylemma \ref{lem:grid-convergence-multidim} functions
\(\onedimvanprojop{j}\), \(j \in \integernumbers\), are continuous and uniformly
bounded by
\begin{equation}
  \label{eq:one-dim-proj-op-uniform-bound-ineq}
  \norm{\onedimvanprojop{j}} \leq \ncover{\myphi} 
  \norminfty{\myphi}
\end{equation}
for all
\begin{math}
  j \in \integernumbers
\end{math}.
\fi
It follows from \mydefs \ref{def:vj-mdim} and
\ref{def:proj-op-one-dim} that operator \(\onedimvanprojop{j}\) is a
projection onto Banach space \(\onedimvanvspace{j}\) for each
\(j \in \integernumbers\).
}
\fi

\if\shortprep0
{
\begin{lemma}
  \label{lem:proj-op-composition-one-dim}
  \(\onedimvanprojop{j} = \onedimvanprojop{j} \circ \onedimvanprojop{j'}\) for all
  \(j, j' \in \integernumbers\),
  \(j' \geq j\).
\end{lemma}
}
\fi

\if\shortprep0
{
\begin{proof}
  Let \(j, j' \in \integernumbers\) and \(j' \geq j\).
  Let \(f \in \vanishingfunc{\realnumbers}{K}\).
  Now
  \if\elsevier0
  {
  \begin{align*}
    \onedimvanprojop{j} ( \onedimvanprojop{j'} f ) & =
    \sum_{k \in \integernumbers}
    \szdualappl{\phidual_{j,k}}
    {\sum_{k' \in \integernumbers}
      \dualappl{\phidual_{j',k'}}{f} \myphi_{j',k'}}
    \myphi_{j,k} \\
    & =
    \sum_{k \in \integernumbers}
    \sum_{k' \in \integernumbers}
    \dualappl{\phidual_{j',k'}}{f}
    \dualappl{\phidual_{j,k}}{\myphi_{j',k'}}
    \myphi_{j,k} \\
    & =
    \sum_{k \in \integernumbers}
    \sum_{k' \in \integernumbers}
    \dualappl{\phidual_{j',k'}}{f}
    \delta_{k',2^{j'-j} k}
    \myphi_{j,k}
    =
    \sum_{k \in \integernumbers}
    \dualappl{\phidual_{j',2^{j'-j} k}}{f}
    \myphi_{j,k} \\
    & =
    \sum_{k \in \integernumbers}
    \dualappl{\phidual_{j,k}}{f}
    \myphi_{j,k}
    = \onedimvanprojop{j} f .
  \end{align*}
  }
  \else
  {
  \begin{align*}
    \onedimvanprojop{j} ( \onedimvanprojop{j'} f ) & =
    \sum_{k \in \integernumbers}
    \szdualappl{\phidual_{j,k}}
    {\sum_{k' \in \integernumbers}
      \dualappl{\phidual_{j',k'}}{f} \myphi_{j',k'}}
    \myphi_{j,k}
    =
    \sum_{k \in \integernumbers}
    \sum_{k' \in \integernumbers}
    \dualappl{\phidual_{j',k'}}{f}
    \dualappl{\phidual_{j,k}}{\myphi_{j',k'}}
    \myphi_{j,k} \\
    & =
    \sum_{k \in \integernumbers}
    \sum_{k' \in \integernumbers}
    \dualappl{\phidual_{j',k'}}{f}
    \delta_{k',2^{j'-j} k}
    \myphi_{j,k}
    =
    \sum_{k \in \integernumbers}
    \dualappl{\phidual_{j',2^{j'-j} k}}{f}
    \myphi_{j,k}
    =
    \sum_{k \in \integernumbers}
    \dualappl{\phidual_{j,k}}{f}
    \myphi_{j,k} \\
    & = \onedimvanprojop{j} f .
  \end{align*}
  }
  \fi
\end{proof}
}
\fi

\if10
{
We get the following theorem as a direct
consequence of \mytheorem \ref{th:upj-ineq}.
}
\fi

\if10
{
\begin{theorem}
  There exists
  \(c_1 \in \positiverealnumbers\)
  and
  \(c_2 \in \positiverealnumbers\)
  so that
  \begin{displaymath}
    \norminfty{f - \mdimvanprojop{1}{j} f}
    \leq
    c_1 \modcont{f}{2^{-j} c_2}
  \end{displaymath}
  for all
  \(f \in \vanishingfunccv{\realnumbers}\)
  and
  \(j \in \integernumbers\).
\end{theorem}
}
\fi

\if\shortprep0
{
We get the following theorem as a direct
consequence of \mytheorem \ref{th:upj-convergence}.
}
\fi

\if\shortprep0
{
\begin{theorem}
  \label{th:one-dim-vanpj-convergence}
  Let \(f \in \vanishingfunccv{\realnumbers}\). Now
  \begin{displaymath}
    \lim_{j \to \infty} \norminfty{f - \onedimvanprojop{j} f} = 0 .
  \end{displaymath}
\end{theorem}
}
\fi

\if10
{
We get the following lemma as a direct consequence
of \mytheorem \ref{th:omega-pj-inequality}.
}
\fi

\if10
{
\begin{lemma}
  \label{lem:one-dim-omega-inequality}
  Let \(t \in \positiverealnumbers\).
  There exists a constant
  \begin{math}
    c_1 \in \positiverealnumbers
  \end{math},
  which may depend on \(t\),
  so that
  \begin{displaymath}
    \forall f \in \vanishingfunccv{\realnumbers} : 
    \modcont{\onedimvanprojop{j} f}{2^{-j} t}
    \leq
    c_1 \modcont{f}{2^{-j}} .
  \end{displaymath}
\end{lemma}
}
\fi

\begin{definition}
  \label{def:wj-one-dim}
  Define
  \begin{math}
    \onedimvanwspace{j} = \left\{ f - \onedimvanprojop{j} f \setsep f \in V_{j + 1} \right\}
  \end{math}
  for all
  \begin{math}
    j \in \integernumbers
  \end{math}.
\end{definition}

\if\shortprep1
{
  \(\onedimvanwspace{j}\) is a closed subspace
  of \(\onedimvanvspace{j + 1}\) for each \(j \in \integernumbers\) and
  \(\onedimvanvspace{j} \cap \onedimvanwspace{j} = \{ 0 \}\)
  for all \(j \in \integernumbers\).
}
\else
{
\begin{lemma}
  \label{lem:one-dim-wj-cs}
  \(\onedimvanwspace{j}\) is a closed subspace
  of \(\onedimvanvspace{j + 1}\) for each \(j \in \integernumbers\) and
  \(\onedimvanvspace{j} \cap \onedimvanwspace{j} = \{ 0 \}\)
  for all \(j \in \integernumbers\).
\end{lemma}
}
\fi

\if\shortprep0
{
\begin{proof}
  If \(f \in \onedimvanvspace{j}\) then \(\onedimvanprojop{j} f = f\). If \(g \in \onedimvanwspace{j}\) and \(g \neq
  0\) then \(g = h - \onedimvanprojop{j} h\) where \(h \in V_{j + 1}\) and \(\onedimvanprojop{j} g =
  0\). So \(\onedimvanprojop{j} g \neq g\) and \(g \not\in \onedimvanvspace{j}\). Hence \(\onedimvanvspace{j} \cap \onedimvanwspace{j}
  = \{ 0 \}\).

  Let \((f_k)_{k=0}^\infty \subset \onedimvanwspace{j}\) be a convergent sequence
  in Banach space \(\onedimvanvspace{j+1}\). Then \(f_k \to g\) as
  \(k \to \infty\), where \(g \in \onedimvanvspace{j+1}\). Since \(\onedimvanprojop{j}\) is
  continuous \(\onedimvanprojop{j} f_k \to \onedimvanprojop{j} g\) as \(k \to \infty\). Now 
  \(\onedimvanprojop{j} f_k
  = 0\) for all \(k \in \naturalnumbers\) and hence \(\onedimvanprojop{j} g =
  0\). Furthermore \(g = g - \onedimvanprojop{j} g \in \onedimvanwspace{j}\). Hence \(\onedimvanwspace{j}\) is closed.
\end{proof}
}
\fi

\begin{definition}
  \label{def:delta-proj-op-one-dim}
  When \(j \in \integernumbers\) define operator
  \( \onedimvandeltaop{j} : \vanishingfunc{\realnumbers}{K} \to \onedimvanwspace{j} \)
  by
  \begin{math}
    \onedimvandeltaop{j} := \onedimvanprojop{j+1} - \onedimvanprojop{j}
  \end{math}.
\end{definition}

It follows that
\begin{math}
  \onedimvanvspace{j+1} = \onedimvanvspace{j} \dsum \onedimvanwspace{j}
\end{math}
for all
\begin{math}
  j \in \integernumbers
\end{math}.
\if\shortprep1
Operator
\begin{math}
  \onedimvandeltaop{j}
\end{math}
is a continuous projection onto
\begin{math}
  \onedimvanwspace{j}
\end{math}
for all
\begin{math}
  j \in \integernumbers
\end{math}.
\else
It follows from \mylemma \ref{lem:proj-op-composition-one-dim} and \myequation
\eqref{eq:one-dim-proj-op-uniform-bound-ineq} that
operator 
\(\onedimvandeltaop{j} : \vanishingfunc{\realnumbers}{K} \to \onedimvanwspace{j} \) is a
continuous projection onto \(\onedimvanwspace{j}\)
\if10
and
\begin{equation}
  \label{eq:one-dim-delta-uniform-bound-ineq}
  \norm{\onedimvandeltaop{j}} \leq 2 \ncover{\myphi} 
  \norminfty{\myphi}
\end{equation}
\fi
for all
\begin{math}
  j \in \integernumbers
\end{math}.
If \(x = v + w \in \onedimvanvspace{j+1}\), \(v \in \onedimvanvspace{j}\), and
\(w \in \onedimvanwspace{j}\) we have \(\onedimvanprojop{j} x = v\) and
\(\onedimvandeltaop{j} x = w\).
Consequently
\begin{equation}
  \label{eq:one-dim-projop-orth-a}
  \forall j' \in \integernumbers, f \in \onedimvanvspace{j'} :
  j' \leq j \implies
  \onedimvandeltaop{j} f = 0 ,
\end{equation}
\begin{equation}
  \label{eq:one-dim-projop-orth-b}
  \forall j' \in \integernumbers, f \in \onedimvanwspace{j'} :
  j' \neq j \implies
  \onedimvandeltaop{j} f = 0 ,
\end{equation}
and
\begin{equation}
  \label{eq:one-dim-projop-orth-c}
  \forall j' \in \integernumbers, f \in \onedimvanwspace{j'} :
  j' \geq j \implies  
  \onedimvanprojop{j} f = 0 .
\end{equation}
\fi

\if10
{
\begin{definition}
  Define function
  \begin{math}
    \vanwjtopisomfunc{1}{j} : \gencospace{\integernumbers}{K}
    \to \onedimvanwspace{j}
  \end{math}
  by
  \begin{displaymath}
    \vanwjtopisom{1}{j}{\seqstyle{a}}
    :=
    \sum_{k \in \integernumbers}
    \seqelem{\seqstyle{a}}{k}
    \psi_{j,k}
  \end{displaymath}
  for all
  \begin{math}
    \seqstyle{a} \in 
    \gencospace{\integernumbers}{K} 
  \end{math}.
\end{definition}
}
\fi

\begin{theorem}
  \label{th:wj-series-representation-one-dim}
  Let \(j \in \integernumbers\).
  \if11
  Then
  \begin{equation}
    \label{eq:wj-expansion-one-dim}
    \onedimvanwspace{j} = \left\{ \sum_{k \in \integernumbers} \seqelem{\seqstyle{a}}{k}
      \psi_{j,k}
      \setsep \seqstyle{a}
      \in \gencospace{\integernumbers}{K} \right\} .
  \end{equation}
  \else
  Then
  \begin{equation}
    \label{eq:wj-expansion-one-dim}
    \onedimvanwspace{j} = \left\{ \sum_{k \in \integernumbers} \seqelem{\seqstyle{a}}{k}
      \psi_{j,k}
      \setsep \seqstyle{a}
      \in \gencospace{\integernumbers}{K} \right\},
  \end{equation}
  function
  \begin{math}
    \vanwjtopisomfunc{1}{j}
  \end{math}
  is a topological isomorphism from
  \begin{math}
    \gencospace{\integernumbers}{K}
  \end{math}
  onto
  \begin{math}
    \onedimvanwspace{j}
  \end{math},
  and
  \begin{equation}
    \label{eq:wj-norm-equivalence}
    \norminfty{\seqstyle{a}}
    \leq \norminfty{\vanwjtopisom{1}{j}{\seqstyle{a}}}
    \leq \ncover{\psi} \norminfty{\psi} \norminfty{\seqstyle{a}}
  \end{equation}
  for all
  \begin{math}
    \seqstyle{a} \in \gencospace{\integernumbers}{K}
  \end{math}
  and for each
  \begin{math}
    j \in \integernumbers
  \end{math}.
  \fi
\end{theorem}

\if11
{
\begin{proof}
  \if\shortprep0
  {
  Use \mylemmas \ref{lem:grid-convergence-multidim}
  and \ref{lem:cospace-top-isomorphism}.
  }
  \fi
  The proof \myequation \eqref{eq:wj-expansion-one-dim} is similar to the
  beginning of the proof of \cite[theorem 2.4]{cl1996}.
\end{proof}
}
\else
{
\begin{proof}
  By \mylemma \ref{lem:grid-convergence-multidim} the series in
  \eqref{eq:wj-expansion-one-dim} converges unconditionally.
  The proof \myequation \eqref{eq:wj-expansion-one-dim} is similar to the
  beginning of the proof of \cite[theorem 2.4]{cl1996}.
  We have
  \begin{math}
    \psi \left( k + \frac{1}{2} \right) = \delta_{k,0}
  \end{math}
  for all \(k \in \integernumbers\).
  Hence by \mylemma \ref{lem:cospace-top-isomorphism}
  vector space \(\onedimvanwspace{j}\) is topologically isomorphic to
  \(\gencospace{\integernumbers}{K}\) and
  \myequation \eqref{eq:wj-norm-equivalence} is true.
\end{proof}
}
\fi

\if\shortprep0
{
It follows from \mytheorem \ref{th:one-dim-vanpj-convergence}
that
\begin{displaymath}
  f = \onedimvanprojop{j_0} f + \sum_{j=j_0}^\infty
   \onedimvandeltaop{j} f
\end{displaymath}
for all
\begin{math}
  f \in \vanishingfunc{\realnumbers}{K}
\end{math}
where \(j_0 \in \integernumbers\).
If
\begin{math}
  f \in \vanishingfunc{\realnumbers}{K}
\end{math}
and
\begin{displaymath}
  f = v + \sum_{j=j_0}^\infty w_j
\end{displaymath}
where \(v \in V_{j_0}\) and \(w_j \in \onedimvanwspace{j}\) for all
\(j \in \integernumbers\), \(j \geq j_0\), 
\if\shortprep1
we have
\else
it follows from \myequations \eqref{eq:one-dim-projop-orth-a},
\eqref{eq:one-dim-projop-orth-b}, and \eqref{eq:one-dim-projop-orth-c}
that
\fi
\(v = \onedimvanprojop{j_0} f\) and
\(\onedimvanwspace{j} = \onedimvandeltaop{j} f\) for all \(j \in \integernumbers\),
\(j \geq j_0\).
}
\fi

\if\shortprep1
{
  We have
  \begin{displaymath}
    \label{eq:deltaop-formula-one-dim}
    \onedimvandeltaop{j} f = \sum_{k \in \integernumbers}
    \dualappl{\psidual_{j,k}}{f} \psi_{j,k}
  \end{displaymath}
  for all \(j \in \integernumbers\) and
  \(f \in \vanishingfunc{\realnumbers}{K}\).
}
\else
{
\begin{lemma}
  \label{lem:deltaop-formula-one-dim}
  We have
  \begin{equation}
    \label{eq:deltaop-formula-one-dim}
    \onedimvandeltaop{j} f = \sum_{k \in \integernumbers}
    \dualappl{\psidual_{j,k}}{f} \psi_{j,k}
  \end{equation}
  for all \(j \in \integernumbers\) and
  \(f \in \vanishingfunc{\realnumbers}{K}\).
  The series in \myequation \eqref{eq:deltaop-formula-one-dim} converges
  unconditionally.
\end{lemma}
}
\fi

\if\shortprep1
{
\if10
{
\begin{proof}
  Use \mylemmas \ref{lem:mdim-vanishing-samples} and
  \ref{lem:grid-convergence-multidim}.
\end{proof}
}
\fi
}
\else
{
\begin{proof}
  Since \(\onedimvandeltaop{j} f \in \onedimvanwspace{j}\) it follows that
  \begin{displaymath}
    \onedimvandeltaop{j} f = \sum_{k \in \integernumbers} a_k \psi_{j,k}
  \end{displaymath}
  where
  \begin{math}
    (a_k)_{k \in \integernumbers} \sequenceof K
  \end{math}.
  Furthermore,
  \begin{math}
    f = \onedimvanprojop{j_0} f + \sum_{j=j_0}^\infty 
    \onedimvandeltaop{j} f
  \end{math}
  and it follows that
  \(\dualappl{\psidual_{j,k}}{f} =
  \dualappl{\psidual_{j,k}}{\onedimvandeltaop{j} f} = a_k\)
  for all \(k \in \integernumbers\).
  As \(\psidual\) is a finite linear combination of functionals
  \(\delta(2 \cdot - k)\), \(k \in \integernumbers\), it follows from
  \mylemmas \ref{lem:mdim-vanishing-samples} and
  \ref{lem:grid-convergence-multidim} that the series in \myequation
  \eqref{eq:deltaop-formula-one-dim} converges unconditionally.
  Hence the theorem is true.
\end{proof}
}
\fi

\begin{definition}
  \label{def:ujspace}
  Define
  \begin{displaymath}
    \uspace{b}{j} =
    \left\{
      \begin{array}{ll}
        \onedimvanvspace{j} ; & b = 0 \\
        \onedimvanwspace{j} ; & b = 1
      \end{array}
    \right.
  \end{displaymath}
  for all \(j \in \integernumbers\) and \(b \in \{ 0, 1 \}\).
\end{definition}

\begin{definition}
  \label{def:genprojoponedim}
  Define operators
  \begin{math}
    \genprojoponedim{b}{j} : \vanishingfunc{\realnumbers}{K} \to
    \uspace{b}{j}
  \end{math}
  by
  \begin{displaymath}
    \genprojoponedim{b}{j} =
    \left\{
      \begin{array}{ll}
        \onedimvanprojop{j} ; & b = 0 \\
        \onedimvandeltaop{j} ; & b = 1
      \end{array}
    \right.
  \end{displaymath}
  for all \(j \in \integernumbers\),
  \(b \in \{ 0, 1 \}\).
\end{definition}

\subsection{Multivariate MRA
\mraprepspace $\vanishingfunc{\rn}{K}$}
\label{sec:mdim-mra}

We shall assume that
\(n \in \positiveintegers\) throughout this subsection.
Unless otherwise stated, we shall assume that 
the same value of \(n\) is used throughout this
subsection.
By \mylemma \ref{lem:vanishing-tp-inclusion} function \(\mdimmothersf{n}\)
belongs to \(\vanishingfunc{\realnumbers^n}{K}\).
Note that many of the series that are pointwise convergent
with \(\ucfunc{\rn}{K}\) converge in the norm in
\(\vanishingfunc{\rn}{K}\). We are also able to represent
many subspaces and operators related to the MRA of
\(\vanishingfunc{\rn}{K}\) as tensor products of the
corresponding one-dimensional cases.

\if\shortprep0
{
\begin{definition}
  \label{def:tensorvspace}
  Define
  \begin{displaymath}
    \mdimtensorvanvspace{n}{j} \defequalns
    \indexedcitp{k=1}{n} \onedimvanvspace{j} , \spaceafter j \in
    \integernumbers .
  \end{displaymath}
\end{definition}
}
\fi

\if\shortprep0
{
The injective tensor product respects subspaces and it follows from
\mytheorem \ref{lem:vanishing-tp-inclusion} that \(\mdimvanvspace{n}{j}\) is a
closed subspace of
\begin{displaymath}
  \indexedcitp{k=1}{n} \vanishingfunc{\realnumbers}{K}
\end{displaymath}
for each \(j \in \integernumbers\).
It follows from
\mytheorem \ref{lem:vanishing-tp-inclusion} that \(\mdimtensorvanvspace{n}{j}\) is a
closed subspace of \(\vanishingfunc{\realnumbers^n}{K}\), too,
for all \(j \in \integernumbers\).
Since \(\onedimvanvspace{j} \closedsubspace \onedimvanvspace{j+1}\) for all
\(j \in \integernumbers\)
it follows that
\begin{math}
  \mdimtensorvanvspace{n}{j} \closedsubspace 
  \mdimtensorvanvspace{n}{j+1}
\end{math}
for all \(n \in \positiveintegers\) and
\(j \in \integernumbers\).
}
\fi

\if\shortprep0
{
\begin{lemma}
  \label{lem:mdim-vj-unconditional-convergence}
  Let \(j \in \integernumbers\) and \(f \in
  \vanishingfunc{\realnumbers^n}{K}\).
  Series
  \begin{displaymath}
    s := \sum_{\firstznvar \in \zn} f \left(
      \frac{\firstznvar}{2^j}
    \right) \mdimsf{n}{j}{\firstznvar}
  \end{displaymath}
  converges unconditionally in
  \(\mdimvanvspace{n}{j}\)
  and
  \begin{math}
    \norminfty{s} \leq
    \ncover{\mdimmothersf{n}} \norminfty{\mdimmothersf{n}}
    \norminfty{f}
  \end{math}.
\end{lemma}
}
\fi

\if\shortprep0
{
\begin{proof}
  This is a consequence of \mylemma \ref{lem:mdim-vanishing-samples}
  and \mylemma \ref{lem:grid-convergence-multidim}.
\end{proof}
}
\fi

\if\shortprep0
{
\begin{theorem}
  \label{th:mdim-vj-unconditional-basis}
  Let \(j \in \integernumbers\).
  The set
  \begin{math}
    \{ \mdimsf{n}{j}{\firstznvar} \setsep \firstznvar \in
    \integernumbers^n \}
  \end{math}
  is an unconditional basis of \(\mdimvanvspace{n}{j}\)
  and
  \begin{math}
    \mdimvanvspace{n}{j} = \mdimtensorvanvspace{n}{j}
  \end{math}.
\end{theorem}
}
\fi

\if\shortprep1
{
\if10
{
\begin{proof}
  Use \mylemmas
  \ref{lem:grid-convergence-multidim},
  \ref{lem:mdim-vanishing-samples},
  and \ref{lem:tp-schauder-basis}.
\end{proof}
}
\fi
}
\else
{
\begin{proof}
  Space \(\mdimtensorvanvspace{n}{j}\) is a closed subspace of
  \(\vfn{n}{K}\). Hence by \mylemma
  \ref{lem:mdim-vj-unconditional-convergence}
  we have
  \begin{math}
    \mdimvanvspace{n}{j} \subset \mdimtensorvanvspace{n}{j}
  \end{math}.
  By \mylemma \ref{lem:tp-schauder-basis}
  the sequence
  \begin{math}
    (\mdimsf{n}{j}{\gensqord{n}{k}})_{k=0}^\infty
  \end{math}
  is a Schauder
  basis of \(\mdimtensorvanvspace{n}{j}\).
  
  Let \(f \in \mdimtensorvanvspace{n}{j}\). Now
  \begin{equation}
    \label{eq:mdim-series-a}
    f = \sum_{k=0}^\infty b_k \mdimsf{n}{j}{\gensqord{n}{k}}
    = \sum_{k=0}^\infty f \left( \frac{\gensqord{n}{k}}{2^j} \right)
    \mdimsf{n}{j}{\gensqord{n}{k}}
  \end{equation}
  We have
  \(\mdimtensorvanvspace{n}{j} \subset 
  \vanishingfunc{\realnumbers^n}{K}\).
  Hence by \mylemma
  \ref{lem:mdim-vj-unconditional-convergence} the
  series in \myequation \eqref{eq:mdim-series-a} converge
  unconditionally.
  Consequently the set
  \begin{math}
    \{ \mdimsf{n}{j}{\firstznvar} \setsep \firstznvar \in
    \integernumbers^n \}
  \end{math}
  is an unconditional basis of \(\mdimtensorvanvspace{n}{j}\).
  As \(f \in \mdimtensorvanvspace{n}{j}\) was arbitrary it follows from \mylemma
  \ref{lem:mdim-vanishing-samples} that
  \begin{math}
    \mdimtensorvanvspace{n}{j} \subset \mdimvanvspace{n}{j}
  \end{math}
  and hence
  \begin{math}
    \mdimvanvspace{n}{j} = \mdimtensorvanvspace{n}{j}
  \end{math}.
\end{proof}
}
\fi

\if\shortprep1
{
  When \(j \in \integernumbers\)
  the set
  \begin{math}
    \{ \mdimsf{n}{j}{\firstznvar} \setsep \firstznvar \in
    \integernumbers^n \}
  \end{math}
  is an unconditional basis of \(\mdimvanvspace{n}{j}\)
  and
  \begin{displaymath}
    \mdimvanvspace{n}{j} \equalns
    \indexedcitp{k=1}{n} \onedimvanvspace{j} .
  \end{displaymath}
Function
}
\else
{
By \mylemma \ref{lem:cospace-top-isomorphism}
function
}
\fi
\begin{math}
  \vanvjtopisomfunc{n}{j} : \gencospace{\zn}{K} \to \mdimvanvspace{n}{j}
\end{math}
defined by
\begin{displaymath}
  \vanvjtopisom{n}{j}{\seqstyle{a}}
  :=
  \sum_{\firstznvar \in \integernumbers} \seqelem{\seqstyle{a}}{\firstznvar}
  \mdimsf{n}{j}{\firstznvar}
\end{displaymath}
for all
\begin{math}
  \seqstyle{a} \in \gencospace{\zn}{K}
\end{math}
\if\shortprep1
{
is a topological isomorphism from
\begin{math}
  \gencospace{\zn}{K}
\end{math}
onto
\begin{math}
  \mdimvanvspace{n}{j}
\end{math}.
}
\else
{
is a topological isomorphism from
\begin{math}
  \gencospace{\zn}{K}
\end{math}
onto
\begin{math}
  \mdimvanvspace{n}{j}
\end{math}
and
\begin{equation}
  \label{eq:tensor-sf-expansion-bound}
  \norminfty{\seqstyle{a}}
  \leq \norminfty{\vanvjtopisom{n}{j}{\seqstyle{a}}}
  \leq \szncover{\mdimmothersf{n}} \norminfty{\mdimmothersf{n}}
  \norminfty{\seqstyle{a}}
\end{equation}
for all
\begin{math}
  \seqstyle{a} \in \gencospace{\zn}{K}
\end{math}
and \(j \in \integernumbers\).
}
\fi
\if\shortprep1
{
As a consequence of \mytheorem \ref{th:intersection}
we get
\begin{displaymath}
  \bigcap_{j \in \integernumbers} \mdimvanvspace{n}{j} = \left\{ 0 \right\} .
\end{displaymath}
}
\else
{
As a consequence of \mytheorem \ref{th:intersection}
we get the following theorem.
}
\fi

\if\shortprep0
{
\begin{theorem}
  \label{th:intersection-mdim}
  We have
  \begin{displaymath}
    \bigcap_{j \in \integernumbers} \mdimvanvspace{n}{j} = \left\{ 0 \right\} .
  \end{displaymath}
\end{theorem}
}
\fi

\begin{definition}
  \label{def:tensorvj-proj-op}
  When \(j \in \integernumbers\) define operator
  \( \mdimvanprojop{n}{j} : \vanishingfunc{\realnumbers^n}{K} \to
  \mdimvanvspace{n}{j}
  \) by
  \begin{displaymath}
    \mdimvanprojop{n}{j} f
    = \sum_{\firstznvar \in \integernumbers^n}
    \dualappl{\mdimdualsf{n}{j}{\firstznvar}}{f}
    \mdimsf{n}{j}{\firstznvar}
  \end{displaymath}
  for all \(f \in \vanishingfunc{\realnumbers^n}{K}\).
\end{definition}

\if\shortprep0
{
It follows from
\if\shortprep1
  \mylemmas \ref{lem:mdim-vanishing-samples} and
  \ref{lem:grid-convergence-multidim}
\else
\mylemma \ref{lem:mdim-vj-unconditional-convergence}
\fi
\(\mdimvanprojop{n}{j}\), \(j \in \integernumbers\) are
operators and uniformly bounded by
\begin{equation}
  \label{eq:mdimprojop-uniform-bound}
  \norm{\mdimvanprojop{n}{j}} \leq
  \szncover{\mdimmothersf{n}} \norminfty{\mdimmothersf{n}}
\end{equation}
for all
\begin{math}
  j \in \integernumbers
\end{math}.
}
\fi
Operator
\begin{math}
  \mdimvanprojop{n}{j}
\end{math}
is a projection onto
\begin{math}
  \mdimvanvspace{n}{j}
\end{math}
for each
\begin{math}
  j \in \integernumbers
\end{math}.

\begin{lemma}
  \label{lem:proj-op-composition-mdim}
  We have \(\mdimvanprojop{n}{j} = \mdimvanprojop{n}{j} \circ \mdimvanprojop{n}{j'}\)
  for all \(j, j' \in \integernumbers\), \(j' \geq j\).
\end{lemma}

\if\shortprep0
{
\begin{proof}
  The proof is similar to the proof of \mylemma
  \ref{lem:proj-op-composition-one-dim}.
\end{proof}
}
\fi

\if\shortprep1
{
  There exist
  \(c_1 \in \positiverealnumbers\)
  and
  \(c_2 \in \positiverealnumbers\)
  so that
  \begin{math}
    \norminfty{f - \mdimvanprojop{n}{j} f}
    \leq
    c_1 \modcont{f}{2^{-j} c_2}
  \end{math}
  for all
  \(f \in \vanishingfunccv{\rn}\)
  and
  \(j \in \integernumbers\).
  Assume that \(g \in \vanishingfunccv{\rn}\). Now
  \begin{displaymath}
    \lim_{j \to \infty} \norminfty{g - \mdimvanprojop{n}{j} g} = 0 .
  \end{displaymath}
  The aforementioned two formulas are direct consequences
  of the corresponding results for
  \(\ucfunccv{\rn}\).
  We give
  now a general result on the tensor products of
  the function space \(\vanishingfunc{\realnumbers}{K}\) with itself.
}
\fi

\if01
{
We get the following theorem as a direct
consequence of \mytheorem \ref{th:upj-ineq}.
}
\fi

\if01
{
\begin{theorem}
  \label{th:vanpj-ineq}
  There exist
  \(c_1 \in \positiverealnumbers\)
  and
  \(c_2 \in \positiverealnumbers\)
  so that
  \begin{displaymath}
    \norminfty{f - \mdimvanprojop{n}{j} f}
    \leq
    c_1 \modcont{f}{2^{-j} c_2}
  \end{displaymath}
  for all
  \(f \in \vanishingfunccv{\rn}\)
  and
  \(j \in \integernumbers\).
\end{theorem}
}
\fi

\if\shortprep0
{
We get the following corollary as a direct
consequence of \mytheorem \ref{th:upj-convergence}.
}
\fi

\if\shortprep0
{
\begin{corollary}
  \label{cor:mdim-vanpj-convergence}
  Let \(f \in \vanishingfunccv{\rn}\). Now
  \begin{displaymath}
    \lim_{j \to \infty} \norminfty{f - \mdimvanprojop{n}{j} f} = 0 .
  \end{displaymath}
\end{corollary}
}
\fi

\if\shortprep0
{
We get the following corollary as a direct consequence
of \mylemma \ref{th:omega-pj-inequality}.
}
\fi

\if\shortprep0
{
\begin{corollary}
  \label{cor:van-omega-inequality}
  Let \(t \in \positiverealnumbers\).
  There exists a constant
  \begin{math}
    c_1 \in \positiverealnumbers
  \end{math},
  which may depend on \(t\),
  so that
  \begin{displaymath}
    \forall f \in \vanishingfunccv{\rn} : 
    \modcont{\mdimvanprojop{n}{j} f}{2^{-j} t}
    \leq c_1 \modcont{f}{2^{-j} \sqrt{n}} .
  \end{displaymath}
\end{corollary}
}
\fi

\if\shortprep0
We prove
now a general result on the tensor products of
the function space \(\vanishingfunc{\realnumbers}{K}\) with itself.
\fi

\begin{theorem}
  \label{th:vanishing-tp-space}
  We have
  \begin{displaymath}
    \indexedcitp{j=1}{n} \vanishingfunc{\realnumbers}{K} \equalns
    \vanishingfunc{\realnumbers^n}{K} .
  \end{displaymath}
\end{theorem}

\if\shortprep0
{
\begin{proof}
  Define
  \begin{displaymath}
    F^{[n]} \defequalns \indexedcitp{j=1}{n} \vanishingfunc{\realnumbers}{K} .
  \end{displaymath}
  By \mylemma \ref{lem:vanishing-tp-inclusion} we have
  \begin{math}
    F^{[n]} \subset \vanishingfunc{\realnumbers^n}{K}
  \end{math}.
  Let \(\myphi\) be the Deslauriers-Dubuc fundamental function of some
  degree \(m\) and
  \(\mdimmothersf{n}\) the \(n\)-dimensional tensor product mother
  scaling function generated by \(\myphi\).
  Let spaces \(\onedimvanvspace{j}\), \(j \in \integernumbers\), belong to the
  interpolating multiresolution analysis generated by \(\myphi\) and
  spaces \(\mdimvanvspace{n}{j}\), \(j \in \integernumbers\), be the
  corresponding tensor product spaces.
  Let \(f \in \vanishingfunc{\realnumbers^n}{K}\). It follows from \mycorollary
  \ref{cor:mdim-vanpj-convergence} that \(\mdimvanprojop{n}{j} f \to f\)
  as \(j \to \infty\) (convergence in the supremum norm).
  Now \(\mdimvanprojop{n}{j} f \in \mdimvanvspace{n}{j} \subset F^{[n]}\)
  for all \(j \in \integernumbers\). Hence
  \begin{math}
    f \in \overline{F^{[n]}}
  \end{math}.
  where we hold \(F^{[n]}\) as a subspace of the Banach space
  \(\vanishingfunc{\realnumbers^n}{K}\).
  Since \(F^{[n]}\) is a Banach space \(\overline{F^{[n]}} \equalns
  F^{[n]}\).
  Hence \(F^{[n]} \equalns \vanishingfunc{\realnumbers^n}{K}\).
\end{proof}
}
\fi

For example,
\begin{math}
  \vanishingfunccv{\realnumbers} \citp \vanishingfunccv{\realnumbers}
  \equalns \vanishingfunccv{\realnumbers^2}
\end{math}.

\begin{definition}
  \label{def:tensorpartialwspace}
  Define
  \begin{displaymath}
    \mdimvanpartialwspace{n}{\zos}{j}
    \defequalns
    \indexedcitp{k=1}{n} \uspace{\cartprodelem{\zos}{k}}{j}
  \end{displaymath}
  where \(j \in \integernumbers\) and \(\zos \in \zeroonesetn\).
\end{definition}

\if\shortprep0
{
Note that
\begin{math}
  \mdimvanpartialwspace{n}{\finitezeroseq{n}}{j}
  \equalns
  \mdimvanvspace{n}{j}
\end{math}.
}
\fi

\begin{definition}
  \label{def:partial-tensor-proj-op}
  Let \(s \in \zeroonesetn\).
  Define operator
  \(\mdimvanpartialprojop{n}{\zos}{j} : \vanishingfunc{\realnumbers}{K} \to
  \mdimvanpartialwspace{n}{\zos}{j}\) by
  \begin{displaymath}
    \mdimvanpartialprojop{n}{\zos}{j} = \indexedopctp{\injtn}{k=1}{n}
    \genprojoponedim{\cartprodelem{\zos}{k}}{j} .
  \end{displaymath}
\end{definition}

\if\shortprep1
{
  Operator
  \(\mdimvanpartialprojop{n}{\zos}{j}\) is a continuous projection onto
  space \(\mdimvanpartialwspace{n}{\zos}{j}\) for each
  \(n \in \positiveintegers\), \(j \in \integernumbers\), and
  \(\zos \in \zeroonesetn\).
  We have
  \begin{displaymath}
    \mdimvanpartialprojop{n}{\zos}{j} f
    = \sum_{\firstznvar \in \integernumbers^n}
    \dualappl{\mdimgendualwavelet{n}{\zos}{j}{\firstznvar}}{f}
    \mdimgenwavelet{n}{\zos}{j}{\firstznvar}
  \end{displaymath}
  for all \(f \in \vanishingfunc{\realnumbers^n}{K}\),
  \(\zos \in \zeroonesetn\), and \(j \in \integernumbers\).
  It follows that
  \begin{equation}
    \label{eq:mdimvanprojop-tp}
    \mdimvanprojop{n}{j} = 
    \mdimvanpartialprojop{n}{\finitezeroseq{n}}{j} =
    \indexedopctp{\injtn}{k=1}{n} \onedimvanprojop{j} .
  \end{equation}
  We also have
  \begin{math}
    \mdimvanpartialprojop{n}{\zos}{j}
    =
    \mdimvanpartialprojop{n}{\zos}{j}
    \circ
    \mdimvanprojop{n}{j'}  
  \end{math}
  for all
  \begin{math}
    j, j' \in \integernumbers
  \end{math},
  \begin{math}
    j' > j
  \end{math}
  and
  \begin{math}
    \zos \in \zeroonesetn
  \end{math}.
}
\fi

\if\shortprep0
By \mylemma \ref{lem:tp-projection} operator
\(\mdimvanpartialprojop{n}{\zos}{j}\) is a continuous projection onto
space \(\mdimvanpartialwspace{n}{\zos}{j}\) for each
\(n \in \positiveintegers\), \(j \in \integernumbers\), and
\(\zos \in \zeroonesetn\).

\begin{lemma}
  \label{lem:deltaop-formula-mdim}
  We have
  \begin{displaymath}
    \mdimvanpartialprojop{n}{\zos}{j} f
    = \sum_{\firstznvar \in \integernumbers^n}
    \dualappl{\mdimgendualwavelet{n}{\zos}{j}{\firstznvar}}{f}
    \mdimgenwavelet{n}{\zos}{j}{\firstznvar}
  \end{displaymath}
  for all \(f \in \vanishingfunc{\realnumbers^n}{K}\),
  \(\zos \in \zeroonesetn\), and \(j \in \integernumbers\).
  Here the series converges unconditionally and
  \begin{math}
    (\dualappl{\mdimgendualwavelet{n}{\zos}{j}{\firstznvar}}
    {f})_{\firstznvar \in \zn} \in \gencospace{\zn}{K}
  \end{math}.
\end{lemma}
\fi

\if\shortprep0
{
\begin{proof}
  Define operators
  \begin{math}
    \toperator{n}{\znstyleprime{s}}{j'} : \vfn{n}{K} \to
    \mdimvanpartialwspace{n}{\znstyleprime{s}}{j'}
  \end{math},
  \begin{math}
    j' \in \integernumbers
  \end{math},
  \begin{math}
    \znstyleprime{s} \in \zeroonesetn
  \end{math},
  by
  \begin{equation}
    \label{eq:toperator}
    \toperator{n}{\znstyleprime{s}}{j'} g
    := \sum_{\firstznvar \in \integernumbers^n}
    \szdualappl{\mdimgendualwavelet{n}{\znstyleprime{s}}{j'}
    {\firstznvar}}{g}
    \mdimgenwavelet{n}{\znstyleprime{s}}{j'}{\firstznvar}
  \end{equation}
  for all
  \begin{math}
    g \in \vfn{n}{K}
  \end{math}.
  Functional \(\mdimmotherdualwavelet{n}{\znstyleprime{s}}\) is a 
  finite linear combination of functionals
  \(\delta(2 \cdot - \firstznvar)\),
  \(\firstznvar \in \zn\). Hence by \mylemma 
  \ref{lem:mdim-vanishing-samples}
  we have
  \begin{math}
    (\dualappl{\mdimgendualwavelet{n}{\znstyleprime{s}}
      {j'}{\firstznvar}}{g})_{\firstznvar \in
      \zn} \in \gencospace{\zn}{K}
  \end{math}
  and there exists \(c \in \positiverealnumbers\) so that
  \begin{math}
    \norminfty{(\dualappl{\mdimgendualwavelet{n}{\znstyleprime{s}}
      {j'}{\firstznvar}}{g})_{\firstznvar \in
      \zn}} \leq c \norminfty{g}
  \end{math}
  for all \(g \in \vfn{n}{K}\).
  By \mylemma \ref{lem:grid-convergence-multidim} the series in \myequation
  \eqref{eq:toperator} converges unconditionally in \(\vfn{n}{K}\) and
  there exists \(c' \in \positiverealnumbers\) so that
  \begin{math}
    \norminfty{\toperator{n}{\znstyleprime{s}}{j'} g} \leq c' 
    \norminfty{g}
  \end{math}
  for all \(g \in \vfn{n}{K}\).

  Let \(g \in \mdimvanvspace{n}{j'}\) for some \(j' \in
  \integernumbers\). Now
  \begin{math}
    g = \sum_{\firstznvar \in \zn} \seqelem{\sqa}{\firstznvar} 
    \mdimsf{n}{j'}{\firstznvar}
  \end{math}
  where 
  \begin{math}
    \seqelem{\sqa}{\firstznvar} \in \gencospace{\zn}{K}
  \end{math}
  and
  \begin{displaymath}
    \mdimvanpartialprojop{n}{\zos}{j} g =
    \sum_{\firstznvar \in \zn}
    \seqelem{\sqa}{\firstznvar}
    \mdimvanpartialprojop{n}{\zos}{j} \mdimsf{n}{j'}{\firstznvar}
    = \sum_{\firstznvar \in \zn}
    \seqelem{\sqa}{\firstznvar}
    \indexedtensorproduct{k=1}{n}
    \genprojoponedim{\cartprodelem{\zos}{k}}{j}
    \myphi_{j',\cartprodelem{\firstznvar}{k}} .
  \end{displaymath}
  \if\longversion0
  {
  Furthermore,
  \begin{align*}
    \indexedtensorproduct{k=1}{n}
    \genprojoponedim{\cartprodelem{\zos}{k}}{j}
    \myphi_{j',\cartprodelem{\firstznvar}{k}}
    & = \indexedtensorproduct{k=1}{n}
    \sum_{l \in \integernumbers}
    \szdualappl{\gendualwaveletonedim{\cartprodelem{\zos}{k}}
    {j}{l}}
      {\myphi_{j',\cartprodelem{\firstznvar}{k}}}
    \genwaveletonedim{\cartprodelem{\zos}{k}}{j}{l} \\
    & = \sum_{\secondznvar \in \zn} \indexedtensorproduct{k=1}{n}
    \szdualappl{\gendualwaveletonedim{\cartprodelem{\zos}{k}}{j}
      {\cartprodelem{\secondznvar}{k}}}
    {\myphi_{j',\cartprodelem{\firstznvar}{k}}}
    \genwaveletonedim{\cartprodelem{\zos}{k}}{j}
      {\cartprodelem{\secondznvar}{k}} \\
    & = \sum_{\secondznvar \in \zn}
    \szdualappl
    {\indexedtensorproduct{k=1}{n}
      \gendualwaveletonedim{\cartprodelem{\zos}{k}}{j}
      {\cartprodelem{\secondznvar}{k}}}
    {\indexedtensorproduct{k=1}{n} 
    \myphi_{j',\cartprodelem{\firstznvar}{k}}}
    \indexedtensorproduct{k=1}{n}
    \genwaveletonedim{\cartprodelem{\zos}{k}}{j}
    {\cartprodelem{\secondznvar}{k}} \\
    & =
    \sum_{\secondznvar \in \zn}
    \szdualappl{\mdimgendualwavelet{n}{\zos}{j}{\secondznvar}}
    {\mdimsf{n}{j'}{\firstznvar}}
    \mdimgenwavelet{n}{\zos}{j}{\secondznvar}
    = \toperator{n}{\zos}{j} \mdimsf{n}{j'}{\firstznvar}
  \end{align*}
  where the second equality follows from \mylemma \ref{lem:tensor-product-of-sums}.
  }
  \else
  {
  We also have
  \begin{displaymath}
    \genprojoponedim{t}{j} \myphi_{j',m}
    = \sum_{l \in \integernumbers}
    \szdualappl{\gendualwaveletonedim{t}{j}{l}}{\myphi_{j',m}}
    \genwaveletonedim{t}{j}{l}
  \end{displaymath}
  for each \(t \in \zerooneset\) and \(m \in \integernumbers\) and
  \begin{align*}
    \indexedtensorproduct{k=1}{n}
    \genprojoponedim{\cartprodelem{\zos}{k}}{j}
    \myphi_{j',\cartprodelem{\firstznvar}{k}}
    & = \indexedtensorproduct{k=1}{n}
    \sum_{l \in \integernumbers}
    \szdualappl{\gendualwaveletonedim{\cartprodelem{\zos}{k}}{j}{l}}
      {\myphi_{j',\cartprodelem{\firstznvar}{k}}}
    \genwaveletonedim{\cartprodelem{\zos}{k}}{j}{l} \\
    & = \sum_{\secondznvar \in \zn} \indexedtensorproduct{k=1}{n}
    \szdualappl{\gendualwaveletonedim{\cartprodelem{\zos}{k}}{j}
      {\cartprodelem{\secondznvar}{k}}}
    {\myphi_{j',\cartprodelem{\firstznvar}{k}}}
    \genwaveletonedim{\cartprodelem{\zos}{k}}{j}
      {\cartprodelem{\secondznvar}{k}} \\
    & = \sum_{\secondznvar \in \zn}
    \szdualappl
    {\indexedtensorproduct{k=1}{n}
      \gendualwaveletonedim{\cartprodelem{\zos}{k}}{j}
      {\cartprodelem{\secondznvar}{k}}}
    {\indexedtensorproduct{k=1}{n}
    \myphi_{j',\cartprodelem{\firstznvar}{k}}}
    \indexedtensorproduct{k=1}{n}
    \genwaveletonedim{\cartprodelem{\zos}{k}}{j}
    {\cartprodelem{\secondznvar}{k}} \\
    & =
    \sum_{\secondznvar \in \zn}
    \szdualappl{\mdimgendualwavelet{n}{\zos}{j}{\secondznvar}}
    {\mdimsf{n}{j'}{\firstznvar}}
    \mdimgenwavelet{n}{\zos}{j}{\secondznvar}
    = \toperator{n}{\zos}{j} \mdimsf{n}{j'}{\firstznvar}
  \end{align*}
  where the second equality follows from \mylemma \ref{lem:tensor-product-of-sums}.
  }
  \fi
  Thus
  \begin{displaymath}
    \mdimvanpartialprojop{n}{\zos}{j} g
    = \sum_{\firstznvar \in \zn} \seqelem{\sqa}{\firstznvar}
    \toperator{n}{\zos}{j} \mdimsf{n}{j'}{\firstznvar}
    = \toperator{n}{\zos}{j} \left( \sum_{\firstznvar \in \zn} 
    \seqelem{\sqa}{\firstznvar}
      \mdimsf{n}{j'}{\firstznvar} \right)
    = \toperator{n}{\zos}{j} g .
  \end{displaymath}

  Let \(f \in \vfn{n}{K}\). Now
  \begin{math}
    \mdimvanprojop{n}{j'} f \to f
  \end{math}
  as
  \begin{math}
     j' \to \infty
  \end{math}
  and
  \begin{math}
    \mdimvanpartialprojop{n}{\zos}{j} (\mdimvanprojop{n}{j'} f)
    = \toperator{n}{\zos}{j} (\mdimvanprojop{n}{j'} f)
  \end{math}
  for all
  \begin{math}
    j' \in \integernumbers    
  \end{math}.
  By continuity of operators
  \begin{math}
    \mdimvanpartialprojop{n}{\zos}{j}
  \end{math}
  and
  \begin{math}
    \toperator{n}{\zos}{j}
  \end{math}
  we have
  \begin{math}
    \mdimvanpartialprojop{n}{\zos}{j} f
    = \toperator{n}{\zos}{j} f
  \end{math}.
\end{proof}
}
\fi

\if\shortprep0
{
It follows that
\begin{displaymath}
  \mdimvanprojop{n}{j} = 
  \mdimvanpartialprojop{n}{\finitezeroseq{n}}{j} =
  \indexedopctp{\injtn}{k=1}{n} \onedimvanprojop{j} .
\end{displaymath}
When \(\zos \in \zeroonesetn\)
it follows from \mylemmas
\ref{lem:grid-convergence-multidim} and
\ref{lem:deltaop-formula-mdim} that 
operators \(\mdimvanpartialprojop{n}{\zos}{j}\),
\(j \in \integernumbers\),
are uniformly bounded by
\begin{displaymath}
  \label{eq:mdimpartialprojop-uniform-bound-ineq}
  \norm{\mdimvanpartialprojop{n}{\zos}{j}}
  \leq
  \szncover{\mdimmotherwavelet{n}{\zos}}
  \norminfty{\mdimmotherwavelet{n}{\zos}}
  \norm{\mdimmotherdualwavelet{n}{\zos}}
\end{displaymath}
for all
\begin{math}
  j \in \integernumbers
\end{math}.
We also have
\begin{math}
  \mdimvanpartialprojop{n}{\zos}{j}
  =
  \mdimvanpartialprojop{n}{\zos}{j}
  \circ
  \mdimvanprojop{n}{j'}  
\end{math}
for all
\begin{math}
  j, j' \in \integernumbers
\end{math},
\begin{math}
  j' > j
\end{math}
and
\begin{math}
  \zos \in \zeroonesetn
\end{math}.
}
\fi

\begin{definition}
  \label{def:tensorwspace}
  Define
  \begin{displaymath}
    \mdimvanwspace{n}{j}
    \defequalns
    \directsumoneindex{\zos \in \zeroonesetnnozero{n}}
      \mdimvanpartialwspace{n}{\zos}{j}
  \end{displaymath}
  for all \(n \in \positiveintegers\) and \(j \in \integernumbers\).
\end{definition}

\if\shortprep0
{
By \mylemma \ref{lem:dsum-closed} we have
\begin{math}
  \mdimvanwspace{n}{j} \closedsubspace \vfn{n}{K}
\end{math}.
}
\fi

\begin{definition}
  \label{def:tensor-delta-proj-op}
  Define operator
  \(\mdimvandeltaop{n}{j} : \vanishingfunc{\realnumbers}{K} \to
  \mdimvanwspace{n}{j}\) by
  \begin{displaymath}
    \mdimvandeltaop{n}{j} = \sum_{\zos \in \zeroonesetnnozero{n}}
    \mdimvanpartialprojop{n}{\zos}{j} .
  \end{displaymath}
\end{definition}

\if\shortprep0
{
\begin{theorem}
  \label{th:mdim-wj-series-rep}
  Let \(j \in \integernumbers\) and
  \(\zos \in \zeroonesetn\).
  Then
  \begin{displaymath}
    \mdimvanpartialwspace{n}{\zos}{j}
    = \left\{ \sum_{\firstznvar \in \zn}
      \seqelem{\seqstyle{a}}{\firstznvar}
      \mdimgenwavelet{n}{\zos}{j}{\firstznvar}
      \setsep
      \seqstyle{a} \in \gencospace{\zn}{K}
      \right\}
  \end{displaymath}
  where the series converges unconditionally.
\end{theorem}
}
\fi

\if\shortprep0
{
\begin{proof}
  Now
  \begin{math}
    \mdimgenwavelet{n}{\zos}{j}{\firstznvar}
    \in \mdimvanpartialwspace{n}{\zos}{j}
  \end{math}
  for each \(\firstznvar \in \zn\).
  Suppose that
  \((a_\firstznvar)_{\firstznvar \in \zn} \in \gencospace{\zn}{K}\).
  By \mylemma \ref{lem:grid-convergence-multidim} the series
  \begin{math}
    t := \sum_{\firstznvar \in \zn} 
    \mdimgenwavelet{n}{\zos}{j}{\firstznvar}
  \end{math}
  converges unconditionally in \(\vfn{n}{K}\). As
  \begin{math}
    \mdimvanpartialwspace{n}{\zos}{j} \closedsubspace \vfn{n}{K}
  \end{math}
  it follows that \(t \in \mdimvanpartialwspace{n}{\zos}{j}\).

  Let \(f \in \mdimvanpartialwspace{n}{\zos}{j}\). Now
  \begin{displaymath}
    f = \mdimvanpartialprojop{n}{\zos}{j} f
    = \sum_{\firstznvar \in \zn}
    \szdualappl{\mdimgendualwavelet{n}{\zos}{j}{\firstznvar}}{f}
      \mdimgenwavelet{n}{\zos}{j}{\firstznvar}
  \end{displaymath}
  and
  \begin{math}
    (
      \dualappl{\mdimgendualwavelet{n}{\zos}{j}{\firstznvar}}{f}
    )_{\firstznvar \in \zn}
    \in \gencospace{\zn}{K}
  \end{math}.
\end{proof}
}
\fi

\if\shortprep1
{
  When \(j \in \integernumbers\) and
  \(\zos \in \zeroonesetn\)
  we have
  \begin{displaymath}
    \mdimvanpartialwspace{n}{\zos}{j}
    = \left\{ \sum_{\firstznvar \in \zn}
      \seqelem{\seqstyle{a}}{\firstznvar}
      \mdimgenwavelet{n}{\zos}{j}{\firstznvar}
      \setsep
      \seqstyle{a} \in \gencospace{\zn}{K}
      \right\} .
  \end{displaymath}
  Function
}
\else
{
  By \mylemma \ref{lem:cospace-top-isomorphism}
  function
}
\fi
\begin{math}
  \vanwsjtopisomfunc{n}{\zos}{j} : \gencospace{\zn}{K} \to 
  \mdimvanpartialwspace{n}{\zos}{j}
\end{math}
defined by
\begin{displaymath}
  \vanwsjtopisom{n}{\zos}{j}{\seqstyle{a}}
  :=
  \sum_{\firstznvar \in \integernumbers} \seqelem{\seqstyle{a}}{\firstznvar}
  \mdimgenwavelet{n}{\zos}{j}{\firstznvar}
\end{displaymath}
for all
\begin{math}
  \seqstyle{a} \in \gencospace{\zn}{K}
\end{math}
\if\shortprep1
{
is a topological isomorphism from
\begin{math}
  \gencospace{\zn}{K}
\end{math}
onto
\begin{math}
  \mdimvanpartialwspace{n}{\zos}{j}
\end{math}.
}
\else
{
is a topological isomorphism from
\begin{math}
  \gencospace{\zn}{K}
\end{math}
onto
\begin{math}
  \mdimvanpartialwspace{n}{\zos}{j}
\end{math}
and
\begin{equation}
  \label{eq:tensorpartialwspace-norm-equiv}
  \norminfty{\seqstyle{a}}
  \leq \norminfty{\vanwsjtopisom{n}{\zos}{j}{\seqstyle{a}}}
  \leq \szncover{\mdimmotherwavelet{n}{\zos}}
  \norminfty{\mdimmotherwavelet{n}{\zos}}
  \norminfty{\seqstyle{a}}
\end{equation}
for all
\begin{math}
  \seqstyle{a} \in \gencospace{\zn}{K}
\end{math},
\begin{math}
  \zos \in \zeroonesetn
\end{math}, and
\begin{math}
  j \in \integernumbers
\end{math}.
}
\fi

\begin{lemma}
  \label{lem:mdim-delta-formula}
  \begin{math}
    \mdimvandeltaop{n}{j} = \mdimvanprojop{n}{j+1} - 
    \mdimvanprojop{n}{j}
  \end{math}.
\end{lemma}

\if\shortprep1
{
\begin{proof}
  Use \mydef \ref{def:delta-proj-op-one-dim}
  and \myequation \eqref{eq:mdimvanprojop-tp}.
\end{proof}
}
\else
{
\begin{proof}
  We have
  \begin{math}
    \mdimvanprojop{n}{j} = 
    \mdimvanpartialprojop{n}{\finitezeroseq{n}}{j}
  \end{math}
  and
  \begin{align*}
    \mdimvanprojop{n}{j+1} & = \indexedopctp{\injtn}{j=1}{n}
    \onedimvanprojop{j+1}
    = \indexedopctp{\injtn}{j=1}{n} (\onedimvanprojop{j} + \onedimvandeltaop{j})
    = \sum_{\zos \in \zeroonesetn} \indexedopctp{\injtn}{j=1}{n}
    \genprojoponedim{\cartprodelem{\zos}{j}}{j} \\
    & = \sum_{\zos \in \zeroonesetn} \mdimvanpartialprojop{n}{\zos}{j}
    .
  \end{align*}
  Consequently the lemma is true.
\end{proof}
}
\fi

\if\shortprep0
{
It follows that the operators \(\mdimvandeltaop{n}{j}\),
\(j \in \integernumbers\), are uniformly bounded by
\begin{equation}
  \label{eq:mdimdeltaprojop-uniform-bound-ineq}
  \norm{\mdimvandeltaop{n}{j}} \leq
  2 \szncover{\mdimmothersf{n}} \norminfty{\mdimmothersf{n}}
\end{equation}
for all
\begin{math}
  j \in \integernumbers
\end{math}.
}
\fi

\begin{lemma}
  \label{lem:mdimpartialprojop-orthogonality}
  Let \(j \in \integernumbers\). Then
  \begin{itemize}
  \item[(i)]
    \begin{math}
      \forall \zos \in \zeroonesetn, \zot \in \zeroonesetn :
      \zos \neq \zot \implies
      \forall f \in 
      \mdimvanpartialwspace{n}{\zos}{j} :
      \mdimvanpartialprojop{n}{\zot}{j} f = 0
    \end{math}
  \item[(ii)]
    \begin{math}
      \forall j' \in \integernumbers,
      f \in \mdimvanvspace{n}{j'} :
      j' \leq j
      \implies \mdimvandeltaop{n}{j} f = 0
    \end{math}
  \item[(iii)]
    \begin{math}
      \forall j' \in \integernumbers,
      f \in \mdimvanwspace{n}{j'} :
      j' \neq j
      \implies \mdimvandeltaop{n}{j} f = 0
    \end{math}
  \item[(iv)]
    \begin{math}
      \forall j' \in \integernumbers,
      f \in \mdimvanwspace{n}{j'} :
      j' \geq j
      \implies \mdimvanprojop{n}{j} f = 0
    \end{math}.
  \end{itemize}
\end{lemma}

\if\shortprep0
{
\begin{proof}
  \mbox{ }
  \begin{itemize}
  \item[(i)]
    Let
    \begin{math}
      \zos, \zot \in \zeroonesetn
    \end{math},
    \begin{math}
      \zos \neq \zot
    \end{math},
    and
    \begin{math}
      f \in \mdimvanpartialwspace{n}{\zos}{j}
    \end{math}.
    Now
    \begin{displaymath}
      f = \sum_{\firstznvar \in \zn} \seqelem{\sqa}{\firstznvar}
      \mdimgenwavelet{n}{\zos}{j}{\firstznvar}
    \end{displaymath}
    where
    \begin{math}
      \sqa \in \gencospace{\zn}{K}
    \end{math}
    and
    \begin{eqnarray*}
      \mdimvanpartialprojop{n}{\zot}{j} f & = &
      \sum_{\secondznvar \in \zn}
      \szdualappl{\mdimgendualwavelet{n}
      {\zot}{j}{\secondznvar}}{f}
      \mdimgenwavelet{n}{\zot}{j}{\secondznvar}
      =
      \sum_{\secondznvar \in \zn} \sum_{\firstznvar \in \zn}
      \seqelem{\sqa}{\firstznvar}
      \szdualappl{\mdimgendualwavelet{n}
      {\zot}{j} 
      {\secondznvar}}
      {\mdimgenwavelet{n}{\zos}{j}{\firstznvar}}
      \mdimgenwavelet{n}{\zot}{j}{\secondznvar} \\
      & = &
      \sum_{\secondznvar \in \zn} \sum_{\firstznvar \in \zn}
      \seqelem{\sqa}{\firstznvar}
      \delta_{\zot,\zos} 
      \delta_{\secondznvar,\firstznvar}
      \mdimgenwavelet{n}{\zot}{j}{\secondznvar}
      = 0
    \end{eqnarray*}
    where the third equality follows from \myequation
    \eqref{eq:mdim-wavelet-orthogonality}.
  \item[(ii)-(iv)]
    The proofs are similar to the proofs in \mylemma
    \ref{lem:mdimpartialuprojop-orthogonality}.
  \end{itemize}
\end{proof}
}
\fi

\if\shortprep1
{
We have
}
\else
{
As
\begin{math}
  \onedimvanvspace{j+1} \equalns \onedimvanvspace{j} \dsum \onedimvanwspace{j}
\end{math}
  it follows from Corollary \ref{cor:distr-law-ucn},
  \mylemma \ref{lem:mdimpartialprojop-orthogonality},
  and \myequations \eqref{eq:one-dim-projop-orth-a},
\eqref{eq:one-dim-projop-orth-b}, and \eqref{eq:one-dim-projop-orth-c}
that
}
\fi
\begin{equation}
  \label{eq:tensorvspace-dsum}
  \mdimvanvspace{n}{j+1} \equalns \mdimvanvspace{n}{j} \dsum
  \mdimvanwspace{n}{j}
\end{equation}
for all \(j \in \integernumbers\).

\if\shortprep1
We have
\else
It follows from \mycorollary \ref{cor:mdim-vanpj-convergence} that
\fi
\begin{equation}
  \label{eq:function-expansion-mdim}
  f = \mdimvanprojop{n}{j_0} f + \sum_{j=j_0}^\infty 
  \mdimvandeltaop{n}{j} f
\end{equation}
for all \(f \in \vanishingfunc{\realnumbers^n}{K}\).
If
\begin{displaymath}
  f = v + \sum_{j=j_0}^\infty w_j
\end{displaymath}
where \(v \in \mdimvanvspace{n}{j_0}\) and
\(w_j \in \mdimvanwspace{n}{j}\) for all \(j \in \integernumbers\),
\(j \geq j_0\),
it follows from \mylemmas \ref{lem:mdim-delta-formula} and
\ref{lem:mdimpartialprojop-orthogonality} that
\(v = \mdimvanprojop{n}{j_0} f\) and
\(w_j = \mdimvandeltaop{n}{j} f\) for all
\(j \in \integernumbers\),
\(j \geq j_0\).

\begin{definition}
  Let \(n \in \positiveintegers\)
  and \(j \in \integernumbers\).
  Let
  \(E \equalns \ucfunc{\rn}{K}\)
  or
  \(E \equalns \vanishingfunc{\rn}{K}\).
  Define
  \begin{displaymath}
    \genvspace{E}{n}{j}
    \defequalns
    \left\{
      \begin{array}{ll}
        \mdimuvspace{n}{j} ; & E = \ucfunc{\rn}{K} \\
        \mdimvanvspace{n}{j} ; & E = \vanishingfunc{\rn}{K} ,
      \end{array}
    \right.
  \end{displaymath}
  \begin{displaymath}
    \genwspace{E}{n}{j}
    \defequalns
    \left\{
      \begin{array}{ll}
        \mdimuwspace{n}{j} ; & E \equalns \ucfunc{\rn}{K} \\
        \mdimvanwspace{n}{j} ; & E \equalns \vanishingfunc{\rn}{K} ,
      \end{array}
    \right.
  \end{displaymath}
  and
  \begin{displaymath}
    \genpartialwspace{E}{n}{\zos}{j}
    \defequalns
    \left\{
      \begin{array}{ll}
        \mdimupartialwspace{n}{\zos}{j} ; & E \equalns \ucfunc{\rn}{K} \\
        \mdimvanpartialwspace{n}{\zos}{j} ; & E \equalns \vanishingfunc{\rn}{K} .
      \end{array}
    \right.
  \end{displaymath}
\end{definition}

As in the case of \(\ucfunc{\rn}{K}\) we get
\begin{displaymath}
  \mdimvanvspace{n}{j_0}
  \dsum \indexeddirectsum{j=j_0}{\infty} \mdimvanwspace{n}{j}
  \neq
  \vanishingfunc{\rn}{K}
\end{displaymath}
when the mother scaling function is Lipschitz continuous.
We also have
\begin{displaymath}
  \vanishingfunc{\realnumbers^{n}}{K}
  \equalns
  \closop \left( \bigcup_{l=j_0}^\infty \left( 
  \mdimvanvspace{n}{j_0} \dsum
    \indexeddirectsum{j=j_0}{l} \mdimvanwspace{n}{j} \right) 
    \right)
\end{displaymath}
for all \(j_0 \in \integernumbers\).

}
\else
{
\input section-mra-vanishing-func-sv.tex
}
\fi

\section{Interpolating Dual MRA}
\label{sec:dual-mra}

We shall have
\(K = \realnumbers\) or \(K = \complexnumbers\)
and \(n \in \positiveintegers\) throughout this
section.
We shall also have
\begin{math}
  E = \ucfunc{\rn}{K}
\end{math}
or
\begin{math}
  E = \vanishingfunc{\rn}{K}
\end{math}.

\subsection{General}

\begin{definition}
  When \(j \in \integernumbers\)
  define
  \begin{eqnarray*}
    \mdimdualvspace{n}{j}
    & \defequalns &
    \left\{
      \sum_{\firstznvar \in \zn}
      \seqelem{\sqd}{\firstznvar}
      \mdimdualsf{n}{j}{\firstznvar}
    \setsep
      \sqd \in \littlelp{1}{\zn}{K}
    \right\} \\
    \norminspace{\dualelem{f}}{\mdimdualvspace{n}{j}}
    & := &
    \norminspace{\dualelem{f}}{\topdual{\ucfunc{\rn}{K}}}
    = \norminspace{\dualelem{f}}{\topdual{\vanishingfunc{\rn}{K}}} ,
    \spaceafter
    \dualelem{f} \in \mdimdualvspace{n}{j} .
  \end{eqnarray*}  
\end{definition}

\if\shortprep1
We may identify
\else
By \mylemma \ref{lem:multidim-delta-isomorphism}
we may identify
\fi
Banach space
\begin{math}
  \mdimdualvspace{n}{j}
\end{math}
for both
\begin{math}
  E = \ucfunc{\rn}{K}
\end{math}
and
\begin{math}
  E = \vanishingfunc{\rn}{K}
\end{math}
for each
\begin{math}
  j \in \integernumbers
\end{math}.
It follows also that function
\begin{math}
  \mdimdualvjiisomfunc{n}{j} : \littlelp{1}{\zn}{K} \to
  \mdimdualvspace{n}{j}
\end{math}
defined by
\begin{displaymath}
  \mdimdualvjiisom{n}{j}{\seqstyle{d}}
  :=
  \sum_{\firstznvar \in \zn} \seqelem{\seqstyle{d}}{\firstznvar}
  \mdimdualsf{n}{j}{\firstznvar}
\end{displaymath}
for all
\begin{math}
  \seqstyle{d} \in \littlelp{1}{\zn}{K}
\end{math}
is an isometric isomorphism from
\begin{math}
  \littlelp{1}{\zn}{K}
\end{math}
onto
\begin{math}
  \mdimdualvspace{n}{j}
\end{math}.
We have
\begin{displaymath}
  \dualelem{f} \in \mdimdualvspace{n}{j}
  \iff
  \dualelem{f}(2 \cdot) \in \mdimdualvspace{n}{j+1}
\end{displaymath}
and
\begin{displaymath}
  \dualelem{f} \in \mdimdualvspace{n}{j}
  \iff
  \dualelem{f}(\cdot - 2^{-j} \firstznvar) \in
  \mdimdualvspace{n}{j} .
\end{displaymath}
for all
\(\dualelem{f} \in \topdual{\ucfunc{\rn}{K}}\),
\(j \in \integernumbers\), and
\(\firstznvar \in \zn\).

\if\shortprep0
{
\begin{lemma}
  \label{lem:dualpartialw-norm}
  Let
  \begin{math}
    E = \ucfunc{\rn}{K}
  \end{math}
  or
  \begin{math}
    E = \vanishingfunc{\rn}{K}
  \end{math}.
  Let
  \begin{math}
    j \in \integernumbers
  \end{math}
  and
  \begin{math}
    \sigma : \naturalnumbers \onto \zn
  \end{math}
  be a bijection.
  Let
  \begin{math}
    \sqd \in \littlelp{1}{\zeroonesetnnozero{n} \times \zn}{K}
  \end{math}
  and
  \begin{math}
    \dualelem{f} \in \topdual{E}
  \end{math},
  \begin{displaymath}
    \dualelem{f} :=
    \sum_{\zos \in \zeroonesetnnozero{n}}
    \sum_{\firstznvar \in \zn}
    \seqelem{\sqd}{\zos, \firstznvar}
    \mdimgendualwavelet{n}{\zos}{j}{\firstznvar} .
  \end{displaymath}
  Then
  \begin{displaymath}
    \nsznorminspace{\dualelem{f}}{\topdual{E}}
    =
    \lim_{p \to \infty}
    \sum_{\secondznvar \in \zn}
    \abs{
      \sum_{\zos \in \zeroonesetnnozero{n}}
      \sum_{q=0}^p
      \seqelem{\sqd}{\zos,\sigma(q)}
      \mdimgendualwaveletfilterelem{n}{\zos}{\secondznvar - 2 \sigma(q)}
    } .
  \end{displaymath}
\end{lemma}
}
\fi

\if\shortprep0
{
\begin{proof}
  Define
  \begin{displaymath}
    \begin{array}{llll}
    \indexeddualelem{g}{p}
    & := &
    \sum_{\zos \in \zeroonesetnnozero{n}}
    \sum_{q = 0}^p
    \seqelem{\sqd}{\zos, \sigma(q)}
    \mdimgendualwavelet{n}{\zos}{j}{\sigma(q)} ,
    & p \in \naturalnumbers
    \\
    I(q, \zos)
    & := &
    \left\{
      \secondznvar \in \zn
      \setsep
      \mdimgendualwaveletfilterelem{n}{\zos}{\secondznvar - 2 \sigma(q)}
      \neq
      0
    \right\} ,
    &
    q \in \naturalnumbers, \; \zos \in \zeroonesetnnozero{n}
    \\
    J(p)
    & := &
    \bigcup_{\zos \in \zeroonesetnnozero{n}}
    \bigcup_{q = 0}^p I(q, \zos) ,
    & p \in \naturalnumbers
    \end{array}    
  \end{displaymath}
  Now
  \begin{displaymath}
    \nsznorm{\dualelem{f}} =
    \lim_{p \to \infty} \norm{\indexeddualelem{g}{p}} 
  \end{displaymath}
  and
  \begin{eqnarray*}
    \indexeddualelem{g}{p}
     & = &
     \sum_{\zos \in \zeroonesetnnozero{n}}
     \sum_{q = 0}^p
     \sum_{\secondznvar \in J(p)}
     \seqelem{\sqd}{\zos, \sigma(q)}
     \mdimgendualwaveletfilterelem{n}{\zos}{\secondznvar - 2 \sigma(q)}
     \mdimsf{n}{j+1}{\secondznvar}
     \\
     & = &
     \sum_{\secondznvar \in J(p)}
     \left(
     \sum_{\zos \in \zeroonesetnnozero{n}}
     \sum_{q = 0}^p
     \seqelem{\sqd}{\zos, \sigma(q)}
     \mdimgendualwaveletfilterelem{n}{\zos}{\secondznvar - 2 \sigma(q)}
     \right)
     \mdimsf{n}{j+1}{\secondznvar}
  \end{eqnarray*}
  from which the lemma follows.
\end{proof}
}
\fi

\begin{definition}
  When \(j \in \integernumbers\)
  and \(\zos \in \zeroonesetn\)
  define
  \begin{eqnarray*}
    \mdimpartialdualwspace{n}{\zos}{j}
    & \defequalns &
    \left\{
      \sum_{\firstznvar \in \zn}
      \seqelem{\sqd}{\firstznvar}
      \mdimgendualwavelet{n}{\zos}{j}{\firstznvar}
    \setsep
      \sqd \in \littlelp{1}{\zn}{K}
    \right\} \\
    \norminspace{\dualelem{f}}{\mdimpartialdualwspace{n}{\zos}{j}}
    & := &
    \norminspace{\dualelem{f}}{\topdual{\ucfunc{\rn}{K}}}
    = \norminspace{\dualelem{f}}{\topdual{\vanishingfunc{\rn}{K}}} ,
    \spaceafter
    \dualelem{f} \in \mdimpartialdualwspace{n}{\zos}{j} .
  \end{eqnarray*}  
\end{definition}

\if\shortprep1
{
  Spaces \(\mdimpartialdualwspace{n}{\zos}{j}\),
  \(\zos \in \zeroonesetn\), \(j \in \integernumbers\),
  are topologically isomorphic to
  \(\littlelp{1}{\zn}{K}\).
}
\else
{
By \mylemma \ref{lem:multidim-dual-isomorphism}
function
  \begin{math}
    \dualwsjtopisomfunc{n}{\zos}{j} : \littlelp{1}{\zn}{K} \to 
    \mdimpartialdualwspace{n}{\zos}{j}
  \end{math}
defined by
  \begin{displaymath}
    \dualwsjtopisom{n}{\zos}{j}{\seqstyle{d}}
    :=
    \sum_{\firstznvar \in \integernumbers}     
    \seqelem{\seqstyle{d}}{\firstznvar}
    \mdimgendualwavelet{n}{\zos}{j}{\firstznvar}
  \end{displaymath}
for all \(\sqd \in \littlelp{1}{\zn}{K}\)
is a topological isomorphism from
  \begin{math}
    \littlelp{1}{\zn}{K}
  \end{math}
onto
  \begin{math}
    \mdimpartialdualwspace{n}{\zos}{j}
  \end{math}
and
  \begin{displaymath}
    \frac{1}{\szncover{\mdimmotherwavelet{n}{\zos}}
      \norminfty{\mdimmotherwavelet{n}{\zos}}}
    \normone{\mathbf{d}}
    \leq \norm{\dualwsjtopisom{n}{\zos}{j}{\sqd}}
    \leq \norm{\mdimmotherdualwavelet{n}{\zos}} \normone{\sqd}
  \end{displaymath}
for all \(\sqd \in \littlelp{1}{\zn}{K}\).
}
\fi
\if\shortprep1
{
  Assume that
  \begin{math}
    j \in \integernumbers
  \end{math},
  \begin{math}
    \sqd \in \littlelp{1}{\zeroonesetnnozero{n} \times \zn}{K}
  \end{math},
  and
  \begin{displaymath}
    \dualelem{f} :=
    \sum_{\zos \in \zeroonesetnnozero{n}}
    \sum_{\firstznvar \in \zn}
    \seqelem{\sqd}{\zos, \firstznvar}
    \mdimgendualwavelet{n}{\zos}{j}{\firstznvar} .
  \end{displaymath}
  Now
  \begin{math}
    \nsznorminspace{\dualelem{f}}{\topdual{\ucfunc{\rn}{K}}}
    =
    \nsznorminspace{\dualelem{f}}{\topdual{\vanishingfunc{\rn}{K}}}
  \end{math}.
  Consequently we may identify the Banach space
  \begin{math}
    \mdimpartialdualwspace{n}{\zos}{j}
  \end{math}
  for both
  \begin{math}
    E = \ucfunc{\rn}{K}
  \end{math}
  and
  \begin{math}
    E = \vanishingfunc{\rn}{K}
  \end{math}
  for each
  \begin{math}
    \zos \in \zeroonesetn
  \end{math}
  and
  \begin{math}
    j \in \integernumbers
  \end{math}.
}
\else
{
By \mylemma \ref{lem:dualpartialw-norm}
we may identify the Banach space
\begin{math}
  \mdimpartialdualwspace{n}{\zos}{j}
\end{math}
for both
\begin{math}
  E = \ucfunc{\rn}{K}
\end{math}
and
\begin{math}
  E = \vanishingfunc{\rn}{K}
\end{math}
for each
\begin{math}
  \zos \in \zeroonesetn
\end{math}
and
\begin{math}
  j \in \integernumbers
\end{math}.
}
\fi

\begin{definition}
  When \(n \in \positiveintegers\)
  and \(j \in \integernumbers\)
  define
  \begin{displaymath}
    \mdimdualprojop{n}{j} \dualelem{f}
    :=
    \sum_{\firstznvar \in \zn}
    \szdualappl{\dualelem{f}}
    {\mdimsf{n}{j}{\firstznvar}}
    \mdimdualsf{n}{j}{\firstznvar}
  \end{displaymath}
  for all
  \(\dualelem{f} \in \topdual{E}\).
\end{definition}

\begin{definition}
  When \(n \in \positiveintegers\),
  \(\zos \in \zeroonesetn\), and
  \(j \in \integernumbers\)
  define
  \begin{displaymath}
    \mdimdualpartialdeltaop{n}{\zos}{j} \dualelem{f}
    :=
    \sum_{\firstznvar \in \zn}
    \szdualappl{\dualelem{f}}
    {\mdimgenwavelet{n}{\zos}{j}{\firstznvar}}
    \mdimgendualwavelet{n}{\zos}{j}{\firstznvar}
  \end{displaymath}
  for all
  \(\dualelem{f} \in \topdual{E}\).
\end{definition}

We have
\begin{math}
  \mdimdualpartialdeltaop{n}{\finitezeroseq{n}}{j}
  =
  \mdimdualprojop{n}{j}
\end{math}
for all
\begin{math}
  j \in \integernumbers
\end{math}.
\if\shortprep0
{
When \(\zos \in \zeroonesetn\)
\if\shortprep0
it follows from \mylemmas \ref{lem:dual-gen-conv} and
\ref{lem:dual-appl-in-l-one} that
\fi
operators \(\mdimdualpartialdeltaop{n}{\zos}{j}\),
\(j \in \integernumbers\),
are uniformly bounded by
\begin{displaymath}
  \norm{\mdimdualpartialdeltaop{n}{\zos}{j}} \leq
  \norm{\mdimmotherdualwavelet{n}{\zos}}
  \norminfty{\mdimmotherwavelet{n}{\zos}}
  \ncover{\mdimmotherwavelet{n}{\zos}}
\end{displaymath}
for all \(j \in \integernumbers\)
and operators \(\mdimdualprojop{n}{j}\),
\(j \in \integernumbers\),
uniformly bounded by
\begin{equation}
  \label{eq:mdimdualprojop-uniform-bound}
  \norm{\mdimdualprojop{n}{j}} \leq
  \norm{\mdimmotherdualsf{n}}
  \norminfty{\mdimmothersf{n}}
  \ncover{\mdimmothersf{n}}
\end{equation}
for all \(j \in \integernumbers\).
}
\fi

\if\shortprep0
{
\givelemmaswithoutproofsn{two}
}
\fi

\begin{lemma}
  \label{lem:dualdeltaop-properties}
  \mbox{ }
  \begin{itemize}
    \item[(i)]
      \begin{math}
        \forall j \in \integernumbers,
        \zos \in \zeroonesetn,
        \zot \in \zeroonesetn :
        \zos \neq \zot \implies
        \forall \dualelem{f}
        \in \mdimpartialdualwspace{n}{\zos}{j} :
        \mdimdualpartialdeltaop{n}{\zot}{j}
        \dualelem{f} = 0
      \end{math}.
    \item[(ii)]
      Operator
      \(\mdimdualpartialdeltaop{n}{\zos}{j}\)
      is a continuous projection \projprep
      \(\topdual{E}\)
      onto
      \(\mdimpartialdualwspace{n}{\zos}{j}\)
      for all \(j \in \integernumbers\)
      and \(\zos \in \zeroonesetn\).
    \item[(iii)]
      \(\mdimdualpartialdeltaop{n}{\zos}{j}
      = \mdimdualpartialdeltaop{n}{\zos}{j}
      \circ \mdimdualprojop{n}{j'}\)
      for all \(j, j' \in \integernumbers\),
      \(j < j'\), and
      \(\zos \in \zeroonesetn\).
  \end{itemize}
\end{lemma}

\if\shortprep1
{
  For example, when \(n = 2\) and \(\zot = (0,1)\)
  we have
  \begin{math}
    \mdimdualpartialdeltaop{n}{\zot}{j} \delta = 0
  \end{math}
  for all
  \begin{math}
    j \in \integernumbers
  \end{math}.
  When \(j \in \integernumbers\),
  \(\zos, \zot \in \zeroonesetn\),
  and \(\zos \neq \zot\)
  we have
  \(
  \mdimpartialdualwspace{n}{\zos}{j}
  \intersection
  \mdimpartialdualwspace{n}{\zot}{j}
  = \zeroset
  \).
}
\else
{

For example, when \(n = 2\) and \(\zot = (0,1)\)
we have
\begin{math}
  \mdimdualpartialdeltaop{n}{\zot}{j} \delta = 0
\end{math}
for all
\begin{math}
  j \in \integernumbers
\end{math}.

\begin{lemma}
  \label{lem:partialdualwspace-intersection}
  Let \(j \in \integernumbers\),
  \(\zos, \zot \in \zeroonesetn\),
  and \(\zos \neq \zot\).
  Then
  \(
  \mdimpartialdualwspace{n}{\zos}{j}
  \intersection
  \mdimpartialdualwspace{n}{\zot}{j}
  = \zeroset
  \).
\end{lemma}

\begin{lemma}
  \label{lem:dual-sum-formula}
  Let \(j \in \integernumbers\) and
  \(\dualelem{f} \in \mdimdualvspace{n}{j+1}\).
  Then
  \begin{displaymath}
    \sum_{\zos \in \zeroonesetn}
    \mdimdualpartialdeltaop{n}{\zos}{j}
    \dualelem{f}
    =
    \dualelem{f} .
  \end{displaymath}
\end{lemma}

\begin{proof}
  Use \mylemmas
  \ref{lem:mdim-wavelet-formulas} and
  \ref{lem:mdim-filter-formula} 
  and
  \myequation \eqref{eq:mdim-dual-wavelet-formula}.
\end{proof}
}
\fi

\begin{definition}
  When \(j \in \integernumbers\)
  define
  \begin{displaymath}
    \mdimdualwspace{n}{j}
    \defequalns
    \directsumoneindex{\zos \in \zeroonesetnnozero{n}}
    \mdimpartialdualwspace{n}{\zos}{j} .
  \end{displaymath}
\end{definition}

\begin{definition}
  When \(j \in \integernumbers\)
  define
  \begin{displaymath}
    \mdimdualdeltaop{n}{j}
    :=
    \sum_{\zos \in \zeroonesetnnozero{n}}
    \mdimdualpartialdeltaop{n}{\zos}{j} .
  \end{displaymath}
\end{definition}

\if\shortprep1
{
  Assume that \(j \in \integernumbers\).
  Now operator \(\mdimdualdeltaop{n}{j}\) is a continuous
  projection of \(\topdual{E}\)
  onto \(\mdimdualwspace{n}{j}\)
  and
  \begin{math}
    \mdimdualdeltaop{n}{j} \dualelem{f} = 0
  \end{math}
  for all
  \(\dualelem{f} \in \mdimdualvspace{n}{j}\).
  We also have
  \begin{math}
    \mdimdualdeltaop{n}{j}
    =
    \mdimdualprojop{n}{j+1}
    -
    \mdimdualprojop{n}{j}
  \end{math}
  and
  \begin{math}
    \mdimdualvspace{n}{j+1}
    \equalns
    \mdimdualvspace{n}{j}
    \dsum
    \mdimdualwspace{n}{j}
  \end{math}.
}
\else
{
\begin{lemma}
  \label{lem:mdimdualdeltaprojop-properties}
  Let \(j \in \integernumbers\).
  Then
  \begin{itemize}
  \item[(i)] Operator \(\mdimdualdeltaop{n}{j}\) is a continuous
    projection of \(\topdual{E}\)
    onto \(\mdimdualwspace{n}{j}\).
  \item[(ii)]
    \begin{math}
      \mdimdualdeltaop{n}{j} \dualelem{f} = 0
    \end{math}
    for all
    \(\dualelem{f} \in \mdimdualvspace{n}{j}\).
  \end{itemize}
\end{lemma}

\begin{proof}
  This is a consequence of \mylemma \ref{lem:dualdeltaop-properties}.
\end{proof}

\begin{lemma}
  When \(j \in \integernumbers\)
  we have
  \begin{math}
    \mdimdualdeltaop{n}{j}
    =
    \mdimdualprojop{n}{j+1}
    -
    \mdimdualprojop{n}{j}
  \end{math} .
\end{lemma}

\begin{proof}
  Use \mylemmas
  \ref{lem:dualdeltaop-properties} (iii)
  and
  \ref{lem:dual-sum-formula}.
\end{proof}

\begin{lemma}
  \label{lem:dual-direct-sum}
  When \(j \in \integernumbers\)
  we have
  \begin{math}
    \mdimdualvspace{n}{j+1}
    \equalns
    \mdimdualvspace{n}{j}
    \dsum
    \mdimdualwspace{n}{j}
  \end{math}.
\end{lemma}

\begin{proof}
  Use \mylemmas
  \ref{lem:dualdeltaop-properties} (ii),
  \ref{lem:partialdualwspace-intersection}, and
  \ref{lem:dual-sum-formula}.
\end{proof}
}
\fi

\if\shortprep1
{
\begin{theorem}
  \label{th:mdimspace-biorth}
  \mbox{ }
  \begin{itemize}
  \item[(i)]
    \begin{math}
      \forall j \in \integernumbers,
      \zos \in \zeroonesetn,
      \zot \in \zeroonesetn :
      \zos \neq \zot \iff
      \mdimpartialdualwspace{n}{\zos}{j} \bot
      \mdimupartialwspace{n}{\zot}{j}
    \end{math}
  \item[(ii)]
    \begin{math}
      \forall j_1, j_2 \in \integernumbers :
      j_1 \leq j_2 \iff
      \mdimdualvspace{n}{j_1} \bot
      \mdimuwspace{n}{j_2}
    \end{math}
  \item[(iii)]
    \begin{math}
      \forall j_1, j_2 \in \integernumbers :
      j_1 \geq j_2 \iff
      \mdimdualwspace{n}{j_1} \bot
      \mdimuvspace{n}{j_2}
    \end{math}
  \item[(iv)]
    \begin{math}
      \forall j_1, j_2 \in \integernumbers :
      j_1 \neq j_2 \iff
      \mdimdualwspace{n}{j_1} \bot
      \mdimuwspace{n}{j_2}
    \end{math}
  \end{itemize}
\end{theorem}
}
\fi

\if\shortprep1
{
\begin{proof}
  Let \(j \in \integernumbers\),
  \(\zos, \zot \in \zeroonesetn\), and
  \(\firstznvar, \secondznvar \in \integernumbers^n\).
  Proposition (i) follows from
  \begin{math}
    \szdualappl{\mdimgendualwavelet{n}{\zos}{j}{\secondznvar}}
    {\mdimgenwavelet{n}{\zot}{j}{\firstznvar}}
    = \delta_{\zos,\zot}
    \delta_{\secondznvar,\firstznvar}
  \end{math}.
  Suppose that \(j_1, j_2 \in \integernumbers\) and \(j_1 \neq j_2\).
  If \(j_1 > j_2\) then \(\mdimuwspace{n}{j_2} \subset
  \mdimuvspace{n}{j_1}\) and \(\mdimdualwspace{n}{j_1} \bot
  \mdimuvspace{n}{j_1}\) and hence
  \(\mdimdualwspace{n}{j_1} \bot \mdimuwspace{n}{j_2}\).
  If \(j_1 < j_2\) then \(\mdimdualwspace{n}{j_1} \subset
  \mdimdualvspace{n}{j_2}\) and \(\mdimdualvspace{n}{j_2} \bot
  \mdimuwspace{n}{j_2}\) and hence
  \(\mdimdualwspace{n}{j_1} \bot \mdimuwspace{n}{j_2}\).
  So proposition (iv) is true.
\end{proof}
}
\fi

\if\shortprep0
{
\begin{theorem}
  \label{th:mdimspace-biorth}
  \mbox{ }
  \begin{itemize}
  \item[(i)]
    \begin{math}
      \forall j \in \integernumbers,
      \zos \in \zeroonesetn,
      \zot \in \zeroonesetn :
      \zos \neq \zot \iff
      \mdimpartialdualwspace{n}{\zos}{j} \bot
      \genpartialwspace{E}{n}{\zot}{j}
    \end{math}
  \item[(ii)]
    \begin{math}
      \forall j_1, j_2 \in \integernumbers :
      j_1 \leq j_2 \iff
      \mdimdualvspace{n}{j_1} \bot
      \genwspace{E}{n}{j_2}
    \end{math}
  \item[(iii)]
    \begin{math}
      \forall j_1, j_2 \in \integernumbers :
      j_1 \geq j_2 \iff
      \mdimdualwspace{n}{j_1} \bot
      \genvspace{E}{n}{j_2}
    \end{math}
  \item[(iv)]
    \begin{math}
      \forall j_1, j_2 \in \integernumbers :
      j_1 \neq j_2 \iff
      \mdimdualwspace{n}{j_1} \bot
      \genwspace{E}{n}{j_2}
    \end{math}
  \end{itemize}
\end{theorem}
}
\fi

\if\shortprep0
{
\begin{proof}
  Let \(j \in \integernumbers\),
  \(\zos, \zot \in \zeroonesetn\), and
  \(\firstznvar, \secondznvar \in \integernumbers^n\).
  Proposition (i) follows from \myequation    
  \eqref{eq:mdim-wavelet-orthogonality}.

  Suppose that \(j_1, j_2 \in \integernumbers\) and \(j_1 \neq j_2\).
  If \(j_1 > j_2\) then \(\genwspace{E}{n}{j_2} \subset
  \genvspace{E}{n}{j_1}\) and \(\mdimdualwspace{n}{j_1} \bot
  \genvspace{E}{n}{j_1}\) and hence
  \(\mdimdualwspace{n}{j_1} \bot \genwspace{E}{n}{j_2}\).
  If \(j_1 < j_2\) then \(\mdimdualwspace{n}{j_1} \subset
  \mdimdualvspace{n}{j_2}\) and \(\mdimdualvspace{n}{j_2} \bot
  \genwspace{E}{n}{j_2}\) and hence
  \(\mdimdualwspace{n}{j_1} \bot \genwspace{E}{n}{j_2}\).
  So proposition (iv) is true.
\end{proof}
}
\fi

\begin{theorem}
  \label{th:dual-weakstar-conv}
  Let \(\dualelem{f} \in \dvfn{n}{K}\).
  Then
  \if\shortprep1
  {
  \begin{math}
    \mdimdualprojop{n}{j} \dualelem{f} \rightarrow \dualelem{f}
  \end{math}
  in the weak-* topology of \(\dvfn{n}{K}\)
  as \(j \to \infty\).
  }
  \else
  {
  \begin{math}
    \mdimdualprojop{n}{j} \dualelem{f} \convweakstar \dualelem{f}
  \end{math}
  as \(j \to \infty\).
  }
  \fi
\end{theorem}

\begin{proof}
  Let \(\dualelem{f} \in \dvfn{n}{K}\).
  Suppose that \(f \in \vfn{n}{K}\). Then
  \begin{align*}
    \szdualappl{\mdimdualprojop{n}{j} 
    \dualelem{f}}{\mdimvanprojop{n}{j} f}
    & = \sum_{\firstznvar \in \zn} 
    \szdualappl{\dualelem{f}}{\mdimsf{n}{j}{\firstznvar}}
    \szdualappl{\mdimdualsf{n}{j}{\firstznvar}}
    {\mdimvanprojop{n}{j} f}
    = \szdualappl{\dualelem{f}}{
      \sum_{\firstznvar \in \zn} 
      \szdualappl{\mdimdualsf{n}{j}{\firstznvar}}{f}
      \mdimsf{n}{j}{\firstznvar}} \\
    & = \szdualappl{\dualelem{f}}{\mdimvanprojop{n}{j} f} .
  \end{align*}
  \if\shortprep0
  {
    By \myequation \eqref{eq:mdimdualprojop-uniform-bound} there exists
  }
  \else
  {
    There exists
  }
  \fi
  \(c \in \positiverealnumbers\) so that 
  \(\norm{\mdimdualprojop{n}{j}} \leq
  c\) for all \(j \in \integernumbers\).
  Consequently
  \begin{eqnarray*}
    \abs{\szdualappl{\dualelem{f}}{f} - 
    \szdualappl{\mdimdualprojop{n}{j}
        \dualelem{f}}{f}}
    & \leq &
    \abs{\szdualappl{\dualelem{f}}{f} - 
    \szdualappl{\dualelem{f}}{\mdimvanprojop{n}{j} f}} \\
    & & + \abs{\szdualappl{\mdimdualprojop{n}{j} 
    \dualelem{f}}{\mdimvanprojop{n}{j} f}
      -
      \szdualappl{\mdimdualprojop{n}{j} \dualelem{f}}{f}} \\
    & = & \abs{\szdualappl{\dualelem{f}}{f - \mdimvanprojop{n}{j} f}}
    + \abs{\szdualappl{\mdimdualprojop{n}{j} \dualelem{f}}{\mdimvanprojop{n}{j} f
        - f}} \\
    & \leq &
    \norm{\dualelem{f}} \norm{f - \mdimvanprojop{n}{j} f} +
    c \norm{\dualelem{f}} \norm{f - \mdimvanprojop{n}{j} f}
    \to 0,
  \end{eqnarray*}
  as \(j \to \infty\).
\end{proof}

For example, let
\begin{displaymath}
  \dualelem{f}(x) :=
  \left\{
  \begin{array}{ll}
    x + 1 ; & x  \in \left[-1, 0\right[ \\
        -x + 1 ; & x \in \left[0, 1 \right[ \\
            0 ; & \textrm{otherwise.}
  \end{array}
  \right.
\end{displaymath}
and
\begin{displaymath}
  f(x) :=
  \left\{
  \begin{array}{ll}
    -2x - 2 ; & x  \in \left[-1, 0\right[ \\
        2x - 2 ; & x \in \left[0, 1 \right[ \\
            0 ; & \textrm{otherwise.}
  \end{array}
  \right.
\end{displaymath}
for all \(x \in \realnumbers\).
Now
\begin{displaymath}
  \szdualappl{\mdimdualprojop{1}{j}
    \dualelem{f}}{f}
  \to
  \szdualappl{\dualelem{f}}{f}
  =
  \int_\realnumbers \dualelem{f}(x) f(x) dx
  = -\frac{4}{3}
\end{displaymath}
as
\begin{math}
  j \to \infty
\end{math}.

\subsection{Tensor Product Representation of the Dual MRA}

\if\shortprep0
{
\begin{lemma}
  \label{lem:dvf-inj-schatten-tp}
  Banach space
  \begin{math}
    \dvfn{n}{K} \cistp \dvf{K}
  \end{math}
  is a closed subspace of
  \begin{math}
    \dvfn{n+1}{K}
  \end{math}.
\end{lemma}
}
\fi

\if\shortprep0
{
\begin{proof}
  \if\shortprep1
  {
  Use \mytheorem \ref{th:vanishing-tp-space} and
  \myequation \eqref{eq:schatten-embedding}.
  }
  \else
  {
  By \mytheorem \ref{th:vanishing-tp-space}
  \if\elsevier0
  {
    \begin{align*}
      \vfn{n}{K} \citp \vf{K}
      & = \left( \indexedcitp{k=1}{n} \vf{K} \right) \citp \vf{K}
      = \indexedcitp{k=1}{n+1} \vf{K} \\
      & = \vfn{n+1}{K} .
    \end{align*}
  }
  \else
  {
    \begin{align*}
      \vfn{n}{K} \citp \vf{K}
      & = \left( \indexedcitp{k=1}{n} \vf{K} \right) \citp \vf{K} \\
      & = \indexedcitp{k=1}{n+1} \vf{K}
      = \vfn{n+1}{K} .
    \end{align*}
  }
  \fi
  Hence it follows from \myequation \eqref{eq:schatten-embedding} that
  \begin{math}
    \dvfn{n}{K} \cistp \dvf{K}
  \end{math}
  is a closed subspace of
  \begin{math}
    \dvfn{n+1}{K}
  \end{math}.
  }
  \fi
\end{proof}
}
\fi

\if\shortprep1
{
  It follows from \myequation \eqref{eq:schatten-embedding} that
  \begin{math}
    \dvfn{n}{K} \cistp \dvf{K}
  \end{math}
  is a closed subspace of
  \begin{math}
    \dvfn{n+1}{K}
  \end{math}.
}
\fi

\begin{definition}
  \label{def:mspace}
  Let \(j \in \integernumbers\).
  Define
  \begin{displaymath}
    \tensorauxdualvspace{n}{j}
    \defequalns
    \indexedinhctp{\topdual{\vanishingfunc{\realnumbers^k}{K}}}
    {k=1}{n}
    \mdimdualvspace{1}{j} .
  \end{displaymath}
\end{definition}

\if\shortprep1
{
Banach space \(\tensorauxdualvspace{n}{j}\) is a closed subspace
of \(\dvfn{n}{K}\).
}
\else
{
The Banach spaces \(\tensorauxdualvspace{n}{j}\) are well defined 
since \(\onedimdualvspace{j}\) is a closed subspace of
\(\topdual{\vanishingfunc{\realnumbers}{K}}\).
Vector space
\(\topdual{\vanishingfunc{\realnumbers^k}{K}} \otimes \onedimdualvspace{j}\)
is a linear
subspace of \(\topdual{\vanishingfunc{\realnumbers^k}{K}}
\otimes \topdual{\vanishingfunc{\realnumbers}{K}}\) for all
\(k \in \positiveintegers\).
It follows from \mylemma \ref{lem:dvf-inj-schatten-tp}
that
\(\topdual{\vanishingfunc{\realnumbers^k}{K}} \otimes \onedimdualvspace{j}\)
is a linear subspace of
\(\topdual{\vanishingfunc{\realnumbers^{k+1}}{K}}\)
for all \(k \in \positiveintegers\).
Consequently
\if10 
the Banach spaces \(\mdimdualvspace{n}{j}\) are well
defined and
\fi 
\(\mdimdualvspace{n}{j}\) is a closed subspace of
\(\topdual{\vanishingfunc{\rn}{K}}\) for all \(n \in \positiveintegers\) and
\(j \in \integernumbers\).
}
\fi
\if\shortprep0
{
Since
\begin{math}
  \onedimdualvspace{j} \closedsubspace \onedimdualvspace{j+1}
\end{math}
for all \(j \in \integernumbers\) it follows that
\begin{math}
  \mdimdualvspace{n}{j} \closedsubspace
  \mdimdualvspace{n}{j+1}
\end{math}
for all \(n \in \positiveintegers\) and \(j \in \integernumbers\).
}
\fi

\if\shortprep1
{
  \begin{definition}
    Define function
    \begin{displaymath}
      \znisomfunc{n} : 
      \indexedcptp{k=1}{n} \littlelp{1}{\integernumbers}{K}
      \to
      \littlelp{1}{\zn}{K}
    \end{displaymath}
    by
    \begin{displaymath}
      \znisom{n}{\ztensornbv{\firstznvar}} := \znatbasvec{\firstznvar} ,
      \spaceafter \firstznvar \in \zn ,
    \end{displaymath}
    and extending by linearity and continuity onto whole
    \begin{math}
      \indexedcptp{k=1}{n} \littlelp{1}{\integernumbers}{K}
    \end{math}.
  \end{definition}

  Function \(\znisomfunc{n}\) is an isometric isomorphism from
  Banach space
  \begin{math}
    \indexedcptp{k=1}{n} \littlelp{1}{\integernumbers}{K}
  \end{math}
  onto Banach space
  \begin{math}
    \littlelp{1}{\zn}{K} .
  \end{math}      
}
\fi

\begin{lemma}
  \label{lem:tensordualvj-proj}
  When \(j \in \integernumbers\)
  \begin{displaymath}
    \indexedcptp{k=1}{n} \onedimdualvspace{j} \equalns \mdimdualvspace{n}{j} .
  \end{displaymath}
\end{lemma}

\begin{proof}
  Let
  \begin{displaymath}
    E \defequalns \indexedcptp{k=1}{n} \onedimdualvspace{j}
  \end{displaymath}
  and
  \begin{displaymath}
    F \defequalns \indexedcptp{k=1}{n} \littlelp{1}{\integernumbers}{K} .
  \end{displaymath}
  Let
  \begin{displaymath}
    \alpha := \indexedopcptp{k=1}{n} \onedimdualvjiisomfunc{j} 
  \end{displaymath}
  Function \(\alpha\) is an isometric
  isomorphism from \(F\) onto \(E\).
  \if\shortprep0
  {
  Define function
  \begin{math}
    \znisomfunc{n}
  \end{math}
  as in \mylemma \ref{lem:lone-zn-isomorphism}.
  }
  \fi
  Let
  \begin{math}
    \beta := \alpha \circ (\znisomfunc{n})^{-1}
  \end{math}.
  Now \(\beta\) is an isometric isomorphism from
  \begin{math}
    \littlelp{1}{\zn}{K}
  \end{math}
  onto \(E\) and
  \begin{math}
    \beta(\znatbasvec{\firstznvar}) = \mdimdualsf{n}{j}{\firstznvar}
    = \mdimdualvjiisom{n}{j}{\znatbasvec{\firstznvar}}
  \end{math}
  for all
  \(\firstznvar \in \zn\).
  When
  \begin{math}
    \seqstyle{d} \in \littlelp{1}{\zn}{K}
  \end{math}
  we have
  \begin{math}
    \beta(\seqstyle{d}) = \mdimdualvjiisom{n}{j}{\seqstyle{d}}
  \end{math}
  and
  \begin{math}
    \norm{\beta(\seqstyle{d})}_E
    =
    \norm{\mdimdualvjiisom{n}{j}{\seqstyle{d}}}_{\dvfn{n}{K}}
  \end{math}.
\end{proof}

\begin{lemma}
  \label{lem:tensordualvj-repr}
  When \(j \in \integernumbers\)
  \begin{displaymath}
    \mdimdualvspace{n}{j} \equalns
    \tensorauxdualvspace{n}{j} \equalns
    \indexedcptp{k=1}{n} \mdimdualvspace{1}{j}
    \equalns
    \indexedctp{\injtn^s}{k=1}{n} \mdimdualvspace{1}{j} .
  \end{displaymath}
\end{lemma}

\begin{proof}
  Use induction by \(n\), metric approximation property of \(l^1\), and
  \cite[prop. 7.1]{ryan2002}.
\end{proof}

\begin{definition}
  \label{def:nspace}
  Define
  \begin{displaymath}
    \nspace{n}{\zos}{j}
    \defequalns
    \indexedinhctp{\topdual{\vanishingfunc{\realnumbers^n}{K}}}
    {k=1}{n}
    \mdimpartialdualwspace{1}{\seqelem{\zos}{k}}{j}
  \end{displaymath}
  where \(j \in \integernumbers\) and
  \(\zos \in \zeroonesetn\).
\end{definition}

\if\shortprep0
{
\begin{lemma}
  \label{lem:tensorpartialdualwspace-clos-span}
  Let \(j \in \integernumbers\) and
  \(\zos \in \zeroonesetn\).
  Then
  \begin{displaymath}
    \nspace{n}{\zos}{j} =
    \closop_{\topdual{\vanishingfunc{\realnumbers^n}{K}}}
    \spanop \{ \mdimgendualwavelet{n}{\zos}{j}{\firstznvar}
    \setsep \firstznvar \in
    \integernumbers^n \} .
  \end{displaymath}
\end{lemma}
}
\fi

\if\shortprep0
{
\begin{proof}
  Let
  \begin{displaymath}
    \aspace{n}{\zos}{j} = \closop_{\topdual{\vanishingfunc{\realnumbers^n}{K}}}
    \spanop \{ \mdimgendualwavelet{n}{\zos}{j}{\firstznvar}
    \setsep \firstznvar \in \integernumbers^n \}
  \end{displaymath}
  for all \(n \in \positiveintegers\), \(j \in \integernumbers\), and
  \(\zos \in \zeroonesetn\).
  We have \(\mdimgendualwavelet{n}{\zos}{j}{\firstznvar} \in
  \nspace{n}{\zos}{j}\) for all
  \(\firstznvar \in \integernumbers^n\).
  Hence
  \begin{math}
    \spanop \{ \mdimgendualwavelet{n}{\zos}{j}{\firstznvar}
    \setsep \firstznvar \in
    \integernumbers^n \} \subset 
    \nspace{n}{\zos}{j} .
  \end{math}
  As \(\nspace{n}{\zos}{j}\) is a Banach 
  space it follows that \(\aspace{n}{\zos}{j}\) is a closed 
  subspace of
  \(\nspace{n}{\zos}{j}\).

  When \(n' = 1\) we have
  \begin{math}
    \nspace{1}{\zos}{j} \equalns 
    \mdimdualvspace{1}{j}
  \end{math}
  or
  \begin{math}
    \nspace{1}{\zos}{j} \equalns
    \mdimdualwspace{1}{j}
  \end{math}
  and hence
  \begin{math}
    \nspace{1}{\zos}{j}
    \closedsubspace
    \aspace{1}{\zos}{j}
  \end{math},
  from which it follows that
  \begin{math}
    \nspace{1}{\zos}{j}
    \equalns
    \aspace{1}{\zos}{j}
  \end{math}.
  Suppose that
  \(\nspace{n'}{\zos}{j} = \aspace{n'}{\zos}{j}\)
  for all \(j \in \integernumbers\), \(\zos \in \zeroonesetnprime\)
  for some \(n' \in \positiveintegers\) (induction assumption).
  Let
  \(\dualelem{x} \in 
  \nspace{n'+1}{\znstyleprime{s}}{j}\)
  for some
  \begin{math}
    \znstyleprime{s} \in \cpzerooneset{n'+1}
  \end{math}.
  Let \(h \in \positiverealnumbers\).
  Let
  \begin{math}
    \znstyle{t} := \seqproj{n'}{\znstyleprime{s}}
  \end{math}
  and
  \begin{math}
    u := \seqelem{\znstyleprime{s}}{n'+1}
  \end{math}.
  There exists
  \(\dualelem{y} \in \nspace{n'}{\znstyle{t}}{j}
  \otimes_{(\topdual{\vanishingfunc{\realnumbers^{n'+1}}{K}})}
  \mdimpartialdualwspace{1}{u}{j}\)
  so that
  \(\norm{\dualelem{x}-\dualelem{y}} < \frac{h}{2}\).
  Now
  \begin{displaymath}
    \dualelem{y} = \sum_{k=1}^m \indexeddualelem{w}{k}
    \otimes \indexeddualelem{v}{k}
  \end{displaymath}
  where \(m \in \naturalnumbers\) and
  \(\indexeddualelem{w}{k} \in 
  \nspace{n'}{\znstyle{t}}{j}\),
  \(\indexeddualelem{v}{l} \in \mdimpartialdualwspace{1}{u}{j}\)
  for each \(k \in \setoneton{m}\).
  Let
  \(c := \max \{ \norm{\indexeddualelem{v}{k}}
    \setsep k \in \setoneton{m} \} + 1\).
  There exist \(\indexeddualelem{r}{k} \in \spanop \{
  \mdimgendualwavelet{n'}{\znstyle{t}}{j}{\firstznvar}
  \setsep \firstznvar \in
  \integernumbers^{n'} \}\) so that
  \begin{displaymath}
    \norm{\indexeddualelem{w}{k} - \indexeddualelem{r}{k}} < 
    \frac{h}{2mc}
  \end{displaymath}
  for each
  \begin{math}
    k \in \setoneton{m}
  \end{math}.

  \if\shortprep0
  {
  Now
  \begin{displaymath}
    \indexeddualelem{r}{k} = \sum_{\secondznvar \in J} 
    b_{k,\secondznvar}
    \mdimgendualwavelet{n'}{\znstyle{t}}{j}{\secondznvar}
  \end{displaymath}
  where \(J\) is a finite subset of \(\integernumbers^{n'}\) and
  \(b_{k,\secondznvar} \in K\) for all \(\secondznvar \in J\)
  and \(k \in \setoneton{m}\).
  We also have
  \begin{displaymath}
    \indexeddualelem{v}{k} = \sum_{p \in \integernumbers} d_{k,p}
    \gendualwaveletonedim{u}{j}{p}
  \end{displaymath}
  where \((d_{k,p})_{p=0}^\infty \in 
  \littlelp{1}{\integernumbers}{K}\)
  and \(k \in \setoneton{m}\).
  }
  \fi
  \if\shortprep0
  {
  It follows that
  \begin{displaymath}
    \indexeddualelem{r}{k} \otimes \indexeddualelem{v}{k} =
    \sum_{\secondznvar \in J}
    \sum_{p \in \integernumbers}
    b_{p,\secondznvar} d_{k,p}
    \mdimgendualwavelet{n'+1}{\znstyleprime{s}}{j}
    {\seqcomb{\secondznvar}{p}}
    \in \aspace{n'+1}{\znstyleprime{s}}{j}
  \end{displaymath}
  for each
  \begin{math}
    k \in \setoneton{m}
  \end{math}.
  }
  \else
  {
  We have
  \begin{math}
    \indexeddualelem{r}{k} \otimes \indexeddualelem{v}{k}
    \in \aspace{n'+1}{\znstyleprime{s}}{j}
  \end{math}
  for each
  \begin{math}
    k \in \setoneton{m}
  \end{math}.
  }
  \fi
  Let
  \begin{displaymath}
    \dualelem{z} = \sum_{k=1}^m \indexeddualelem{r}{k}
    \otimes \indexeddualelem{v}{k} .
  \end{displaymath}
  Now \(\dualelem{z} \in \aspace{n'+1}{\zos}{j}\) and
  \begin{displaymath}
    \norm{\dualelem{y} - \dualelem{z}} \leq \sum_{k=1}^m
    \norm{\indexeddualelem{w}{k} - \indexeddualelem{r}{k}}
    \norm{\indexeddualelem{v}{k}}
    < \sum_{p=1}^m \frac{h}{2mc} c = \frac{h}{2} .
  \end{displaymath}
  Hence
  \(\norm{\dualelem{x} - \dualelem{z}}
  \leq \norm{\dualelem{x} - \dualelem{y}}
  + \norm{\dualelem{y} - \dualelem{z}} <
  \frac{h}{2} + \frac{h}{2} = h\).
  \if\shortprep0
  {
  Number \(h > 0\) was arbitrary and hence
  \(\dualelem{x} \in \aspace{n'+1}{\znstyleprime{s}}{j}\).
  Therefore the proposition is true for \(n'+1\) and consequently 
  for
  all \(n \in \positiveintegers\).
  }
  \else
  {
  Thus
  \(\dualelem{x} \in \aspace{n'+1}{\znstyleprime{s}}{j}\).
  }
  \fi
\end{proof}
}
\fi

\begin{definition}
  \label{def:mdimaltpartialdualprojop}
  Let \(j \in \integernumbers\) and
  \(\zos \in \zeroonesetn\). Define
  \if\shortprep1
  operator
  \else
  function
  \fi
  \begin{math}
    \mdimaltpartialdualprojop{n}{j}{\zos} : \mdimdualvspace{n}{j+1} \to
    \nspace{n}{\zos}{j}
  \end{math}
  by
  \begin{displaymath}
    \mdimaltpartialdualprojop{n}{j}{\zos} =
    \indexedinhopctp{\topdual{\vanishingfunc{\realnumbers^k}{K}}}
    {\topdual{\vanishingfunc{\realnumbers^k}{K}}}{k=1}{n}
    \left( \mdimpartialdualprojop{1}{\seqelem{\zos}{k}}{j} \vert 
    \mdimdualvspace{1}{j+1}
    \right)
  \end{displaymath}
\end{definition}

\if\shortprep0
{
\begin{theorem}
  \label{th:mdimaltpartialdualprojop-well-defined}
  Functions \(\mdimaltpartialdualprojop{n}{j}{\zos}\),
  \(j \in \integernumbers\),
  \(\zos \in \zeroonesetn\)
  are well defined, linear,
  and continuous.
\end{theorem}
}
\fi

\if\shortprep0
{
\begin{proof}
  Let \(j \in \integernumbers\) and
  \(\zos \in \zeroonesetn\).
  Functions \(\mdimvanpartialdualwspace{1}{b}{j}\),
  \(b \in \{ 0, 1 \}\), are linear and
  continuous.
  Let \(P_k = \mdimvanpartialdualprojop{1}{\seqelem{\zos}{k}}{j}
  \vert 
  \onedimvandualvspace{j+1}\) for
  \(l = 1, \ldots, n\). Define the operators \(S_k\) and \(T_k\) as in
  \mydef \ref{def:indexed-inh-op-tp}.
  Now \(T_1 = S_1 = P_1\) is continuous.
  Suppose that \(T_k\) is continuous for some
  \(k \in \setoneton{n - 1}\).
  Operator \(T_k\) is a linear operator from Banach space
  \(\mdimdualvspace{k}{j+1}\) into Banach space
  \begin{displaymath}
    \indexedinhctp{\topdual{\vanishingfunc{\realnumbers^m}{K}}}{m=1}{k}
    \mdimvanpartialdualwspace{1}{\cartprodelem{\zos}{m}}{j} \equalns
    \nspace{k}{\seqproj{k}{\zos}}{j} .
  \end{displaymath}
  Operator \(P_{k+1}\) is a continuous linear operator from
  \(\onedimvandualvspace{j+1}\) into 
  \(\mdimvanpartialdualwspace{1}{\cartprodelem{\zos}{k+1}}{j}\).
  Subspace \(\mdimdualvspace{k}{j+1}\) is topologically complemented
  in \(\topdual{\vanishingfunc{\realnumbers^k}{K}}\) and subspace
  \(\onedimvandualvspace{j+1}\) is topologically complemented in
  \(\topdual{\vanishingfunc{\realnumbers}{K}}\).
  By \mylemma \ref{lem:operator-norm-tp-operator-continuity} operator
  \(S_{k+1} = T_k \otimes P_{k+1}\) is a continuous linear operator
  from normed vector space
  \begin{math}
    \mdimdualvspace{k}{j+1}
    \otimes_{(\topdual{\vanishingfunc{\realnumbers^{k+1}}{K}})}
    \onedimvandualvspace{j+1}
  \end{math}
  into normed vector space
  \begin{math}
    \nspace{k}{\seqproj{k}{\zos}}{j}
    \otimes_{(\topdual{\vanishingfunc{\realnumbers^{k+1}}{K}})}
    \mdimvanpartialdualwspace{1}{\cartprodelem{\zos}{k+1}}{j}
  \end{math}.
  Hence \(T_{k+1}\) is a continuous linear operator from Banach space
  \(\mdimdualvspace{k+1}{j+1}\) into Banach space
  \(\nspace{k+1}{\seqproj{k+1}{\zos}}{j}\).
  It follows that \(\mdimaltpartialdualprojop{n}{j}{\zos} = T_n\) is
  well defined
  and it is a continuous linear operator from  Banach space
  \(\mdimdualvspace{n}{j+1}\) into Banach space
  \(\nspace{n}{\zos}{j}\).
\end{proof}
\fi

\if\shortprep0
{
\begin{lemma}
  \label{lem:mdimaltpartialdualprojop-series-formula}
  Let \(j \in \integernumbers\), and
  \(\zos \in \zeroonesetn\).
  Then
  \begin{equation}
    \label{eq:mdimaltpartialdualprojop-series-formula}
    \mdimaltpartialdualprojop{n}{j}{\zos} \dualelem{f}
    = \sum_{\firstznvar \in \integernumbers^n}
    \dualappl{\dualelem{f}}
    {\mdimgenwavelet{n}{\zos}{j}{\firstznvar}}
    \mdimgendualwavelet{n}{\zos}{j}{\firstznvar}
  \end{equation}
  for all \(\dualelem{f} \in \mdimdualvspace{n}{j+1}\).
  The series in \myequation 
  \eqref{eq:mdimaltpartialdualprojop-series-formula}
  converges absolutely for all \(\dualelem{f} \in 
  \mdimdualvspace{n}{j+1}\).
  \if\shortprep0
  When \(\zos \in \zeroonesetn\)
  operators \(\mdimaltpartialdualprojop{n}{j}{\zos}\),
  \(j \in \integernumbers\),
  are uniformly bounded by
  \begin{displaymath}
    \norm{\mdimaltpartialdualprojop{n}{j}{\zos}}
    \leq \ncover{\mdimmotherwavelet{n}{\zos}}
    \norminfty{\mdimmotherwavelet{n}{\zos}}
    \norm{\mdimmotherdualwavelet{n}{\zos}}
  \end{displaymath}
  for all
  \begin{math}
    j \in \integernumbers
  \end{math}.
  \fi
\end{lemma}
}
\fi

\if\shortprep0
{
\begin{proof}
  Define linear operators
  \begin{math}
    \dualtoperator{n}{\zosp}{j'} : \dvfn{n}{K}
    \to \mdimvanpartialdualwspace{n}{\zosp}{j'}
  \end{math},
  \begin{math}
    j' \in \integernumbers
  \end{math},
  \begin{math}
    \zosp \in \zeroonesetn
  \end{math},
  by
  \begin{equation}
    \label{eq:mapdp-a}
    \dualtoperator{n}{\zosp}{j'} \dualelem{g}
    := \sum_{\firstznvar \in 
    \integernumbers^n}
    \szdualappl{\dualelem{g}}
    {\mdimgenwavelet{n}{\zosp}{j'}{\firstznvar}}
    \mdimgendualwavelet{n}{\zosp}{j'}{\firstznvar}
  \end{equation}
  for all
  \begin{math}
    \dualelem{g} \in \dvfn{n}{K}
  \end{math}.
  By \mylemmas \ref{lem:dual-appl-in-l-one} and
  \ref{lem:dual-gen-conv} the series in \myequation
  \eqref{eq:mapdp-a} converges absolutely,
  operators \(\dualtoperator{n}{\zosp}{j'}\)
  \if\shortprep1
  {
  are well-defined and continuous.
  }
  \else
  {
  are well-defined and continuous, and
  \begin{equation}
    \label{eq:toperator-uniform-bound-ineq}
    \norm{\dualtoperator{n}{\zosp}{j'}}
    \leq \ncover{\mdimmotherwavelet{n}{\zosp}}
    \norminfty{\mdimmotherwavelet{n}{\zosp}}
    \norm{\mdimmotherdualwavelet{n}{\zosp}}
  \end{equation}
  for all
  \begin{math}
    j' \in \integernumbers
  \end{math}
  and
  \begin{math}
    \zosp \in \zeroonesetn
  \end{math}.
  }
  \fi
  
  Let \(\dualelem{f} \in \mdimdualvspace{n}{j+1}\). By \mylemma
  \ref{lem:tensordualvj-repr}
  \begin{displaymath}
    \dualelem{f} 
    = \sum_{\firstznvar \in \integernumbers^n}
    \seqelem{\sqa}{\firstznvar}
    \mdimdualsf{n}{j+1}{\firstznvar}
  \end{displaymath}
  where
  \begin{math}
    \sqa \in
    \littlelp{1}{\zn}{K}
  \end{math}
  and the series converges absolutely.
  Let \(\sigma : \naturalnumbers \onto \integernumbers^n\) be a bijection.
  Let
  \begin{displaymath}
    \indexeddualelem{f}{m} := \sum_{k=0}^m \seqelem{\sqa}{\sigma(k)}
    \mdimdualsf{n}{j+1}{\sigma(k)}
  \end{displaymath}
  for all \(m \in \naturalnumbers\).
  Now \(\indexeddualelem{f}{m} \to \dualelem{f}\) as
  \(m \to \infty\).
  Furthermore,
  \if\elsevier0
  {
  \begin{align}
    \nonumber
    \mdimaltpartialdualprojop{n}{j}{\zos}
    \indexeddualelem{f}{m}
    & = \sum_{k=0}^m
    \seqelem{\sqa}{\sigma(k)} \indexedtensorproduct{l=1}{n} \left(
      \genprojoponedim{\cartprodelem{\zos}{l}}{j}
      \mdimdualsf{n}{j+1}{\sigma(k)}
    \right) \\
    \nonumber
    & = \sum_{k=0}^m
    \seqelem{\sqa}{\sigma(k)} \indexedtensorproduct{l=1}{n} \left(
      \sum_{p \in \integernumbers}
      \szdualappl{\phidual_{j+1,\cartprodelem{\sigma(k)}{l}}}
        {\genwaveletonedim{\cartprodelem{\zos}{l}}{j}{p}}
      \gendualwaveletonedim{\cartprodelem{\zos}{l}}{j}{p}
    \right)
    \\
    \label{eq:mapdp-b}
    & = \sum_{k=0}^m
    \seqelem{\sqa}{\sigma(k)}
    \sum_{p_1 \in \integernumbers} \cdots \sum_{p_n \in \integernumbers}
    \szdualappl{\mdimdualsf{n}{j+1}{\sigma(k)}}
      {\mdimgenwavelet{n}{\zos}{j}{(p_1,\ldots,p_n)}}
    \mdimgendualwavelet{n}{\zos}{j}{(p_1,\ldots,p_n)}
  \end{align}
  }
  \else
  {
  \begin{align}
    \nonumber
    \mdimaltpartialdualprojop{n}{j}{\zos}
    \indexeddualelem{f}{m}
    & = \sum_{k=0}^m
    \seqelem{\sqa}{\sigma(k)} \indexedtensorproduct{l=1}{n} \left(
      \genprojoponedim{\cartprodelem{\zos}{l}}{j}
      \mdimdualsf{n}{j+1}{\sigma(k)}
    \right) \\
    \if\shortprep0
    {
    \nonumber
    & = \sum_{k=0}^m
    \seqelem{\sqa}{\sigma(k)} \indexedtensorproduct{l=1}{n} \left(
      \sum_{p \in \integernumbers}
      \szdualappl{\phidual_{j+1,\cartprodelem{\sigma(k)}{l}}}
        {\genwaveletonedim{\cartprodelem{\zos}{l}}{j}{p}}
      \gendualwaveletonedim{\cartprodelem{\zos}{l}}{j}{p}
    \right)
    \\
    }
    \fi
    \label{eq:mapdp-b}
    & = \sum_{k=0}^m
    \seqelem{\sqa}{\sigma(k)}
    \sum_{p_1 \in \integernumbers} \cdots \sum_{p_n \in \integernumbers}
    \szdualappl{\mdimdualsf{n}{j+1}{\sigma(k)}}
      {\mdimgenwavelet{n}{\zos}{j}{(p_1,\ldots,p_n)}}
    \mdimgendualwavelet{n}{\zos}{j}{(p_1,\ldots,p_n)}
  \end{align}
  }
  \fi
  Since
  \begin{math}
    \mdimgenwavelet{n}{\zos}{j}{\firstznvar} =
    \mdimmotherwavelet{n}{\zos}(2^j \cdot - \firstznvar)
  \end{math}
  and functions \(\mdimmotherwavelet{n}{\zos}\) are compactly supported
  the series in formula \eqref{eq:mapdp-b} has only finite number of
  nonzero terms.
  Hence
  \if\elsevier0
  {
  \begin{align*}
    \mdimaltpartialdualprojop{n}{j}{\zos} \indexeddualelem{f}{m}
    & = \sum_{k=0}^m
    \seqelem{\sqa}{\sigma(k)}
    \sum_{\firstznvar \in \integernumbers^n}
    \szdualappl{\mdimdualsf{n}{j+1}{\sigma(k)}}
    {\mdimgenwavelet{n}{\zos}{j}{\firstznvar}}
    \mdimgendualwavelet{n}{\zos}{j}{\firstznvar} \\
    & =
    \sum_{\firstznvar \in \integernumbers^n}
    \szdualappl{\sum_{k=0}^m
      \seqelem{\sqa}{\sigma(k)}
      \mdimdualsf{n}{j+1}{\sigma(k)}}
      {\mdimgenwavelet{n}{\zos}{j}{\firstznvar}}
    \mdimgendualwavelet{n}{\zos}{j}{\firstznvar}
    \\
    & =
    \sum_{\firstznvar \in \integernumbers^n}
    \szdualappl{\indexeddualelem{f}{m}}
    {\mdimgenwavelet{n}{\zos}{j}{\firstznvar}}
    \mdimgendualwavelet{n}{\zos}{j}{\firstznvar}
    = \dualtoperator{n}{\zos}{j} \indexeddualelem{f}{m} .
  \end{align*}
  }
  \else
  {
  \begin{align*}
    \mdimaltpartialdualprojop{n}{j}{\zos} \indexeddualelem{f}{m}
    & = \sum_{k=0}^m
    \seqelem{\sqa}{\sigma(k)}
    \sum_{\firstznvar \in \integernumbers^n}
    \szdualappl{\mdimdualsf{n}{j+1}{\sigma(k)}}
    {\mdimgenwavelet{n}{\zos}{j}{\firstznvar}}
    \mdimgendualwavelet{n}{\zos}{j}{\firstznvar} \\
    \if\shortprep0
    {
    & =
    \sum_{\firstznvar \in \integernumbers^n}
    \szdualappl{\sum_{k=0}^m
      \seqelem{\sqa}{\sigma(k)}
      \mdimdualsf{n}{j+1}{\sigma(k)}}
      {\mdimgenwavelet{n}{\zos}{j}{\firstznvar}}
    \mdimgendualwavelet{n}{\zos}{j}{\firstznvar} \\
    }
    \fi
    & =
    \sum_{\firstznvar \in \integernumbers^n}
    \szdualappl{\indexeddualelem{f}{m}}
    {\mdimgenwavelet{n}{\zos}{j}{\firstznvar}}
    \mdimgendualwavelet{n}{\zos}{j}{\firstznvar}
    = \dualtoperator{n}{\zos}{j} \indexeddualelem{f}{m} .
  \end{align*}
  }
  \fi
  Since operators \(\mdimaltpartialdualprojop{n}{j}{\zos}\) and
  \(\dualtoperator{n}{\zos}{j}\) are continuous we
  \if\shortprep0
  {
  have
  \begin{displaymath}
    \mdimaltpartialdualprojop{n}{j}{\zos} \dualelem{f}
    = \lim_{m \to \infty} \mdimaltpartialdualprojop{n}{j}{\zos} 
    \indexeddualelem{f}{m}
    = \lim_{m \to \infty} \dualtoperator{n}{\zos}{j} 
    \indexeddualelem{f}{m}
    = \dualtoperator{n}{\zos}{j} \dualelem{f} .
  \end{displaymath}
  }
  \else
  {
  get
  \begin{math}
    \mdimaltpartialdualprojop{n}{j}{\zos}
    =
    \dualtoperator{n}{\zos}{j}
  \end{math}.
  }
  \fi
\end{proof}
}
\fi

\if\shortprep0
{
\begin{lemma}
  \label{lem:mdimaltpartialdualprojop-properties}
  Let \(j \in \integernumbers\).
  Then
  \begin{itemize}
  \item[(i)]
    \begin{math}
      \forall \zos \in \zeroonesetn, \dualelem{f} \in
      \nspace{n}{\zos}{j} : \;
      \mdimaltpartialdualprojop{n}{j}{\zos} \dualelem{f} = 
      \dualelem{f}
    \end{math}
  \item[(ii)]
    \begin{math}
      \forall \zos \in \zeroonesetn, \zot \in \zeroonesetn,
      \dualelem{f} \in
      \nspace{n}{\zos}{j} :
      \zos \neq \zot \implies 
      \mdimaltpartialdualprojop{n}{j}{\zot}
      \dualelem{f} = 0
    \end{math}
  \item[(iii)]
    \begin{math}
      \forall \zos \in \zeroonesetnnozero{n}, \dualelem{f} \in
      \mdimdualvspace{n}{j} :
      \mdimaltpartialdualprojop{n}{j}{\zos} \dualelem{f} = 0
    \end{math}.
  \item[(iv)]
    Operator \(\mdimaltpartialdualprojop{n}{j}{\zos}\) is a 
    projection \projprep
    \(\mdimdualvspace{n+1}{j+1}\) onto
    \(\nspace{n}{\zos}{j}\)
    for each \(\zos \in \zeroonesetn\).
  \end{itemize}
\end{lemma}
}
\fi

\if\shortprep1
{
\if10
{
\begin{proof}
    Let \(\zos, \zot \in \zeroonesetn\) and
    \(\dualelem{f} \in \nspace{n}{\zos}{j}\).
    By \mylemma \ref{lem:tensorpartialdualwspace-clos-span}
    there exists a sequence
    \begin{math}
      (\indexeddualelem{f}{m})_{m=0}^\infty \subset
      \spanop \left\{ \mdimgendualwavelet{n}{\zos}{j}{\firstznvar} 
      \setsep
        \firstznvar \in \integernumbers^n \right\}
    \end{math}
    so that \(\indexeddualelem{f}{m} \to \dualelem{f}\) as \(m \to 
    \infty\).
    Now
    \begin{displaymath}
      \indexeddualelem{f}{m} = \sum_{\firstznvar \in J_m} 
      a_{m,\firstznvar}
        \mdimgendualwavelet{n}{\zos}{j}{\firstznvar}
    \end{displaymath}
    where \(J_m\) is a finite subset of \(\integernumbers^n\) and
    \(a_{m,\firstznvar} \in K\) for each \(\firstznvar \in J_m\)
    for all \(m \in \naturalnumbers\).
    Furthermore,
    \begin{align*}
      \mdimaltpartialdualprojop{n}{j}{\zot} \indexeddualelem{f}{m}
      & = \sum_{\secondznvar \in \integernumbers^n}
      \szdualappl{\sum_{\firstznvar \in J_m} a_{m,\firstznvar}
        \mdimgendualwavelet{n}{\zos}{j}{\firstznvar}}
          {\mdimgenwavelet{n}{\zot}{j}{\secondznvar}}
      \mdimgendualwavelet{n}{\zot}{j}{\secondznvar}
      = \sum_{\secondznvar \in \integernumbers^n}
      \left( \sum_{\firstznvar \in J_m} a_{m,\firstznvar}
        \delta_{\zos,\zot} \delta_{\firstznvar,\secondznvar} \right)
      \mdimgendualwavelet{n}{\zot}{j}{\secondznvar} \\
      & = \delta_{\zos,\zot}
      \sum_{\secondznvar \in \integernumbers^n}
      a_{m,\secondznvar}
      \mdimgendualwavelet{n}{\zot}{j}{\secondznvar}
      = \delta_{\zos,\zot} \indexeddualelem{f}{m} .
    \end{align*}
    Since \(\mdimaltpartialdualprojop{n}{j}{\zot}\) is continuous it
    follows that
    \begin{displaymath}
      \mdimaltpartialdualprojop{n}{j}{\zot} \dualelem{f}
      = \lim_{m \to \infty} \mdimaltpartialdualprojop{n}{j}{\zos}
        \indexeddualelem{f}{m}
      = \delta_{\zos,\zot} \lim_{m \to \infty} \indexeddualelem{f}{m}
      = \delta_{\zos,\zot} \dualelem{f} .
    \end{displaymath}
\end{proof}
}
\fi
}
\else
{
\begin{proof}
  \mbox{ }
  \begin{itemize}
  \item[(i) and (ii)]
    Let \(\zos, \zot \in \zeroonesetn\) and
    \(\dualelem{f} \in \nspace{n}{\zos}{j}\).
    By \mylemma \ref{lem:tensorpartialdualwspace-clos-span}
    there exists a sequence
    \begin{math}
      (\indexeddualelem{f}{m})_{m=0}^\infty \subset
      \spanop \left\{ \mdimgendualwavelet{n}{\zos}{j}{\firstznvar} 
      \setsep
        \firstznvar \in \integernumbers^n \right\}
    \end{math}
    so that \(\indexeddualelem{f}{m} \to \dualelem{f}\) as \(m \to 
    \infty\)
    (strong convergence).
    Now
    \begin{displaymath}
      \indexeddualelem{f}{m} = \sum_{\firstznvar \in J_m} 
      a_{m,\firstznvar}
        \mdimgendualwavelet{n}{\zos}{j}{\firstznvar}
    \end{displaymath}
    where \(J_m\) is a finite subset of \(\integernumbers^n\) and
    \(a_{m,\firstznvar} \in K\) for each \(\firstznvar \in J_m\)
    for all \(m \in \naturalnumbers\).
    Furthermore,
    \begin{align*}
      \mdimaltpartialdualprojop{n}{j}{\zot} \indexeddualelem{f}{m}
      & = \sum_{\secondznvar \in \integernumbers^n}
      \szdualappl{\sum_{\firstznvar \in J_m} a_{m,\firstznvar}
        \mdimgendualwavelet{n}{\zos}{j}{\firstznvar}}
          {\mdimgenwavelet{n}{\zot}{j}{\secondznvar}}
      \mdimgendualwavelet{n}{\zot}{j}{\secondznvar} \\
      & = \sum_{\secondznvar \in \integernumbers^n}
      \left( \sum_{\firstznvar \in J_m} a_{m,\firstznvar}
        \delta_{\zos,\zot} \delta_{\firstznvar,\secondznvar} \right)
      \mdimgendualwavelet{n}{\zot}{j}{\secondznvar} \\
      & = \delta_{\zos,\zot}
      \sum_{\secondznvar \in \integernumbers^n}
      a_{m,\secondznvar}
      \mdimgendualwavelet{n}{\zot}{j}{\secondznvar}
      = \delta_{\zos,\zot} \indexeddualelem{f}{m} .
    \end{align*}
    Since \(\mdimaltpartialdualprojop{n}{j}{\zot}\) is continuous it
    follows that
    \begin{displaymath}
      \mdimaltpartialdualprojop{n}{j}{\zot} \dualelem{f}
      = \lim_{m \to \infty} \mdimaltpartialdualprojop{n}{j}{\zos}
        \indexeddualelem{f}{m}
      = \delta_{\zos,\zot} \lim_{m \to \infty} \indexeddualelem{f}{m}
      = \delta_{\zos,\zot} \dualelem{f} .
    \end{displaymath}
    Thus both (i) and (ii) are true.
  \item[(iii)]
    This is a consequence of (ii).
  \item[(iv)]
    The range of operator
    \(\mdimaltpartialdualprojop{n}{j}{\zos}\)
    is \(\nspace{n}{\zos}{j}\).
    Hence (i) implies (iv).
  \end{itemize}
\end{proof}
}
\fi

\if\shortprep1
{
  Assume that \(j \in \integernumbers\), and
  \(\zos \in \zeroonesetn\).
  Now
  \begin{equation}
    \label{eq:mdimaltpartialdualprojop-series-formula}
    \mdimaltpartialdualprojop{n}{j}{\zos} \dualelem{f}
    = \sum_{\firstznvar \in \integernumbers^n}
    \dualappl{\dualelem{f}}
    {\mdimgenwavelet{n}{\zos}{j}{\firstznvar}}
    \mdimgendualwavelet{n}{\zos}{j}{\firstznvar}
  \end{equation}
  for all \(\dualelem{f} \in \mdimdualvspace{n}{j+1}\).
  The series in \myequation 
  \eqref{eq:mdimaltpartialdualprojop-series-formula}
  converges absolutely for all \(\dualelem{f} \in 
  \mdimdualvspace{n}{j+1}\).
  We also have
  \begin{math}
    \mdimpartialdualwspace{n}{\zos}{j}
    \equalns
    \nspace{n}{\zos}{j}
  \end{math}
  and
  \begin{math}
    \mdimdualpartialdeltaop{n}{\zos}{j}
    =
    \mdimaltpartialdualprojop{n}{j}{\zos}
    \circ
    \mdimdualprojop{n}{j+1}
  \end{math}.
}
\else
{
\begin{theorem}
  \label{th:tensorpartialdualwspace-abs-basis}
  Let \(j \in \integernumbers\) and
  \(\zos \in \zeroonesetn\).
  Then
  \begin{itemize}
  \item[(i)] 
    The set
    \begin{math}
      \{ \mdimgendualwavelet{n}{\zos}{j}{\firstznvar} \setsep
      \firstznvar \in \integernumbers^n \}
    \end{math}
    is an absolutely convergent basis of Banach space
    \(\nspace{n}{\zos}{j}\).
  \item[(ii)] 
    \begin{math}
      \mdimpartialdualwspace{n}{\zos}{j}
      \equalns
      \nspace{n}{\zos}{j}
    \end{math}.
  \end{itemize}
\end{theorem}
}
\fi

\if\shortprep0
{
\begin{proof}
  Since \(\nspace{n}{\zos}{j}\) is a closed 
  subspace of
  \(\dvfn{n}{K}\) it follows that
  \begin{math}
    \mdimpartialdualwspace{n}{\zos}{j} \subset \nspace{n}{\zos}{j}
  \end{math}.
  
  Let \(\dualelem{g} \in \nspace{n}{\zos}{j}\).
  By \mylemmas \ref{lem:mdimaltpartialdualprojop-series-formula}
  and \ref{lem:mdimaltpartialdualprojop-properties}
  \begin{equation}
    \label{eq:tpdw-abs-basis-a}
    \dualelem{g} = \mdimaltpartialdualprojop{n}{j}{\zos}
    \dualelem{g}
    = \sum_{\firstznvar \in \integernumbers^n}     
    \dualappl{\dualelem{g}}
    {\mdimgenwavelet{n}{\zos}{j}{\firstznvar}}
    \mdimgendualwavelet{n}{\zos}{j}{\firstznvar} .
  \end{equation}
  where the series converges absolutely.
  Thus
  \begin{math}
    \{ \mdimgendualwavelet{n}{\zos}{j}{\firstznvar} \setsep
    \firstznvar \in \integernumbers^n \}
  \end{math}
  is an absolutely convergent basis of
  \(\nspace{n}{\zos}{j}\) so (i) is true.
  
  It follows from \mylemma \ref{lem:dual-appl-in-l-one} that
  \begin{math}
    (\dualappl{\dualelem{g}}
    {\mdimgenwavelet{n}{\zos}{j}{\firstznvar}})_{\firstznvar \in
      \integernumbers^n}
    \in l^1(\integernumbers^n, K)
  \end{math}.
  As \(\dualelem{g} \in \nspace{n}{\zos}{j}\) was 
  arbitrary we have
  \begin{math}
    \nspace{n}{\zos}{j} \closedsubspace \mdimpartialdualwspace{n}{\zos}{j}
  \end{math}.
  Hence
  \begin{math}
    \mdimpartialdualwspace{n}{\zos}{j} =
    \nspace{n}{\zos}{j}
  \end{math} and (ii) is true.
\end{proof}
}
\fi


\if\shortprep0
{
\begin{lemma}
  \label{lem:mdimpartialdualprojop-properties}
  Let
  \begin{math}
    j \in \integernumbers
  \end{math}
  and
  \begin{math}
      \zos \in \zeroonesetn
  \end{math}.
  Then
  \begin{math}
    \mdimdualpartialdeltaop{n}{\zos}{j}
    =
    \mdimaltpartialdualprojop{n}{j}{\zos}
    \circ
    \mdimdualprojop{n}{j+1}
  \end{math}.
\end{lemma}
}
\fi

\if\shortprep0
{
\begin{proof}
  Let
  \begin{math}
    \dualelem{f} \in \topdual{E}
  \end{math}.
  By \mylemma \ref{lem:dualdeltaop-properties} (ii) we have
  \begin{math}
    \mdimdualprojop{n}{j+1} \dualelem{f} \in \mdimdualvspace{n}{j+1}
  \end{math}.
  It follows by \mylemma \ref{lem:mdimaltpartialdualprojop-series-formula}
  that
  \begin{displaymath}
    \mdimaltpartialdualprojop{n}{j}{\zos}(
    \mdimdualprojop{n}{j+1}(
      \dualelem{f}
    )
    )
    =
    \sum_{\firstznvar \in \zn}
    \szdualappl{\mdimdualprojop{n}{j+1} \dualelem{f}}
    {\mdimgenwavelet{n}{\zos}{j}{\firstznvar}}
    \mdimgendualwavelet{n}{\zos}{j}{\firstznvar}
    =
    \mdimpartialdualprojop{n}{\zos}{j}(
    \mdimdualprojop{n}{j+1}(
      \dualelem{f}
    )
    ) .   
  \end{displaymath}
  By \mylemma \ref{lem:dualdeltaop-properties} (iii) we get
  \begin{math}
    \mdimaltpartialdualprojop{n}{j}{\zos}(
    \mdimdualprojop{n}{j+1}(
      \dualelem{f}
    )
    )
    =
    \mdimpartialdualprojop{n}{\zos}{j} \dualelem{f}
  \end{math}.
\end{proof}
}
\fi

\section{Besov Space Norm Equivalence}
\label{sec:besov-norm-equivalence}

\if\shortprep0
{
\subsection{Results from Donoho \cite{donoho1992}}

Donoho has derived norm equivalences for the Besov and Triebel-Lizorkin spaces
in the one-dimensional case using the interpolating wavelet 
expansion (in our notation and norming)
\begin{displaymath}
  f = \sum_{k \in \integernumbers} a_k \myphi_{j_0,k} +
  \sum_{j \geq j_0} \sum_{k \in \integernumbers} b_{j,k} \psi_{j,k}
\end{displaymath}
for an arbitrary function \(f \in \vanishingfunccv{\realnumbers}\)
\cite{donoho1992}. Numbers \(R \in \positiverealnumbers\) and
\(D \in \naturalnumbers\) are defined so that the
mother scaling function \(\myphi\) is H\"older continuous of order \(R\)
and the collection of formal sums \(\sum_k a_k \myphi(t - k)\)
contains all polynomials of degree \(D\).
Define
\begin{math}
  \seqstyle{a} = (a_k)_{k \in \integernumbers}
\end{math}
and
\begin{math}
  \mathbf{b}_j = (b_{j,k})_{k \in \integernumbers}
\end{math}
where \(j \in \integernumbers\), \(j \geq j_0\).
The equivalent norm for the Besov space is
\begin{displaymath}
  \norm{f} = 2^{-\frac{j_0}{2}}
  \norm{\seqstyle{a}}_{l^p}
  + \norm{\left(2^{\left(\sigma-\frac{1}{p}\right) j}
      \norm{\mathbf{b}_j}_{l^p} \right)_{j \geq j_0}}_{l^q} ,
\end{displaymath}
where \(\min \{ R, D \} > \sigma > 1/p\) and \(p, q \in ]0, \infty]\).
}
\fi

\subsection{Norm Equivalence for the Besov Spaces in the $n$-dimensional Case}

This derivation is based on the corresponding one-dimensional
derivation in \cite{donoho1992}. The cases \(p < 1\) or \(q < 1\)
yielding quasi-Banach spaces \(B^\sigma_{p,q}\) are not discussed
in this article.
We assume that
\if01
{
\begin{math}
  K = \realnumbers
\end{math}
or
\begin{math}
  K = \complexnumbers
\end{math}
and
}
\fi
\begin{math}
  n \in \positiveintegers
\end{math}
throughout this section.

We give first some definitions similar to those in \cite{meyer1992}
related to orthonormal wavelets.

\begin{definition}
Let \(\bphi : \rn \rightarrow K\) be
a mother scaling function of an orthonormal wavelet family.
When
\begin{math}
  j \in \integernumbers
\end{math},
\begin{math}
  \firstznvar \in \integernumbers^n
\end{math}, and
\begin{math}
  f \in \biglpcv{\infty}{\rn}
\end{math}
define
\begin{displaymath}
  \bar{\alpha}_{j,\firstznvar}(f) := 2^{nj} \int_{\rnx \in \rn}
  \bar{\myphi}^*(2^j \rnx - \firstznvar) f(\rnx) d\tau .
\end{displaymath}
\end{definition}

The spaces
\begin{math}
  \tensorvp{p}{j}
\end{math}
that are defined in \cite{meyer1992}
are denoted by
\begin{math}
  \orthvp{p}{j}
\end{math}
in this document.
Space
\begin{math}
  \orthvp{p}{j}
\end{math}
is a closed subspace of
\begin{math}
  \biglpcv{p}{\rn}
\end{math}
for each
\begin{math}
  j \in \integernumbers
\end{math}
and
\begin{math}
  p \in [ 1, \infty ]
\end{math}
at least when the mother scaling function \(\bphi\) is continuous and
compactly supported.

\begin{definition}
Let \(\bphi : \rn \rightarrow \complexnumbers\) be a scaling
function of an orthonormal wavelet family.
Let \(p \in [ 1, \infty ]\)
and \(j \in \integernumbers\).
Define operator \(\orthprojp{p}{j} : \biglpcv{p}{\rn}
\rightarrow \orthvp{p}{j}\) by
\begin{displaymath}
  (\orthprojp{p}{j} f)(\rnx)
  :=
  \sum_{\firstznvar \in \integernumbers^n}
  \bar{\alpha}_{j,\firstznvar} (f) \bar{\myphi}(2^j \rnx - \firstznvar)
\end{displaymath}
for all
\begin{math}
  f \in \biglpcv{p}{\rn}
\end{math}
and
\begin{math}
  \rnx \in \realnumbers^n  
\end{math}.
Define also
\begin{math}
  \orthdeltaprojp{p}{j} := \orthprojp{p}{j+1} - \orthprojp{p}{j}
\end{math}.
\end{definition}

When \(\bphi\) is a compactly supported and continuous function
operator \(\orthprojp{p}{j}\) is a continuous linear projection
\projprep \(\biglpcv{p}{\rn}\) onto \(\orthvp{p}{j}\) for each
\(j \in \integernumbers\) and \(p \in [ 1, \infty ]\).

\if\shortprep0
{
\begin{definition}
  When \(m \in \naturalnumbers\)
  define
  \begin{displaymath}
    \indexsetequal{n}{m}
    :=
    \left\{
      \alpha \in \naturalnumbers^n
      \setsep
      \normone{\alpha} = m
    \right\} .
  \end{displaymath}
\end{definition}
}
\fi

\if\shortprep0
{
\begin{definition}
  \label{def:holder-coefficient}
  When \(m \in \naturalnumbers\), \(r \in [0, 1[\), and
  \(f \in \contdiffspace{m}{\rn}\)
  define
  \begin{displaymath}
    \holdercoeff{f}{m}{r}
    :=
    \max_{\alpha \in \indexsetequal{n}{m}}
      \sup \left\{ \frac{\abs{\left(\mdimderiv{\alpha}{f}\right)\left(\rnx\right)
        -\left(\mdimderiv{\alpha}{f}\right)\left(\rny\right)}}
      {\norm{\rnx-\rny}^r}
      \setsep
      \rnx, \rny \in \rn \land \rnx \neq \rny \right\} .
  \end{displaymath}
\end{definition}
}
\fi

\if\shortprep0
{
\begin{definition}
Let \(j \in \integernumbers\).
Define
\[
S_j f := \left( f \left( \frac{\firstznvar}{2^j} \right) 
\right)_{\firstznvar \in \zn}
\]
for all functions \(f : \realnumbers^n \rightarrow \complexnumbers\).
\end{definition}
}
\fi

\if\shortprep0
{
\begin{definition}
  Let \(j \in \integernumbers\)
  and \(\rnb \in \rn\).
  Define
  \[
  S_{j,\rnb} f := \left( f \left( \frac{\firstznvar}{2^j}
  + \rnb \right)
  \right)_{\firstznvar \in \zn}
  \]
  for all functions \(f : \realnumbers^n \rightarrow \complexnumbers\).
\end{definition}
}
\fi

\begin{definition}
Let
\begin{math}
  p \in [ 1, \infty ]
\end{math},
\begin{math}
  q \in [ 1, \infty ]
\end{math},
and
\begin{math}
  \sigma \in \positiverealnumbers
\end{math}.
Let \(j_0 \in \integernumbers\)
and \(\br \in \positiveintegers\), \(\br > \sigma\). 
Let \(\bphi : \rn \rightarrow \complexnumbers\) be a mother scaling
function of an
\(\br\)-regular orthonormal MRA \mraprepspace \(\biglpcv{2}{\rn}\) and
\(\orthprojp{p}{j}\) and \(\orthdeltaprojp{p}{j}\),
\(j \in \integernumbers\), be the corresponding projection operators.
When
\(f \in \biglpcv{p}{\rn}\)
define
\begin{displaymath}
  \besovorthwaveletnorm{n}{j_0}{\sigma}{p}{q}{f}
  := \norminspace{\orthprojp{p}{j_0} f}{\biglpcv{p}{\rn}} +
  \norminspace{\bbfh}
  {\littlelpcv{q}{\naturalnumbers + j_0}}
\end{displaymath}
where
\begin{align*}
\bar{h}_j & := 2^{j\sigma}
\norminspace{\orthdeltaprojp{p}{j} f}{\biglpcv{p}{\rn}}, \spaceafter j \in \naturalnumbers + j_0 \\
\bbfh & := (\bar{h}_j)_{j=j_0}^\infty .
\end{align*}
\end{definition}

Norm
\begin{math}
  \besovorthwaveletnorm{n}{j_0}{\sigma}{p}{q}{\cdot}
\end{math}
is an equivalent norm for the Besov space
\begin{math}
  \besovspace{\sigma}{p}{q}{\rn}
\end{math}
and characterizes
\begin{math}
  \besovspace{\sigma}{p}{q}{\rn}
\end{math}
on
\begin{math}
  \borelfunc{\rn}{\complexnumbers}
\end{math}.
This has been proved in
\cite[chapter 2.9 proposition 4]{meyer1992}.

\begin{definition}
Let
\begin{math}
  p \in [ 1, \infty ]
\end{math},
\begin{math}
  q \in [ 1, \infty ]
\end{math},
and
\begin{math}
  \sigma \in \positiverealnumbers
\end{math}.
Let \(j_0 \in \integernumbers\).
Let \(\mdimuprojop{n}{j}\) and \(\mdimudeltaop{n}{j}\), \(j
\in \integernumbers\), be the projection operators belonging to a
compactly supported interpolating tensor product MRA \mraprepspace
\(\ucfunccv{\rn}\). When
\(f \in \ucfunccv{\rn}\) define
\begin{displaymath}
  \besovinterpwaveletnorm{n}{j_0}{\sigma}{p}{q}{f}
  := \norminspace{\mdimuprojop{n}{j_0} f}{\biglpcv{p}{\rn}} +
  \norminspace{\mathbf{h}}{\littlelpcv{q}{\naturalnumbers + j_0}}
\end{displaymath}
where
\begin{align*}
  h_j & := 2^{j\sigma}
  \norminspace{\mdimudeltaop{n}{j} f}{\biglpcv{p}{\rn}},
  \spaceafter j \in \naturalnumbers + j_0 \\
  \mathbf{h} & := (h_j)_{j=j_0}^\infty .
\end{align*}
\end{definition}

\begin{definition}
  Let
  \begin{math}
    p \in [ 1, \infty ]
  \end{math},
  \begin{math}
    q \in [ 1, \infty ]
  \end{math},
  and
  \begin{math}
    \sigma \in \positiverealnumbers
  \end{math}.
  Let \(j_0 \in \integernumbers\).
  Let \(\mdimdualsf{n}{j}{\firstznvar}\) and
  \(\mdimgendualwavelet{n}{\zos}{j}{\firstznvar}\),
  \(j \in \integernumbers\),
  \(\zos \in \zeroonesetnnozero{n}\), \(\firstznvar \in \zn\),
  be the dual
  scaling functions and dual wavelets belonging to a compactly supported
  interpolating tensor product MRA \mraprepspace \(\ucfunccv{\rn}\). 
  When
  \(f \in \ucfunccv{\rn}\)
  define
  \begin{eqnarray*}
    \besovinterpcoeffnorm{n}{j_0}{\sigma}{p}{q}{f}
    & := &      
    \norm{\left(\szdualappl{\mdimdualsf{n}{j_0}{\firstznvar}}
    {f}\right)_{\firstznvar \in \zn}}_p \\
    & & +
    \norm{\left(2^{\left(\sigma - \frac{n}{p}\right) j}
    \norm{\left( 
    \szdualappl{\mdimgendualwavelet{n}{\zos}{j}{\secondznvar}}
    {f}
    \right)_{{\zos \in \zeroonesetnnozero{n}},\secondznvar \in \zn}}_p 
    \right)_{j = j_0}^\infty}_q .
  \end{eqnarray*}
\end{definition}

\begin{definition}
  When \(j \in \integernumbers\) and
  \(p \in [1, \infty]\)
  define
  \begin{displaymath}
    \mdimuvpspace{n}{p}{j} :=
    \left\{
      \rnx \in \rn \mapsto
      \sum_{\firstznvar \in \zn}
      \seqelem{\sqa}{\firstznvar}
      \mdimsf{n}{j}{\firstznvar}(\rnx)
      \setsep
      \sqa \in \littlelpcv{p}{\zn}
    \right\}
  \end{displaymath}
  and
  \begin{math}
    \norminspace{f}{\mdimuvpspace{n}{p}{j}} :=
    \norminspace{f}{\biglpcv{p}{\rn}}
  \end{math}
  for all
  \begin{math}
    f \in \mdimuvpspace{n}{p}{j}
  \end{math}.
\end{definition}

\begin{definition}
  Let \(j \in \integernumbers\) and \(p \in [ 1, \infty ]\).
  When \(T \in \linopautom{\biglpcv{p}{\rn}}\)
  define
  \begin{displaymath}
    \normintensorvp{n}{j}{p}{T} :=
    \norminspace{(\restrictfunc{T}{\mdimuvpspace{n}{p}{j}})}
    {\contlinop{\mdimuvpspace{n}{p}{j}}{\biglpcv{p}{\rn}}} .
  \end{displaymath}
\end{definition}

\begin{definition}
  Let \(j \in \integernumbers\) and \(p \in [ 1, \infty ]\).
  When \(T \in \linopautom{\ucfunccv{\rn}}\) and
  \begin{math}
    \setimage{T}{\orthvp{p}{j}} \subset \biglpcv{p}{\rn}
  \end{math}
  define
  \begin{displaymath}
    \norminorthvp{p}{j}{T} :=
    \norminspace{(\restrictfunc{T}{\orthvp{p}{j}})}
    	{\contlinop{\orthvp{p}{j}}{\biglpcv{p}{\rn}}} .
  \end{displaymath}
\end{definition}

\if\shortprep0
{
\begin{definition}
  When \(j \in \integernumbers\) and
  \(p \in [ 1, \infty]\) define
  \begin{displaymath}
    \mdimuwspace{n}{j}(p) :=
    \left\{
      \rnx \in \rn \mapsto
      \sum_{\zos \in \zeroonesetnnozero{n}}
      \sum_{\firstznvar \in \zn}
      \seqelem{\sqa}{\zos,\firstznvar}
      \mdimgenwavelet{n}{\zos}{j}{\firstznvar}(\rnx)
      \setsep
      \sqa \in
      \szlittlelpcv{p}{\zeroonesetnnozero{n} \times \zn}
    \right\}
  \end{displaymath}
  and
  \begin{math}
    \norminspace{f}{\mdimuwspace{n}{j}(p)}
    :=
    \norminspace{f}{\biglpcv{p}{\rn}}
  \end{math}
  for all
  \begin{math}
    f \in \mdimuwspace{n}{j}(p)
  \end{math}.
  When \(j \in \integernumbers\), \(\zos \in \zeroonesetn\), and
  \(p \in [ 1, \infty]\) define
  \begin{displaymath}
    \mdimupartialwspace{n}{\zos}{j}(p) :=
    \left\{
      \rnx \in \rn \mapsto
      \sum_{\firstznvar \in \zn}
      \seqelem{\sqa}{\firstznvar}
      \mdimgenwavelet{n}{\zos}{j}{\firstznvar}(\rnx)
      \setsep
      \sqa \in
      \szlittlelpcv{p}{\zn}
    \right\}
  \end{displaymath}
  and
  \begin{math}
    \norminspace{f}{\mdimupartialwspace{n}{\zos}{j}(p)}
    :=
    \norminspace{f}{\biglpcv{p}{\rn}}
  \end{math}
  for all
  \begin{math}
    f \in \mdimupartialwspace{n}{\zos}{j}(p)
  \end{math}.
\end{definition}
}
\fi

\begin{definition}
  When \(p \in [1, \infty]\) define
  \begin{displaymath}
    \cinterp{p} :=
    \left\{
      \begin{array}{ll}
      \frac{p}{p-1} ; \spaceafter & p \in ] 1, \infty [ \\
      1 ; \spaceafter & p = 1 \lor p = \infty .
      \end{array}
    \right.
  \end{displaymath}
\end{definition}

\if\shortprep0
{
\givelemmawithoutproof
}
\fi

\if\shortprep0
{
\begin{lemma}
  \label{lem:orth-coeff-formula}
  Let \(\bphi\) be a mother scaling function
  of an orthonormal MRA of \(\biglpcv{2}{\rn}\)
  spanning polynomials of degree 0.
  Let \(f \in \cbfunccv{\rn}\).
  Suppose that \(j \in \integernumbers\)
  and \(\firstznvar \in \zn\).
  Then
  \begin{displaymath}
    \abs{\orthcoefffunctionalvalue{j}{\firstznvar}{f}
      - \left( \frac{\firstznvar}{2^j} \right)}
    \leq
    \normone{\bphi} \szmodcont{f}{2^{-j} \supportradius{\bphi}}
  \end{displaymath}
\end{lemma}
}
\fi

\if\shortprep0
{
  \begin{proof}
    By the polynomial span we have
    \begin{displaymath}
      f \left( \frac{\firstznvar}{2^j} \right) =
      2^{nj} \int_{\rnx \in\rn} \bphi^*(2^j\rnx - \firstznvar) 
      f \left( \frac{\firstznvar}{2^j} \right) .
    \end{displaymath}
    Let
    \begin{math}
      M := \closedball{\rn}{2^{-j} \firstznvar}{2^{-j} \supportradius{\bphi}}
    \end{math}.
    Now
    \begin{eqnarray*}
      \abs{\orthcoefffunctionalvalue{j}{\firstznvar}{f}
        - f \left( \frac{\firstznvar}{2^j} \right)}
      & = &
      \abs{2^{nj} \int_{\rnx \in \rn} \bphi^* \left( 2^j \rnx -
      \firstznvar \right)
      \left( f(\rnx) - f \left( \frac{\firstznvar}{2^j} \right)
      \right)
      d \tau} \\
      & \leq &
      2^{nj} \int_{\rnx \in M}
      \abs{\bphi^* \left( 2^j \rnx -
      \firstznvar \right)}
      \abs{\left( f(\rnx) - f \left( \frac{\firstznvar}{2^j} \right)
        \right)}
      d \tau .
    \end{eqnarray*}
    As
    \begin{displaymath}
      \abs{\left( f(\rnx) - f \left( \frac{\firstznvar}{2^j} \right)
        \right)}
      \leq
      \modcont{f}{2^{-j} \supportradius{\bphi}}
    \end{displaymath}
    for all
    \begin{math}
      \rnx \in M
    \end{math}
    it follows that
    \begin{eqnarray*}
      \abs{\orthcoefffunctionalvalue{j}{\firstznvar}{f}
        - f \left( \frac{\firstznvar}{2^j} \right)}
      & \leq &
      \modcont{f}{2^{-j} \supportradius{\bphi}}
      \cdot
      2^{nj}
      \int_{\rnx \in M}
      \abs{\bphi^* \left( 2^j \rnx -
        \firstznvar \right)}
      d\tau \\
      & = &
      2^{nj} \cdot 2^{-nj} \normone{\bphi^*}
      \modcont{f}{2^{-j} \supportradius{\bphi}} \\
      & = &
      \normone{\bphi^*}
      \modcont{f}{2^{-j} \supportradius{\bphi}} .
    \end{eqnarray*}
  \end{proof}
}
\fi


\if\shortprep0
{
\begin{lemma}
  \label{lem:orth-proj-op-complement-infty-norm}
  Let \(\bphi\) be a mother scaling function
  of an orthonormal MRA of \(\biglpcv{2}{\rn}\)
  spanning polynomials of degree 0.
  Then
  \begin{displaymath}
    \norminfty{f - \bbigp^{(\infty)}_j f} \leq (1 + 
    \normone{\bphi})
    \norminfty{\bphi} \ncover{\bphi} \omega(f; 
    \supportradius{\bphi} \cdot 2^{-j})    
  \end{displaymath}
  for all \(j \in \integernumbers\) and
  \(f \in \cbfunccv{\rn}\).
\end{lemma}
}
\fi

\if\shortprep0
{
\begin{proof}
  As
  \begin{displaymath}
    \sum_{\firstznvar \in \zn} \bphi(\rny - \firstznvar)
    = 1, \spaceafter \rny \in \rn,
  \end{displaymath}
  we have
  \begin{align*}
    \left( f - \bbigp^{(\infty)}_j \right) \left( \rnx \right)
    & =
    \sum_{\firstznvar \in \zn} \left( f(\rnx) -
    \orthcoefffunctionalvalue{j}{\firstznvar}{f} \right)
    \bphi(2^j\rnx-\firstznvar) \\
    & =
    \sum_{\firstznvar \in \itrans{\bphi}{2^j\rnx}} \left( f(\rnx) -
    \orthcoefffunctionalvalue{j}{\firstznvar}{f} \right)
    \bphi(2^j\rnx-\firstznvar)
  \end{align*}
  for all \(\rnx \in \rn\).
  Consequently
  \begin{equation}
    \label{eq:alpha-formula}
    \abs{\left( f - \bbigp^{(\infty)}_j \right) \left( \rnx \right)}
    \leq
    \sum_{\firstznvar \in \itrans{\bphi}{2^j\rnx}}
    \abs{f(\rnx) - \orthcoefffunctionalvalue{j}{\firstznvar}{f}}
    \norminfty{\bphi} 
  \end{equation}
  for all \(\rnx \in \rn\).
  Let \(\rnx_1 \in \rn\) and \(\firstznvar_1 \in \itrans{\bphi}{2^j\rnx_1}\).
  By \mylemma \ref{lem:orth-coeff-formula} we have
  \begin{eqnarray}
    \nonumber
    \abs{f(\rnx_1) - \orthcoefffunctionalvalue{j}{\firstznvar_1}{f}}
    & \leq &
    \abs{f(\rnx_1) - f \left( \frac{\firstznvar_1}{2^j} \right)}
    +
    \abs{f \left( \frac{\firstznvar_1}{2^j} \right) -
      \orthcoefffunctionalvalue{j}{\firstznvar_1}{f}} \\
    \label{eq:supp-formula}
    & \leq &
    \modcont{f}{2^{-j} \supportradius{\bphi}}
    +
    \normone{\bphi} \modcont{f}{2^{-j} \supportradius{\bphi}} .
  \end{eqnarray}
  It follows from \myequations \eqref{eq:alpha-formula} and
  \eqref{eq:supp-formula} that
  \begin{displaymath}
    \abs{\left( f - \bbigp^{(\infty)}_j \right) \left( \rnx_1 \right)}
    \leq
    \ncover{\bphi} \norminfty{\bphi} (1 + \normone{\bphi})
    \modcont{f}{2^{-j} \supportradius{\bphi}} .
  \end{displaymath}
\end{proof}
}
\fi

\if\shortprep0
{
See also Jackson's inequality
\cite[proposition 9.6]{wojtaszczyk1997}.
\givelemmaswithoutproofsn{three}
}
\fi

\begin{lemma}
  \label{lem:off-diagonal-op}
  Suppose that \(q \in [ 1, \infty ] \) and \(t \in \positiverealnumbers\). 
  Let \(a(j,j') := 2^{-\abs{j'-j} t}\) for all \(j, j' \in \naturalnumbers\).
  Define
  \[
  A \mathbf{b} := \left( \sum_{j'=0}^\infty a(j,j') b_{j'}
  \right)_{j \in \naturalnumbers}
  \]
  for all \(\mathbf{b} \in \littlelpcv{q}{\naturalnumbers}\).
  Then \(A \in \contlinopautom{\littlelpcv{q}{\naturalnumbers}}\).
\end{lemma}

\if\longversion1
{
\begin{proof}
  Cases \(q = 1\) and \(q = \infty\) can be proved by
  starting from the expression of
  \begin{math}
    A \mathbf{b}
  \end{math}
  where
  \begin{math}
    \mathbf{b} \in \littlelpcv{q}{\naturalnumbers}
  \end{math}.
  When
  \begin{math}
    q \in ] 1, \infty [
  \end{math}
  the proof is based on interpolation of Banach spaces, see
  \cite{bs1988}.
  Let \(\nu\) be the counting measure on \(\naturalnumbers\).
  Operator \(A\) is admissible with respect to couple
  \((l^1, l^\infty)\).
  As \(l^p = L^p(\naturalnumbers, \nu)\)
  for all \(p \in ]1, \infty[\)
  it follows by \cite[corollary IV.1.8]{bs1988} that
  \(A\) is a continuous linear
  operator from \(l^q\) into \(l^q\).
\end{proof}
}
\else
{
\begin{proof}
  Cases \(q = 1\) and \(q = \infty\) can be proved by
  starting from the expression of
  \begin{math}
    A \mathbf{b}
  \end{math}
  where
  \begin{math}
    \mathbf{b} \in \littlelpcv{q}{\naturalnumbers}
  \end{math}.
  When
  \begin{math}
    q \in ] 1, \infty [
  \end{math}
  the proof is based on interpolation of Banach spaces, see
  \cite{bs1988}.
\end{proof}
}
\fi

\if\shortprep0
{
\begin{lemma}
  \label{lem:orth-coeff-sequence}
  Let
  \begin{displaymath}
    A_j f := (\orthcoefffunctionalvalue{j}{\firstznvar}{f})_{\firstznvar \in \zn}
  \end{displaymath}
  for all
  \(f \in \biglpcv{1}{\rn} \cup \biglpcv{\infty}{\rn}\)
  and \(j \in \integernumbers\).
  Then
  \(A_j f \in \littlelpcv{p}{\zn}\)
  for all
  \(f \in \biglpcv{p}{\rn}\),
  \(p \in [ 1, \infty ]\),
  and \(j \in \integernumbers\).
  Furthermore, there exists constant \(c_1 \in \positiverealnumbers\)
  so that
  \begin{displaymath}   
    \norminspace{A_j}{\contlinop{\biglpcv{p}{\rn}}
    {\littlelpcv{p}{\zn}}}
    \leq
    2^{\frac{nj}{p}} c_1
  \end{displaymath}
  for all
  \(p \in [ 1, \infty ]\)
  and \(j \in \integernumbers\).  
\end{lemma}
}
\fi

\if\shortprep0
{
\begin{lemma}
  \label{lem:sampling}
  Let \(p \in [1, \infty]\),
  \(v \in \cscfunccv{\rn}\),
  and \(\rnb \in \rn\).
  There exists constant \(c_1 \in \positiverealnumbers\)
  so that
  for all \(j \in \integernumbers\)
  and for all functions
  \begin{displaymath}
    f(\rnx) := \sum_{\firstznvar \in \zn}
    \seqelem{\sqa}{\firstznvar} v(2^j \rnx - \firstznvar)
  \end{displaymath}
  where \(\sqa \in \seqset{\zn}{\complexnumbers}\)
  we have
  \begin{displaymath}
    \norm{S_{j,\rnb} f}_p \leq 2^{\frac{nj}{p}} c_1 \norm{f}_p .
  \end{displaymath}
\end{lemma}
}
\fi

\if\shortprep0
{
\begin{lemma}
  \label{lem:sequence-inclusion}
  Suppose that \(v \in \cscfunccv{\rn}\) and
  \(v(\firstznvar) = \delta_{\firstznvar,0}\)
  for all \(\firstznvar \in \zn\).
  Let \(j \in \integernumbers\) and
  \(\sqa \in \seqset{\zn}{\complexnumbers}\).
  Let
  \begin{equation}
    \label{eq:function}
    f(\rnx) := \sum_{\firstznvar \in \zn}
    \seqelem{\sqa}{\firstznvar} v(2^j \rnx - \firstznvar)
  \end{equation}
  for all \(\rnx \in \rn\)
  and assume that series \eqref{eq:function} converges absolutely
  for each \(\rnx \in \rn\).
  Furthermore, assume that \(f \in \biglpcv{p}{\rn}\)
  for some \(p \in [ 1, \infty ]\).
  Then
  \begin{math}
    \sqa \in \littlelpcv{p}{\zn}
  \end{math}.
\end{lemma}
}
\fi

\if\shortprep0
{
\begin{proof}
  Use \mylemma \ref{lem:sampling}.
\end{proof}
}
\fi

\if\shortprep0
{
\givelemmawithoutproof
}
\fi

\if\shortprep0
{
\begin{lemma}
  \label{lem:iota-interp}
  Suppose that
  \begin{math}
    v_k \in \biglpcv{1}{\rn} \intersection \biglpcv{\infty}{\rn}
  \end{math}
  for each \(k \in \naturalnumbers\).
  Define
  \begin{equation}
    \label{eq:iota-a}
    \iota_{p}(\sqa) := \rnx \mapsto \sum_{k \in \naturalnumbers}
    \seqelem{\sqa}{k} v_k(\rnx)
  \end{equation}
  for all \(\sqa \in l^p\) and \(p \in [1, \infty]\).
  When \(p = 1\) or \(p = \infty\)
  assume that the series in \eqref{eq:iota-a}
  converges absolutely for all \(\sqa \in l^p\)
  and \(\rnx \in \rn\).
  Then \(\iota_{p}\) is an operator
  from \(l^p\) into \(\biglpcv{p}{\rn}\)
  for each \(p \in ] 1, \infty [\).
\end{lemma}
}
\fi

\if\shortprep0
{
\begin{lemma}
  \label{lem:gen-interp}
  Assume that the following conditions are true.
  \begin{itemize}
    \item[(A1)]
  	  \(I\) is a countably infinite set.
  	\item[(A2)]
  	  \begin{math}
  	    \forall \alpha \in I, j \in \integernumbers : v_{j,\alpha} \in \cscfunccv{\rn}
  	  \end{math}
  	\item[(A3)]
  	  Either of the following conditions is true:
  	  \begin{itemize}
  	    \item[(A3.1)]
  	      We have
  	      \(m(\alpha) \in \positiveintegers\),
  	      \(r_{j,\alpha,k} \in \complexnumbers\),
  	      and
  	      \(\rns_{j,\alpha,k} \in \rn\)
  	      for all
  	      \(\alpha \in I\),
  	      \(j \in \integernumbers\),
  	      and
  	      \(k \in \setoneton{m(\alpha)}\).
  	      Furthermore,
  	      \begin{displaymath}
  	        \dualelem{v}_{j,\alpha}
  	        =
  	        \sum_{k=1}^{m(\alpha)}
  	        r_{j,\alpha,k}
  	        \delta\left( \cdot - \rns_{j,\alpha,k} \right)
  	      \end{displaymath}
  	      for all
  	      \(\alpha \in I\)
  	      and
  	      \(j \in \integernumbers\).
  	    \item[(A3.2)]
  	      We have
  	      \(\dualelem{v}_{j,\alpha} \in \cscfunccv{\rn}\)
  	      for all
  	      \(\alpha \in I\)
  	      and
  	      \(j \in \integernumbers\).
  	      The set
  	      \begin{math}
  	        \left\{ \beta \in I \setsep \suppop \dualelem{v}_{j,\alpha}
  	          \intersection \suppop v_{j,\beta} \right\}
  	      \end{math}
  	      is finite for each
  	      \(\alpha \in I\)
  	      and
  	      \(j \in \integernumbers\).
  	  \end{itemize}
  	\item[(A4)]
  	  \begin{math}
  	    \forall j \in \integernumbers, \alpha \in I, \beta \in I :
  	    \dualappl{\dualelem{v}_{j,\alpha}}{v_{j,\beta}} = \delta_{\alpha,\beta}
		  \end{math}
		\item[(A5)]
		  The series
		  \begin{displaymath}
		    \sum_{\alpha \in I} \seqelem{\sqa}{\alpha} v_{j,\alpha}(\rnx)
		  \end{displaymath}
		  converges absolutely for each
		  \(\sqa \in \littlelpcv{\infty}{I}\),
		  \(j \in \integernumbers\),
		  and \(\rnx \in \rn\).
  \end{itemize}
  Define
  \begin{displaymath}
    \left( \iota_{p,j} \left(\sqa\right) \right) \left( \rnx \right)
    := \sum_{\alpha \in I} \seqelem{\sqa}{\alpha}
     v_{j,\alpha}\left(\rnx\right)
  \end{displaymath}
  for all
  \(\sqa \in \littlelpcv{p}{I}\),
  \(j \in \integernumbers\),
  \(\rnx \in \rn\), and
  \(p \in \{1, \infty\}\).
  Define also
  \begin{math}
    A_{p,j} := \setimage{\iota_{p,j}}{\littlelpcv{p}{I}}
  \end{math}
  for all \(p \in \{1,\infty\}\),
  \(j \in \integernumbers\)
  and
  \begin{math}
    \norminspace{f}{A_{p,j}} := \norm{f}_p
  \end{math}
  for all \(p \in \{1,\infty\}\),
  \(j \in \integernumbers\),
  and \(f \in A_{p,j}\).
  Assume further that the following conditions are true.
  \begin{itemize}
		\item[(A6)]
		  \(\iota_{1.j}\) is a topological isomorphism from
		  \(\littlelpcv{1}{I}\) onto \(A_{1,j}\)
		  for each \(j \in \integernumbers\).
		\item[(A7)]
		  \(\iota_{\infty,j}\) is a topological isomorphism from
		  \(\littlelpcv{\infty}{I}\) onto \(A_{\infty,j}\)
		  for each \(j \in \integernumbers\).
		\item[(A8)]
		  \(A_{1,j} \closedsubspace \biglpcv{1}{\rn}\)
		  and
		  \(A_{\infty,j} \closedsubspace \biglpcv{\infty}{\rn}\)
		  for each \(j \in \integernumbers\).
		\item[(A9)]
		  \begin{math}
		    \existsdef c_1 \in \positiverealnumbers :
		    \forall j \in \integernumbers :
		    \norm{\iota_{1,j}} \leq 2^{-nj} c_1
		  \end{math}
		\item[(A10)]
		  \begin{math}
		    \existsdef c_2 \in \positiverealnumbers :
		    \forall j \in \integernumbers :
		    \norm{(\iota_{1,j})^{-1}} \leq 2^{nj} c_2
		  \end{math}
		\item[(A11)]
		  \begin{math}
		    \existsdef c_3 \in \positiverealnumbers :
		    \forall j \in \integernumbers :
		    \norm{\iota_{\infty,j}} \leq c_3
		  \end{math}
		\item[(A12)]
		  \begin{math}
		    \existsdef c_4 \in \positiverealnumbers :
		    \forall j \in \integernumbers :
		    \norm{(\iota_{\infty,j})^{-1}} \leq c_4
		  \end{math}
  \end{itemize}
  Then
  \begin{itemize}
    \item[(i)]
      Set
      \begin{displaymath}
        \left( \iota_{p,j} \left(\sqa\right) \right) \left( \rnx \right)
        := \sum_{\alpha \in I} \seqelem{\sqa}{\alpha} v_{j,\alpha}\left(\rnx\right)
      \end{displaymath}
      for all
      \(\sqa \in \littlelpcv{p}{I}\),
      \(j \in \integernumbers\),
      \(\rnx \in \rn\), and
      \(p \in ]1, \infty[\).
      Then
      \(\iota_{p,j}\) is an operator
      from \(\littlelpcv{p}{I}\) into
      \(\biglpcv{p}{\rn}\).
    \item[(ii)]
      Define
      \begin{math}
        A_{p,j} := \setimage{\iota_{p,j}}{\littlelpcv{p}{I}}
      \end{math}
      for all \(p \in ]1,\infty[\) and
      \(j \in \integernumbers\).
      Now
      \begin{displaymath}
        A_{p,j} \equaltvs
        \kinterp{1-\frac{1}{p}}{p}{A_{1,j}}{A_{\infty,j}}
      \end{displaymath}
      for all \(p \in ]1,\infty[\) and
      \(j \in \integernumbers\).
    \item[(iii)]
      Function \(\iota_{p,j}\) is a topological isomorphism
      from \(\littlelpcv{p}{I}\) onto \(A_{p,j}\)
      for each
      \(p \in [ 1, \infty ]\)
      and \(j \in \integernumbers\).
    \item[(iv)]
      \begin{math}
        \exists c_5 \in \positiverealnumbers :
        \forall p \in ]1, \infty], j \in \integernumbers :
        \norm{\iota_{p,j}} \leq
        \frac{p}{p-1}
        2^{-\frac{nj}{p}} c_5
      \end{math}
    \item[(v)]
      \begin{math}
        \exists c_6 \in \positiverealnumbers :
        \forall p \in ]1, \infty], j \in \integernumbers :
        \norm{\left(\iota_{p,j}\right)^{-1}} \leq
        \frac{p}{p-1}
        2^{\frac{nj}{p}} c_6
      \end{math}.
  \end{itemize}
\end{lemma}
}
\fi

\if\shortprep0
{
\begin{proof}
  \mbox{ }
  \begin{itemize}
  \item[(i)]
    This is a consequence of \mylemma \ref{lem:iota-interp}.
  \item[(ii)]
    Let \(p_1 \in ]1, \infty[\) and
    \(j_1 \in \integernumbers\).
    Let
    \begin{displaymath}
      \eta :=
      \szkinterp{1 - \frac{1}{p_1}}{p_1}{\iota_{1,j_1}}{\iota_{\infty,j_1}}
    \end{displaymath}
    
    Suppose first that \(f \in A_{p_1,j_1}\).
    Now \(f = \iota_{p_1,j_1}(\sqa)\) for some
    \(\sqa \in \littlelpcv{p}{I}\).
    As
    \begin{displaymath}
      l^{p_1} \equaltvs
      \szkinterp{1-\frac{1}{p_1}}{p_1}{l^1}{l^\infty}
    \end{displaymath}
    we have
    \(\sqa = \sqb + \sqc\)
    for some
    \(\sqb \in \littlelpcv{1}{I}\)
    and
    \(\sqc \in \littlelpcv{\infty}{I}\).
    It follows that
    \begin{displaymath}
      f = \iota_{1,j_1}(\sqb) + \iota_{\infty,j_1}(\sqc)
      = \eta(\sqa)
      \in
      \szkinterp{1-\frac{1}{p_1}}{p_1}{A_{1,j_1}}{A_{\infty,j_1}} .
    \end{displaymath}
    
    Suppose then that
    \begin{displaymath}
      g \in \szkinterp{1-\frac{1}{p_1}}{p_1}{A_{1,j_1}}{A_{\infty,j_1}} .
    \end{displaymath}
    Now  \(g = w + z\) for some \(w \in A_{1,j_1}\) and
    \(z \in A_{\infty,j_1}\).
    Furthermore,
    \(w = \iota_{1,j_1}(\seqstyler)\)
    for some
    \(\seqstyler \in \littlelpcv{1}{I}\)
    and
    \(z = \iota_{\infty,j_1}(\sqs)\)
    for some
    \(\sqs \in \littlelpcv{\infty}{I}\).
    We have
    \begin{displaymath}
      g(\rnx)
      =
      \sum_{\alpha \in I}
      \left(
        \seqelem{\seqstyler}{\alpha}
        +
        \seqelem{\sqs}{\alpha}
      \right)
      v_{j_1,\alpha}\left(\rnx\right)
    \end{displaymath}
    for all \(\rnx \in \rn\) and it follows that
    \(g = \eta(\seqstyler + \sqs)\).
    Hence
    \begin{displaymath}
      \eta^{-1}\left( g \right) = \seqstyler + \sqs
      \in
      \szkinterp{1-\frac{1}{p_1}}{p_1}
      {\littlelpcv{1}{I}}{\littlelpcv{\infty}{I}}
    \end{displaymath}
    from which it follows that
    \(\seqstyler + \sqs \in \littlelpcv{p}{I}\).
    Thus
    \(g = \iota_{p_1,j_1}(\seqstyler + \sqs) \in A_{p_1,j_1}\).
    
  \item[(iii)]
  Let \(\sqb \in \littlelpcv{p}{I} \setminus \zeroset\).
  There exists \(\alpha_1 \in I\) so that
  \(\seqelem{\sqb}{\alpha_1} \neq 0\).
  If (A3.1) is true we have
  \begin{math}
    \dualappl{\dualelem{v}_{j,\alpha_1}}{\iota_{p_1,j_1}(\sqb)}
    = \seqelem{\sqb}{\alpha_1} \neq 0
  \end{math}.
  If (A3.2) is true we have
  \begin{displaymath}
    J := \{ \beta \in I \setsep
    \suppop \dualelem{v}_{j_1,\alpha_1}
    \intersection
    \suppop v_{j_1,\beta}
    \neq
    \zeroset \}
  \end{displaymath}
  and it follows that
  \begin{eqnarray*}
    \szdualappl{\dualelem{v}_{j_1,\alpha_1}}
    {\iota_{p_1,j_1}(\sqb)}
    & = &
    \int_{\rnx \in \rn}
    \dualelem{v}_{j_1,\alpha_1}\left(\rnx\right)
    \left(
      \sum_{\beta \in J} \seqelem{\sqb}{\beta}
      v_{j,\beta}
      \left( \rnx \right)
    \right)
    d \tau \\
    & = &
    \sum_{\beta \in J} \seqelem{\sqb}{\beta}
    \szdualappl{\dualelem{v}_{j_1,\alpha_1}}{v_{j,\beta}} \\
    & = &
    \seqelem{\sqb}{\alpha_1} \neq 0 .
  \end{eqnarray*}
  Consequently
  \(\iota_{p_1,j_1}(\sqb) \neq 0\).
  Hence \(\iota_{p_1,j_1}\) is an injection.

  By \mylemma \ref{lem:iota-interp} function \(\iota_{p_1,j_1}\)
  is continuous. By the Inverse Mapping Theorem
  function \(\iota_{p_1,j_1}^{-1}\) is continuous.
  Thus \(\iota_{p_1,j_1}\) is a topological isomorphism
  from Banach space \(\littlelpcv{p_1}{I}\) onto Banach space
  \(A_{p_1,j_1}\).
  
  \item[(iv)]
  Let \(p_1 \in ] 1, \infty [\) and \(j_1 \in \integernumbers\).
  Now by \mylemma \ref{lem:lp-oper-interp}
  \begin{displaymath}
    \norm{\iota_{p_1,j_1}}
    \leq
    \frac{p}{p-1}
    \norm{\iota_{1,j_1}}^{\frac{1}{p_1}}
    \norm{\iota_{\infty,j_1}}^{1-\frac{1}{p_1}}
    \leq c_5 \frac{p}{p-1} \cdot 2^{-\frac{nj_1}{p_1}}
  \end{displaymath}
  where
  \begin{math}
    c_5 := c_1^{\frac{1}{p_1}} c_3^{1-\frac{1}{p_1}}
  \end{math}.

  \item[(v)]
  Let \(p_1 \in ] 1, \infty [\) and \(j_1 \in \integernumbers\).
  Now by \mylemma \ref{lem:lp-oper-interp}
  \begin{displaymath}
    \norm{(\iota_{p_1,j_1})^{-1}}
    \leq
    \frac{p}{p-1}
    \norm{(\iota_{1,j_1})^{-1}}^{\frac{1}{p_1}}
    \norm{(\iota_{\infty,j_1})^{-1}}^{1-\frac{1}{p_1}}
    \leq c_6 \frac{p}{p-1} \cdot 2^{\frac{nj_1}{p_1}}
  \end{displaymath}
  where
  \begin{math}
    c_6 := c_2^{\frac{1}{p_1}} c_4^{1-\frac{1}{p_1}}
  \end{math}.
  \end{itemize}
\end{proof}
}
\fi

\if\shortprep0
{
\begin{lemma}
  \label{lem:lone-top-isomorphism}
  Let
  \begin{math}
    f \in \cscfunccv{\rn}
  \end{math},
  \begin{math}
    \sqb \in \rn
  \end{math}
  and suppose that
  \begin{math}
    f(\firstznvar + \sqb) = \delta_{\firstznvar,0}
  \end{math}
  for all \(\firstznvar \in \zn\).
  Define
  \begin{displaymath}
    \left(\iota_j \left(\sqa\right)\right)\left(\rnx\right)
    :=
    \sum_{\firstznvar \in \zn} \seqelem{\sqa}{\firstznvar}
    f(2^j \rnx - \firstznvar)
  \end{displaymath}
  for all
  \begin{math}
    \sqa \in \littlelpcv{1}{\zn}
  \end{math}
  and
  \begin{math}
    j \in \integernumbers
  \end{math}.
  Define also
  \begin{displaymath}
    A_j := \left\{ \iota_j(\sqa) \setsep \sqa \in
    \littlelpcv{1}{\zn}\right\}
  \end{displaymath}
  for all
  \begin{math}
    j \in \integernumbers
  \end{math}
  and
  \begin{math}
    \norminspace{g}{A_j} := \norminspace{g}{\biglpcv{1}{\rn}}
  \end{math}
  for all
  \begin{math}
    g \in A_j
  \end{math}
  and
  \begin{math}
    j \in \integernumbers
  \end{math}.
  Then
  \begin{itemize}
    \item[(i)]
      \(\iota_j\) is a topological isomorphism from
      \(\littlelpcv{1}{\zn}\) onto \(A_j\) for each
      \(j \in \integernumbers\).
    \item[(ii)]
      \begin{math}
        \exists c_1 \in \positiverealnumbers :
        \forall j \in \integernumbers :
        \norm{\iota_j} \leq c_1 \cdot 2^{-nj}
      \end{math}
    \item[(iii)]
      \begin{math}
        \exists c_2 \in \positiverealnumbers :
        \forall j \in \integernumbers :
        \norm{\left(\iota_j\right)^{-1}} \leq c_2 \cdot 2^{nj}
      \end{math}
  \end{itemize}
\end{lemma}
}
\fi

\if\shortprep0
{
\begin{proof}
  Use \mylemmas \ref{lem:change-int-and-sum-order}
  and \ref{lem:sampling}.
\end{proof}
}
\fi

\if\shortprep0
{
\begin{lemma}
  \label{lem:mdimpartialuwspace-interp}
  Let
  \begin{displaymath}
    \left( \iota_{p,\zos,j} \left( \sqa \right) \right)
    \left( \rnx \right)
    :=
    \sum_{\firstznvar \in \zn}
    \seqelem{\sqa}{\firstznvar}
    \mdimgenwavelet{n}{\zos}{j}{\firstznvar}(\rnx)
  \end{displaymath}
  for all
  \(p \in [1, \infty]\),
  \(\zos \in \zeroonesetn\),
  \(j \in \integernumbers\),
  \(\sqa \in \littlelpcv{p}{\zn}\),
  and \(\rnx \in \rn\).
  Then
  \begin{itemize}
    \item[(i)]
    We have
    \begin{displaymath}
    \mdimupartialwspace{n}{\zos}{j}(p)
    \equaltvs
    \left( \mdimupartialwspace{n}{\zos}{j}(1),
      \mdimupartialwspace{n}{\zos}{j}(\infty)
    \right)
    \end{displaymath}
    for all
    \(p \in ] 1, \infty [\),
    \(\zos \in \zeroonesetn\), and
    \(j \in \integernumbers\).
    \item[(ii)]
      Function \(\iota_{p,\zos,j}\) is a topological
      isomorphism from \(\littlelpcv{p}{\zn}\)
      onto \(\mdimupartialwspace{n}{\zos}{j}(p)\)
      for all
      \(p \in ] 1, \infty [\),
      \(\zos \in \zeroonesetn\), and
      \(j \in \integernumbers\).
    \item[(iii)]
      \begin{math}
        \exists c_1 \in \positiverealnumbers :
        \forall p \in [ 1, \infty],
        j \in \integernumbers,
        \zos \in \zeroonesetn :
        \norm{\iota_{p,\zos,j}}
        \leq
        c_1 \cdot \cinterp{p} \cdot 2^{-\frac{nj}{p}}
      \end{math}
    \item[(iv)]
      \begin{math}
        \exists c_2 \in \positiverealnumbers :
        \forall p \in [ 1, \infty],
        j \in \integernumbers,
        \zos \in \zeroonesetn :
        \norm{\left(\iota_{p,\zos,j}\right)^{-1}}
        \leq
        c_2 \cdot \cinterp{p} \cdot 2^{\frac{nj}{p}}
      \end{math}
  \end{itemize}
\end{lemma}
}
\fi

\if\shortprep0
{
\begin{proof}
  Use \mylemmas  
  \ref{lem:linftyspace-top-isomorphism},
  \ref{lem:uc-convergence},
  \ref{lem:gen-interp}, and
  \ref{lem:lone-top-isomorphism}.
\end{proof}
}
\fi

\if\shortprep0
{
\begin{lemma}
  \label{lem:orthvp-interp}
  Let \(\orthvp{p}{j}\) belong to an orthonormal
  MRA \mraprepdomain \(\rn\).
  Define functions \(\iota_{p,j}\) by
  \begin{displaymath}
    \left( \iota_{p,j} \left( \sqa \right) \right)
    \left( \rnx \right)
    :=
    \sum_{\firstznvar \in \zn}
    \seqelem{\sqa}{\firstznvar}
    \bphi_{j,\firstznvar} \left( \rnx \right)    
  \end{displaymath}
  for all
  \(\sqa \in \littlelpcv{p}{\zn}\),
  \(\rnx \in \rn\),
  \(p \in [1, \infty]\), and
  \(j \in \integernumbers\).
  Then
  \begin{itemize}
    \item[(i)]
      \begin{math}
        \forall j \in \integernumbers, p \in ] 1, \infty [ :
        \orthvp{p}{j} \equaltvs
        \left( \orthvp{1}{j}, \orthvp{\infty}{j}
        \right)_{1-\frac{1}{p},p}
      \end{math}
    \item[(ii)]
      Function \(\iota_{p,j}\) is a topological isomorphism
      from \(\littlelpcv{p}{\zn}\) onto \(\orthvp{p}{j}\)
      for each \(j \in \integernumbers\) and
      \(p \in [1, \infty]\).
    \item[(iii)]
      \begin{math}
        \exists c_1 \in \positiverealnumbers :
        \forall j \in \integernumbers, p \in [ 1, \infty ] :
        \norm{\iota_{p,j}} \leq
        c_1 \cdot \cinterp{p} \cdot 2^{-\frac{nj}{p}}
      \end{math}
    \item[(iv)]
      \begin{math}
        \exists c_2 \in \positiverealnumbers :
        \forall j \in \integernumbers, p \in [ 1, \infty ] :
        \norm{\left(\iota_{p,j}\right)^{-1}}
        \leq c_2 \cdot \cinterp{p} \cdot 2^{\frac{nj}{p}}
      \end{math}
  \end{itemize}
\end{lemma}
}
\fi

\if\shortprep0
{
\begin{proof}
  Use \mylemmas \ref{lem:change-int-and-sum-order},
  \ref{lem:uc-convergence}, and \ref{lem:gen-interp}.
\end{proof}
}
\fi

\if\shortprep0
{
\begin{lemma}
  \label{lem:mdimuwspace-interp}
  Define functions \(\iota_{p,j}\) by
  \begin{displaymath}
    \left( \iota_{p,j} \left( \sqa \right) \right)
    \left( \rnx \right)
    :=
    \sum_{\zos \in \zeroonesetnnozero{n}}
    \sum_{\firstznvar \in \zn}
    \seqelem{\sqa}{\zos,\firstznvar}
    \mdimgenwavelet{n}{\zos}{j}{\firstznvar} \left( \rnx \right)    
  \end{displaymath}
  for all
  \(\sqa \in \littlelpcv{p}{\zeroonesetnnozero{n} \times \zn}\),
  \(\rnx \in \rn\),
  \(p \in [1, \infty]\), and
  \(j \in \integernumbers\).
  Then
  \begin{itemize}
    \item[(i)]
      \begin{math}
        \forall j \in \integernumbers, p \in ] 1, \infty [ :
        \left( \mdimuwspace{n}{j}(1), \mdimuwspace{n}{j}(\infty)
        \right)_{1-\frac{1}{p},p}
        \equaltvs
        \mdimuwspace{n}{j}(p)
        \closedsubspace
        \biglpcv{p}{\rn}
      \end{math}
    \item[(ii)]
      Function \(\iota_{p,j}\) is a topological isomorphism
      from \(\littlelpcv{p}{\zeroonesetnnozero{n} \times \zn}\) 
      onto \(\mdimuwspace{n}{j}(p)\)
      for each \(j \in \integernumbers\) and
      \(p \in [1, \infty]\).
    \item[(iii)]
      \begin{math}
        \exists c_1 \in \positiverealnumbers :
        \forall j \in \integernumbers, p \in [ 1, \infty ] :
        \norm{\iota_{p,j}}
        \leq
        c_1 \cdot \cinterp{p} \cdot 2^{-\frac{nj}{p}}
      \end{math}
    \item[(iv)]
      \begin{math}
        \exists c_2 \in \positiverealnumbers :
        \forall j \in \integernumbers, p \in [ 1, \infty ] :
        \norm{\left(\iota_{p,j}\right)^{-1}}
        \leq
        c_2 \cdot \cinterp{p} \cdot 2^{\frac{nj}{p}}
      \end{math}
  \end{itemize}
\end{lemma}
}
\fi

\if\shortprep0
{
\begin{proof}
  Use \mylemmas
  \ref{lem:uw-top-isomorphism},
  \ref{lem:sampling}, and
  \ref{lem:gen-interp}
  and
  \myequation \eqref{eq:mdim-dual-wavelet-formula}.
\end{proof}
}
\fi

\if\shortprep0
{
\begin{lemma}
  \label{lem:uprojop-in-orthvp}
  Let \(p \in [1, \infty]\).
  Then
  \begin{displaymath}
    \exists c_1 \in \positiverealnumbers :
    \forall j \in \integernumbers, j' \in \integernumbers :
    j' \geq j \implies
    \norminorthvp{p}{j'}{\mdimuprojop{n}{j}}
    \leq
    c_1 \cdot 2^{n(j'-j)/p} .
  \end{displaymath}
\end{lemma}
}
\fi

\if\shortprep1
{
\if10
{
\begin{proof}
  Use \mylemmas
  \ref{lem:change-int-and-sum-order} and
  \ref{lem:lp-oper-interp}.
\end{proof}
}
\fi
}
\else
{
\begin{proof}
  Use \mylemmas
  \ref{lem:change-int-and-sum-order},
  \ref{lem:lp-oper-interp},
  \ref{lem:sampling}, and
  \ref{lem:orthvp-interp}.
\end{proof}
}
\fi

\if\shortprep0
{
We have \cite{meyer1992}
\begin{equation}
  \label{eq:orthprojop-uniform-bound}
  \forall p \in [1, \infty] :
  \exists c_1 \in \positiverealnumbers :
  \forall j \in \integernumbers :
  \norm{\orthprojp{p}{j}} \leq c_1 .
\end{equation}
}
\fi

\if\shortprep0
{
\begin{lemma}
  \label{lem:uprojop-inclusion}
  Let \(j \in \integernumbers\),
  \(s \in \positiverealnumbers\),
  \(p \in [1, \infty]\), and
  \(q \in [1, \infty]\).
  Then
  \begin{math}
    \mdimuprojop{n}{j} f \in \mdimuvpspace{n}{p}{j}
  \end{math}
  for all
  \begin{math}
    f \in \besovspace{s}{p}{q}{\rn} \intersection \ucfunccv{\rn}
  \end{math}.
\end{lemma}
}
\fi

\if\shortprep0
{
\begin{proof}
  We have
  \begin{displaymath}
    f = \orthprojp{p}{j} f
    + \sum_{j'=j}^\infty \orthdeltaprojp{p}{j'} f
    = \orthprojp{\infty}{j} f
    + \sum_{j'=j}^\infty \orthdeltaprojp{\infty}{j'} f
  \end{displaymath}
  and
  \begin{displaymath}
    \norm{\mdimuprojop{n}{j} f}_p
    \leq
    \norm{\mdimuprojop{n}{j} \orthprojp{p}{j} f}_p
    + \sum_{j'=j}^\infty
    \norm{\mdimuprojop{n}{j} \orthdeltaprojp{p}{j'} f}_p
    \in \nonnegrealnumbers .
  \end{displaymath}
  Consequently
  \begin{math}
    \mdimuprojop{n}{j} f \in \biglpcv{p}{\rn}
  \end{math}.
  By \mylemma \ref{lem:sequence-inclusion} we have
  \begin{displaymath}
    \left( f \left( \frac{\firstznvar}{2^j} \right)
    \right)_{\firstznvar \in \zn}
    \in \littlelpcv{p}{\zn} .
  \end{displaymath}
  Thus
  \begin{math}
    \mdimuprojop{n}{j} f \in \mdimuvpspace{n}{p}{j}
  \end{math}.  
\end{proof}
}
\fi

\if\shortprep0
{
\givelemmaswithoutproofsn{two}
}
\fi

\if\shortprep0
{
\begin{lemma}
  \label{lem:holder-coeff}
  Let \(m \in \naturalnumbers\),
  \(a \in ]0, 1]\), and
  \(v \in \cscfunccv{\rn}\).
  Then there exists \(c_1 \in \positiverealnumbers\)
  so that for all \(\sqb \in \littlelpcv{\infty}{\zn}\),
  \(j \in \integernumbers\), and
  \begin{displaymath}
    f = \rnx \in \rn \mapsto
    \sum_{\firstznvar \in \zn}
    \seqelem{\sqb}{\firstznvar}
    v(2^j \rnx - \firstznvar)
  \end{displaymath}
  we have
  \begin{displaymath}
    \holdercoeff{f}{m}{a}
    \leq
    c_1 \cdot
    2^{j(m+a)}
    \holdercoeff{v}{m}{a}
    \norminfty{\sqb} .
  \end{displaymath}
\end{lemma}
}
\fi

\if\shortprep0
{
\begin{lemma}
  \label{lem:taylor-residual}
  Let \(s \in \positiverealnumbers \setminus
  \positiveintegers\), \(s > 1\), and
  \(f \in \holderspace{s}{\rn}\).
  Let \(m := \floor{s}\) and \(a := s - m\).
  Let \(\indexedseqstyle{x}{0} \in \rn\), \(p_m\)
  be the \(m\)th degree Taylor polynomial at point
  \(\indexedseqstyle{x}{0}\), and \(r_m := f - p_m\).
  Then
  \begin{displaymath}
    \abs{r_m(\rnx)}
    \leq
    \frac{n}{(m-1)!}
    \holdercoeff{f}{m}{a}
    \norm{\rnx - \indexedseqstyle{x}{0}}^s
  \end{displaymath}
  for all \(\rnx \in \rn\).
\end{lemma}
}
\fi

\if\shortprep0
{
See also \cite[chapter 10]{grossman} and \cite[chapter V.5]{bs1988}.
}
\fi

\if\shortprep0
{
\begin{lemma}
  \label{lem:i-minus-pj-bound}
  Assume that the following conditions are true.
  \begin{itemize}
    \item[(B1)]
      \(s \in \positiverealnumbers\)
      and
      \(v \in \cscfunccv{\rn} \intersection
      \holderspace{s}{\rn}\)
    \item[(B2)]
      We have \(\dualw \in \topdual{\ucfunccv{\rn}}\)
      and \(w \in \cscfunccv{\rn}\).
    \item[(B3)]
      Either of the following conditions is true:
      \begin{itemize}
        \item[(B3.1)]
          We have
          \begin{displaymath}
            \dualw = \sum_{k=1}^m r_k \delta \left( \cdot -
            \indexeddualelem{s}{k} \right)
          \end{displaymath}
          where \(m \in \positiveintegers\) and
          \(r_k \in \complexnumbers\).
        \item[(B3.2)]
          \begin{math}
            \dualw \in \cscfunccv{\rn}
          \end{math}.
      \end{itemize}
    \item[(B4)]
      We have
      \begin{eqnarray*}
        \suppop v & \subset & \closedball{\rn}{0}{\rho_1} ,
        \spaceafter \rho_1 \in \positiverealnumbers , \\
        \suppop w & \subset & \closedball{\rn}{0}{\rho_2} ,
        \spaceafter \rho_2 \in \positiverealnumbers , \\
        \suppop \dualw & \subset & \closedball{\rn}{0}{\rho_3} ,
        \spaceafter \rho_3 \in \positiverealnumbers .
      \end{eqnarray*}
    \item[(B5)]
      We have
      \begin{eqnarray*}
        v_{j,\firstznvar} & = & v(2^j \cdot - \firstznvar) \\
        w_{j,\firstznvar} & = & w(2^j \cdot - \firstznvar) \\
        \indexeddualelem{w}{j,\firstznvar} & = &
        2^{nj} \dualw (2^j \cdot - \firstznvar)
      \end{eqnarray*}
      where
      \(\forall j \in \integernumbers\)
      and
      \(\firstznvar \in \zn\).
    \item[(B6)]
      \begin{math}
        \forall \firstznvar \in \zn, \secondznvar \in \zn,
        j \in \integernumbers :
        \szdualappl{\indexeddualelem{w}{j,\firstznvar}}
        {w_{j,\secondznvar}}
        = \delta_{\firstznvar,\secondznvar}
      \end{math}.
    \item[(B7)]
      Define
      \begin{math}
        M_{j,\firstznvar} :=
        \closedball{\rn}{2^{-j} \firstznvar}
        {2^{-j} \rho_3}
      \end{math}
      for all \(j \in \integernumbers\) and
      \(\firstznvar \in \zn\).
      We have
      \begin{math}
        \indexeddualelem{z}{j,\firstznvar}
        \in
        \topdual{\contfunccv{M_{j,\firstznvar}}}
      \end{math}
      for all \(j \in \integernumbers\),
      \(\firstznvar \in \zn\),
      and
      \begin{displaymath}
        \szdualappl{\indexeddualelem{w}{j,\firstznvar}}{f}
        =
        \szdualappl{\indexeddualelem{z}{j,\firstznvar}}
        {\restrictfunc{f}{M_{j,\firstznvar}}}
      \end{displaymath}
      for all \(j \in \integernumbers\),
      \(\firstznvar \in \zn\),
      and \(f \in \ucfunccv{\rn}\).
    \item[(B8)]
      When \(j \in \integernumbers\)
      define linear function \(P_j\) by
      \begin{displaymath}
        \left( P_j f \right) \left( \rnx  \right)
        :=
        \sum_{\firstznvar \in \zn}
        \szdualappl{\indexeddualelem{z}{j,\firstznvar}}
        {\restrictfunc{f}{M_{j,\firstznvar}}}
        w_{j,\firstznvar} \left( \rnx \right)
      \end{displaymath} 
      for all \(\rnx \in \rn\),
      \(j \in \integernumbers\), and
      \(f \in \contfunccv{\rn}\).
    \item[(B9)]
      Define
      \begin{eqnarray*}
        & & A_{p,j} :=
        \left\{
          \rnx \in \rn \mapsto
          \sum_{\firstznvar \in \zn}
          \seqelem{\sqa}{\firstznvar}
          v_{j,\firstznvar}(\rnx)
          \setsep
          \sqa \in \littlelpcv{p}{\zn}
        \right\} \\
        & &
        \norminspace{f}{A_{p,j}}
        :=
        \norm{f}_p
      \end{eqnarray*}
      for all
      \(p \in [1,\infty]\)
      and
      \(j \in \integernumbers\).
    \item[(B10)]
      Pair \((\dualw, w)\) spans all polynomials
      of \(n\) variables and of degree at most
      \(\ceil{s} - 1\).
  \end{itemize}
  Then
  \begin{eqnarray*}
    & & \exists c_1 \in \positiverealnumbers :
    \forall p \in [ 1, \infty ],
    j \in \integernumbers,
    j' \in \integernumbers : \\
    & & j' \leq j \implies
    \sup \left\{
      \norm{(I - P_j)f}_p
      \setsep
      f \in A_{p,j'}, \norm{f}_p \leq 1
    \right\}
    \leq c_1 \cdot \cinterp{p} \cdot 2^{(j'-j)s} .
  \end{eqnarray*}
\end{lemma}
}
\fi

\if\shortprep0
{
\begin{proof}
  Suppose first that \(p = \infty\).
  Let \(j_1, j'_1 \in \integernumbers\) and \(j'_1 \leq j_1\).
  Let \(g \in A_{\infty,j'_1}\) and
  \begin{displaymath}
    g(\rnx) := \sum_{\firstznvar \in \zn}
    \seqelem{\sqb}{\firstznvar}
    v_{j'_1,\firstznvar}(\rnx)
  \end{displaymath}
  for all \(\rnx \in \rn\).
  Let \(\rnxo \in \rn\).
  Let
  \begin{math}
    J := \{ \firstznvar \in \zn \setsep
    w_{j_1,\firstznvar}(\rnxo) \neq 0
    \}
  \end{math}
  and
  \begin{math}
    S_1 := \closedball{\rn}{\rnxo}{2^{-j_1}(\rho_2 + \rho_3)}
  \end{math}.
  Now \(\card{J} \leq \ncover{w}\) and
  \(M_{j_1,\firstznvar} \subset S_1\) for all
  \(\firstznvar \in J\). 
  Suppose that \(s \in \positiverealnumbers \setminus
  \positiveintegers\).
  Let \(m := \floor{s}\) and \(a := s - m\).
  Let \(t\) be the \(m\)th degree Taylor polynomial of
  \(g\) at \(\rnxo\) and \(r := g - t\).
  It follows from the polynomial span of \((\dualw, w)\)
  that
  \begin{math}
    (I - P_{j_1}) g = r - P_{j_1} r
  \end{math}.
  Hence
  \begin{math}
    ((I - P_{j_1})g)(\rnxo)
    = r(\rnxo) - (P_{j_1} r)(\rnxo)
    = - (P_{j_1} r)(\rnxo)
  \end{math}.
  Let \(\secondznvar \in J\) and
  \(\rny \in M_{j_1,\secondznvar}\).
  Now \(\norm{\rny - \rnxo} \leq 2^{-j_1}(\rho_2 + \rho_3)\).
  If \(m = 0\) we have
  \begin{math}
    \abs{r(\rny)} = \abs{g(\rny) - g(\rnxo)}
    \leq \holdercoeff{g}{0}{a} \norm{\rny - \rnxo}^s
    \leq c_2 \holdercoeff{g}{0}{a} \cdot 2^{-js}
  \end{math}.
  If \(m  > 0\) it follows from \mylemma \ref{lem:taylor-residual}
  that
  \begin{math}
    \abs{r(\rny)} =
    \leq c_3 \holdercoeff{g}{m}{a} \norm{\rny - \rnxo}^s
    \leq c_4 \holdercoeff{g}{m}{a} \cdot 2^{-js}
  \end{math}.
  By \mylemma \ref{lem:holder-coeff} we have
  \begin{math}
    \holdercoeff{g}{m}{a} \leq
    c_5 \cdot 2^{j's} \holdercoeff{v}{m}{a}
    \norminfty{\sqb}
  \end{math}.
  We also have
  \begin{math}
    \norminfty{\sqb} \leq c_6 \norminfty{g}
  \end{math}
  where \(c_6\) does not depend on \(\sqb\) or \(g\).
  Thus
  \begin{math}
    \abs{\dualappl{\indexeddualelem{z}{j_1,\secondznvar}}
    {\restrictfunc{r}{M_{j_1,\secondznvar}}}}
    \leq
    c_8 \holdercoeff{v}{m}{a} \cdot 2^{(j'_1-j_1)s}
    \norminfty{\sqb}
    \leq
    c_9 \cdot 2^{(j'_1-j_1)s} \norminfty{g}
  \end{math}.
  Hence
  \begin{align*}
    \abs{\left(\left(I - P_{j_1} \right) g \right)
    \left( \rnxo \right)}
    & \leq \sum_{\firstznvar \in J}
    \abs{\dualappl{\indexeddualelem{z}{j_1,\firstznvar}}
    {\restrictfunc{r}{M_{j_1,\firstznvar}}}}
    \abs{w_{j_1,\firstznvar}\left(\rnxo\right)} \\
    & \leq
    \ncover{w} \norminfty{w} \cdot c_9 \cdot
    2^{(j'_1 - j_1)s} \norminfty{g} .
  \end{align*}
  
  Suppose then that \(s \in \positiveintegers\).
  Let \(t\) be the \(m-1\) degree Taylor polynomial
  of \(g\) at point \(\rnxo\) and \(r := g - t\).
  The result follows from
  \begin{displaymath}
    r(\rny) =
    \frac{1}{m!}
    \sum_{i_1,\ldots,i_m=1}^n
    f_{i_1,\ldots,i_m}\left(\rnxo + c\left(\rny\right)
    \left(\rny - \rnxo\right)\right)
    \prod_{l=1}^m
    \left(
      \cartprodelem{\rny}{i_l} - \cartprodelem{\rnxo}{i_l}
    \right) .
  \end{displaymath}
  
  We have
  \begin{equation}
    \label{eq:supp}
    \suppop (I - P_j) v_{j',\firstznvar}
    \subset
    \closedball{\rn}{2^{-j'}\firstznvar}
    {2^{-j} \rho_2 + 2^{-j'} (\rho_1 + \rho_3)}
  \end{equation}
  for all \(j, j' \in \integernumbers\)
  and \(j' \leq j\).
  The result in case \(p = 1\) follows from
  \myequation \eqref{eq:supp}, \mylemma
  \ref{lem:lone-top-isomorphism} and from case
  \(p = \infty\).
    
  When \(p \in ] 1, \infty [\)
  the result follows from \mylemma \ref{lem:gen-interp}
  and Banach space interpolation.
\end{proof}
}
\fi

\if\shortprep0
{
\givelemmawithoutproof
}
\fi

\if\shortprep0
{
\begin{lemma}
  Let \(\sigma, s \in \positiverealnumbers\),
  \(s > \sigma\),
  \(p \in [ 1, \infty ]\),
  and
  \(q \in [ 1, \infty ]\).
  Suppose that
  \(\mdimmothersf{n} \in \holderspace{s}{\rn}\).
  Then
  \(\mdimmothersf{n} \in \besovspace{\sigma}{p}{q}{\rn}\).  
\end{lemma}
}
\fi

\begin{theorem}
  \label{th:tensor-and-orth-func-besov-norm-equiv}
  Let \(p \in [ 1, \infty ]\),
  \(q \in [ 1, \infty ]\),
  \(j_0 \in \integernumbers\),
  \(\sigma \in \positiverealnumbers\),
  and \(n / p < \sigma < r_0\).
  Let \(\mdimmothersf{n}\) be a scaling function of a compactly
  supported interpolating tensor product MRA \mraprepspace
  \(\ucfunccv{\rn}\) and
  \(\bphi : \rn \rightarrow \complexnumbers\) a compactly supported and 
  continuous mother scaling function of a
  \((\floor{\sigma} + 1)\)-regular
  orthonormal wavelet family.
  Suppose that
  \(\mdimmothersf{n} \in \holderspace{r_0}{\rn}\).
  Suppose also that
  \((\mdimmotherdualsf{n}, \mdimmothersf{n})\)
  and \((2^{nj} \bphi^*, \bphi)\)
  span all polynomials of degree at most
  \(\ceil{\sigma} - 1\).
  Let \(B\) be the normed space
  \begin{math}
    \besovspace{\sigma}{p}{q}{\rn} \intersection \ucfunccv{\rn}
  \end{math}
  equipped with some Besov space norm.
  Then
  \(\besovinterpwaveletnorm{n}{j_0}{\sigma}{p}{q}{\cdot}\)
  and
  \(\besovorthwaveletnorm{n}{j_0}{\sigma}{p}{q}{\cdot}\)
  are equivalent norms on the vector space
  \(B\)
  and norm
  \(\besovinterpwaveletnorm{n}{j_0}{\sigma}{p}{q}{\cdot}\)
  characterizes \(B\) on \(\ucfunccv{\rn}\).
\end{theorem}

\begin{proof}
  Let \(\br_0 := \floor{\sigma} + 1\).
  We prove first that
  \begin{math}
    \besovinterpwaveletnorm{n}{j_0}{\sigma}{p}{q}{f}
    \leq
    c \besovorthwaveletnorm{n}{j_0}{\sigma}{p}{q}{f}
  \end{math}
  for all 
  \(f \in B\)
  and for some \(c \in \positiverealnumbers\).
  Define \(q'\) by
  \begin{displaymath}
    \frac{1}{q'} + \frac{1}{q} = 1 .
  \end{displaymath}
  \if\shortprep1
  {
    There exists \(c_1 \in \positiverealnumbers\)
    so that
  }
  \else
  {
    By \mylemma \ref{lem:uprojop-in-orthvp}
    there exists \(c_1 \in \positiverealnumbers\)
    so that
  }
  \fi
  \begin{equation}
    \label{eq:proj-in-orthvp}
    \norminorthvp{p}{j'}{\mdimuprojop{n}{j}} \leq c_1 \cdot 2^{n(j'-j)/p}
  \end{equation}
  for all
  \(j, j' \in \integernumbers\), \(j' \geq j\).
  \if\shortprep1
  {
    There exists \(c_2 \in \positiverealnumbers\)
    so that
  }
  \else
  {
    By \mylemma \ref{lem:i-minus-pj-bound}
    there exists \(c_2 \in \positiverealnumbers\)
    so that
  }
  \fi
  \begin{equation}
    \label{eq:complement-in-orthvp}
    \norminorthvp{p}{j'}{I - \mdimuprojop{n}{j}} \leq c_2 \cdot
    \cinterp{p} \cdot 2^{(j'-j) \br_0}
  \end{equation}
  for all
  \(j, j' \in \integernumbers\), \(j' \leq j\).
  \if\shortprep1
  {
    These two results can be proved by Banach space interpolation between cases
    \(p=1\) and \(p=\infty\).
  }
  \fi
  Let \(f \in B\).
  We have
  \begin{eqnarray*}
    \norm{\mdimuprojop{n}{j_0} f}_p & = &
    \norm{\mdimuprojop{n}{j_0} \left( \orthprojp{p}{j_0} f + \sum_{j=j_0}^\infty
      \orthdeltaprojp{p}{j} f \right)}_p \\
    & \leq &
    \norminorthvp{p}{j_0}{\mdimuprojop {n}{j_0}}
    \norm{\orthprojp{p}{j_0} f}_p
    +
    \sum_{j=j_0}^{\infty}
    \norminorthvp{p}{j+1}{\mdimuprojop {n}{j_0}}
    \norm{\orthdeltaprojp{p}{j} f}_p \\
    & \leq &
    c_3 \norm{\orthprojp{p}{j_0} f}_p
    +
    c_3 \cdot 2^{(n-nj_0)/p}
    \sum_{j=j_0}^{\infty} 2^{-j(\sigma - \frac{n}{p})}
    \seqelem{\bbfh}{j}
  \end{eqnarray*}
  where
  \begin{math}
    \bbfh := ( 2^{j\sigma} \nsznorm{\orthdeltaprojp{p}{j} f}_p)_{j=j_0}^\infty
  \end{math}.
  Let
  \begin{displaymath}
    c_4 := c_3 \cdot 2^{(n-nj_0)/p} \cdot \norm{\left( 2^{-j'(\sigma -
      \frac{n}{p})} \right)_{j'=j_0}^\infty}_{q'} .
  \end{displaymath}
  Now
  \begin{displaymath}
    \norm{\mdimuprojop{n}{j_0} f}_p \leq
    c_3 \norm{\orthprojp{p}{j_0} f}_p +
    c_4 \norm{\bbfh}_q
    \leq \max \left\{ c_3, c_4 \right\} \cdot \besovorthwaveletnorm{n}{j_0}{\sigma}{p}{q}{f}
  \end{displaymath}
  where \(c_3\) and \(c_4\) do not depend on \(f\).
  We have
  \begin{eqnarray*}
    \norm{\mdimudeltaop{n}{j} f}_p & = &
    \norm{\mdimudeltaop{n}{j} \left( \orthprojp{p}{j_0} f + \sum_{j'=j_0}^\infty
      \orthdeltaprojp{p}{j'} f \right)}_p \\
    & \leq &
    \norminorthvp{p}{j_0}{\mdimudeltaop {n}{j}}
    \norm{\orthprojp{p}{j_0} f}_p
    +
    \sum_{j'=j_0}^{\infty}
    \norminorthvp{p}{j'+1}{\mdimudeltaop {n}{j}}
    \norm{\orthdeltaprojp{p}{j'} f}_p
  \end{eqnarray*}
  for all \(j \in \integernumbers\), \(j \geq j_0\).
  It follows from \myequation \eqref{eq:complement-in-orthvp} that
  \begin{displaymath}
    2^{j\sigma} \norminorthvp{p}{j_0}{\mdimudeltaop {n}{j}}
    \norm{\orthprojp{p}{j_0} f}_p
    \leq c_5 \cdot 2^{-j(\br_0 - \sigma)} \norm{\orthprojp{p}{j_0} f}_p
  \end{displaymath}
  for all \(j \in \integernumbers\), \(j \geq j_0\),
  where \(c_5\) is independent of \(f\).
  Consequently
  \if\ijwmip1
  {
  \begin{eqnarray}
    \nonumber
    \norm{\left( 2^{j\sigma} \norminorthvp{p}{j_0}{\mdimudeltaop {n}{j}}
    \norm{\orthprojp{p}{j_0} f}_p \right)_{j=j_0}^\infty}_q
    & \leq &
    c_5 \norm{\left( 2^{-j(\br_0 - \sigma)} \right)_{j=j_0}^\infty}_q
    \norm{\orthprojp{p}{j_0} f}_p \\
    \label{eq:ineq-i-o-a}
    & = & c_6 \norm{\orthprojp{p}{j_0} f}_p
    \leq c_6 \besovorthwaveletnorm{n}{j_0}{\sigma}{p}{q}{f} .
  \end{eqnarray}
  }
  \else
  {
  \begin{align}
    \nonumber
    & \norm{\left( 2^{j\sigma} \norminorthvp{p}{j_0}{\mdimudeltaop {n}{j}}
    \norm{\orthprojp{p}{j_0} f}_p \right)_{j=j_0}^\infty}_q \\
    \nonumber
    & \leq
    c_5 \norm{\left( 2^{-j(\br_0 - \sigma)} \right)_{j=j_0}^\infty}_q
    \norm{\orthprojp{p}{j_0} f}_p \\
    \label{eq:ineq-i-o-a}
    & = c_6 \norm{\orthprojp{p}{j_0} f}_p
    \leq c_6 \besovorthwaveletnorm{n}{j_0}{\sigma}{p}{q}{f} .
  \end{align}
  }
  \fi
  When \(j' < j\) it follows from \myequation
  \eqref{eq:complement-in-orthvp} that
  \begin{equation}
    \label{eq:delta-a}
    2^{j\sigma} \norminorthvp{p}{j'+1}{\mdimudeltaop {n}{j}}
    \norm{\orthdeltaprojp{p}{j'} f}_p
    \leq
    c_7 \cdot 2^{-\abs{j'-j}(\br_0 - \sigma)} \seqelem{\bbfh}{j'}
  \end{equation}
  When \(j' \geq j\) it follows from \myequation
  \eqref{eq:proj-in-orthvp} that
  \begin{equation}
    \label{eq:delta-b}
    2^{j\sigma} \norminorthvp{p}{j'+1}{\mdimudeltaop {n}{j}}
    \norm{\orthdeltaprojp{p}{j'} f}_p
    \leq
    c_8 \cdot 2^{-\abs{j'-j}(\sigma - \frac{n}{p})}
    \seqelem{\bbfh}{j'} .
  \end{equation}
  Let
  \begin{displaymath}
    w := \min \left\{ \sigma - \frac{n}{p}, \br_0 - \sigma \right\}
  \end{displaymath}
  and \(c_9 := \max \{ c_7, c_8 \}\).
  Now by \myequations \eqref{eq:delta-a} and \eqref{eq:delta-b}
  \begin{eqnarray*}
    a & := & \norm{\left( 2^{j\sigma} \sum_{j'=j_0}^{\infty}
      \norminorthvp{p}{j'+1}{\mdimudeltaop {n}{j}}
      \norm{\orthdeltaprojp{p}{j'} f}_p
    \right)_{j=j_0}^\infty}_q \\
    & \leq &
    c_9 \norm{\left( \sum_{j'=j_0}^\infty 2^{-\abs{j'-j}w}
      \seqelem{\bbfh}{j'} \right)_{j=j_0}^\infty}_q .
  \end{eqnarray*}
  It follows from \mylemma \ref{lem:off-diagonal-op} that
  \begin{equation}
    \label{eq:ineq-i-o-b}
    a \leq c_{10} \norm{\bbfh}_q
    \leq c_{10} \besovorthwaveletnorm{n}{j_0}{\sigma}{p}{q}{f} .
  \end{equation}
  where \(c_{10}\) depends only on \(w\) and \(j_0\).
  By \myequations \eqref{eq:ineq-i-o-a} and \eqref{eq:ineq-i-o-b} we have
  \begin{math}
    \besovinterpwaveletnorm{n}{j_0}{\sigma}{p}{q}{f}
    \leq
    c_{10} \besovorthwaveletnorm{n}{j_0}{\sigma}{p}{q}{f}
  \end{math}

  We prove then that
  \begin{math}
    \besovorthwaveletnorm{n}{j_0}{\sigma}{p}{q}{g}
    \leq
    c \besovinterpwaveletnorm{n}{j_0}{\sigma}{p}{q}{g}
  \end{math}
  for all 
  \(g \in B\)
  and for some \(c \in \positiverealnumbers\).
  Now
  \begin{eqnarray*}
    \norm{\mdimuprojop{n}{j_0} g}_p
    + \sum_{j=j_0}^\infty \norm{\mdimudeltaop{n}{j} g}_p
    & = &
    \norm{\mdimuprojop{n}{j_0} g}_p
    + \sum_{j=j_0}^\infty \norm{2^{-j\sigma} \cdot 2^{j\sigma}
      \mdimudeltaop{n}{j} g}_p \\
    & \leq &
    \norm{\mdimuprojop{n}{j_0} g}_p
    +
    \norm{\left( 2^{j\sigma}
      \norm{\mdimudeltaop{n}{j} g}_p \right)_{j=j_0}^\infty}_q
    \sum_{j=j_0}^\infty 2^{-j\sigma} \\
    & \leq & c_{11} \besovinterpwaveletnorm{n}{j_0}{\sigma}{p}{q}{g}
    \leq c_{12} \besovorthwaveletnorm{n}{j_0}{\sigma}{p}{q}{g}
    \in \nonnegrealnumbers
  \end{eqnarray*}
  where \(c_{11}\) and \(c_{12}\) do not depend on \(g\).
  Consequently
  \begin{displaymath}
    \mdimuprojop{n}{j_0} g + \sum_{j=j_0}^m \mdimudeltaop{n}{j} g
    \to
    g
  \end{displaymath}
  in \(\biglpcv{p}{\rn}\) as \(m \to \infty\).
  We have
  \begin{eqnarray*}
    \norm{\orthprojp{p}{j_0} g}_p & = &
    \norm{\orthprojp{p}{j_0} \left( \mdimuprojop{n}{j_0} g + \sum_{j=j_0}^\infty
      \mdimudeltaop{n}{j} g \right)}_p \\
    & \leq &
    \norm{\orthprojp{p}{j_0}}
    \norm{\mdimuprojop{n}{j_0} g}_p
    +
    \sum_{j'=j_0}^{\infty}
    \norm{\orthprojp{p}{j_0}}
    \norm{\mdimudeltaop{n}{j} g}_p \\
    & \leq &
    c_{13} \norm{\mdimuprojop{n}{j_0} g}_p
    +
    c_{13}
    \sum_{j=j_0}^{\infty} 2^{-j\sigma}
    \seqelem{\mathbf{h}}{j} \\
    & \leq &
    c_{13} \norm{\mdimuprojop{n}{j_0} g}_p
    +
    c_{13}
    \norm{\left( 2^{-j\sigma} \right)_{j=j_0}^\infty}_{q'}
    \norm{\mathbf{h}}_q \\
    & \leq &
    c_{14} \besovinterpwaveletnorm{n}{j_0}{\sigma}{p}{q}{g}
  \end{eqnarray*}
  where
  \begin{math}
    \mathbf{h} := (2^{j\sigma}\nsznorm{\mdimudeltaop{n}{j} g}_p)_{j=j_=}^\infty
  \end{math}.
  \if\shortprep1
  {
    There exists \(c_{15} \in \positiverealnumbers\)
    so that
  }
  \else
  {
    By \mylemma \ref{lem:i-minus-pj-bound}
    there exists \(c_{15} \in \positiverealnumbers\)
    so that
  }
  \fi
  \begin{equation}
    \label{eq:complement-in-tensorvp}
    \normintensorvp{n}{j'}{p}{I - \orthprojp{p}{j}} \leq c_{15} \cdot
    \cinterp{p} \cdot 2^{(j'-j) \sigma}
  \end{equation}
  for all
  \(j, j' \in \integernumbers\), \(j' \leq j\).
  \if\shortprep1
  {
    This can be proved by Banach space interpolation between cases
    \(p=1\) and \(p=\infty\).
  }
  \fi
  Furthermore,
  \begin{eqnarray*}
    \norm{\orthdeltaprojp{p}{j} g}_p & = &
    \norm{\orthdeltaprojp{p}{j} \left( \mdimuprojop{n}{j_0} g + \sum_{j'=j_0}^\infty
      \mdimudeltaop{n}{j'} g \right)}_p \\
    & \leq &
    \norm{\orthdeltaprojp{p}{j}}
    \norm{\mdimuprojop{n}{j_0} g}_p
    +
    \sum_{j'=j_0}^{\infty}
    \norm{\orthdeltaprojp{p}{j}}
    \norm{\mdimudeltaop{n}{j'} g}_p .
  \end{eqnarray*}
  It follows from \myequation \eqref{eq:complement-in-tensorvp} that
  \begin{displaymath}
    2^{j\sigma} 
    \norm{\orthdeltaprojp{p}{j}}
    \norm{\mdimuprojop{n}{j_0} g}_p
    \leq c_{16} \cdot 2^{-j(r_0-\sigma)}
    \norm{\mdimuprojop{n}{j_0} g}_p
  \end{displaymath}
  for all \(j \geq j_0\) where \(c_{16}\) is independent of \(g\).
  Consequently
  \begin{eqnarray}
    \nonumber
    \norm{\left( 2^{j\sigma} 
      \norm{\orthdeltaprojp{p}{j}}
      \norm{\mdimuprojop{n}{j_0} g}_p
      \right)_{j=j_0}^\infty}_q
    & \leq &
    c_{16} \norm{\left( 2^{-j(r_0-\sigma)} \right)_{j=j_0}^\infty}_q
    \norm{\mdimuprojop{n}{j_0} g}_p \\
    & = &
    \label{eq:ineq-o-i-a}
    c_{17}
    \norm{\mdimuprojop{n}{j_0} g}_p
    \leq c_{17} \besovinterpwaveletnorm{n}{j_0}{\sigma}{p}{q}{g} .
  \end{eqnarray}
  When \(j' < j\) it follows from \myequation
  \eqref{eq:complement-in-tensorvp} that
  \begin{equation}
    \label{eq:delta-c}
    2^{j\sigma} \norm{\orthdeltaprojp{p}{j}}
    \norm{\mdimudeltaop{n}{j'} g}_p
    \leq
    c_{18} \cdot 2^{-\abs{j'-j}(r_0 - \sigma)}
    \seqelem{\mathbf{h}}{j'} .
  \end{equation}
  When \(j' \geq j\) it follows from uniform boundedness of the
  projection operators that
  \begin{equation}
    \label{eq:delta-d}
    2^{j\sigma} \norm{\orthdeltaprojp{p}{j}}
    \norm{\mdimudeltaop{n}{j'} g}_p
    \leq
    c_{19} \cdot 2^{-\abs{j'-j}\sigma} \seqelem{\mathbf{h}}{j'} .
  \end{equation}
  Let
  \begin{displaymath}
    z := \min \left\{ \sigma, r_0 - \sigma \right\} .
  \end{displaymath}
  Now by \myequations \eqref{eq:delta-c} and \eqref{eq:delta-d} we
  have
  \begin{eqnarray*}
    b & := & \norm{\left( \sum_{j'=j_0}^\infty
    2^{j\sigma} \norm{\orthdeltaprojp{p}{j}}
    \norm{\mdimudeltaop{n}{j'} g}_p
    \right)_{j=j_0}^\infty}_q \\
    & \leq &
    c_{20} \norm{\left( \sum_{j'=j_0}^\infty 2^{-\abs{j'-j}z}
      \seqelem{\mathbf{h}}{j'} \right)_{j=j_0}^\infty}_q .
  \end{eqnarray*}
  It follows from \mylemma \ref{lem:off-diagonal-op} that
  \begin{equation}
    \label{eq:ineq-o-i-b}
    b \leq c_{21} \norm{\mathbf{h}}_q
    \leq c_{21} \besovinterpwaveletnorm{n}{j_0}{\sigma}{p}{q}{f} .
  \end{equation}
  where \(c_{21}\) depends only on \(z\) and \(j_0\).
  By \myequations \eqref{eq:ineq-o-i-a} and \eqref{eq:ineq-o-i-b} we have
  \begin{math}
    \besovorthwaveletnorm{n}{j_0}{\sigma}{p}{q}{g}
    \leq
    c_{22}
    \besovinterpwaveletnorm{n}{j_0}{\sigma}{p}{q}{g}
  \end{math}.

  Suppose then that
  \begin{math}
    \besovorthwaveletnorm{n}{j_0}{\sigma}{p}{q}{g}
    = \infty
  \end{math}.
  Let
  \begin{displaymath}
    A \seqstyle{b} := \left( \sum_{j'=j_0}^\infty 2^{-\abs{j'-j}z}
      \seqelem{\seqstyle{b}}{j'} \right)_{j=j_0}^\infty ,
      \spaceafter
      \seqstyle{b} \in \complexnumbers^{\naturalnumbers + j_0} .
  \end{displaymath}
  Now
  \begin{displaymath}
    \infty = \norm{\left( 2^{j\sigma} \norm{\orthdeltaprojp{p}{j} g}_p
      \right)_{j=j_0}^\infty}_q
    \leq
    c_{17} \norm{\mdimuprojop{n}{j_0} g}_p
    + c_{20} \norm{A\mathbf{h}}_q .
  \end{displaymath}
  By the continuity of
  \(\restrictfunc{A}{\littlelpcv{q}{\naturalnumbers + j_0}}\) we must
  have \(\norm{\mathbf{h}}_q = \infty\). Consequently
  \if\ijwmip1
  \begin{math}
    \besovinterpwaveletnorm{n}{j_0}{\sigma}{p}{q}{g} = \infty
  \end{math}.
  \else
  \begin{displaymath}
    \besovinterpwaveletnorm{n}{j_0}{\sigma}{p}{q}{g} = \infty .
  \end{displaymath}
  \fi
\end{proof}

It is possible to construct an orthonormal compactly supported
Daubechies wavelet family
with arbitrary high regularity \(\br \in \positiverealnumbers\)
and polynomial span (number of vanishing moments)
\(\bar{d} \in \naturalnumbers\)
\cite{daubechies1988}. Consequently, when
\(j_0 \in \integernumbers\),
\(\sigma \in \positiverealnumbers\),
\(p \in [ 1, \infty ]\),
\(q \in [ 1, \infty ]\),
\(n / p < \sigma < r_0\),
\(\mdimmothersf{n} \in \holderspace{r_0}{\rn}\), and
\((\mdimmotherdualsf{n}, \mdimmothersf{n})\)
spans all polynomials of degree at most
\(\ceil{\sigma} - 1\) function \(\mdimmothersf{n}\) defines an
equivalent norm \(\besovinterpwaveletnorm{n}{j_0}{\sigma}{p}{q}{\cdot}\)
on Banach space \(B\).

\begin{theorem}
  \label{th:tensor-func-and-coeff-norm-equiv}
  Let \(p \in [ 1, \infty ]\),
  \(q \in [ 1, \infty ]\),
  \(\sigma \in \positiverealnumbers\),
  \(j_0 \in \integernumbers\),
  \(r_0 \in \positiverealnumbers\),
  and \(n / p < \sigma < r_0\).  Let
  \(\mdimmothersf{n}\) be a scaling function of a compactly 
  supported interpolating tensor product MRA
  \mraprepspace 
  \(\ucfunccv{\rn}\) and suppose that
  \(\mdimmothersf{n} \in \holderspace{r_0}{\rn}\).
  Let \(B\) be the normed space
  \begin{math}
    \besovspace{\sigma}{p}{q}{\rn} \intersection \ucfunccv{\rn}
  \end{math}
  equipped with some Besov space norm.
  Then
  \(\besovinterpwaveletnorm{n}{j_0}{\sigma}{p}{q}{\cdot}\)
  and
  \(\besovinterpcoeffnorm{n}{j_0}{\sigma}{p}{q}{\cdot}\)
  are equivalent to the norm of \(B\) and
  they characterize \(B\) on \(\ucfunccv{\rn}\).
\end{theorem}

\begin{proof}
  Define
  \begin{displaymath}
    \sqa_j(f) := \left( \szdualappl{\mdimgendualwavelet{n}{\zos}{j}{\firstznvar}}{f}
      \right)_{\zos \in \zeroonesetnnozero{n},\firstznvar \in \zn} 
  \end{displaymath}
  where
  \(f \in B\) and \(j \in \integernumbers\).
  \if\shortprep1
  {
    By applying Banach space interpolation to cases \(p = 1\) and \(p =
    \infty\) we get
  }
  \else
  {
    By \mylemmas
    \ref{lem:mdimpartialuwspace-interp} and
    \ref{lem:mdimuwspace-interp}
    we get
  }
  \fi
  \begin{eqnarray}
    \nonumber
    c_1 \cdot 2^{-\frac{nj_0}{p}}
    \norm{\left( \szdualappl{\mdimdualsf{n}{j_0}{\firstznvar}}{f}
      \right)_{\firstznvar \in \zn}}_p
    & \leq &
    \norm{\mdimuprojop{n}{j_0} f} \\
    \label{eq:coeff-pj0}
    & \leq &
    c_2 \cdot 2^{-\frac{nj_0}{p}}
    \norm{\left( \szdualappl{\mdimdualsf{n}{j_0}{\firstznvar}}{f}
      \right)_{\firstznvar \in \zn}}_p
  \end{eqnarray}
  and
  \begin{displaymath}
    c_3 \cdot 2^{-\frac{nj}{p}}
    \norm{\sqa_j(f)}_p
    \leq
    \norm{\mdimudeltaop{n}{j} f}_p
    \leq
    c_4 \cdot 2^{-\frac{nj}{p}}
    \norm{\sqa_j(f)}_p
  \end{displaymath}
  for all \(j \in \naturalnumbers + j_0\) and
  \(f \in B\).
  Consequently
  \begin{eqnarray}
    \nonumber
    c_3 \norm{\left( 2^{j(\sigma-\frac{n}{p})} \norm{\sqa_j(f)}_p \right)_{j=j_0}^\infty}_q
    & \leq &
    \norm{\left( 2^{j\sigma} \norm{\mdimudeltaop{n}{j} f}_p \right)_{j=j_0}^\infty}_q \\
    \label{eq:coeff-qj}
    & \leq &
    c_4 \norm{\left( 2^{j(\sigma-\frac{n}{p})} \norm{\sqa_j(f)}_p \right)_{j=j_0}^\infty}_q
  \end{eqnarray}
  for all
  \(f \in B\).

  Let
  \begin{math}
    c_5 := \min \{ c_1 \cdot 2^{-nj_0/p}, c_3\}
  \end{math}
  and
  \begin{math}
    c_6 := \max \{ c_2 \cdot 2^{-nj_0/p}, c_4\}
  \end{math}.
  It follows from \myequations \eqref{eq:coeff-pj0} and
  \eqref{eq:coeff-qj} that
  \begin{displaymath}
    c_5 \besovinterpcoeffnorm{n}{j_0}{\sigma}{p}{q}{f}
    \leq
    \besovinterpwaveletnorm{n}{j_0}{\sigma}{p}{q}{f}  
    \leq
    c_6 \besovinterpcoeffnorm{n}{j_0}{\sigma}{p}{q}{f}
  \end{displaymath}
  for all
  \(f \in B\).
\end{proof}

It is possible to construct Deslauriers-Dubuc wavelet families with
arbitrary high regularity \(r \in \positiverealnumbers\) (i.e.
\(\mdimmothersf{n} \in \holderspace{r}{\rn}\)) and polynomial span
\(d \in \naturalnumbers\) \cite{donoho1992}.
As a direct consequence of \mytheorems
\ref{th:tensor-and-orth-func-besov-norm-equiv} and
\ref{th:tensor-func-and-coeff-norm-equiv}
norms
\(\besovinterpcoeffnorm{n}{j_0}{\sigma}{p}{q}{\cdot}\)
and
\(\besovinterpwaveletnorm{n}{j_0}{\sigma}{p}{q}{\cdot}\)
are equivalent to the restriction of some norm of
Besov space \(\besovspace{\sigma}{p}{q}{\rn}\)
onto the vector space
\(\besovspace{\sigma}{p}{q}{\rn} \intersection \ucfunccv{\rn}\)
and these norms characterize \(\besovspace{\sigma}{p}{q}{\rn}
\intersection \ucfunccv{\rn}\)
on \(\ucfunccv{\rn}\).

\if01
{
\begin{lemma}
  \label{lem:tensor-and-orth-func-besov-norm-equiv}
  Let \(p \in [ 1, \infty ]\),
  \(q \in [ 1, \infty ]\),
  \(\sigma \in \positiverealnumbers\),
  and \(n / p < \sigma < r_0\).
  Let \(\mdimmothersf{n}\) be a scaling function of a compactly
  supported interpolating tensor product MRA \mraprepspace
  \(\ucfunccv{\rn}\) and
  \(\bphi : \rn \rightarrow K\) a compactly supported and 
  continuous mother scaling function of a
  \((\floor{\sigma} + 1)\)-regular
  orthonormal wavelet family.
  Suppose that
  \(\mdimmothersf{n} \in \holderspace{r_0}{\rn}\)
  and
  \((\mdimmotherdualsf{n}, \mdimmothersf{n})\)
  spans all polynomials of degree at most
  \(\ceil{\sigma} - 1\).
  Then
  \(\besovinterpwaveletnorm{n}{j_0}{\sigma}{p}{q}{\cdot}\)
  and
  \(\besovorthwaveletnorm{n}{j_0}{\sigma}{p}{q}{\cdot}\)
  are equivalent norms on the vector space
  \(\besovspace{\sigma}{p}{q}{\rn} \intersection \ucfunccv{\rn}\).
\end{lemma}

\begin{proof}
  We have \(\bphi \in \holderspace{\br_0}{\rn}\)
  where \(\sigma < \br_0 < \floor{\sigma} + 1\).
  By \cite[section 2.6 corollary]{meyer1992}
  \((2^{nj} \bphi^*, \bphi)\)
  spans all polynomials of degree at most
  \(\floor{\sigma} + 1\).
  The proof is based on \mylemmas
  \ref{lem:conv-in-different-norms},
  \ref{lem:orth-proj-op-complement-infty-norm},
  \ref{lem:off-diagonal-op},
  \ref{lem:uprojop-in-orthvp},
  \ref{lem:uprojop-inclusion}, and
  \ref{lem:i-minus-pj-bound}
  and \myequation
  \eqref{eq:orthprojop-uniform-bound}.
\end{proof}

\begin{lemma}
  \label{lem:tensor-func-and-coeff-norm-equiv}
  Let \(p \in [ 1, \infty ]\),
  \(q \in [ 1, \infty ]\),
  \(\sigma \in \positiverealnumbers\),
  \(j_0 \in \integernumbers\),
  \(r_0 \in \positiverealnumbers\),
  and \(n / p < \sigma < r_0\).  Let
  \(\mdimmothersf{n}\) be a scaling function of a compactly 
  supported interpolating tensor product MRA
  \mraprepspace 
  \(\ucfunccv{\rn}\) and suppose that
  \(\mdimmothersf{n} \in \besovspace{\sigma}{p}{q}{\rn}\) and
  \(\mdimmothersf{n} \in \holderspace{r_0}{\rn}\).  Then
  \(\besovinterpwaveletnorm{n}{j_0}{\sigma}{p}{q}{\cdot}\)
  and
  \(\besovinterpcoeffnorm{n}{j_0}{\sigma}{p}{q}{\cdot}\)
  are equivalent norms on the 
  vector space
  \(\besovspace{\sigma}{p}{q}{\rn} \intersection \ucfunccv{\rn}\).
\end{lemma}

\begin{proof}
  Use \mylemmas
  \ref{lem:mdimpartialuwspace-interp} and
  \ref{lem:mdimuwspace-interp}.
\end{proof}

\begin{theorem}
  \label{th:tensor-besov-norm-equiv}
  Let \(p \in [ 1, \infty ]\),
  \(q \in [ 1, \infty ]\), \(\sigma \in ]0, 1[\),
  \(r_0 \in \positiverealnumbers\),
  and
  \(n / p < \sigma < r_0\).  Let
  \(\mdimmothersf{n} : \rn \rightarrow K\) be a scaling function of 
  a compactly supported
  interpolating tensor product MRA \mraprepspace
  \(\ucfunccv{\rn}\) and 
  suppose that
  \(\mdimmothersf{n} \in \holderspace{r_0}{\rn}\).
  Suppose also that
  \((\mdimmotherdualsf{n}, \mdimmothersf{n})\)
  spans all polynomials of degree at most
  \(\ceil{\sigma} - 1\).
  Then
  \(\besovinterpcoeffnorm{n}{j_0}{\sigma}{p}{q}{\cdot}\)
  and
  \(\besovinterpwaveletnorm{n}{j_0}{\sigma}{p}{q}{\cdot}\)
  are equivalent to the restriction of some norm of
  Besov space \(\besovspace{\sigma}{p}{q}{\rn}\)
  onto the vector space
  \(\besovspace{\sigma}{p}{q}{\rn} \intersection \ucfunccv{\rn}\).
\end{theorem}

\begin{proof}
  The theorem follows from \mylemmas
  \ref{lem:tensor-and-orth-func-besov-norm-equiv} and
  \ref{lem:tensor-func-and-coeff-norm-equiv}.
\end{proof}
}
\fi

\subsection{Consequences of the Besov Space Norm Equivalence}

In particular, the Besov space norm equivalence holds for
H\"older spaces \(\holderspace{\sigma}{\rn}\) for
\(\sigma \in \positiverealnumbers \setminus \integernumbers\).
When \(n \in \positiveintegers\),
\(\sigma \in \positiverealnumbers\),
and \(q \in [1, \infty]\)
we have
\begin{math}
  \besovspace{\sigma}{\infty}{q}{\realnumbers^n}
  \subsetset
  \ucfunccv{\realnumbers^n}
\end{math}
\cite[prop. 2.3.2/2(i) and eq. (2.3.5/1)]{triebel1983}.

\begin{definition}
  \label{def:bmiset}
  When \(n \in \positiveintegers\) and \(j_0 \in \integernumbers\) 
  define
  \begin{displaymath}
    \bmiset{n}{j_0}
    :=
    \left\{ (\finitezeroseq{n}, j_0, \firstznvar) \setsep \firstznvar \in \zn \right\}
    \cup
    \bigcup_{j=j_0}^\infty
    \bigcup_{\zos \in \zeroonesetnnozero{n}}
    \left\{ (\zos, j, \firstznvar) \setsep \firstznvar \in \zn \right\} .
  \end{displaymath}
\end{definition}

\if\shortprep1
{
  As the Banach space
  \begin{math}
    \besovspace{\sigma}{\infty}{\infty}{\rn}
    \equaltvs \zygmundspace{\sigma}{\rn}
  \end{math}
  is isomorphic to \(l^\infty\)
  \cite[section 2.5.5. page 87]{triebel1983}
  the set
  \begin{math}
    \{ \psi_\alpha \setsep \alpha \in \bmiset{n}{j_0} \}
  \end{math}
  cannot be a Schauder basis of
  \begin{math}
    \besovspace{\sigma}{\infty}{\infty}{\rn}
  \end{math}
  with any summing order.

  \begin{remark}
  Suppose that
  \begin{math}
    n \in \positiveintegers
  \end{math},
  \begin{math}
    j_0 \in \integernumbers
  \end{math},
  \begin{math}
    \sigma \in \positiverealnumbers
  \end{math}, and
  \begin{math}
    \mdimmothersf{n} \in \holderspace{r}{\rn}
  \end{math}
  where \(r \in \positiverealnumbers\), \(r > \sigma\).
  There exists \(r_1 \in \positiveintegers\) so that
  \begin{math}
    \suppop \mdimmothersf{n} \in \closedball{\rn}{0}{r_1}
  \end{math}.
  Consider
  \begin{displaymath}
    f :=
    \sum_{k=0}^\infty 2^{-(j_0+k)\sigma}
    \mdimgenwavelet{n}{\onetonnbv{n}{1}}{j_0+k}{\eta(k)}
  \end{displaymath}
  where
  \begin{eqnarray*}
    c & := & 2 r_1 + \ceil{2^{j_0}} (2\ceil{\sigma}+1) \\
    \eta(k) & := & c \cdot 2^{k} k \onetonnbv{n}{1}
  \end{eqnarray*}
  for all \(k \in \integernumbers\), \(k \geq j_0\),
  and the series converges in
  \begin{math}
    \vanishingfunccv{\rn}
  \end{math}.
  Let \(m := \ceil{r}\)
  and
  \begin{displaymath}
    g_k := 2^{-(j_0+k)\sigma}
    \mdimgenwavelet{n}{\onetonnbv{n}{1}}{j_0+k}{\eta(k)}
  \end{displaymath}
  for all
  \begin{math}
    k \in \naturalnumbers
  \end{math}.
  As
  \begin{math}
    \mdimmotherwavelet{n}{\onetonnbv{n}{1}}
    \in \besovspace{\sigma}{\infty}{\infty}{\rn}
  \end{math}
  it follows that
  there exists \(c_1 \in \positiverealnumbers\) so that
  \begin{equation}
    \label{eq:mod-cont-a}
    \szgenmodcont{m}{\infty}{g_k}{t}
    \leq
    c_1 t^\sigma
  \end{equation}
  for all \(k \in \naturalnumbers\) and \(t \in ] 0, 1 [\).
  Define
  \begin{math}
    S_{j,\secondznvar}(y)
    := \closedball{\rn}{2^{-j} \secondznvar}{2^{-j} r_1 + y}
  \end{math}.
  where \(j \in \integernumbers\),
  \(\secondznvar \in \zn\),
  and \(y \in \nonnegrealnumbers\).
  Suppose that \(t_1 \in ]0, 1[\).
  We have
  \begin{equation}
    \label{eq:supp-b}
    \suppopst \funcdiff{m}{\rnh}{g_k}
    \subset
    S_{j_0+k,\eta(k)}(m)
  \end{equation}
  for all
  \(\rnh \in \closedball{\rn}{0}{t_1}\)
  and \(k \in \naturalnumbers\).
  Furthermore,
  \begin{equation}
    \label{eq:empty-intersection}
    S_{j_0+k,\eta(k)}(m)
    \intersection
    S_{j_0+l,\eta(l)}(m)
    =
    \emptyset
  \end{equation}
  for all \(k, l \in \naturalnumbers\)
  and \(k \neq l\).
  Using \myequations \eqref{eq:supp-b} and \eqref{eq:empty-intersection}
  we obtain
  \begin{equation}
    \label{eq:mod-cont-b}
    \genmodcont{m}{\infty}{f}{t_1}
    \leq
    \sup \{ \genmodcont{m}{\infty}{g_k}{t_1} \setsep
    k \in \naturalnumbers\} .
  \end{equation}
  We have
  \(f \in \biglpcv{\infty}{\rn}\)
  and by
  \myequations \eqref{eq:mod-cont-a} and \eqref{eq:mod-cont-b}
  \(
  \genmodcont{m}{\infty}{f}{t_1}
  \leq c_2 t_1^\sigma
  \).
  Thus \(f \in \besovspace{\sigma}{\infty}{\infty}{\rn}\).

  Suppose that
  \begin{displaymath}
    f = \sum_{k=0}^\infty \seqelem{\sqd}{\iota(k)} \psi_{\iota(k)}
  \end{displaymath}
  for some sequence
  \(\sqd \in \seqset{\bmiset{n}{j_0}}{\complexnumbers}\)
  where the series converges in 
  \(\besovspace{\sigma}{\infty}{\infty}{\rn}\).
  Define
  \begin{displaymath}
    \xi_b :=
    \sum_{k = 0}^b
    \seqelem{\sqd}{\iota(k)}
    \psi_{\iota(k)}
  \end{displaymath}
  where \(b \in \naturalnumbers\).
  Let \(b_1 \in \naturalnumbers\).
  Define
  \begin{displaymath}
    j_1 := \max \{ j(\iota(k)) \in \integernumbers \setsep
    k \in \setzeroton{b_1}
    \}
    + 1 .
  \end{displaymath}
  There exists \(k_1 \in \naturalnumbers\)
  so that \(\iota(k_1) = \alpha(j_1 - j_0)\).
  Now \(j_1 \neq j(\iota(k))\) for all \(k \in \setzeroton{b_1}\)
  and consequently \(\iota(k_1) \neq \iota(k)\) for all
  \(k \in \setzeroton{b_1}\). It follows that \(k_1 > b_1\).
  Let \(k_2 := j_1 - j_0\).
  Now
  \begin{displaymath}
    \iota(k_1) = (\onetonnbv{n}{1}, j_1, c \cdot 2^{k_2} k_2 \onetonnbv{n}{1}) .
  \end{displaymath}
  We have
  \begin{displaymath}
    \abs{\dualappl{\psidual_{\iota(k_1)}}{f - \xi_{b_1}}}
    \leq
    \norminfty{\left( \dualappl{\mdimgendualwavelet{n}{\zos}{j_1}{\firstznvar}}
      {f - \xi_{b_1}}
      \right)_{\zos \in \zeroonesetnnozero{n},\firstznvar \in \zn}}
  \end{displaymath}
  and
  \begin{math}
    \abs{\dualappl{\mdimmotherdualwavelet{n}{\iota(k_1)}}{f - \xi_{b_1}}}
    = 2^{-(j_0+k_2)\sigma}
  \end{math}.
  Hence
  \begin{eqnarray*}
    1
    & \leq &
    2^{(j_0+k_2)\sigma}
    \norminfty{\left( \dualappl{\mdimgendualwavelet{n}{\zos}{j_1}{\firstznvar}}
      {f - \xi_{b_1}}
      \right)_{\zos \in \zeroonesetnnozero{n},\firstznvar \in \zn}} \\
    & \leq &
    \norminfty{
      \left(
      2^{j \sigma}
      \norminfty{\left( \dualappl{\mdimgendualwavelet{n}{\zos}{j}{\firstznvar}}
        {f - \xi_{b_1}}
        \right)_{\zos \in \zeroonesetnnozero{n},\firstznvar \in \zn}}
      \right)_{j=j_0}^\infty
    } \\
    & \leq &
    \besovinterpcoeffnorm{n}{j_0}{\sigma}{\infty}{\infty}{f - \xi_{b_1}} .
  \end{eqnarray*}
  Thus \(\xi_b \not\to f\) in Banach space
  \(\besovspace{\sigma}{\infty}{\infty}{\rn}\)
  as \(b \to \infty\).
  \end{remark}
}
\fi

\if\shortprep0
{
\begin{theorem}
  \label{th:wavelet-basis}
  Let \(n \in \positiveintegers\),
  \(\sigma \in \positiverealnumbers\), and
  \(j_0 \in \integernumbers\).
  Let \(\mdimmothersf{n}\) be a mother scaling function of a
  compactly supported tensor product MRA \mraprepspace
  \(\ucfunccv{\rn}\). Suppose that
  \(\mdimmothersf{n} \in \holderspace{r}{\rn}\)
  for some \(r \in \positiverealnumbers\), \(r > \sigma\).
  Then the sequence
  \((\psi_{\iota(k)})_{k=0}^\infty\) is not a 
  Schauder basis of
  \begin{math}
    \besovspace{\sigma}{\infty}{\infty}{\rn}
    \equaltvs \zygmundspace{\sigma}{\rn}
  \end{math}
  with any summing order \(\iota\) where \(\iota\) is a 
  bijection from
  \(\naturalnumbers\) onto
  \(\bmiset{n}{j_0}\).
\end{theorem}
}
\fi

\if\shortprep1
{
\if10
{
\begin{proof}
  There exists \(r_1 \in \positiveintegers\) so that
  \begin{math}
    \suppop \mdimmothersf{n} \in \closedball{\rn}{0}{r_1}
  \end{math}.
  Let
  \(c := 2 r_1 + \ceil{2^{j_0}} (2\ceil{\sigma}+1)\)
  and \(m := \ceil{r}\).
  Let
  \begin{eqnarray*}
    \eta(k) & := & c \cdot 2^{k} k \onetonnbv{n}{1} \\
    g_k & := & 2^{-(j_0+k)\sigma}
    \mdimgenwavelet{n}{\onetonnbv{n}{1}}{j_0+k}{\eta(k)}
  \end{eqnarray*}
  where \(k \in \naturalnumbers\) and
  \begin{displaymath}
    f :=
    \sum_{k=0}^\infty g_k
  \end{displaymath}
  where the series converges in
  \begin{math}
    \vanishingfunccv{\rn}
  \end{math}.
  It can be shown that
  \begin{math}
    f \in \besovspace{\sigma}{\infty}{\infty}{\rn}
  \end{math}
  by proving that
  \begin{math}
    \contbesovnorm{\sigma}{\infty}{\infty}{\rn}{m}{f}
    <
    \infty
  \end{math}.
  On the other hand, it can be shown that
  \begin{math}
    f
  \end{math}
  is not the sum of the series
  \begin{displaymath}
    \sum_{k=0}^\infty \seqelem{\sqd}{k} \iota(k)
  \end{displaymath}
  in
  \begin{math}
    \besovspace{\sigma}{\infty}{\infty}{\rn}
  \end{math}
  for any bijection
  \begin{math}
    \iota : \naturalnumbers \onto \bmiset{n}{j_0}
  \end{math}
  and any sequence
  \begin{math}
    \sqd \in \seqset{\naturalnumbers}{\complexnumbers}
  \end{math}.
\end{proof}
}
\fi
}
\else
{
\begin{proof}
  Define
  \begin{math}
    j(\sqa) = j'
  \end{math}
  for all \(\sqa = (\zosp, j', \firstznvarp)\),
  \(\zosp \in \zeroonesetn\),
  \(j' \in \integernumbers\), and
  \(\firstznvarp \in \zn\).
    There exists \(r_1 \in \positiveintegers\) so that
    \begin{math}
      \suppop \mdimmothersf{n} \in \closedball{\rn}{0}{r_1}
    \end{math}.
    Let
    \(c := 2 r_1 + \ceil{2^{j_0}} (2\ceil{\sigma}+1)\)
    and \(m := \ceil{r}\).
    Let
    \begin{eqnarray*}
      \eta(k) & := & c \cdot 2^{k} k \onetonnbv{n}{1} \\
      \alpha(k) & := &
      (\onetonnbv{n}{1},{j_0+k},\eta(k)) \\
      a_k & := &
      \psi_{\alpha(k)} \\
      g_k & := & 2^{-(j_0+k)\sigma} a_k \\
      z & := & \mdimmotherwavelet{n}{\onetonnbv{n}{1}}
    \end{eqnarray*}
    where \(k \in \naturalnumbers\)
    and
    \begin{displaymath}
      f_l := \sum_{k=0}^l g_k
    \end{displaymath}
    where \(l \in \naturalnumbers\).
    Let
    \(f\)
    be the limit of sequence \((f_l)_{l=0}^\infty\) in
    \(\vanishingfunccv{\rn}\).
    We have
    \begin{displaymath}
      \sum_{k=0}^\infty \norminfty{g_k}
      = \norminfty{z} \sum_{k=0}^\infty
      2^{-(j_0+k)s}
      \in \positiverealnumbers
    \end{displaymath}
    and hence series \(\sum_{k=0}^\infty g_k\)
    converges absolutely in Banach space \(\vanishingfunccv{\rn}\).
    We also have \(a_k, g_k \in \besovspace{\sigma}{\infty}{\infty}{\rn}\)
    for all \(k \in \naturalnumbers\).
    Define
    \begin{math}
      S_{j,\secondznvar}(y)
      := \closedball{\rn}{2^{-j} \secondznvar}{2^{-j} r_1 + y}
    \end{math}.
    where \(j \in \integernumbers\),
    \(\secondznvar \in \zn\),
    and \(y \in \nonnegrealnumbers\).
  Suppose that \(t_1 \in ]0, 1[\).
  We have
  \begin{equation}
    \label{eq:supp-b}
    \suppopst \funcdiff{m}{\rnh}{g_k}
    =
    \suppopst \funcdiff{m}{\rnh}{a_k}
    \subset
    S_{j_0+k,\eta(k)}(m)
  \end{equation}
  for all
  \(\rnh \in \closedball{\rn}{0}{t_1}\)
  and \(k \in \naturalnumbers\).
  Furthermore,
  \begin{equation}
    \label{eq:empty-intersection}
    S_{j_0+k,\eta(k)}(m)
    \intersection
    S_{j_0+l,\eta(l)}(m)
    =
    \emptyset
  \end{equation}
  for all \(k, l \in \naturalnumbers\)
  and \(k \neq l\).
  Using \myequations \eqref{eq:supp-b} and \eqref{eq:empty-intersection}
  we obtain
  \begin{equation}
    \label{eq:mod-cont}
    \genmodcont{m}{\infty}{f}{t_1}
    \leq
    \sup \{ \genmodcont{m}{\infty}{g_k}{t_1} \setsep
    k \in \naturalnumbers\} .
  \end{equation}
  Furthermore,
  \(\genmodcont{m}{\infty}{a_k}{t_1}
  = \genmodcont{m}{\infty}{z}{2^{j_0+k} t_1}\)
  for all \(k \in \naturalnumbers\).
  Consequently
  \begin{displaymath}
    \genmodcont{m}{\infty}{z}{2^{j_0+k} t_1}
    \leq
    c_2
    \cdot
    \left( 2^{j_0+k} t_1 \right)^\sigma
  \end{displaymath}
  where
  \(
  c_2 := \contbesovnormrl{\sigma}{\infty}{\infty}{\rn}{m}{z}
  \).
  We also have
  \begin{displaymath}
    \genmodcont{m}{\infty}{g_k}{t_1}
    = 2^{-\left(j_0+k\right) \sigma}
    \genmodcont{m}{\infty}{a_k}{t_1}
    \leq c_2 t_1^\sigma .
  \end{displaymath}
  We have
  \(f \in \biglpcv{\infty}{\rn}\)
  and by
  \myequation \eqref{eq:mod-cont}
  \(
  \genmodcont{m}{\infty}{f}{t_1}
  \leq c_2 t_1^\sigma
  \).
  Thus \(f \in \besovspace{\sigma}{\infty}{\infty}{\rn}\).
  
  Suppose that \((\psi_{\iota(k)})_{j=j_0}^\infty\)
  would be a Schauder basis of
  \(\besovspace{\sigma}{\infty}{\infty}{\rn}\).
  Then
  \begin{displaymath}
    f = \sum_{k=0}^\infty \seqelem{\sqd}{\iota(k)} \psi_{\iota(k)}
  \end{displaymath}
  for some sequence
  \(\sqd \in \seqset{\bmiset{n}{j_0}}{\complexnumbers}\).
  Define
  \begin{displaymath}
    \xi_b :=
    \sum_{k = 0}^b
    \seqelem{\sqd}{\iota(k)}
    \psi_{\iota(k)}
  \end{displaymath}
  where \(\beta \in \naturalnumbers\).
  Let \(b_1 \in \naturalnumbers\).
  Define
  \begin{displaymath}
    j_1 := \max \{ j(\iota(k)) \in \integernumbers \setsep
      k \in \setzeroton{b_1}
    \}
    + 1 .
  \end{displaymath}
  There exists \(k_1 \in \naturalnumbers\)
  so that \(\iota(k_1) = \alpha(j_1 - j_0)\).
  Now \(j_1 \neq j(\iota(k))\) for all \(k \in \setzeroton{b_1}\)
  and consequently \(\iota(k_1) \neq \iota(k)\) for all
  \(k \in \setzeroton{b_1}\). It follows that \(k_1 > b_1\).
  Let \(k_2 := j_1 - j_0\).
  Now
  \begin{displaymath}
    \iota(k_1) = (\onetonnbv{n}{1}, j_1, c \cdot 2^{k_2} k_2 \onetonnbv{n}{1}) .
  \end{displaymath}
  We have
  \begin{displaymath}
  \abs{\dualappl{\psidual_{\iota(k_1)}}{f - \xi_{b_1}}}
  \leq
  \norminfty{\left( \dualappl{\mdimgendualwavelet{n}{\zos}{j_1}{\firstznvar}}
  {f - \xi_{b_1}}
  \right)_{\zos \in \zeroonesetnnozero{n},\firstznvar \in \zn}}
  \end{displaymath}
  and
  \begin{math}
    \abs{\dualappl{\mdimmotherdualwavelet{n}{\iota(k_1)}}{f - \xi_{b_1}}}
    = 2^{-(j_0+k_2)\sigma}
  \end{math}.
  Hence
  \begin{eqnarray*}
    1
    & \leq &
    2^{(j_0+k_2)\sigma}
    \norminfty{\left( \dualappl{\mdimgendualwavelet{n}{\zos}{j_1}{\firstznvar}}
    {f - \xi_{b_1}}
    \right)_{\zos \in \zeroonesetnnozero{n},\firstznvar \in \zn}} \\
    & \leq &
    \norminfty{
    \left(
      2^{j \sigma}
      \norminfty{\left( \dualappl{\mdimgendualwavelet{n}{\zos}{j}{\firstznvar}}
      {f - \xi_{b_1}}
      \right)_{\zos \in \zeroonesetnnozero{n},\firstznvar \in \zn}}
      \right)_{j=j_0}^\infty
    } \\
    & \leq &
    \besovinterpcoeffnorm{n}{j_0}{\sigma}{\infty}{\infty}{f - \xi_{b_1}} .
  \end{eqnarray*}
  Thus \(\xi_b \not\to f\) in Banach space
  \(\besovspace{\sigma}{\infty}{\infty}{\rn}\)
  as \(b \to \infty\).
  Consequently
  \((\psi_{\iota(k)})_{j=j_0}^\infty\)
  is not a Schauder basis of
  \(\besovspace{\sigma}{\infty}{\infty}{\rn}\).
\end{proof}
}
\fi

\begin{theorem}
  \label{th:wavelet-unc-basis}
  Let \(n \in \positiveintegers\),
  \(\sigma \in \positiverealnumbers\),
  \(q \in [1, \infty[\),
  and
  \(j_0 \in \integernumbers\).
  Let \(\mdimmothersf{n}\) be a mother scaling function of a
  compactly supported tensor product MRA \mraprepspace
  \(\vanishingfunccv{\rn}\). Suppose that
  \(\mdimmothersf{n} \in \holderspace{r}{\rn}\)
  for some \(r \in \positiverealnumbers\), \(r > \sigma\).
  Then \(\{\psi_\alpha \setsep \alpha \in \bmiset{n}{j_0}\}\)
  is an unconditional basis of Banach space
  \(\besovspace{\sigma}{\infty}{q}{\rn} \intersection
  \vanishingfunccv{\rn}\) equipped with a norm of
  \(\besovspace{\sigma}{\infty}{q}{\rn}\)
  and the coefficient functional corresponding to basis vector
  \(\psi_\alpha\) is \(\psidual_\alpha\) for each
  \(\alpha \in \bmiset{n}{j_0}\).
\end{theorem}

\begin{proof}
  Let
  \begin{eqnarray*}
    S_j(m) & := & \{ (\finitezeroseq{n}, j, \firstznvar)
    \setsep \firstznvar \in \zn,
    \normtwo{\firstznvar} \leq m \},
    \spaceafter j \in \integernumbers, m \in \naturalnumbers \\
    D_{\zos,j}(m) & := & \{ (\zos, j, \firstznvar)
    \setsep \firstznvar \in \zn,
    \normtwo{\firstznvar} \leq m \},
    \spaceafter
    \zos \in \zeroonesetn,
    j \in \integernumbers, m \in \naturalnumbers \\
    D_j(m) & := &
    \bigcup_{\zos \in \zeroonesetnnozero{n}}
    D_{\zos,j}(m),
    \spaceafter j \in \integernumbers, m \in \naturalnumbers .
  \end{eqnarray*}
  Let
  \begin{eqnarray*}
    m_1 & := & \min \{
      m \in \naturalnumbers
      \setsep
      \forall \firstznvar \in \zn :
      (
        \norminfty{\firstznvar} > m
        \implies
        \forall \zos \in \zeroonesetn :
        \mdimgendualwaveletfilterelem{n}{\zos}{\firstznvar}
        = 0
      )
    \} \\
    m_2 & := & \ceil{m_1 \sqrt{n}} \\
    A(k) & := & S_{j_0}(2^{j_0} k)
    \cup
    \bigcup_{l = j_0}^{j_0 + k - 1}
    D_l(m_2 + 2^{j_0 + k} k),
    \spaceafter
    k \in \positiveintegers .
  \end{eqnarray*}
  Define
  \begin{math}
    \zos(\sqa) = \zosp
  \end{math},
  \begin{math}
    j(\sqa) = j'
  \end{math}, and
  \begin{math}
    \firstznvar(\sqa) = \firstznvarp
  \end{math}
  for all \(\sqa = (\zosp, j', \firstznvarp)\),
  \(\zosp \in \zeroonesetn\),
  \(j' \in \integernumbers\), and
  \(\firstznvarp \in \zn\).
  Let \(f \in \besovspace{\sigma}{\infty}{q}{\rn} \intersection
  \vanishingfunccv{\rn}\).
  Let \(\eta : \naturalnumbers \onto \bmiset{n}{j_0}\)
  be a bijection.
  Let
  \begin{math}
    \beta_\alpha
    :=
    \dualappl{\psidual_\alpha}{f}
  \end{math}
  for all
  \begin{math}
    \alpha \in \bmiset{n}{j_0}
  \end{math}.
  Define
  \begin{displaymath}
    \beta^{(m)}_\alpha :=
    \left\{
      \begin{array}{ll}
        \dualappl{\psidual_\alpha}{f} ; \spaceafter
        & \alpha \not\in \setimage{\eta}{\setzeroton{m}}  \\
        0 ; & \textrm{otherwise}
      \end{array}
    \right.
  \end{displaymath}
  for all
  \(m \in \naturalnumbers\)
  and
  \(\alpha \in \bmiset{n}{j_0}\).
  Define
  \begin{displaymath}
    \xi_m :=
    \sum_{\alpha \in \bmiset{n}{j_0}}
    \beta^{(m)}_\alpha \psi_\alpha
  \end{displaymath}
  for all
  \(m \in \naturalnumbers\).
  There exists \(c_1 \in \positiverealnumbers\)
  so that
  \begin{equation}
    \label{eq:functional-bound}
    \abs{\dualappl{\psidual_\alpha}{f}}
    \leq
    c_1 \norminfty{f}
  \end{equation}
  for all
  \(f \in \vanishingfunccv{\rn}\)
  and
  \(\alpha \in \bmiset{n}{j_0}\).

  Let \(h \in \positiverealnumbers\).
  Choose \(j_1 \in \integernumbers\),
  \(j_1 > \max \{ j_0, 0 \}\)
  so that
  \begin{equation}
    \label{eq:tail-bound}
    \norm{
    \left(
      2^{j \sigma}
      \norm{
        \left(
          \szdualappl{\mdimgendualwavelet{n}{\zos}{j}{\firstznvar}}
          {f}
        \right)_{
          \zos \in \zeroonesetnnozero{n},
          \firstznvar \in \zn
        }
      }_\infty
    \right)_{j = j_1}^\infty
    }_q
    < \frac{h}{4} .
  \end{equation}
  Choose \(r_1 \in \positiverealnumbers\) so that
  \begin{equation}
    \label{eq:vanishing-bound}
    \sup \{
      \abs{f(\rnx)}
      \setsep
      \rnx \in \rn,
      \norm{\rnx} \geq r_1
    \}
    <
    \frac{h}{2^{n + 2 + j_1 \sigma}(j_1 - j_0) c_1}
  \end{equation}
  and let
  \begin{math}
    m_3 := \max \{ j_1 - j_0, \ceil{r_1} \}
  \end{math}.
  
  \if\shortprep1
  {
  Choose \(m_4 \in \positiveintegers\) so that
  \(A(m_3) \subset \setimage{\eta}{\setzeroton{m_4}}\).
  Let
  \begin{math}
    l \in \naturalnumbers
  \end{math},
  \begin{math}
    l > m_4
  \end{math}.
  If
  \begin{math}
    \zos(\eta(l)) = \finitezeroseq{n}
  \end{math}
  we have
  \begin{equation}
    \label{eq:scf-bound}
    \norminfty{(\beta^{(m_4)}_{\finitezeroseq{n}, j_0, 
    \firstznvar})_
    {\firstznvar \in \zn}} < \frac{h}{4} .
  \end{equation}
  If
  \begin{math}
    \zos(\eta(l)) \neq \finitezeroseq{n}
  \end{math}
  we have
  \begin{equation}
    \label{eq:head-bound}
    \sznorminspace{\left( 2^{j \sigma}
    \norminfty{\left( \beta^{(m_4)}_{\zos,j,\firstznvar}\right)
    _{\zos \in \zeroonesetnnozero{n},
    \firstznvar \in \zn}}
    \right)_{j = j_0}^{j_1 - 1}}{\littlelpcv{q}{j_0 + 
    \naturalnumbers}}
    < \frac{1}{4} h .
  \end{equation}
  Thus by \myequations \eqref{eq:scf-bound} and \eqref{eq:head-bound}
  we get
  \begin{equation}
    \label{eq:wavelet-bound}
    \sznorminspace{
      \left(
      2^{j \sigma}
      \norm{
        \left(
          \beta^{(m_4)}_{\zos,j,\firstznvar}
        \right)_
        {\zos \in \zeroonesetnnozero{n}, \firstznvar \in \zn}
      }_\infty
      \right)_{j = j_0}^\infty
    }{\littlelpcv{q}{j_0 + \naturalnumbers}}
    < \frac{h}{2} .
  \end{equation}
  By \myequations \eqref{eq:scf-bound} and \eqref{eq:wavelet-bound}
  we have
  \begin{math}
    \besovinterpcoeffnorm{n}{j_0}{\sigma}{\infty}{q}{\xi_{m_4}}  
    < h
  \end{math}.
  It follows that
  \begin{displaymath}
    s_m := \sum_{l = 0}^m
    \beta_{\eta(l)} \psi_{\eta(l)}
    \to g
  \end{displaymath}
  as \(m \to \infty\)
  where the series converges in Banach space
  \begin{math}
    \besovspace{\sigma}{\infty}{q}{\rn}
  \end{math}.
  We also have
  \begin{math}
    s_m \to g
  \end{math}
  as
  \begin{math}
    m \to \infty
  \end{math}
  in Banach space
  \begin{math}
    \vanishingfunccv{\rn}
  \end{math}.
  If we had
  \begin{math}
    g \neq f
  \end{math}
  we would have
  \begin{math}
    \dualappl{\psidual_\gamma}{f - g}
    = \dualappl{\psidual_\gamma}{f} -
    \dualappl{\psidual_\gamma}{g}
    \neq 0
  \end{math}
  for some \(\gamma \in \bmiset{n}{j_0}\).
  Now
  \begin{math}
    \dualappl{\psidual_\gamma}{g}
    = \beta_\gamma
    = \dualappl{\psidual_\gamma}{f}
  \end{math},
  which is a contradiction.
  Hence \(f = g\).
  }
  \else
  {
  Choose \(m_4 \in \positiveintegers\) so that
  \(A(m_3) \subset \setimage{\eta}{\setzeroton{m_4}}\).
  Suppose that \(j_2 \in \integernumbers\),
  \(j_0 \leq j_2 < j_1\).
  Let \(l \in \integernumbers\), \(l > m_4\), and
  \(j(\eta(l)) = j_2\).
  If \(\zos(\eta(l)) = \finitezeroseq{n}\)
  we have \(j(\eta(l)) = j_0\) and
  \begin{math}
    \normtwo{\firstznvar(\eta(l))} > 2^{j_0} m_3
    \geq 2^{j_0} r_1
  \end{math}
  from which it follows that
  \begin{displaymath}
    \abs{\szdualappl{\psidual_{\eta(l)}}{f}}
    =
    \abs{f \left( \frac{\firstznvar(\eta(l))}{2^{j^0}} \right)}
    < \frac{h}{4} .
  \end{displaymath}
  Consequently
  \begin{equation}
    \label{eq:scf-bound}
    \norminfty{(\beta^{(m_4)}_{\finitezeroseq{n}, j_0, 
    \firstznvar})_
    {\firstznvar \in \zn}} < \frac{h}{4} .
  \end{equation}
  Suppose then that
  \(\zos(\eta(l)) \neq \finitezeroseq{n}\).
  Now
  \(\normtwo{\firstznvar(\eta(l))} >
  m_2 + 2^{j_0 + m_3} m_3\)
  from which it follows that
  \begin{math}
    \forall \secondznvar \in \zn :
    \norminfty{\secondznvar} \leq m_1
    \implies
    \normtwo{\firstznvar(\eta(l)) + \secondznvar}
    \geq
    2^{j_0 + m_3} m_3
  \end{math}.
  We have \(2^{j_0 + m_3} > 2^{j_2}\)
  and
  \begin{math}
    \forall \secondznvar \in \zn :
    \norminfty{\secondznvar} \leq m_1
    \implies
    \normtwo{\firstznvar(\eta(l)) + \secondznvar}
    >
    2^{j_2} m_3 
    \geq
    2^{j_2} r_1
  \end{math}
  By \myequations \eqref{eq:functional-bound}
  and \eqref{eq:vanishing-bound}
  we have
  \begin{displaymath}
    \abs{\dualappl{\psidual_{\eta(l)}}{f}}
    <
    \frac{h}{2^{n + 2 + j_1 \sigma} (j_1 - j_0)} .
  \end{displaymath}
\if\longversion0
{
  Hence
  \begin{displaymath}
    2^{j_2 \sigma} \norminfty{
    \left( \beta^{(m_4)}_{\zos(\eta(l)),j_2,
    \firstznvar}
    \right)_{\firstznvar \in \zn}}
    <
    \frac{h}{2^{n + 2} (j_1 - j_0)} .
  \end{displaymath}
  Consequently
  \begin{equation}
    \label{eq:head-bound}
    \sznorminspace{\left( 2^{j \sigma}
    \norminfty{\left( \beta^{(m_4)}_{\zos,j,\firstznvar}\right)
    _{\zos \in \zeroonesetnnozero{n},
    \firstznvar \in \zn}}
    \right)_{j = j_0}^{j_1 - 1}}{\littlelpcv{q}{j_0 + 
    \naturalnumbers}}
    < \frac{1}{4} h .
  \end{equation}
}
\else
{
  Hence
  \begin{displaymath}
    \norminfty{\left( \beta^{(m_4)}_{\zos(\eta(l)),j_2,
    \firstznvar}
    \right)_{\firstznvar \in \zn}}
    <
    \frac{h}{2^{n + 2 + j_1 \sigma} (j_1 - j_0)} .
  \end{displaymath}
  from which it follows that
  \begin{displaymath}
    2^{j_2 \sigma} \norminfty{
    \left( \beta^{(m_4)}_{\zos(\eta(l)),j_2,
    \firstznvar}
    \right)_{\firstznvar \in \zn}}
    <
    \frac{h}{2^{n + 2} (j_1 - j_0)} .
  \end{displaymath}
  Consequently
  \begin{eqnarray}
    \nonumber
    & & \sznorminspace{\left( 2^{j \sigma}
    \norminfty{\left( \beta^{(m_4)}_{\zos,j,\firstznvar}\right)
    _{\zos \in \zeroonesetnnozero{n},
    \firstznvar \in \zn}}
    \right)_{j = j_0}^{j_1 - 1}}{\littlelpcv{q}{j_0 + 
    \naturalnumbers}} \\
    \nonumber
    & & <
    \left( \left(
      \frac{h}{2^{n + 2}(j_1 - j_0)}
    \right)^q
    (j_1 - j_0) \right)^{\frac{1}{q}} \\
    \label{eq:head-bound}
    & & < \frac{1}{4} h .
  \end{eqnarray}
}
\fi
  Thus by \myequations \eqref{eq:tail-bound}
  and \eqref{eq:head-bound}
\if\longversion0
{
 \begin{displaymath}
    \sznorminspace{
      \left(
      2^{j \sigma}
      \norm{
        \left(
          \beta^{(m_4)}_{\zos,j,\firstznvar}
        \right)_
        {\zos \in \zeroonesetnnozero{n}, \firstznvar \in \zn}
      }_\infty
      \right)_{j = j_0}^\infty
    }{\littlelpcv{q}{j_0 + \naturalnumbers}}^q \\
    < \left( \frac{h}{4} \right)^q +
    \left( \frac{h}{4} \right)^q
    = \frac{2}{4^q} h^q .
  \end{displaymath}
}
\else
{
  \begin{eqnarray*}
    & & \sznorminspace{
      \left(
      2^{j \sigma}
      \norm{
        \left(
          \beta^{(m_4)}_{\zos,j,\firstznvar}
        \right)_
        {\zos \in \zeroonesetnnozero{n}, \firstznvar \in \zn}
      }_\infty
      \right)_{j = j_0}^\infty
    }{\littlelpcv{q}{j_0 + \naturalnumbers}}^q \\
    & & =
    \sznorminspace{
      \left(
      2^{j \sigma}
      \norm{
        \left(
          \beta^{(m_4)}_{\zos,j,\firstznvar}
        \right)_
        {\zos \in \zeroonesetnnozero{n}, \firstznvar \in \zn}
      }_\infty
      \right)_{j = j_0}^{j_1 - 1}
    }{\littlelpcv{q}{\{j_0,\ldots,j_1 - 1\}}}^q \\
    & & +
    \sznorminspace{
      \left(
      2^{j \sigma}
      \norm{
        \left(
          \beta^{(m_4)}_{\zos,j,\firstznvar}
        \right)_
        {\zos \in \zeroonesetnnozero{n}, \firstznvar \in \zn}
      }_\infty
      \right)_{j = j_1}^\infty
    }{\littlelpcv{q}{j_1 + \naturalnumbers}}^q \\
    & < & \left( \frac{h}{4} \right)^q +
    \left( \frac{h}{4} \right)^q
    = \frac{2}{4^q} h^q .
  \end{eqnarray*}
}
\fi
  It follows that
  \begin{equation}
    \label{eq:wavelet-bound}
    \sznorminspace{
      \left(
      2^{j \sigma}
      \norm{
        \left(
          \beta^{(m_4)}_{\zos,j,\firstznvar}
        \right)_
        {\zos \in \zeroonesetnnozero{n}, \firstznvar \in \zn}
      }_\infty
      \right)_{j = j_0}^\infty
    }{\littlelpcv{q}{j_0 + \naturalnumbers}}
    < \frac{h}{2} .
  \end{equation}
  By \myequations \eqref{eq:scf-bound} and \eqref{eq:wavelet-bound}
  we have
  \begin{math}
    \besovinterpcoeffnorm{n}{j_0}{\sigma}{\infty}{q}{\xi_{m_4}}  
    < h
  \end{math}.
  Hence
  \begin{math}
    \besovinterpcoeffnorm{n}{j_0}{\sigma}{\infty}{q}{\xi_m}  
    \to 0
  \end{math}
  as \(m \to \infty\).
  It follows that
  \begin{displaymath}
    \sznorminspace{
    \sum_{l = 0}^\infty
      \beta^{(m)}_{\eta(l)} \psi_{\eta(l)}
    }
    {\besovspace{\sigma}{\infty}{q}{\rn}}
    =
    \sznorminspace{
    \sum_{l = m + 1}^\infty
      \beta_{\eta(l)} \psi_{\eta(l)}
    }
    {\besovspace{\sigma}{\infty}{q}{\rn}}
    \to 0
  \end{displaymath}
  as \(m \to \infty\).
  Using
  \mydef \ref{def:contbesovnorm}
  we get
  \begin{math}
    \norminspace{\xi_m}{\vanishingfunccv{\rn}}
    \to 0
  \end{math}
  as \(m \to \infty\).
  We also have
  \begin{displaymath}
    s_m := \sum_{l = 0}^m
    \beta_{\eta(l)} \psi_{\eta(l)}
    \to g \in \besovspace{\sigma}{\infty}{q}{\rn}
  \end{displaymath}
  from which it follows that
  \begin{math}
    \norminspace{g - s_m}{\besovspace{\sigma}{\infty}{q}{\rn}}
    \to 0
  \end{math}
  (use some equivalent norm of the Besov space)
  as \(m \to \infty\).
  It follows from \mydef \ref{def:contbesovnorm}
  that
  \begin{math}
    \norminspace{g - s_m}{\vanishingfunccv{\rn}}
    \to 0
  \end{math}
  as \(m \to \infty\).
  Hence \(g \in \vanishingfunccv{\rn}\)
  and
  \(s_m \to g\) in \(\vanishingfunccv{\rn}\)
  as \(m \to \infty\).
  If we had \(g \neq f\) we would have
  \begin{math}
    \dualappl{\psidual_\gamma}{f - g}
    = \dualappl{\psidual_\gamma}{f} -
    \dualappl{\psidual_\gamma}{g}
    \neq 0
  \end{math}
  for some \(\gamma \in \bmiset{n}{j_0}\).
  Now \(\gamma = \eta(l_0)\) for some
  \(l_0 \in \naturalnumbers\) and
  \(\dualappl{\psidual_{\eta(l_0)}}{g}
  = \beta_{\eta(l_0)}
  = \dualappl{\psidual_{\eta(l_0)}}{f}\),
  which is a contradiction.
  Hence \(f = g\).
  }
  \fi
\end{proof}

\begin{theorem}
  \label{th:wavelet-basis-pointwise}
  Let \(n \in \positiveintegers\) and
  \(j_0 \in \integernumbers\).
  Let \(\sigma \in \positiverealnumbers\),
  \(p \in [ 1, \infty ]\), \(q \in [ 1, \infty ]\), and \(\sigma > n / p\).
  Let \(\mdimmothersf{n}\) be a mother scaling function of a
  compactly supported tensor product MRA \mraprepspace
  \(\ucfunccv{\rn}\). Suppose that
  \(\mdimmothersf{n} \in \holderspace{r}{\rn}\)
  for some \(r \in \positiverealnumbers\), \(r > \sigma\).
  Let \(f \in \besovspace{\sigma}{p}{q}{\rn}\).
  Then
  \begin{displaymath}
    f(\rnx) = \sum_{\alpha \in \bmiset{n}{j_0}}
    \dualappl{\psidual_\alpha}{f}
    \psi_\alpha(\rnx)
  \end{displaymath}
  for all \(\rnx \in \rn\) and
  the series above converges absolutely
  for each \(\rnx \in \rn\).
\end{theorem}

\begin{proof}
  Let
  \begin{displaymath}
    a := \norm{
      \left(
        2^{\left(\sigma - \frac{n}{p}\right) j}
        \norm{
          \left(
          \dualappl{\mdimgendualwavelet{n}{\zos}{j}{\firstznvar}}{f}
          \right)_{\zos \in \zeroonesetnnozero{n}, \firstznvar \in \zn}
        }_p
      \right)_{j=j_0}^\infty
    }_q .
  \end{displaymath}
  Now
  \begin{equation}
    \label{eq:coeff-bound}
    \abs{\dualappl{\mdimgendualwavelet{n}{\zos}{j}{\firstznvar}}{f}}
    \leq
    2^{-\left( \sigma - \frac{n}{p} \right) j} a
  \end{equation}
  for all
  \begin{math}
    \zos \in \zeroonesetnnozero{n}
  \end{math},
  \begin{math}
    j \in \naturalnumbers + j_0
  \end{math},
  and
  \begin{math}
    \firstznvar \in \zn
  \end{math}.
  Let
  \begin{displaymath}
    A_0(\rnx) :=
    \left\{
      (\finitezeroseq{n}, j_0, \firstznvar)
      \setsep
      \rnx \in \suppop \mdimsf{n}{j_0}{\firstznvar},
      \firstznvar \in \zn
    \right\}
  \end{displaymath}
  for each \(\rnx \in \rn\)
  and
  \begin{displaymath}
    A_k(\rnx) :=
    \left\{
      (\zos, j_0 + k - 1, \firstznvar)
      \setsep
      \rnx \in \suppop \mdimgenwavelet{n}{\zos}{j_0 + k - 1}{\firstznvar},
      \zos \in \zeroonesetnnozero{n},
      \firstznvar \in \zn
    \right\}
  \end{displaymath}
  for each \(\rnx \in \rn\) and
  \(k \in \positiveintegers\).
  Let \(\eta : \naturalnumbers \onto \bmiset{n}{j_0}\)
  be a bijection and define
  \begin{equation}
    \label{eq:beta-def}
    \beta^{(m)}_\alpha
    :=
    \left\{
      \begin{array}{ll}
        \dualappl{\psidual_\alpha}{f} ; &
        \alpha \notin \setimage{\eta}{\setzeroton{m}} \\
        0 ; &
        \alpha \in \setimage{\eta}{\setzeroton{m}} \\
      \end{array}
    \right.
  \end{equation}
  for each \(m \in \naturalnumbers\) and
  \(\alpha \in \bmiset{n}{j_0}\).
  Let
  \begin{displaymath}
    g_m(\rnx)
    :=
    \sum_{k=0}^m
    \dualappl{\psidual_{\eta(k)}}{f}
    \psi_{\eta(k)}(\rnx)
  \end{displaymath}
  for all \(\rnx \in \rn\) and
  \(m \in \naturalnumbers\).
  It follows from \mytheorem \ref{th:upj-convergence} that
  \begin{displaymath}
    f(\rnx)
    =
    \sum_{k=0}^\infty
    \sum_{\alpha \in A_k(\rnx)}
    \dualappl{\psidual_\alpha}{f}
    \psi_\alpha(\rnx)
  \end{displaymath}
  for all \(\rnx \in \rn\).
  Define \(c_1 := \max \{ \nsznorminfty{\mdimmotherwavelet{n}{\zos}}
    \setsep \zos \in \zeroonesetn \}\).
  Let \(\rny \in \rn\) and \(h \in \positiverealnumbers\).
  Define \(m_1 := \max \{ \card{A_k(\rnx)} \setsep \rnx \in \rn\} \).  
\if\longversion1
{
  Now \(m_1 \in \positiveintegers\).
}
\fi
  Let
  \begin{displaymath}
    j_1
    :=
    \max
    \left\{
      j_0,
      \szceil{\left(\sigma - \frac{n}{p}\right)^{-1} \log_2
      \frac{c_1 m_1 a}
      {\left(1 - 2^{-\left(\sigma - \frac{n}{p}\right)}\right) h}}
    \right\} .
  \end{displaymath}
\if\longversion1
{
  Now
  \begin{displaymath}
    2^{-\left(\sigma - \frac{n}{p}\right) j_1} a
    \leq
    \frac{(1 - 2^{-\left(\sigma - \frac{n}{p}\right)})h}{c_1 m_1} .
  \end{displaymath}
}
\fi
  Choose \(m_2 \in \naturalnumbers\) so that
  \(A_k (\rny) \subset \setimage{\eta}{\setzeroton{m_2}}\)
  for all \(k \in \setzeroton{j_1 - j_0}\).
  Suppose that \(m \in \naturalnumbers\), \(m > m_2\),
  and \(\alpha_0 \in A_{j_0+k}\)
  for some \(k_0 \in \setzeroton{j_1 - j_0}\).
\if\shortprep1
{
  By \myequations \eqref{eq:coeff-bound} and \eqref{eq:beta-def} we have
}
\else
{
  Now \(\alpha_0 \in \setimage{\eta}{\setzeroton{m_2}}
  \subset \setimage{\eta}{\setzeroton{m}}\).
  By \myequation \eqref{eq:beta-def} we have
  \begin{math}
    \beta^{(m)}_{\alpha_0} = 0
  \end{math}.
  Hence by \myequation \eqref{eq:coeff-bound}
}
\fi
\if\shortprep1
{
  \begin{eqnarray*}
    \abs{f(\rny)-g_m(\rny)}
    & \leq &
    c_1 \sum_{k=0}^\infty
    \sum_{\alpha \in A_k(\rny)}
    \abs{\beta^{(m)}_\alpha}
    =
    c_1 \sum_{k=j_1-j_0+1}^\infty
    \sum_{\alpha \in A_k(\rny)}
    \abs{\beta^{(m)}_\alpha} \\
    & \leq &
    c_1 m_1 a \cdot 2^{-j_1 \left(\sigma-\frac{n}{p}\right)}
    \frac{1}{1 - 2^{-\left(\sigma-\frac{n}{p}\right)}}
    < h
  \end{eqnarray*}
}
\else
{
  \begin{eqnarray*}
    \abs{f(\rny)-g_m(\rny)}
    & \leq &
    c_1 \sum_{k=0}^\infty
    \sum_{\alpha \in A_k(\rny)}
    \abs{\beta^{(m)}_\alpha}
    =
    c_1 \sum_{k=j_1-j_0+1}^\infty
    \sum_{\alpha \in A_k(\rny)}
    \abs{\beta^{(m)}_\alpha} \\
    & \leq &
    c_1 \sum_{k=j_1-j_0+1}^ \infty
    m_1 \cdot 2^{-(j_0 + k - 1) \left(\sigma - \frac{n}{p}\right)} \cdot a \\
    & = &
    c_1 m_1 a \cdot 2^{-j_1 \left(\sigma-\frac{n}{p}\right)}
    \frac{1}{1 - 2^{-\left(\sigma-\frac{n}{p}\right)}} \\
    & \leq & h
  \end{eqnarray*}
}
\fi
  Thus
  \begin{math}
    g_m(\rnx) \to f(\rnx)
  \end{math}
  as \(m \to \infty\)
  for all \(\rnx \in \rn\).
\end{proof}

\if\elsevier0
{
\bibliographystyle{abbrv}
}
\else
{
\bibliographystyle{elsarticle-num}
}
\fi

\bibliography{interp-mra}

\end{document}